\documentclass[letterpaper, 11pt,  reqno]{amsart}

\usepackage[margin=1.2in,marginparwidth=1.5cm, marginparsep=0.5cm]{geometry}

\usepackage{amsmath, mathrsfs, amssymb,amscd,amsthm,amsxtra, esint}
\usepackage{tikz} 
\usepackage{color}

\usepackage[implicit=true]{hyperref}

\makeatletter
\AtBeginDocument{%
  \@ifundefined{SK@@ref}{}{%
    \ifx\SK@ref\@empty\else
      \renewcommand*\NR@ref@showkeys[1]{\SK@\SK@@ref{#1}}%
    \fi
  }%
}
\makeatother

\usepackage{cases}

\usepackage{upgreek}

\allowdisplaybreaks[2]

\definecolor{gr}{rgb}   {0.,   0.69,   0.23 }
\definecolor{bl}{rgb}   {0.,   0.5,   1. }
\definecolor{mg}{rgb}   {0.85,  0.,    0.85}
\definecolor{yl}{rgb}   {0.8,  0.7,   0.}
\definecolor{or}{rgb}  {0.7,0.2,0.2}

\usetikzlibrary{shapes.misc}
\usetikzlibrary{shapes.symbols}
\usetikzlibrary{shapes.geometric}

\tikzset{
	ddot/.style={circle,fill=white,draw=black,inner sep=0pt,minimum size=0.8mm},
	>=stealth,
	}

\tikzset{
	ddot2/.style={circle,fill=black,draw=black,inner sep=0pt,minimum size=0.8mm},
	>=stealth,
	}

\newtheorem{theorem}{Theorem} [section]

\newtheorem{lemma}[theorem]{Lemma}
\newtheorem{proposition}[theorem]{Proposition}
\newtheorem{remark}[theorem]{Remark}

\newtheorem{definition}[theorem]{Definition}
\newtheorem{corollary}[theorem]{Corollary}

\newtheorem{oldtheorem}{Theorem}

\DeclareMathOperator*{\intt}{\int}

\DeclareMathOperator{\med}{med}

\DeclareMathOperator{\HS}{HS}
\DeclareMathOperator{\com}{com}

\DeclareMathOperator{\Id}{Id}
\DeclareMathOperator{\sgn}{sgn}

\DeclareMathOperator{\Law}{Law}
\DeclareMathOperator{\Ker}{Ker}

\newcommand{\1}{\hspace{0.5mm}\text{I}\hspace{0.5mm}}
\newcommand{\II}{\text{I \hspace{-2.8mm} I} }
\newcommand{\I}{\mathcal{I}}

\newcommand{\noi}{\noindent}
\newcommand{\Z}{\mathbb{Z}}
\newcommand{\R}{\mathbb{R}}
\newcommand{\C}{\mathbb{C}}
\newcommand{\T}{\mathbb{T}}
\newcommand{\bul}{\bullet}

\let\Re=\undefined\DeclareMathOperator*{\Re}{Re}
\let\Im=\undefined\DeclareMathOperator*{\Im}{Im}

\let\P= \undefined
\newcommand{\P}{\mathbf{P}}

\newcommand{\PP}{\mathbb{P}}

\newcommand{\E}{\mathbb{E}}

\renewcommand{\H}{\mathcal{H}}

\newcommand{\CC}{\mathcal{C}}

\renewcommand{\L}{\mathcal{L}}

\newcommand{\F}{\mathcal{F}}

\newcommand{\NB}{\mathbb{N}}
\newcommand{\low}{\textup{low}}
\newcommand{\high}{\textup{high}}

\newcommand{\al}{\alpha}
\newcommand{\be}{\beta}
\newcommand{\dl}{\delta}
\newcommand{\updl}{\updelta}

\newcommand{\Dl}{\Delta}
\newcommand{\eps}{\varepsilon}

\newcommand{\g}{\gamma}
\newcommand{\G}{\Gamma}
\newcommand{\ld}{\lambda}
\newcommand{\Ld}{\Lambda}
\newcommand{\s}{\sigma}
\newcommand{\Si}{\Sigma}
\newcommand{\ft}{\widehat}

\newcommand{\wt}{\widetilde}
\newcommand{\cj}{\overline}
\newcommand{\dx}{\partial_x}
\newcommand{\dt}{\partial_t}
\newcommand{\dd}{\partial}

\newcommand{\ta}{\theta}

\renewcommand{\l}{\ell}
\renewcommand{\o}{\omega}
\renewcommand{\O}{\Omega}

\newcommand{\bi}{\bfseries\itshape} 

\newcommand{\Gdl}{\mathcal{G}_{\dl} }

\newcommand{\Gd}{\wt{\mathcal{G}}_\dl}

\newcommand{\les}{\lesssim}
\newcommand{\ges}{\gtrsim}

\newcommand{\jb}[1]
{\langle #1 \rangle}

\newcommand{\ind}{\mathbf 1}

\renewcommand{\S}{\mathcal{S}}

\newcommand{\M}{\mathcal{M}}

\newcommand{\N}{\mathbb{N}}

\newcommand{\NN}{\mathcal{N}}
\newcommand{\cL}{\mathcal{L}}
\newcommand{\cC}{\mathcal{C}}

\newcommand{\cX}{\mathcal{X}}
\newcommand{\cY}{\mathcal{Y}}

\newcommand{\cD}{\mathscr{D}}

\newcommand{\W}{\mathcal{W}}

\newcommand{\Lip}{\mathrm{Lip}}
\newcommand{\uw}{U^w}
\newcommand{\uu}{\mathbf{u}}
\newcommand{\vv}{\mathbf{v}}
\newcommand{\z}{\zeta}

\newcommand{\Ta}{\Theta}

\newcommand{\sub}{\substack}

\newcommand{\BO}{\text{\rm BO} }
\newcommand{\KDV}{\text{\rm KdV} }
\newcommand{\MKDV}{\text{\rm mKdV} }

\newcommand{\too}{\longrightarrow}

\newtheorem*{ackno}{Acknowledgements}

\numberwithin{equation}{section}
\numberwithin{theorem}{section}

\makeatletter
\@namedef{subjclassname@2020}{\textup{2020} Mathematics Subject Classification}
\makeatother

\begin{document}
\baselineskip = 14pt

\title[Nonlinear PDEs with modulated dispersion II]
{Nonlinear PDEs with modulated dispersion II:\\
Korteweg-de Vries equation}

\author[K.~Chouk, M.~Gubinelli, G.~Li, J.~Li, and T.~Oh]
{Khalil Chouk, Massimiliano Gubinelli, Guopeng Li, Jiawei Li,  and Tadahiro Oh}

\address{Khalil Chouk\\
\begin{itemize}
\item[-]
 Universit\'e Paris Dauphine, Place du Mar\'echal De Lattre De Tassigny 75775, Paris cedex 16, France,  
\item[-] 
 Humboldt-Universitaet zu Berlin
Institut f\"ur Mathematik
10099 Berlin, Germany,
\item[-] 
School of Mathematics\\
The University of Edinburgh, and the Maxwell Institute for the Mathematical Sciences. James Clerk Maxwell Building, The King's Buildings, Peter Guthrie Tait Road, Edinburgh, EH9 3FD, United Kingdom
\end{itemize}
}

\email{khalil.chouk@gmail.com}

\address{
Massimiliano Gubinelli\\
\begin{itemize}
\item[-]
 Universit\'e Paris Dauphine, Place du Mar\'echal De Lattre De Tassigny 75775, Paris cedex 16, France,
\item[-]    
Institut Universitaire de France, France, 
\item[-]
Hausdorff Center for Mathematics \&  Institut f\"ur Angewandte Mathematik\\
 Universit\"at Bonn\\
Endenicher Allee 60\\
D-53115 Bonn\\
Germany,
\item[-]
Mathematical Institute\\ University of Oxford\\ United Kingdom
\end{itemize}
}

\email{gubinelli@maths.ox.ac.uk}

\address{
Guopeng Li, 
\begin{itemize}
\item[-]
School of Mathematics and Statistics, Beijing Institute of Technology,
Beijing 100081, China\\
\item[-]
School of Mathematics\\
The University of Edinburgh\\
and The Maxwell Institute for the Mathematical Sciences\\
James Clerk Maxwell Building\\
The King's Buildings\\
Peter Guthrie Tait Road\\
Edinburgh\\ 
EH9 3FD\\
 United Kingdom
\end{itemize}
}

\email{guopeng.li@bit.edu.cn}

\address{Jiawei Li, School of Mathematics\\
The University of Edinburgh\\
and The Maxwell Institute for the Mathematical Sciences\\
James Clerk Maxwell Building\\
The King's Buildings\\
Peter Guthrie Tait Road\\
Edinburgh\\ 
EH9 3FD\\
 United Kingdom}

\email{jiawei.li@ed.ac.uk}

%
%

%
\address{
Tadahiro Oh
\begin{itemize}
\item[-]
School of Mathematics and Statistics, Beijing Institute of Technology,
Beijing 100081, China,
\item[-]
School of Mathematics\\
The University of Edinburgh\\
and The Maxwell Institute for the Mathematical Sciences\\
James Clerk Maxwell Building\\
The King's Buildings\\
Peter Guthrie Tait Road\\
Edinburgh\\ 
EH9 3FD\\
 United Kingdom
\end{itemize}
}

\email{hiro.oh@ed.ac.uk}



\begin{abstract}

We study dispersive equations with a time non-homogeneous modulation acting on the linear dispersion term. As  primary models, we consider the Korteweg-de~Vries equation (KdV) and related equations such as the Benjamin-Ono equation (BO) and the intermediate long wave equation (ILW), imposing certain irregularity conditions on the time non-homogeneous modulation. In this work, we establish phenomena called {\it regularization by noise} in three-folds: (i)~When the modulation is sufficiently irregular, we show that the modulated KdV on both the circle and the real line is locally well-posed in the regime where the (unmodulated) KdV equation is known to be ill-posed. In particular, given {\it any} $s \in \mathbb R$, we show that the modulated KdV on the circle with a sufficiently irregular modulation is locally well-posed in $H^s(\mathbb T)$. For example, this result implies that if the modulation is given by a fractional Brownian motion with Hurst index $0 < H < \frac 23$, the modulated KdV on the circle is locally well-posed in $H^{-\frac{1}{2H}+\varepsilon}(\mathbb{T})$ for any $\varepsilon > 0$. Moreover, by adapting the $I$-method to the current modulated setting, we prove global well-posedness of the modulated KdV in negative Sobolev spaces.  (ii)~It is known that  (semilinear) dispersive equations such as KdV and BO exhibit  quasilinear nature below certain  regularity thresholds.
We show that sufficiently irregular modulations make the modulated versions of these equations semilinear by establishing their local well-posedness by a contraction argument, providing local Lipschitz continuity of the solution map. (iii)~We also prove nonlinear smoothing for these modulated equations, where we show that a gain of regularity of the nonlinear part becomes (arbitrarily) larger for more irregular modulations.

As applications of our approach, we also include the following examples: (iv)~deep-water and shallow-water convergence of the modulated (scaled) ILW in arbitrary negative Sobolev spaces. (v)~the modulated KdV on the circle with white noise initial data. (vi)~stochastic modulated KdV on the circle with an additive space-time white noise and singular additive noises of arbitrarily low regularity.

\end{abstract}

%
\maketitle

\tableofcontents

%
%
%

\newpage

\section{Introduction}\label{SEC:1}

\subsection{Modulated dispersive equations}
In this paper, we study  a modulated dispersive equation of the following form:
\begin{align}
\begin{cases}
\dt u +  L u \cdot \dt w=  \NN(u)  \\ 
u|_{t = 0} = u_0,
\end{cases}
\qquad ( t, x) \in \R_+ \times \M,
\label{Maineq}
\end{align}

\noi
where
$\M = \R$ or $  \T = \R/ (2\pi \Z)$,\footnote{By convention, we endow
$\T$ with the normalized Lebesgue measure $ dx_\T =  (2\pi)^{-1}dx$
such that we do not need to carry factors involving $2\pi$.}
 $u :\R_+\times \M \to \R$ or $\C$ is the unknown,
the operator  $L$ denotes  an (unbounded) linear dispersive operator, 
 and $\NN(u)$ denotes the nonlinearity.
Here, 
 $w:\R_+\to\R$ is an arbitrary continuous function of time, 
called a {\it modulation}, 
which will play a central role in our analysis.
Previously, 
de Bouard and Debussche \cite{DD} and Debussche and Tsutsumi~\cite{DT} 
used stochastic calculus to study 
the modulated nonlinear Schr\"odinger equations (NLS) with the modulation $w$ given by a Brownian motion;
see also \cite{DR2, Ste}.
For a general modulation $w$, 
such an approach based on stochastic calculus, 
interpreting~\eqref{Maineq}
as an Ito or Stratonovich stochastic PDE
is not available, and our goal in this paper is to develop a {\it pathwise} approach 
to study the modulated dispersive equation~\eqref{Maineq}, 
exploiting {\it irregularity} of the modulation function $w$;
see Definition \ref{DEF:ir}.
In the previous work
\cite{CG1}, the first two authors developed a pathwise
approach to study 
well-posedness of the modulated NLS. 
In this paper, we study
the following modulated 
 Korteweg-de Vries equation (KdV)
 on both the circle and the real line:
\begin{equation}
\label{kdv1}
\dt u+  \dx^3 u \cdot \dt w =\dx u^2
\end{equation}

\noi
as our primary model example. 
The (unmodulated) KdV equation:
\begin{equation}
\dt u+  \dx^3 u  =\dx u^2, 
\label{kdv2}
\end{equation}

\noi
derived as a model for uni-directional shallow water waves, 
is of fundamental importance
 in mathematics and has been studied extensively
 from both the applied and theoretical viewpoints.
The modulated KdV \eqref{kdv1} naturally appears 
as a model 
for
weakly nonlinear long waves in an inhomogeneous waveguide;
see \cite{CMG, HZ}, where 
the modulation $w$ is taken to be periodic but not differentiable.
There is, however, no analytical study  on the modulated KdV~\eqref{kdv1}
and we hope to make an important first step in this direction.
We also study the related modulated dispersive equations, 
appearing in the study of fluid motion, 
such as 

\smallskip
\begin{itemize}
\item modulated mKdV equation \eqref{mkdv1}, 

\smallskip
\item 
modulated Benjamin-Ono equation \eqref{BO}, 

\smallskip
\item 
modulated intermediate long wave equation \eqref{ILW1}, 

\smallskip
\item 
stochastic modulated KdV equation with an additive forcing \eqref{skdv1}.

\end{itemize}

Our main motivation to study the modulated dispersive equation of 
form \eqref{Maineq} 
is 
to underline the regularization effect of a non-homogeneous time modulation,
in the spirit of the  work of Flandoli, Priola, and the second author on the stochastic transport equation; 
see~\cite{FGP}. 
The so-called {\it regularization by noise}
phenomenon
 has been observed for different models:

\smallskip

\begin{itemize}
\item
many results in the SDE case since the 70's
\cite{Zv, Ver, KR, Davie, DF, 
CatellierGubinelli, DFRV},

\smallskip
\item 
stochastic parabolic equations
\cite{GP, GLT, DGT12, BM}, 
linear transport equations with rough drifts~\cite{FGP}, 
scalar conservation laws \cite{CGess},  
Hamilton-Jacobi equations \cite{LS},

\smallskip
\item 
regularization by non-random noises:
\cite{GG0, FHLN, RT}.

\end{itemize}

\noi
See \cite{RT} for further references.
We refer interested readers to  survey works \cite{Flan, Gess}.

We point out, however, that  rather few results 
on regularization by noise are known
for  dispersive equations, mainly limited to random initial data\,/\,additive noises of super-critical regularity;
see Remark~\ref{REM:random1}.
In this paper, we carry out pathwise analysis to establish several forms of  regularization by noise
for canonical dispersive equations
such as the KdV and BO equations.
See the subsequent subsections
in this introduction
for concrete examples.

The presence of a modulation $w$ in \eqref{Maineq} prevents  
a straightforward use of the Fourier restriction norm method
due to Bourgain \cite{BO93}, 
thus forcing us  
 to develop an alternative approach based on nonlinear Young integration as in~\cite{CG1}.  
See Section \ref{SEC:Young}, 
where we go over a general theory of nonlinear Young integrals
and well-posedness of a Young differential equation.
We point out that 
the nonlinear Young integral approach
we develop here
has 
some  similarity 
with 
the refinement of the Fourier restriction norm method
due to Koch and Tataru~\cite{KT}, 
where an `endpoint Young integral' via the $U^2$-$V^2$ duality
plays a crucial role.
Lastly, we mention recent works
\cite{Tanaka, Robert1, Robert2, Robert3, Robert4}, 
on pathwise well-posedness of  various modulated dispersive equations.
See Remark~\ref{REM:Tristan}
for a further discussion.

\subsection{Irregularity of the modulation}

At this stage, 
the modulated dispersive equation~\eqref{Maineq} 
is merely formal, since the time derivative of the modulation  $w$ does not exist in general.
One of our primary interests
is to  consider the  case where $w$ is a sample path of a stochastic process
such as a Brownian motion and, more generally,  a fractional Brownian motion. 

To bypass this problem, 
 we will interpret the equation~\eqref{Maineq} in terms of the 
 Duhamel formulation (= mild formulation):
\begin{equation}
u(t) = U^w(t) u_0 +  U^w({t}) \int_0^t  U^w({t'})^{-1}\NN( u(t') ) dt', 
\label{mild1}
\end{equation}

\noi
where $U^w(t) = e^{- w(t) L }$ denotes the modulated linear propagator
and we impose $w(0) = 0$ such that  
 $U^w(0) = \Id$.\footnote{The normalization $w(0) = 0$
 is not an additional restriction since only the time derivative $\dt w$
 appears in the modulated equation \eqref{Maineq}.}
Note that the Duhamel formulation \eqref{mild1}
only involves the modulation $w$ but not its derivative, 
allowing us to study the problem
even when $w$ is only continuous.
The aim of this paper is to study 
the Duhamel formulation \eqref{mild1} under 
the assumption on the  \emph{irregularity} of the modulation~$w$ 
introduced in~\cite{CatellierGubinelli, CG1};
see Definition \ref{DEF:ir}.
In particular,
for the modulated KdV \eqref{kdv1}, 
we will show that if $w$ is sufficiently irregular,
then \eqref{mild1} is well-posed even in Sobolev spaces
of very low regularities.

The next definition concerns a particular notion of irregularity of the modulation $w$, 
introduced in~\cite{CatellierGubinelli, CG1}, 
which  will play a fundamental role  in our analysis.

\begin{definition}
\label{DEF:ir}
\rm

Let $\rho>0$ and  $0 < \g < 1$.
Given $T > 0$, we say that a function $w\in C([0,T];\R)$ is $(\rho,\g)$-irregular 
on the time interval $[0, T]$  if we have
\begin{align}
\|\Phi^w\|_{  \W^{\rho,\g}_T} 
:= \sup_{a\in \R} \sup_{0\leq r < t\leq T} \langle a \rangle^\rho \frac{|\Phi^w_{t,r}(a)|}{|t-r|^\g} 
< \infty, 
\label{rho1}
\end{align}

\noi
where 
\begin{align}
\Phi^w_{t,r}(a)
=\int_r^t e^{i  a w(t') } d t'.
\label{rho2}
\end{align}

\noi
We say that $w$ is 
$(\rho,\g)$-irregular 
on $\R_+$  if it is 
$(\rho,\g)$-irregular 
on $[0, T]$ for each finite $T > 0$.

\end{definition}

For further study on the notion of irregularity, 
we refer readers to 
 recent works \cite{GG, RT}
which appeared after the first version of the current paper.

In order to study \eqref{mild1} with an irregular modulation $w$, 
we  need to exploit the temporal regularity of a solution~$u$.
More precisely, 
we will exploit the temporal regularity 
of the 
{\it modulated interaction representation} 
$\uu$ of the unknown $u$ defined by 
\begin{align}
\uu(t)=\uw(t)^{-1}u(t).
\label{int1} 
\end{align}

\noi
When $w(t) = t$
(i.e.~for usual (unmodulated) dispersive PDEs), the (standard) interaction representation $e^{tL} u(t)$
has played a fundamental  role
in the Fourier restriction norm method
\cite{BO93, KT}
(see also \cite{GTV}, 
where an analogy
to the 
 interaction representation in quantum mechanics was made)
 and the normal form method 
 in both the deterministic and probabilistic settings
 \cite{BIT, KO, GKO, OTz, OST, OW}.
In general, there are two advantages
in working with (modulated) interaction representations:
(i)~it allows us to exploit temporal regularity of
the interaction representation,
while, at the level of the original unknown, positive temporal regularity normally comes 
with a loss of spatial derivative
(precisely due to 
the presence of 
the modulated linear propagator
$U^w(t)$
 in~\eqref{mild1})
 and 
(ii)~it allows us to show multilinear dispersive smoothing 
in a more direct manner.

In terms of the modulated interaction representation $\uu$, 
the Duhamel formulation \eqref{mild1}
is written as 
\begin{equation}
\uu(t) = u_0 + \int_0 ^t \uw(t')^{-1}\NN(\uw(t') \uu(t'))dt'.
\label{mild2}
\end{equation}

\noi
In this paper, we study the fixed point problem \eqref{mild2}
by assuming  a certain regularity condition 
on 
the modulated interaction representation $\uu$.
A particular subspace of $C([0,T];H^s(\M))$ will play a special role in our analysis, 
which we introduce in the next definition.

\begin{definition}
\label{DEF:con1}
\rm

Let  $\M =\R$ or $\T$
and 
$s \in \R$.
 Suppose that 
 $w\in C(I;\R)$
 for an interval $I \subset \R_+$
 is $(\rho,\g)$-irregular for some $\rho >0$ and $\frac12< \g < 1$.
 Then, we
denote by 
\begin{align*}
 \cD_{w}^s(I\times \M) \subset C(I;H^s(\M))
\end{align*}

\noi
the space of {\it paths controlled by $w$}, 
which consists of 
paths $u \in C(I;H^s(\M))$ 
such that 
their modulated interaction representations 
$\uu(t)=\uw(t)^{-1}u(t)$
belong
to  $\CC^{\g}(I;H^s(\M))$; see~\eqref{Ho2a}
for the definition of the $\CC^\g$-norm.

We also define the class
$\dot  \cD_{w}^s(I\times \M) $
in a similar manner.
Namely, 
we say that 
 $u \in C(I; \dot H^s(\M))$ 
belongs to 
$\dot  \cD_{w}^s(I\times \M) $
if its  modulated interaction representation
$\uu$ 
belongs
to  $\CC^{\g}(I; \dot H^s(\M))$.

\end{definition}

We conclude this subsection by presenting an example
of irregular modulations.
In~\cite{CatellierGubinelli}, 
the authors proved that a fractional Brownian motion is a $(\rho,\g)$-irregular function.

\begin{oldtheorem}
\label{THM:A}

Let $\{W_t\}_{t\in \R_+}$ be a fractional Brownian motion 
of Hurst index $H\in(0,1)$.
Then, 
for any $\rho < \frac{1}{2H}$,  
there exists $\frac 12 < \g < 1$
such that,  with probability one,  the sample paths of 
$W$ are $(\rho,\g)$-irregular on $\R_+$.

\end{oldtheorem}

Theorem \ref{THM:A}
shows that there exist continuous paths which are $(\rho, \g)$-irregular for arbitrarily large $\rho$. 
In a later work~\cite{GG}, it was shown  that $(\rho,\g)$-irregularity is a generic property of H\"older functions of sufficiently low regularity; see \cite[p.\,2418 and Theorem 3.1]{GG}.

\begin{oldtheorem}
\label{THM:B}
Let  $d \ge 1$.
Given  any $\delta \in (0,1)$, almost every $\dl$-H\"older continuous  function $w \in C^\delta([0,1];\R^d)$ 
 is $(\rho,\g)$-irregular for any $\rho < \frac1{2\delta}$ with  some $\gamma = \gamma(\rho) \in (\frac 12,1)$.
\end{oldtheorem}

In this statement, ``almost every'' is to be understood according to the notion of {\it prevalence}; see~\cite{GG} for details and the  references therein.

As observed in~\cite{CatellierGubinelli,CG1}, the irregularity of the path $w$ allows 
us to obtain a regularizing effect in the ODE context and 
to prove pathwise well-posedness for 
the modulated  NLS.
In the present paper, we show that irregular modulations 
induce strong regularization effects 
in the case of  KdV
and related  equations.
In order to treat 
 irregular modulations in the sense of Definition~\ref{DEF:ir},
we will use a nonlinear version of Young integration theory introduced in~\cite{CG1}
to give a proper meaning to the integral term in~\eqref{mild2};
see Section~\ref{SEC:Young}, 
where we provide an extensive review on the nonlinear Young integration theory
for non-experts.
These ideas, essentially coming from the controlled  path theory~\cite{Gub04}, were already used in the context of the (unmodulated)  KdV \eqref{kdv2} on the circle  in the work~\cite{GubinelliKdV} 
by the second author, 
where the equation has been studied in Sobolev spaces (and more general Fourier--Lebesgue spaces) 
of negative regularity, 
without relying on the Fourier restriction norm method.
We also point out that the refinement of the Fourier restriction norm method
due to Koch and Tataru~\cite{KT}
is also based on interpreting the Duhamel integral (at the level of the interaction representation)
as an `endpoint Young integral' via the $U^2$-$V^2$ duality.
We refer interested readers to the lecture notes \cite{Koch};
see also \cite{HHK}.
Indeed, this point of view was used in 
recent works
\cite{Robert2, Robert3}
to study modulated dispersive PDEs;
see Remark \ref{REM:Tristan}.
Lastly, we  mention recent works \cite{CLO2, COZ,  CLOO}
on 
pathwise well-posedness of stochastic dispersive PDEs
with multiplicative noises, 
where 
the Young\,/\,rough integration theory
(via the sewing lemma)
are combined with
the Fourier restriction norm method, 
where the random tensor estimate, 
originally introduced in the context 
of random variables by Deng, Nahmod, and Yue \cite{DNY3}
and later adapted  
to multiple stochastic integrals with respect to (fractional) Brownian motions \cite{OWZ, CLO2}, 
plays a crucial role in establishing a pathwise bound
on the effect of multiplicative stochastic forcing; 
see also \cite{CGLLO2, SLO}.

\subsection{Well-posedness results}

In this subsection, 
we state our main results on well-posedness
of the modulated dispersive equations mentioned above.
For simplicity of the presentation, 
we impose
the mean-zero assumption\footnote{For the modulated KdV, BO, and  ILW equations, 
if initial data has non-zero mean $\al_0$, then the following Galilean transformation:
$u(t, x)\mapsto u (t, x - 2\al_0t)-\al_0$, 
along with the conservation of the spatial mean (which can be verified
by arguing as in Section \ref{SEC:GWP1}), 
transforms
these equations with a non-zero mean into the mean-zero versions.
Note that, while  the mean-zero assumption is not necessary
for the modulated mKdV \eqref{mkdv1},
we need to apply a gauge transform and consider its renormalized version \eqref{mkdv2}.
}
 on solutions to
these modulated dispersive equations on the circle
in the remaining part of this paper.

Our primary model example  is the modulated KdV equation \eqref{kdv1}
on both the circle and the real line.
Before stating our main result on the modulated KdV \eqref{kdv1}, 
let us go over the known well-posedness (and ill-posedness) results
on the (unmodulated) KdV~\eqref{kdv2}.
KdV  has attracted wide attention 
from both the applied and theoretical communities
and its well-posedness has been studied
extensively on both the circle and the real line.
In~\cite{BO93}, 
Bourgain introduced the so-called
 Fourier restriction norm method
and proved local and global well-posedness of \eqref{kdv2}
in $L^2(\M)$, $\M = \T$ and $\R$,
which was
improved
by the subsequent works
\cite{KPV96, CKSTT03, Ki09, G09}
(see also \cite{BK})
and now KdV \eqref{kdv2}
is known to be 
globally well-posed
in $H^s(\T)$ for $s \ge - \frac 12$
and in $H^{s}(\R)$ for $s \ge -\frac 34$.
In particular, in these works, 
local well-posedness was shown by a contraction argument, 
exhibiting semilinear nature of the equation;
in particular, the solution map$: u_0\in H^s(\M)\mapsto u\in C([-T, T]; H^s(\M))$
is smooth (in fact, real analytic).\footnote{On the circle, 
the smoothness of the solution map holds
only after  restricting our attention to the subspace 
$H^s_\al(\T)$
of $H^s(\T)$, 
consisting of functions with spatial mean $\al \in \R$.
Recall that, in the periodic setting,  we restrict our attention to $H^s_0(\T)$
in this paper.}
In the following, we refer local well-posedness
by a contraction argument 
as {\it semilinear local well-posedness}.
On the other hand, it is known that the solution map
fails to be locally 
uniformly continuous below the aforementioned regularity thresholds
\cite{CCT}; see also~\cite{BO97}.
Namely,  the KdV equation~\eqref{kdv2} exhibits 
{\it quasilinear} nature in the low regularity setting
(i.e.~below $s = -\frac 12$ on $\T$
and $s=- \frac 34$ on $\R$).

There are also well-posedness results of KdV based on its complete integrability;
see
\cite{KT1,  KV19}
for the sharp global well-posedness of KdV in $H^{-1}(\M)$, $\M = \T$ and $\R$
(see  \cite{M11, M12}
for ill-posedness below $s = -1$).
We note that 
 the time modulation effect 
breaks  the integrability of the KdV equation, 
and thus we aim to develop  a Fourier analytic approach
adapted to the current modulated setting.

\medskip

We now state our main well-posedness
results
for  the modulated KdV \eqref{kdv1} on the circle.

\begin{theorem}\label{THM:1}
Given $\rho \ge\frac 12$,  $\frac12< \g < 1$, and $T> 0$, 
let  $w$ be $(\rho,\g)$-irregular on $[0, T]$ in the sense of Definition~\ref{DEF:ir}.

\smallskip

\noi
\textup{(i)}
Suppose that $\rho \ge \frac 12$ and $s\in \R$
satisfy one of the following conditions\textup{:}
\begin{align}
\begin{split}
\textup{(i.a)} &\ \  \tfrac 12 \le \rho \le \tfrac 34 \quad \text{and} \quad 
s > \tfrac 32 - 3 \rho, \\
\textup{(i.b)} &\ \  \rho > \tfrac 34\quad  \text{and} \quad  s \ge - \rho.
\end{split}
\label{reg1}
\end{align}

\noi
Then, the modulated KdV equation \eqref{kdv1}
on  $\T$
is semilinearly locally well-posed in $H^s(\T)$.

\smallskip

\noi
\textup{(ii)}
In addition, suppose that 
  $w$ is $(\rho,\g)$-irregular on $\R_+$
and  that $\rho > \frac 12$ and $s\in \R$
satisfy one of the following conditions\textup{:}
\begin{align}
\begin{split}
\textup{(ii.a)} &\ \ 
\tfrac 12 < \rho \le \tfrac 3{4\g}
 \quad \text{and} \quad 
s > \tfrac {3- 6\rho}{6-4\g}, \\
\textup{(ii.b)} &\ \ 
\tfrac 3{4 \g}< \rho < \tfrac 32
\quad  \text{and} \quad  s \ge - \rho,\\
\textup{(ii.c)} &\ \ 
  \rho \ge \tfrac 32
\quad  \text{and} \quad 
s  > -\tfrac 32.
\end{split}
\label{s1}
\end{align}

\noi
Then, 
 the modulated KdV equation \eqref{kdv1}
on $\T$
is globally  well-posed in $H^s(\T)$.

\end{theorem}

\begin{remark}\label{REM:pers1}\rm

(i)
Uniqueness of a solution $u$ constructed in Theorem \ref{THM:1}\,(i)
holds 
 in the class $\cD_{w}^s([0, \tau] \times \T)$ 
 defined in Definition \ref{DEF:con1}, 
 where 
 $0 < \tau \le  T$ denotes the local existence time.
A similar comment applies to the other modulated equations
considered in this paper.

\smallskip

\noi 
(ii)
The solution map:
 \[(u_0, w) \in H^s(\T) 
\times  \W^{\rho,\g}_T \longmapsto \uu \in \cC^\g([0, \tau]; H^s(\T))\]

\noi
is locally Lipschitz continuous;
see Remarks \ref{REM:main2}\,(iii) and \ref{REM:conv1}.
A similar comment applies to the other modulated equations
considered in this paper. 

\smallskip

\noi
(iii) All of our local well-posedness results
in this paper
come with {\it persistence of regularity}.
Suppose that  a modulated equation is locally well-posed
in $H^{s}(\M)$ for some $s \in \R$
with local existence time\footnote{Clearly, $\tau$ also depends
on the modulation $w$ but we suppress its dependence on $w$
in this discussion.}
 $\tau = \tau(\|u_0\|_{H^{s}}) \in (0, T]$.
If, in addition, we have $u_0 \in H^{s_0}(\M)$
for some $s_0 > s$, 
then the corresponding solution exists
with the same local existence time $\tau$, 
depending only on the $H^s$-norm of the initial data $u_0$.

\end{remark}

As a corollary to Theorem \ref{THM:1}\,(i)
with Theorem \ref{THM:A}, we have
the following well-posedness
result for the modulated KdV \eqref{kdv1} on $\T$, 
where a modulation is given by a fractional Brownian motion.
For simplicity of the presentation, 
we only state its local well-posedness.

\begin{corollary}\label{COR:1}
Let $\{W_t\}_{t\in \R_+}$ be a fractional Brownian motion 
of Hurst index $H\in(0,1)$.
Suppose that $s \in \R$ satisfies 
\begin{align*}
\textup{(a)} &\ \  \tfrac 23 \le H  <1 \quad \text{and} \quad 
s > \tfrac 32 - \tfrac{3}{2H} , \\
\textup{(b)} &\ \  0 < H < \tfrac 23\quad  \text{and} \quad  s > - \tfrac 1{2H}.
\end{align*}

\noi
Then, 
with probability one, 
the modulated KdV equation \eqref{kdv1}
on  $\T$
with the modulation given by 
a sample path of the
 fractional Brownian motion 
 $\{W_t\}_{t\in \R_+}$
is semilinearly locally well-posed in $H^s(\T)$.

\end{corollary}

\begin{remark}\rm
Similar results hold for each of the dispersive equations considered in this paper, with the appropriate choice of exponents, when the modulation function is given by a fractional Brownian motion, or,  thanks to Theorem \ref{THM:B}, when the modulation function is given by a generic (in the sense of prevalence) H\"older function.
\end{remark}

Next, 
we state our main well-posedness
results
for  the modulated KdV \eqref{kdv1} on the real line.

\begin{theorem}
\label{THM:2}
Given $\rho > \frac 12$,  $\frac12< \g < 1$, and $T> 0$, 
let  $w$ be $(\rho,\g)$-irregular on $[0, T]$ in the sense of Definition~\ref{DEF:ir}.

\smallskip

\noi
\textup{(i)}
Suppose that $\rho >  \frac 12$ and $s\in \R$
satisfy one of the following conditions\textup{:}
\begin{align}
\begin{split}
\textup{(i.a)} &\ \  \tfrac 12 < \rho \le \tfrac 34 \quad \text{and} \quad 
s > \tfrac 32 - 3 \rho, \\
\textup{(i.b)} &\ \  \tfrac 34 <  \rho <  \tfrac 32 \quad \text{and} \quad 
s \ge - \rho, \\
\textup{(i.c)} &\ \  \rho \ge  \tfrac 32\quad  \text{and} \quad  s > -\tfrac 32.
\end{split}
\label{regR1}
\end{align}

\noi
Then, the modulated KdV equation \eqref{kdv1}
on  $\R$
is semilinearly  locally well-posed in $H^s(\R)$.

\smallskip

\noi
\textup{(ii)}
In addition, suppose that 
  $w$ is $(\rho,\g)$-irregular on $\R_+$
and that $\rho > \frac 12$ and $s\in \R$
satisfy~\eqref{s1}.
Then, 
 the modulated KdV equation \eqref{kdv1}
on $\R$
is globally  well-posed in $H^s(\R)$.

\end{theorem}

Theorems \ref{THM:1} and \ref{THM:2}
show
that an irregular modulation provides a \emph{regularization effect} on the KdV equation, 
allowing us to prove well-posedness in the regime, where the unmodulated 
KdV equation \eqref{kdv2} is ill-posed
 (namely, $s < -1$).
We also mention that our construction of a solution is based on a contraction argument,
yielding local Lipschitz continuous dependence.
Such a semilinear property fails for 
 the unmodulated 
KdV equation \eqref{kdv2}
in $H^s(\T)$ for $s <  - \frac 12$
and in $H^{s}(\R)$ for $s <  -\frac 34$, 
thus exhibiting a {\it semilinearization phenomenon}.
Last but not least,  in the periodic case, 
Theorem \ref{THM:1} shows 
a very strong regularization effect;
given {\it any} $s \in \R$, 
the modulated KdV \eqref{kdv1}
on  the circle
with a sufficiently irregular modulation
is locally well-posed in $H^s(\T)$; see also Theorem~\ref{THM:4}
for an analogous result on the modulated BO and  ILW equations.
Theorems~\ref{THM:1} and~\ref{THM:2}
present  the first 
results 
on regularization by noise
in 
the study of nonlinear dispersive equations with rough (deterministic) initial data.
See also Remark \ref{REM:random1}.
In the next subsection, 
we also present another kind of regularization-by-noise phenomenon
in the context of 
{\it nonlinear smoothing};
see Theorem \ref{THM:G1}.

\begin{remark}\rm
Contrary to the periodic case, the local well-posedness result of the modulated KdV \eqref{kdv1}
on the real line 
(Theorem \ref{THM:2}\,(i)) holds only 
in the scaling-subcritical regime
$s > s_\text{crit} = -\frac 32$, even if we take $\rho$ to be large.
This is due to the fact that frequencies can be arbitrarily small on the real line case, 
which prevents us from 
exploiting the irregularity of the modulation function $w$
in a certain frequency configuration.
See Remark~\ref{REM:low1}, 
where we show
the essential sharpness of the condition $s > - \frac 32$ for $\rho \ge \frac 32$ in \eqref{regR1}
(modulo the endpoint $s = -\frac 32$)
by establishing ``mild'' ill-posedness
of the modulated KdV \eqref{kdv1} on the real line for $s < - \frac 32$.

\end{remark}

Let us point  out that our results are limited to 
sufficiently irregular modulations 
and they do not provide any sharp information about the dependence of solutions on modulations, especially about the transition between smooth and irregular modulations. 
Indeed, a bit surprisingly, in the current modulated context,  the application of 
 controlled path techniques (including nonlinear Young integration theory via the sewing lemma) is easier if the modulation is very irregular. 
In fact, by considering
sufficiently irregular modulations, 
we have avoided the need of  
a second-order controlled expansion which was necessary in
the work~\cite{GubinelliKdV}, where the second author developed a controlled path
approach to study the (unmodulated) KdV equation \eqref{kdv2}.
It remains   an open problem to fill the gap between regular and  irregular modulations.

Our main strategy for proving local well-posedness is to make sense of the integral term in~\eqref{mild2}
as a {\it nonlinear Young integral}
and prove local well-posedness of  the equation \eqref{mild2}
via a contraction argument, 
just as in the previous work \cite{CG1}
by the first two authors.
While our approach is equivalent
to that presented in \cite{CG1}, 
we first construct the integral in \eqref{mild2}
as a nonlinear Young integral
$\I^X(\uu)$ with a given driver $X$
in an abstract setting
via the sewing lemma
(Lemma \ref{LEM:sew})
from 
the theory of controlled rough paths~\cite{Gub04}
developed by the second author; see Lemma \ref{LEM:int1}.
We made a particular effort to 
present this part 
in a pedestrian manner, accessible
to non-experts.
We then study well-posedness and related properties of the following
Young differential equation (YDE):
\begin{equation}
\uu = u_0 + \I^X(\uu)
\label{YDE0}
\end{equation}

\noi
in our functional framework adapted to Sobolev spaces, 
under appropriate assumptions 
on the driver $X$; 
see Proposition \ref{PROP:main}.
In Sections~\ref{SEC:LWP} and~\ref{SEC:BO}, 
we verify 
that such assumptions are indeed satisfied
for each of the modulated dispersive equations
we consider in this paper.

The controlled rough path theory, 
introduced 
as an alternative formulation of the rough path theory
due to Lyons
\cite{Lyons1, Lyons2}, 
has been applied and extended further in a variety of contexts; 
for example, 

\smallskip
\begin{itemize}
\item
well-posedness of stochastic parabolic  equations \cite{GT10, DGT12}, 

\smallskip
\item
Hairer's work on well-posedness for the KPZ equation \cite{Hairer1},

\smallskip
\item
Hairer's theory of regularity structures, 
 dealing with very singular parabolic PDEs such as the parabolic  
 $\Phi_3^4$-model \cite{Hairer2};

\smallskip
\item
paracontrolled calculus \cite{GIP, CK, MW}
which has also been applied 
to study  nonlinear wave equations (NLW)
in recent years
\cite{GKO2, OOT1, Bring, OOT2, BDNY}.

\end{itemize}

\noi
Our approach to study modulated dispersive equations, 
exhibiting various regularization-by-noise effects
is an interesting and useful addition to this list.

All our techniques are deterministic and they provide novel results even in the stochastic context; for example,   we can take the modulation $w$ to be the sample path of a fractional Brownian motion
as we have seen in Corollary \ref{COR:1}. 
Even for a Brownian motion, our results on KdV, mKdV, BO, and ILW are novel  to our knowledge.
See Subsections~\ref{SUBSEC:1.5} and~\ref{SUBSEC:1.6}
for applications with random initial data and stochastic forcing.
We also note that 
YDEs
allow for a straightforward Euler approximation scheme
to discretize in time; see, for example,  \cite[Subsection 2.3]{CG1}.

\medskip

In Section \ref{SEC:I}, we prove the global well-posedness results
in Theorems \ref{THM:1}\,(ii) and \ref{THM:2}\,(ii)
by the so-called $I$-method (= method of almost conservation laws)
introduced in \cite{CKSTT03}.
In particular, we prove these results
by 
combining a key crucial commutator estimate (Proposition~\ref{PROP:com}), 
a common  tool in an application 
of the $I$-method, with the sewing lemma (Lemma~\ref{LEM:sew}).
Such a combination of tools and ideas from dispersive PDEs and the controlled rough path theory is novel in the literature.
We also mention a recent work \cite{GKOT}
for 
a highly non-trivial application of 
the $I$-method in the context of the stochastic NLW.
The $I$-method relies on the scaling property
of a given equation (see \eqref{scaling1}) and thus is limited to the scaling-subcritical regime
$s > s_\text{crit} = -\frac 32$ in the KdV case.
This is the reason why the restriction $s > -\frac 32$ appears in Theorems~\ref{THM:1}\,(ii)
and~\ref{THM:2}\,(ii).\footnote{In this paper, 
we apply the classical KdV scaling 
\eqref{scaling1} to $u$, while we apply  \eqref{scaling2}
to the modulation~$w$.
In fact, by adjusting the scaling for $w$, 
we can have 
the modulated KdV (on $\R$) invariant under
a non-KdV scaling  applied to  $u$.
In a recent work \cite{GGLO}, 
the second and fifth
authors with D.\,Greco and S.\,Liu
made use of this observation and 
proved  that, given  any $s \in \R$, 
the modulated KdV on the circle is 
globally  well-posed
in $H^s(\T)$ , provided that the modulation
is sufficiently irregular.

}
On the other hand, it is 
intriguing to see the same scaling-subcritical restriction appears
in  local well-posedness on the real line (Theorem \ref{THM:2}\,(i)).
Our implementation of the $I$-method in the current modulated context is of the first order.
It would be of interest to see if we can improve the range of regularity
by adding a ``correction'' term
and proceed with the second order $I$-method;
see, for example,  \cite{CKSTT03} for a higher order implementation of the $I$-method
in the unmodulated setting.

\begin{remark}\label{REM:random1}
\rm
In \cite{BO94, BO96}, Bourgain initiated the study 
of nonlinear dispersive equations with random initial data of low regularity.
Over the last fifteen years, 
this study has been particularly active
\cite{BO96, BT08, CO, BT3, BOP2, Poc, OP};
see \cite{BOP4} for a (slightly outdated) survey on this topic.
See also more recent works
\cite{GKO2, DNY2, OOT1, Bring, OOT2, BDNY}, including those
on stochastic dispersive PDEs with additive noises
(which can be easily adapted to random initial data).
In these works, 
non-deterministic viewpoints
allowed us to go beyond the limit of deterministic analysis, 
and thus, 
broadly speaking, 
we can  regard them as 
regularization-by-noise phenomena
(by considering random initial data).

We also note that 
a  noise can sometimes  worsen the behavior of an equation; 
see, 
for example, 
\cite{de_bouard_effect_2002, de_bouard_blow-up_2005}
for the construction of 
finite-time blowup  for the stochastic NLS.

\end{remark}

\medskip

Let us now consider other modulated dispersive equations.
For simplicity of presentation, 
we restrict our attention to the periodic case.
We first consider the following
  modulated
 modified KdV equation (mKdV) on $\T$:
\begin{equation}
\label{mkdv1}
\dt u+ \dx^3 u \cdot \dt w=\dx u^3.
\end{equation}

\noi
As in the unmodulated setting \cite{BO93}, 
the equation \eqref{mkdv1} on the circle is not suitable 
for well-posedness studies.
Using the $L^2$-conservation,\footnote{In the current modulated setting, 
we need to justify the $L^2$-conservation; see Section \ref{SEC:GWP1}.} 
the gauge transform: $x \mapsto x -3\|u_0\|_{L^2}^2 t$
converts~\eqref{mkdv1}
into the  following renormalized modulated mKdV:
\begin{equation}
\label{mkdv2}
\dt u+   \dx^3 u \cdot \dt w=3   \big(u^2- \|u\|_{L^2}^2\big)\dx u, 
\end{equation}

\noi
where certain resonant interactions have been removed;
see \eqref{mK0x} and \eqref{mK0}.

\begin{theorem}
\label{THM:3}
Given $\rho \ge \frac 12$,  $\frac12< \g < 1$, and $T> 0$, 
let  $w$ be $(\rho,\g)$-irregular on $[0, T]$ in the sense of Definition~\ref{DEF:ir}.
Then, for $s\geq \frac 12$, 
the modulated \textup{(}renormalized\textup{)} mKdV equation~\eqref{mkdv2}
on $\T$
is semilinearly  locally well-posed in $H^s(\T)$.

\end{theorem}

At the level of the original modulated mKdV \eqref{mkdv1}, 
local well-posedness holds
in a class $H^s(\T) \cap \{ \|f \|_{L^2} = \mu\}$
for any prescribed value of  $\mu \in \R_+$.

\begin{remark}\label{REM:mkdv1}\rm

In the unmodulated setting, 
mKdV on the circle is known to be semilinearly locally well-posed in 
$H^s(\T)$ for $s \ge \frac 12$
\cite{BO93, KO}, 
while 
 the solution map
fails to be locally 
uniformly continuous below $s = \frac 12$
\cite{CCT, TT}; see also \cite{BO97}.
See \cite{KT2, Forlano}
for well-posedness results via the complete integrability of the equation.

Unlike the KdV case, 
Theorem \ref{THM:3}
does not provide any improvement as compared to the unmodulated setting.
This is due to the fact that 
the regularity restriction $s \ge \frac 12$  in Theorem~\ref{THM:3}
comes from  the resonant part of the nonlinearity, 
where the modulation plays no role;
see the proof of Proposition \ref{PROP:mkdv1}.
By adapting  the approach in  \cite{TT, NTT}
(see also \cite{OTzW, BLLZ})
to weaken the resonant interaction,
one may be able to improve Theorem \ref{THM:3} for some $s  < \frac 12$, 
but we do not pursue this issue here.
See also Remarks \ref{REM:Tristan} and \ref{REM:mkdv2}.
We also mention the recent works
 \cite{C1, C2}
for a further discussion on the renormalization of the complex-valued mKdV.

\end{remark}


Next, we consider the following 
modulated Benjamin-Ono equation (BO) on the circle:
\begin{equation}
\label{BO}
\dt u-    \H\dx^2 u \cdot \dt w=\dx u^2,
\end{equation}

\noi
where $\H$ denotes the
Hilbert transform
with the Fourier multiplier $- \ind_{n\ne 0}\cdot i \sgn (n)$.
Given $0 < \dl < \infty$, we also consider the 
modulated 
intermediate long wave equation (ILW) on the circle: 
\begin{equation}
\label{ILW1}
\dt u-    \Gdl\dx^2 u \cdot \dt w =\dx u^2, 
\end{equation}

\noi
where the operator~$\Gdl$ is given 
as the following Fourier multiplier operator:
\begin{align}
\ft{\Gdl f}(n) =
\ft{\Gdl}(n) \ft f(n)
:=
- \ind_{n\ne 0}\cdot i \bigg(\coth(\dl n)  - \frac{1}{\dl n}\bigg) \ft{f}(n), \quad n\in\Z.
\label{GG1}
\end{align}

\noi
In the unmodulated setting (i.e.~$w(t) = t$), 
the ILW equation  
models the propagation of an internal wave at the interface of a stratified fluid of
finite depth $\dl>0$,
and the unknown~$u$ denotes the amplitude of the internal wave at the interface.
Moreover, it serves as an ``intermediate'' equation of finite depth $0 < \dl < \infty$,
providing a natural connection
between the BO equation, modeling fluid of infinite depth  ($\dl = \infty$),
and the  KdV equation \eqref{kdv2}, modeling shallow water ($\dl = 0$).
See \cite{S19, KS21} for an overview of the subject and the references therein.

In considering 
 the shallow-water limit ($\dl\to0$), 
we need to magnify the amplitude $u$
of  the internal wave at the interface
to observe any non-trivial limiting behavior.
For this purpose, we introduce the following scaling transform \cite{ABFS}:
\noi
\begin{align}
v(t,x) = \tfrac{3}{ \dl} u\big(\tfrac{3}{ \dl} t,x\big),
\label{scaling0}
\end{align}

\noi
which converts \eqref{ILW1}
into the following modulated scaled ILW equation,
where the rescaled modulation
$\widetilde w_\dl(t)
= \frac{\dl}{3}w\big(\frac{3}{\dl} t \big)$
is denoted again by $w$:
\begin{equation}
\label{ILW2}
\dt v-    \wt \Gdl\dx^2 v \cdot \dt w =\dx v^2,
\end{equation}

\noi
where 
\begin{align}
 \Gd = \frac{3}{\dl} \Gdl.
\label{GG2}
\end{align}

\noi
In studying shallow-water convergence, 
we fix the modulation $w$ in the modulated scaled ILW~\eqref{ILW2}, 
independent of $\dl$.

We first state well-posedness results
of these modulated equations.

\begin{theorem}
\label{THM:4}
Given $\rho \ge  1$, $\frac12< \g < 1$, and $T > 0$, 
let  $w$ be $(\rho,\g)$-irregular on $[0, T]$ in the sense of Definition~\ref{DEF:ir}.

\smallskip

\noi
\textup{(i)}
Suppose that $\rho \ge 1$ and $s\in \R$
satisfy one of the following conditions\textup{:}
\begin{align}
\begin{split}
\textup{(i.a)} &\ \   \rho =  1 \quad \text{and} \quad s > - \tfrac 12, \\
\textup{(i.b)} &\ \  \rho >  1 \quad  \text{and} \quad  s \ge - \tfrac 12 \rho.
\end{split}
\label{regBO1}
\end{align}

\noi
Then, the modulated BO equation \eqref{BO}
on  $\T$
is semilinearly  locally well-posed in $H^s(\T)$.
Given $0 < \dl < \infty$, 
the modulated 
ILW \eqref{ILW1}
and
the modulated 
scaled ILW \eqref{ILW2}
on  $\T$
are also semilinearly  locally well-posed in $H^s(\T)$.

\smallskip

\noi
\textup{(ii)}
In addition, suppose that 
  $w$ is $(\rho,\g)$-irregular on $\R_+$.
Then, the modulated BO equation~\eqref{BO} 
on $\T$ is globally well-posed in $H^s(\T)$ for any $s \ge 0$.
Similarly, 
given $0 < \dl < \infty$, 
the modulated 
ILW \eqref{ILW1}
and
the modulated 
scaled ILW \eqref{ILW2}
on $\T$ are globally well-posed in $H^s(\T)$ for any $s \ge 0$.

\end{theorem}

In the unmodulated setting, 
Fourier analytic approaches
allow us to show well-posedness of 
BO and the (scaled) ILW  in $L^2(\T)$
\cite{Moli1, Moli2, MP, CLOP};
see also \cite{Deng15}.
Furthermore, 
by the complete integrability of the equation, 
the (unmodulated) BO equation 
is known to be globally well-posed
in the entire scaling-subcritical regime, 
namely in $H^s(\T)$ for any $s > - \frac 12$, 
while it is ill-posed for $s \le -\frac 12$
\cite{GKT23, KLV}.\footnote{This global well-posedness
for $s > - \frac 12$ was recently extended to ILW via a perturbative
argument; see~\cite{GL}.}
See also \cite{CFLOP}
for the study of ILW in negative Sobolev spaces.
In the current modulated setting, 
Theorem \ref{THM:4}
shows that, 
 given  any $s \in \R$, 
the modulated BO, ILW, and scaled ILW
on  the circle
with sufficiently irregular modulations
are locally well-posed in $H^s(\T)$, 
exhibiting a strong regularization-by-noise phenomenon.

\begin{remark}\label{REM:BO}\rm

In the unmodulated setting, 
BO and ILW on the circle
are known to behave 
quasilinearly
on $H^s(\T)$ for $s < 0$.
More precisely, 
semilinear local well-posedness
of these equations
on $H^s(\T)$
fails for  $s < 0$
(even if we restrict our attention to the subclass $H^s_0(\T)$, 
consisting of mean-zero functions);
see \cite[Theorem 1.2]{Moli2}.
See also 
\cite[second paragraph on p.\,640]{Moli2}.
Thus, Theorem~\ref{THM:4}
establishes 
a semilinearization phenomenon on the modulated BO and  ILW equations
in negative Sobolev spaces.

Given $\al \in \R$, let 
$H^s_\al(\T)$
denote the subspace of  $H^s(\T)$
consisting of functions with spatial mean $\al$.
It is worthwhile to note that 
while
Molinet~\cite[Theorem~1.2]{Moli2} showed that the solution map
for the unmodulated BO on $\T$ 
is  analytic
in $H^s_\al(\T)$
for any $s \ge 0$ and $\al \in \R$, 
 there is {\it no} known well-posedness result 
for the unmodulated BO  
in $H^s_\al(\T)$, $s \ge 0$, 
  via a contraction argument
at this point,
which is in sharp contrast to 
the semilinear local well-posedness of 
the modulated BO on $\T$, stated in 
Theorem \ref{THM:4}.
We also mention a recent work~\cite{GKT24}, 
where certain quasilinear nature of the unmodulated BO on $\T$
was manifested.

Lastly, 
we point out  that 
 semilinear local well-posedness 
of the unmodulated  BO and ILW on the real line  
fails in $H^s(\R)$ for {\it any}  
 $s \in \R$; see~\cite{MST, KTz2, CHLO}.
Namely, on the real line, these (unmodulated) equations 
manifest their quasilinear nature in 
a stronger manner.
Thus, it would be of interest to investigate 
if an analogue
of 
Theorem \ref{THM:4} holds on the real line
(which would then 
show that, 
given any $s \in \R$, 
irregular modulations 
semilinearize  these equations).
We, however, do not pursue this issue further in this paper.

\end{remark}

Lastly, we consider deep-water convergence ($\dl \to \infty$)
and shallow-water convergence ($\dl \to 0$)
of the modulated (scaled) ILW.
In the unmodulated setting,  
the convergence  issue
 for the (scaled) ILW
 has been studied in 
 \cite{ABFS, Li24, LOZ,  CLOP,  FLZ, CLO, GL, HKV, CLOZ, CHLO}
from both deterministic and probabilistic viewpoints.
The following theorem establishes
deep-water 
and shallow-water convergence results
for the modulated (scaled) ILW.

\begin{theorem}
\label{THM:conv}
Given $\rho >   1$,  $\frac12< \g < 1$, and $T> 0$, 
let  $w$ be $(\rho,\g)$-irregular on $[0, T]$ in the sense of Definition~\ref{DEF:ir}.

\smallskip

\noi
\textup{(i) (deep-water convergence).}
Let 
$s > - \frac 12 \rho$.
Given $u_0 \in H^s (\T)$,
let $u^\dl$, $0<   \dl < \infty$,  and $u$ denote the local-in-time solutions
to the modulated ILW equation \eqref{ILW1}
and the modulated BO equation \eqref{BO}
with $u^\dl|_{t = 0} = u|_{t = 0} = u_0$, respectively.
Then, as $\dl \to \infty$, 
$u^\dl$ converges to $u$ in 
$C([0, \tau]; H^s(\T))$, 
where $\tau \in (0,  T]$ denotes the local existence time of the modulated BO equation
\eqref{BO}.
Moreover, if 
  $w$ is $(\rho,\g)$-irregular on $\R_+$ 
and 
 $s \ge 0 $, 
then the deep-water convergence holds globally in time.

\smallskip

\noi
\textup{(ii) (shallow-water convergence).}
Let 
$s > - \frac 12 \rho$.
Given $u_0 \in H^s (\T)$,
let $v^\dl$, $0<   \dl < \infty$,  and $u$ denote the local-in-time solutions
to the modulated scaled ILW equation \eqref{ILW2}
and the modulated KdV equation \eqref{kdv1}
with $v^\dl|_{t = 0} = u|_{t = 0} = u_0$, respectively.
Then, as $\dl \to 0$, 
$v^\dl$ converges to $u$ in 
$C([0, \tau]; H^s(\T))$, 
where $\tau \in (0, T]$ denotes the local existence time of the modulated KdV equation
\eqref{kdv1}.
Moreover, if 
  $w$ is $(\rho,\g)$-irregular on $\R_+$ 
and 
 $s \ge 0 $, 
then the shallow-water convergence holds globally in time.

\end{theorem}

It follows from 
the discussion in 
Subsections \ref{SUBSEC:BO3}
and \ref{SUBSEC:BO4}
together with 
Lemma \ref{LEM:OBS1}\,(iv), 
Remark \ref{REM:OBS2}, 
and 
Proposition \ref{PROP:main}\,(iv)
that, at the level of the modulated interaction representations, 
convergence holds
in  $\CC^\g([0, \tau]; H^s(\M))$.
For example, 
$\uu^\dl(t) = e^{-  w(t) \Gdl \dx^2} u^\dl(t)$
converges
to $\uu(t) = e^{ -w(t) \H \dx^2} u(t)$
in  $\CC^\g([0, \tau]; H^s(\M))$
as $\dl \to \infty$.

\begin{remark}\label{REM:Tristan}\rm
(i) After the appearance of the first version of this paper, 
there have been 
interesting well-posedness results on 
 modulated dispersive equations
\cite{DR2, G22, Ste, Tanaka, Robert1, Robert2, Robert3, Robert4}.
In particular, 
in recent preprints \cite{Robert2, Robert3}, 
Robert 
studied pathwise  well-posedness of the modulated dispersive equation \eqref{Maineq}
such as the modulated (fractional) NLS and modulated (fractional) KdV, 
%
%
by working within the framework
of 
 the refined Fourier restriction norm method
 involving the $U^2$- and $V^2$-spaces
(namely,  without the use of the sewing lemma).
In~\cite{Robert3}, he proved 
local well-posedness of 
 the following modulated fractional KdV on $\T$:
\[\dt u+    \dx|\dx|^\be u \cdot \dt w=\dx u^2, \]

\noi
 which essentially 
agrees with 
Theorem~\ref{THM:1}\,(i)
for the modulated KdV 
\eqref{kdv1} on $\T$ (corresponding to the $\be = 2$ case)
and 
Theorem \ref{THM:4}\,(i)
for the modulated BO \eqref{BO}, the modulated ILW \eqref{ILW1}, 
and the modulated scaled ILW \eqref{ILW2} on $\T$
(corresponding to the $\be = 1$ case).

In  \cite{Robert3}, 
Robert employed the so-called  short-time Fourier restriction norm method
(introduced in~\cite{KT}; see also \cite{CCT2, IKT, OW2}),
which is often used to  weaken resonant interactions, 
and proved local well-posedness of the (real-valued) 
renormalized modulated mKdV \eqref{mkdv2} in $H^s(\T)$
for any $s \in \R$, provided that $\rho = \rho(s) \gg 1$
is sufficiently large, 
thus improving Theorem \ref{THM:3}
when $\rho$ is large.
We point out that the solution map constructed in 
\cite{Robert3} is merely continuous, which is consistent
with the failure of semilinear local well-posedness 
in $H^s(\T)$ for $s < \frac 12$
discussed in Remark \ref{REM:mkdv1}.

As mentioned above, 
 our strategy in this paper is 
 to  study  the YDE \eqref{YDE0}
by constructing 
 $\I^X(\uu)$ as a nonlinear Young integral 
  via the sewing lemma
(without relying on the (refined) Fourier restriction norm method).
 One of the benefits of our approach 
is that we can allow the temporal regularity $\al > 0$ of
a solution $\uu$ to \eqref{YDE0} to be less than $\frac 12$
as long as $\al + \g > 1$;
see Lemma~\ref{LEM:int1}.\footnote{In a recent preprint \cite{GLLO}, 
we implemented a normal form approach to the modulated KdV \eqref{kdv1}
and the modulated BO \eqref{BO}, where we constructed 
a solution with temporal regularity $\al = 0$.}
This in particular allows us to 
study well-posedness
of stochastic modulated dispersive PDEs with additive forcing.
See 
Subsection \ref{SUBSEC:1.6}
and Section \ref{SEC:SK}
for a further discussion.

\smallskip

\noi
(ii)
In a recent preprint
\cite{CGLLO2}, with A.\,Chapouto, 
we studied pathwise well-posedness of stochastic modulated dispersive equations with 
multiplicative noises and established a new regularization-by-noise phenomenon by exploiting the nonlinear interaction between the unknown and the noise. For example, for the stochastic 
modulated KdV with a multiplicative fractional-in-time noise in the Young regime, we showed that irregularity of the modulation induces smoothing on the stochastic convolution in a pathwise manner, where a gain of spatial regularity becomes (arbitrarily) larger for more irregular modulations.

In a recent preprint \cite{GLLO}, 
we implemented a  normal form method to study modulated dispersive PDEs,
which provides an alternative way to give a meaning to the 
integral term in~\eqref{mild2} (and 
$ \I^X(\uu)$ in \eqref{YDE0})
without relying on the sewing lemma (or the Fourier restriction norm method 
as in  \cite{Robert2, Robert3}).
Furthermore, our normal form approach in \cite{GLLO}
does not require any positive temporal regularity
of the modulated interaction representation
$\uu$ defined in \eqref{int1}, allowing us to prove 
unconditional uniqueness of solutions
to modulated dispersive PDEs; see \cite{GLLO}
for further details.

\end{remark}

\subsection{Nonlinear smoothing \& Galerkin approximation}

Our approach to  modulated dispersive equations on the circle
also yields
(i)~nonlinear smoothing 
and (ii)~convergence of the Galerkin approximation.

Given $N \in \N$, 
we consider 
the following truncated modulated KdV on the circle:
\begin{equation}
\label{kdv1x}
\dt u^N+  \dx^3 u^N \cdot \dt w =\dx \P_N\big((\P_Nu^N)^2\big), 
\end{equation}

\noi
where   $\P_N$ denotes  the Dirichlet projector onto the (spatial) frequencies 
$\{|n| \leq N\}$; see \eqref{Diri1}.
Then, a slight modification of the proof of Theorem \ref{THM:1}\,(i)
yields  the following results
for the modulated KdV \eqref{kdv1} on the circle.

\begin{theorem}\label{THM:G1}
Given $\rho \ge\frac 12$,  $\frac12< \g < 1$, and $T> 0$, 
let  $w$ be $(\rho,\g)$-irregular on $[0, T]$ in the sense of Definition~\ref{DEF:ir}.

\smallskip

\noi
\textup{(i) (nonlinear smoothing).}
Suppose that $\rho \ge \frac 12$ and $s\in \R$
satisfy 
 \eqref{reg1}.
In addition, 
suppose that $s_0 > s$ satisfies 
 one of the following conditions\textup{:}
\begin{align}
\begin{split}
\textup{(i.a)} &\ \  
0 \le s + \rho \le  \tfrac 12
 \quad \text{and} \quad 
s_0 < 2s +  3\rho - \tfrac 32, \\
\textup{(i.b)} &\ \  
s + \rho >  \tfrac 12 
\quad  \text{and} \quad s_0 \le s + 2\rho - 1.
\end{split}
\label{reg3}
\end{align}

\noi
Given $u_0 \in H^s(\T)$, 
let $u
\in C([0, \tau]; H^s(\T))$ be the  solution 
to the modulated KdV equation~\eqref{kdv1}
on $\T$
with $u|_{t = 0} = u_0$, 
where 
$\tau \in (0,  T]$ denotes the local existence time.
Then, we have 
\begin{align*}
 u - e^{-w(t) \dx^3}u_0 \in C([0, \tau]; H^{s_0}(\T)).
\end{align*}

\smallskip

\noi
\textup{(ii) (Galerkin approximation).}
Suppose that $\rho> \frac 12$ and $s\in \R$
satisfy one of the following conditions\textup{:}
\begin{align}
\begin{split}
\textup{(ii.a)} &\ \  \tfrac 12 < \rho \le  \tfrac 34 \quad \text{and} \quad 
s > \tfrac 32 - 3 \rho, \\
\textup{(ii.b)} &\ \  \rho > \tfrac 34\quad  \text{and} \quad  s > - \rho.
\end{split}
\label{reg3a}
\end{align}

\noi
Then, given $u_0 \in H^s(\T)$, 
the solution $u^N$ to the truncated modulated KdV equation~\eqref{kdv1x}
on $\T$
with $u^N|_{t = 0} = u_0$
converges to the solution $u$
to the modulated KdV equation \eqref{kdv1}
on $\T$ with $u|_{t = 0} = u_0$
in $C([0, \tau]; H^s(\T))$, 
where $\tau \in (0, T]$ denotes the local existence time of the modulated KdV equation
\eqref{kdv1}.
Moreover, 
given any $r > 0$, 
the rate of convergence of $u^N$ to $u$ is
uniform in $u_0 \in B_r$ 
where
$B_r$ denotes the ball 
in $H^s(\T)$
of radius $r > 0$ centered at the origin.

\end{theorem}

We point out that
our argument does not yield 
nonlinear smoothing or convergence 
of the Galerkin approximation 
for the modulated KdV on the real line;
see Remark \ref{REM:R1}.

In the case of the (unmodulated) KdV \eqref{kdv2} on $\T$, 
nonlinear smoothing 
was shown 
for $s > - \frac 12$ and $s_0 < \min(3s + 1, s + 1)$
by the normal form method
\cite{ETz}.
We prove  Theorem~\ref{THM:G1}\,(i)
by a slight modification of the contraction argument
needed to prove local well-posedness
and thus we do not need to rely on 
 the normal form method.
 Furthermore, 
 the regularity gain on the nonlinear part
$ u - e^{-w(t) \dx^3}u_0$
can be (arbitrarily) large
for irregular modulations,
which is in sharp contrast to the unmodulated case. 
We also point out that 
the uniformity claim in 
Theorem \ref{THM:G1}\,(ii)
fails for 
 the (unmodulated) KdV \eqref{kdv2} on $\T$
 (see \cite[Theorem 1.1]{CKSTT3}), 
 which shows another example of regularization by noise
 for the modulated KdV \eqref{kdv1} on the circle.

\medskip

Next, we state analogous results for the modulated mKdV \eqref{mkdv1} on the circle.

\begin{theorem}\label{THM:G2}
Given $\rho \ge\frac 12$,  $\frac12< \g < 1$, and $T> 0$, 
let  $w$ be $(\rho,\g)$-irregular on $[0, T]$ in the sense of Definition~\ref{DEF:ir}.

\smallskip

\noi
\textup{(i) (nonlinear smoothing).}
Suppose that 
$s_0 > s \ge \frac 12$ satisfies 
\begin{align}
 s_0 \le \min(s + 2 \min(\rho, s) - 1, 3s - 1).
\label{regM2}
\end{align}

\noi
Given $u_0 \in H^s(\T)$, 
let $u
\in C([0, \tau]; H^s(\T))$ be the  solution 
to the modulated mKdV equation~\eqref{mkdv1}
on $\T$
with $u|_{t = 0} = u_0$, 
where 
$\tau \in (0,  T]$ denotes the local existence time.
Then, we have 
\begin{align*}
 u - e^{-w(t) \dx^3 + 3t\|u_0\|_{L^2}^2\dx}u_0 \in C([0, \tau]; H^{s_0}(\T)).
\end{align*}

\smallskip

\noi
\textup{(ii) (Galerkin approximation).}
Let $\rho > \frac 12$ and 
$s > \frac 12$.
Then, given $u_0 \in H^s(\T)$, 
the solution $u^N$ to the truncated modulated mKdV equation 
on $\T$\textup{:}
\begin{equation*}
\begin{cases}
\dt u^N+  \dx^3 u^N \cdot \dt w =\dx \P_N\big((\P_Nu^N)^3\big)\\
u^N|_{t = 0} = u_0
\end{cases}
\end{equation*}

\noi
converges to the solution $u$
to the modulated mKdV equation \eqref{mkdv1}
on $\T$
with $u|_{t = 0} = u_0$
in $C([0, \tau]; H^s(\T))$, 
where $\tau  \in (0, T]$ denotes the local existence time of the modulated mKdV equation
\eqref{mkdv1}.
%
%

\end{theorem}

\begin{remark}\label{REM:mkdv2}\rm
Contrary to Theorem \ref{THM:G1}\,(i) for
the modulated KdV \eqref{kdv1} (and also Theorem~\ref{THM:G3}\,(i) for 
the modulated BO \eqref{BO}), 
the gain of regularity 
in nonlinear smoothing for the modulated mKdV \eqref{mkdv1}
is bounded by $3s - 1$ even if we take $\rho \gg1$.
This is (once again) due to the resonant part of the nonlinearity, where
the modulation plays no role; 
see the proof of Proposition \ref{PROP:mkdv1}.
See also Remark \ref{REM:mkdv1}.

\end{remark}

Lastly, we state the corresponding results
for the modulated BO and the modulated (scaled) ILW
on the circle.

\begin{theorem}\label{THM:G3}
Given $\rho > 1$,  $\frac12< \g < 1$, and $T> 0$, 
let  $w$ be $(\rho,\g)$-irregular on $[0, T]$ in the sense of Definition~\ref{DEF:ir}.

\smallskip

\noi
\textup{(i) (nonlinear smoothing).}
Suppose that $\rho > 1$ and $s\in \R$
satisfy \eqref{regBO1}.
In addition, suppose that $s_0 > s$ satisfies
\begin{align}
\begin{split}
\textup{(i.a)} &\ \
0 \le s < s_0 \le \rho - 1,\\
\textup{(i.b)} &\ \
 s <0
\quad 
\text{and}\quad 
s_0 \le s + \rho - 1.
\end{split}
\label{regBO2}
\end{align}

\noi
Given $u_0 \in H^s(\T)$, 
let $u
\in C([0, \tau]; H^s(\T))$ be the  solution 
to the modulated BO equation~\eqref{BO}
on $\T$
with $u|_{t = 0} = u_0$, 
where 
$\tau \in (0,  T]$ denotes the local existence time.
Then, we have 
\begin{align*}
 u - e^{w(t) \H\dx^2}u_0 \in C([0, \tau]; H^{s_0}(\T)).
\end{align*}

\smallskip

\noi
\textup{(ii) (Galerkin approximation).}
Suppose that 
$ \rho >  1$
and $ 
s > - \tfrac 12\rho$.
Then, given $u_0 \in H^s(\T)$, 
the solution $u^N$ to the truncated modulated BO equation 
on $\T$\textup{:}
\begin{equation*}
\begin{cases}
\dt u^N- \H  \dx^2 u^N \cdot \dt w =\dx \P_N\big((\P_Nu^N)^2\big)\\
u^N|_{t = 0} = u_0
\end{cases}
\end{equation*}

\noi
converges to the solution $u$
to the modulated BO equation \eqref{BO} on $\T$
with $u|_{t = 0} = u_0$
in $C([0, \tau]; H^s(\T))$, 
where $\tau \in (0, T]$ denotes the local existence time of the modulated BO equation~\eqref{BO}.
Moreover, 
given any $r > 0$, 
the rate of convergence of $u^N$ to $u$ is
uniform in $u_0 \in B_r \subset H^s(\T)$.

\smallskip

\noi
\textup{(iii)}
Let $0<  \dl < \infty$.
Then, 
the conclusions in Part (i) and (ii) 
also hold
for the modulated
ILW equation \eqref{ILW1}
and
 the modulated
scaled ILW equation \eqref{ILW2}
on $\T$.

\end{theorem}

\subsection{Invariance of the white noise under the modulated KdV}
\label{SUBSEC:1.5}

In Theorem \ref{THM:1}\,(ii), 
we established global well-posedness
of the periodic modulated KdV \eqref{kdv1}
 via a deterministic tool (the $I$-method).
In this subsection, 
we use a probabilistic tool 
to globalize solutions to 
 the periodic modulated KdV~\eqref{kdv1},
 which allows us to construct global-in-time solutions
 in a regime not covered by 
Theorem \ref{THM:1}\,(ii).

Let $\mu$ be  
a  white noise\footnote{As it is customary in the literature, 
with a slight abuse of notation, 
we use the term `white noise' to refer to both the distribution-valued random variable $u_0^\o$ in \eqref{series1}
and its law $\mu = \Law (u_0^\o)$, when there is no confusion.
Here, $\Law(X)$ denotes the law of a random variable $X$.
For clarity, we may refer to $\mu = \Law (u_0^\o)$ as the white noise measure.} 
on the circle $\T$ with formal density:
\begin{align}
d\mu = Z^{-1} \exp\bigg(- \frac 12 \int_\T u^2 dx_\T\bigg) du.
\label{white1}
\end{align}

\noi
More precisely, $\mu$ is the induced probability
measure
under the map:\footnote{For simplicity of
the presentation, we assume that the white noise has spatial mean zero on $\T$.} 
\begin{align}
\o \in \O \longmapsto
u_0^\o = \sum_{n \in \Z^*} g_n(\o) e_n, 
\label{series1}
\end{align}

\noi
where $\Z^* = \Z \setminus \{0\}$, 
$e_n(x) = e^{ i n  x}$,  and 
$\{g_n \}_{n \in \Z^*}$ is a family of independent standard 
complex-valued Gaussian random variables\footnote{Note that $\text{Var}(g_n) = 1$.} 
on a probability space $(\O, \F, \PP)$
conditioned that $g_{-n} = \cj{g_n}$, $n \in \Z^*$.
It is easy to see that 
$u_0^\o$ in  \eqref{series1} belongs almost surely to $H^s(\T) \setminus H^{-\frac 12} (\T)$, $s < - \frac 12$.

In particular, 
Theorem~\ref{THM:1}\,(i)
states that, for $\rho > \frac 23$, 
the modulated KdV  \eqref{kdv1} on $\T$
is locally well-posed almost surely with respect to the white noise initial data.
On the other hand, 
the  regularity threshold 
in~\eqref{s1} for (deterministic)
global well-posedness tends to $- \frac 14$
as $\rho \to \frac 23$ and $\g \to \frac 12$.
More precisely, 
when $ \rho < 1 - \frac 13 \g$, 
the regularity threshold in \eqref{s1} is above $- \frac 12$
and thus 
we can not apply Theorem \ref{THM:1}\,(ii)
to globalize local-in-time solutions with the  white noise initial data.

By using a probabilistic argument, 
we establish the following global-in-time result
for 
the periodic modulated KdV equation \eqref{kdv1}
with the white noise initial data, 
which complements the deterministic
global well-posedness result in Theorem \ref{THM:1}\,(ii).

\begin{theorem}\label{THM:white1}
Given $\rho > \frac 23 $ and  $\frac12< \g < 1$, 
let  $w$ be $(\rho,\g)$-irregular on $\R_+$ in the sense of Definition~\ref{DEF:ir}.
Then, the modulated KdV equation \eqref{kdv1}
on the circle $\T$
is almost surely globally well-posed
with the white noise initial data\footnote{As it is common
in the study of dispersive PDEs with random initial data, 
we use $u_0^\o$ to denote both the white noise
and a sample path of the white noise but the meaning is clear from the context.}
$u_0^\o$ in \eqref{series1}, 
and the white noise $\mu$ is invariant under the resulting dynamics.

\end{theorem}

Theorem \ref{THM:white1}
extends 
invariance of the white noise under the KdV dynamics
from the unmodulated setting (see \cite{QV, Oh1, Oh2, Oh10, OQV, KMV})
to the current modulated setting, 
thus establishing the first result 
on the construction of invariant measures
for modulated dispersive PDEs.
We also mention 
that 
Theorem \ref{THM:white1}
plays a key role in the globalization
argument for
the stochastic modulated KdV \eqref{skdv1}
with an additive space-time white noise
 (Theorem~\ref{THM:skdv1}\,(ii))
 considered in the next subsection.

In view of 
Theorem \ref{THM:G1} on 
the convergence of the truncated modulated KdV dynamics~\eqref{kdv1x}
to the modulated KdV dynamics \eqref{kdv1}, 
Theorem \ref{THM:white1}
follows from a standard application of 
 Bourgain's invariant measure argument \cite{BO94, BO96}
(whose details we omit; see~\cite{BT1} for 
details of  Bourgain's invariant measure argument), 
once we prove invariance of the white noise
for the truncated modulated KdV~\eqref{kdv1x}; see Proposition \ref{PROP:white2}.

Theorem \ref{THM:4}
states that 
for $\rho > 1$, 
 the modulated BO \eqref{BO}, 
the modulated ILW~\eqref{ILW1}, 
and the modulated scaled ILW \eqref{ILW2}
are locally well-posed 
in $H^s(\T)$, $s > -\frac 12 \rho$, 
containing the support of the white noise.
Moreover, 
for the same range of $\rho$ and $s$, 
Theorem \ref{THM:G3}\,(ii)
guarantees convergence of the Galerkin approximation.
Hence,
a slight modification of a proof of Theorem \ref{THM:white1}
yields
 the following result
on almost sure global well-posedness
of these equations with the white noise initial data
and invariance of the white noise.
We omit details.

\begin{theorem}\label{THM:white2}

Given $\rho > 1 $ and  $\frac12< \g < 1$, 
let  $w$ be $(\rho,\g)$-irregular on $\R_+$ in the sense of Definition~\ref{DEF:ir}.

\smallskip

\noi
\textup{(i)}
The modulated BO equation \eqref{BO}
on the circle $\T$
is almost surely globally well-posed
with the white noise initial data
$u_0^\o$ in \eqref{series1}, 
and the white noise $\mu$ is invariant under the resulting dynamics.
The same result also 
holds
for the modulated ILW equation \eqref{ILW1}
and the modulated scaled ILW equation \eqref{ILW2} on the circle.

\smallskip

\noi
\textup{(ii)}
Let $\{\dl_j \}_{j \in \N}$ be an increasing sequence of positive numbers, tending to $\infty$
as $j \to \infty$.
Let $u^{\dl_j}$, $j \in \N$,  and $u$ denote the  solutions
to the modulated ILW equation \eqref{ILW1}
and the modulated BO equation \eqref{BO}
with the white noise initial data
$u_0^\o$ in \eqref{series1}, respectively.
Then, as $j \to \infty$ \textup{(}and hence $\dl_j \to \infty$\textup{)}, 
$u^{\dl_j}$ converges  to $u$ in 
$C(\R_+; H^{- \frac 12 - \eps}(\T))$, $\eps > 0$, 
almost surely.
Here, 
$C(\R_+; H^{- \frac 12 - \eps}(\T))$
is endowed with the compact-open topology.

\smallskip

\noi
\textup{(iii)}
Let $\{\dl_j \}_{j \in \N}$ be a decreasing sequence of positive numbers, tending to $0$
as $j \to \infty$.
Let $v^{\dl_j}$, $j \in \N$,  and $u$ denote the  solutions
to the modulated scaled ILW equation~\eqref{ILW2}
and the modulated KdV equation \eqref{kdv1}
with the white noise initial data 
$u_0^\o$ in~\eqref{series1}, respectively.
Then, as $j \to \infty$ \textup{(}and hence $\dl_j \to 0$\textup{)}, 
$v^{\dl_j}$ converges  to $u$ in 
$C(\R_+; H^{-\frac 12 - \eps}(\T))$, $\eps > 0$, 
almost surely.

\end{theorem}

\begin{remark}\rm
In Theorem \ref{THM:white2}\,(ii), 
we took a discrete sequence of depth parameters $\dl_j$, $j \in \N$,  
such that, in view of  
Theorem \ref{THM:white2}\,(i), 
there exists a set of full probability, 
guaranteeing global existence of 
$u^{\dl_j}$ for any  $j \in \N$.

In the case of  a continuous parameter $\dl \in (0, \infty)$, 
we have  almost sure asymptotic  global-in-time deep-water convergence
in the following sense.
It follows from Theorem \ref{THM:white2}\,(i)
that there exists  $\Si \subset \O$ with $\PP(\Si) = 1$
such that for each $\o \in \Si$, 
the solution $u$ 
to the modulated BO~\eqref{BO}
with $u|_{t= 0} = u_0^\o$
exists globally in time.
Given any $T \gg 1$, by iterating the local-in-time argument
(as in the proof of Theorem \ref{THM:conv}\,(i)), 
we see that there exists $\dl_ 0 = \dl_0(\o, T) \gg 1$
such that, for any $\dl > \dl_0$,  
the solution $u^{\dl}$ 
to the modulated ILW~\eqref{ILW1}
with $u^{\dl}|_{t= 0} = u_0^\o$
exists on the time interval $[0, T]$
and, moreover, $u^\dl$ converges to $u$ on $[0, T]$ as $\dl \to \infty$.

By making use of 
a bi-parameter Kolmogorov continuity criterion
\cite[Lemma 2.1]{Baldi}, 
it may be possible to prove full deep-water convergence
along a continuous parameter $\dl \in (0, \infty)$.
We, however, do not pursue this issue here.
See \cite{Zine} for such an argument
in the context of 
the Smoluchowski-Kramers approximation for singular stochastic
wave equations.

Similar comments apply to the shallow-water convergence
considered in 
 Theorem \ref{THM:white2}\,(iii).

\end{remark}

\subsection{Stochastic modulated KdV}
\label{SUBSEC:1.6}

As the last example, 
we consider 
the following stochastic modulated KdV equation on the circle
with an additive white-in-time noise:
\begin{equation}
\label{skdv1}
\dt u+  \dx^3 u \cdot \dt w =\dx u^2 + \phi \xi, 
\end{equation}

\noi
where $\xi$ denotes a (Gaussian) space-time white noise 
on $\R_+\times \T$
whose space-time covariance is (formally) given by 
\begin{align*}
 \E[ \xi(t_1, x_1)\xi(t_2, x_2) ] = \dl (t_1 - t_2) \dl(x_1 - x_2) 
\end{align*} 

\noi
for $t_1, t_2 \in \R_+$  and $x_1, x_2 \in \T$ 
with $\dl$ denoting the Dirac delta function.
Here, 
$\phi$ is a Hilbert-Schmidt operator from $L^2_0(\T) = \P_{\ne 0}L^2(\T)$ to 
$H_0^s(\T)= \P_{\ne 0}H^s(\T)$
for some $s \in \R$, 
where $\P_{\ne 0}$ denotes the projection onto non-zero frequencies.
Namely, 
$L^2_0(\T)$
and $H_0^s(\T)$
are subspaces of $L^2(\T)$ and $H^s(\T)$, 
consisting of mean-zero functions.

Given $\s \ge 0$,\footnote{By convention, we have $X\equiv 0$ when $\s = 0$.
Namely, $\mu_0 = \dl_0$, where $\dl_0$ is the Dirac delta distribution at the trivial function.} we say that a distribution-valued random variable $X_\s$ on $\T$ (and its law, denoted by $\mu_\s$)
is a (spatial) white noise on~$\T$ with variance~$\s$ if
\begin{align}
 \mu_{\s} = \Law(X_\s) = \Law (\sqrt \s\, u_0^\o),
 \label{mu1}
\end{align}

\noi
where $u_0^\o$ is the white noise (with variance 1) in \eqref{series1}.

Let us recall the notion of 
 an {\it evolution system of measures} \cite{DR, DD}
 in a somewhat formal manner.
Let $\Phi_{t_1, t_2}= \Phi_{t_1, t_2}^\o$, $t_1 \ge t_2 \ge 0$,  
be a solution map for a given  (random) dynamical system, 
sending the data $\varphi$ at time $t_2$ to the solution $\Phi_{t_1, t_2}\varphi$ at time $t_1$.
Then, we define the transition semigroup $P_{t_1, t_2}$
by 
\begin{align}
 P_{t_1, t_2}F(\varphi) = \E[ F(\Phi_{t_1, t_2}^\o \varphi)]
\label{trans1}
\end{align} 

\noi
for a bounded measurable function $F$ on the phase space $\mathcal{X}$.
Then, we say that\footnote{Strictly speaking, an evolution system of measures
is the mapping $t \in \R_+ \mapsto \rho_t \in \mathcal{P}(\M)$, 
where $\mathcal{P}(\M)$ denotes the family of probability measures on $\M$.
However, we simply refer to the family $\{\rho_t\}_{t \in \R_+}$ of measures
as an evolution system of measures.} 
a family $\{\rho_t\}_{t \in \R_+}$ of probability measures on $\mathcal{X}$
is an evolution system of measures indexed by $\R_+$ if
\begin{align}
\int_\mathcal{X}  F(\varphi) \rho_{t_1}(d\varphi)= 
 \int_\mathcal{X} P_{t_1, t_2} F(\varphi) \rho_{t_2}(d\varphi)
 \label{evo1}
\end{align}

\noi
for any bounded continuous function $F$ on $\mathcal{X}$ and $t_1 \ge t_2 \ge 0$.
Note that \eqref{evo1} is equivalent to 
\begin{align}
 \rho_{t_1} = P^*_{t_1, t_2} \rho_{t_2} 
 \label{evo2}
\end{align}

\noi
for any $t_1 \ge t_2 \ge 0$.
If there exists an invariant measure $\rho$, 
then by setting $\rho_t = \rho$, $t \in \R_+$, 
the family $\{\rho_t\}_{t\in \R_+}$ is obviously
an evolution system of measures.
It is in this sense that the notion of an evolution system
of measures is a generalization of the notion of an invariant measure.

We now state our main result on the stochastic modulated KdV \eqref{skdv1}.

\begin{theorem}\label{THM:skdv1}
\noi
\textup{(i)}
Given $\rho \ge\frac 12$,  $\frac12< \g < 1$, and $T> 0$, 
let  $w$ be $(\rho,\g)$-irregular on $[0, T]$ in the sense of Definition~\ref{DEF:ir}.
Suppose that  $s \in \R$ satisfies~\eqref{reg1}
and that $\phi \in \HS(L^2_0(\T); H^s_0(\T))$.
Then, 
 the stochastic modulated KdV equation~\eqref{skdv1}
is locally well-posed in $H^s(\T)$
and the modulated interaction representation $\uu(t) = \uw(t)^{-1}u(t)$
of the solution $u$
belongs to 
$\cC^\al([0, \tau]; H^s(\T))$
for any $ 0 < \al < \frac 12$, 
where $\tau > 0$ denotes the almost surely positive local existence time.

In particular, when $\rho > \frac 23$, 
 the stochastic modulated KdV equation~\eqref{skdv1}
 with an additive space-time white noise \textup{(}namely, $\phi = \Id$\textup{)}
is locally well-posed in $H^s(\T)$.

\smallskip

\noi
\textup{(ii)}
Given $\rho > \frac 23$ and   $\frac12< \g < 1$, 
let  $w$ be $(\rho,\g)$-irregular on $\R_+$.
Then, 
the stochastic modulated KdV equation \eqref{skdv1}
on $\T$
 with an additive space-time white noise
 \textup{(}namely, $\phi = \Id$\textup{)}
 is globally well-posed with 
 the  \textup{(}spatial\textup{)}
 white noise initial data $u_0^\o$ in~\eqref{series1}, 
 independent of the space-time white noise forcing  $\xi$.
   Moreover, for any $t\ge 0$, 
 we have
\begin{align}
\Law (u(t)) = \mu_{1+t}.
\label{th1}
\end{align}

\noi
Namely,  the family  $\{ \mu_{1+t}\}_{t\in \R_+}$
is an evolution system of measures
for the stochastic modulated KdV equation  \eqref{skdv1} with the white noise initial data $u_0^\o$ in~\eqref{series1}.

\end{theorem}

When $\rho > \frac 12$, 
convergence of the Galerkin approximation (see \eqref{KdV9}) holds
in the current stochastic setting;
see Lemma~\ref{LEM:Ga1}.

The well-posedness issue of the stochastic KdV with an additive forcing:
\begin{align}
\dt u + \dx^3 u + u \dx u = \phi \xi
\label{KdV2}
\end{align}

\noi
 has 
been studied   on both the real line and on the torus;
see \cite{DDT2, DDT1, Oh2, OQS}.
In particular, in~\cite{OQS}, 
Quastel, Sosoe, and the fifth author
proved an analogue of Theorem \ref{THM:skdv1}\,(ii)
for~\eqref{KdV2}
with $\phi = \Id$
by implementing 
 a variant of Bourgain's invariant measure argument~\cite{BO94}
in the context of an evolution system of measures.
The local well-posedness claim 
in Theorem~\ref{THM:skdv1}\,(i)
follows from a slight modification 
of the proof of Theorem \ref{THM:1}\,(i).
On the other hand, 
Theorem~\ref{THM:skdv1}\,(ii)
follows from 
a straightforward adaptation of  the argument developed in \cite{OQS} to the current problem.
See also \cite[Corollary 1.3\,(ii)]{OQS}
for a global existence result of \eqref{KdV2}
with initial data of the form 
$w_0 + u_0^\o$ for $w_0 \in L^2(\T)$ based on the ideas in~\cite{OQ}, 
which also applies to the stochastic modulated KdV \eqref{skdv1}
with $\phi = \Id$.

The Duhamel formulation of  \eqref{skdv1} is given by 
\begin{align*}
u(t) & = U^w(t) u_0 +  U^w({t}) \int_0^t  (U^w({t'}))^{-1}\NN( u(t') ) dt'\\
& \quad + U^w({t}) \int_0^t  (U^w({t'}))^{-1} \phi \xi( dt'), 
\end{align*}

\noi
where 
$\uw (t)=e^{-   w(t)\dx^3}  $, 
$\NN(u) = \dx  u^2 $, 
and 
 the last term on the right-hand side 
represents
the effect of the stochastic forcing $\phi\xi$ in \eqref{skdv1}.
As in the deterministic case, 
we study the following Duhamel formulation satisfied by 
the modulated interaction representation $\uu(t) = \uw(t)^{-1}u(t)$:
\begin{equation}
\uu(t) = u_0 + \int_0 ^t \uw(t')^{-1}\NN(\uw(t') \uu(t'))dt'
+ \Psi(t), 
\label{mild4}
\end{equation}

\noi
where $\Psi$ is 
a solution to the following stochastic equation:
\begin{equation*}
\dt \Psi = \uw(t)^{-1} \phi \xi, \qquad \Psi(0) = 0,
\end{equation*}

\noi
formally given by the following stochastic integral:
\begin{align}
 \Psi(t) =  \int_0^t  (U^w({t'}))^{-1} \phi \xi( dt').
\label{sconv1}
\end{align}

\noi
See 
\eqref{sconv2} for a precise definition of $\Psi$.
Note that, due to the roughness-in-time of the noise, 
the stochastic term $\Psi$
has temporal regularity $\al < \frac 12$
(see 
Lemma \ref{LEM:sto1})
and thus we expect that a solution $\uu$
to \eqref{mild4} 
has temporal regularity $\al < \frac 12$.
We point out that 
 our nonlinear Young integral theory
(Proposition \ref{PROP:young1})
and a general local well-posedness
result 
for the YDE~\eqref{YDE0}
are adapted to the case 
with temporal regularity $\al < \frac 12$
(as long as $\al + \g > 1$; see the proof of Proposition \ref{PROP:main}).
Hence, once we establish regularity properties
of the stochastic term $\Psi$ (Lemma \ref{LEM:sto1}), 
Theorem \ref{THM:skdv1}\,(i) follows
from a minor modification of the proof of 
Proposition~\ref{PROP:main}.
See Section \ref{SEC:SK}
for details.

\medskip

Next, we briefly discuss a proof of 
the almost sure global well-posedness result in Theorem~\ref{THM:skdv1}\,(ii), 
based on the argument in \cite{OQS}.
We first 
 provide a heuristic argument for the claim \eqref{th1}
 by viewing the stochastic modulated KdV \eqref{skdv1}
 as a
superposition of the deterministic modulated KdV~\eqref{kdv1}
and the stochastic flow:
\begin{align}
\dt  u = \xi
\label{KdV3a}
\end{align}

\noi
(at the level of infinitesimal generators).
While 
Theorem \ref{THM:white1}
states that 
the white noise (with any variance) is invariant under \eqref{kdv1}, 
it is easy to see that 
the stochastic flow \eqref{KdV3a} with a white noise initial data (with any variance)
increases the variance by the length of the time interval under consideration.
Then, at a formal level, 
the claim \eqref{th1} follows
from these observations 
 together with 
the 
Lie-Trotter product formula \cite[Section VIII.8]{RS}:
\begin{align}
e^{t (A+ B)} = \lim_{n \to \infty} \big[ e^{\frac tn A}e^{\frac tn B}\big]^n
\label{tro}
\end{align}

\noi
(which holds, for example, 
for finite-dimensional matrices $A, B$).
Obviously, 
the Lie-Trotter product formula \eqref{tro}
is not directly applicable to our problem, 
and we justify the heuristic argument above
by  the Galerkin approximation (Lemma \ref{LEM:Ga1}).
See Section \ref{SEC:SK}
for further details.

As in Subsection \ref{SUBSEC:1.5}, 
by a straightforward modification
of the proof of Theorem \ref{THM:skdv1}, 
we obtain
 the corresponding  results
for the following  stochastic modulated BO equation on the circle:
\begin{equation}
\label{SBO1}
\dt u-\H  \dx^2 u \cdot \dt w =\dx u^2 + \phi \xi.
\end{equation}

\begin{theorem}\label{THM:skdv2}
\noi
\textup{(i)}
Given $\rho \ge 1$,  $\frac12< \g < 1$, and $T> 0$, 
let  $w$ be $(\rho,\g)$-irregular on $[0, T]$ in the sense of Definition~\ref{DEF:ir}.
Suppose that  $s \in \R$ satisfies~\eqref{regBO1}
and that $\phi \in \HS(L^2_0(\T); H^s_0(\T))$.
Then, 
 the stochastic modulated BO equation~\eqref{SBO1}
is locally well-posed in $H^s(\T)$.
In particular, when $\rho > 1$, 
 the stochastic modulated BO equation~\eqref{SBO1}
 with an additive space-time white noise \textup{(}namely, $\phi = \Id$\textup{)}
is locally well-posed in $H^s(\T)$.

\smallskip

\noi
\textup{(ii)}
Given $\rho > 1$ and   $\frac12< \g< 1$, 
let  $w$ be $(\rho,\g)$-irregular on $\R_+$.
Then, 
the stochastic modulated BO equation \eqref{SBO1}
on $\T$
 with an additive space-time white noise
 \textup{(}namely, $\phi = \Id$\textup{)}
 is globally well-posed with 
 the  \textup{(}spatial\textup{)}
 white noise initial data $u_0^\o$ in~\eqref{series1}, 
 independent of the space-time white noise forcing  $\xi$.
   Moreover,   the family  $\{ \mu_{1+t}\}_{t\in \R_+}$
is an evolution system of measures
for the stochastic modulated BO equation~\eqref{SBO1} with the white noise initial data $u_0^\o$ in~\eqref{series1}.

\end{theorem}

The same results
also hold for the 
 stochastic modulated ILW equation on the circle:
\begin{equation}
\dt u-\Gdl  \dx^2 u \cdot \dt w =\dx u^2 + \phi \xi
\label{SILW1}
\end{equation}

\noi
and 
the 
 stochastic modulated scaled ILW equation on the circle:
\begin{equation}
\dt v-    \wt \Gdl\dx^2 v \cdot \dt w =\dx v^2
+ \phi \xi.
\label{SILW2}
\end{equation}

\noi
In particular, 
deep-water and shallow-water convergence
analogous to those in 
Theorems~\ref{THM:conv}
and \ref{THM:white2}
also holds for \eqref{SILW1} and \eqref{SILW2}.
We omit details.

\subsection{Organization}

In Section
\ref{SEC:2}, 
we introduce basic notations and function spaces.
We also introduce various classes of drivers
such as $\cX^{s, \g}_k ([0, T]\times \M)$ in \eqref{X1}
and \eqref{X1z}
which play a crucial role in the remaining part of this paper.
In Section~\ref{SEC:Young}, 
after providing simple criteria
for checking regularity of drivers (Lemma \ref{LEM:OBS1}), 
we go over the basic theory
for nonlinear Young integrals
(Lemma \ref{LEM:int1} and Proposition \ref{PROP:young1}).
We then present a general local well-posedness
result for the YDE \eqref{YDE0}
(Proposition \ref{PROP:main}).
In Section~\ref{SEC:LWP}, 
by establishing regularity properties
of the associated drivers, 
we prove local well-posedness and related properties
for the modulated KdV-type equations.
In Section~\ref{SEC:BO}, 
we establish local well-posedness and related properties
for the modulated BO and ILW  equations.
In Section~\ref{SEC:GWP1}, 
we prove the conservation of the $L^2$-norm
for the modulated dispersive equations considered in this paper
along with their $L^2$-global well-posedness
in the regime where they are locally well-posed in $L^2(\M)$.
In Section~\ref{SEC:I}, we adapt the $I$-method
to the current modulated setting.
By combining a crucial  commutator estimate
(Proposition \ref{PROP:com})
with the sewing lemma (Lemma~\ref{LEM:sew}), 
we prove global well-posedness of the modulated
KdV \eqref{kdv1} in negative Sobolev spaces
(Theorems~\ref{THM:1}\,(ii) and \ref{THM:2}\,(ii)).
In Section~\ref{SEC:WN}, we provide a brief discussion 
for
the proof of Theorem~\ref{THM:white1}
on invariance of the white noise under the modulated KdV on the circle.
Finally, 
in Section~\ref{SEC:SK}, 
we study the stochastic modulated KdV
\eqref{skdv1} on the circle and prove Theorem \ref{THM:skdv1}.

\section{Notations and preliminary tools}
\label{SEC:2}

\subsection{Basic notations}
Let $A\les B$ denote an estimate of the form $A\leq CB$ for some constant $C>0$. We write $A\sim B$ if $A\les B$ and $B\les A$, while $A\ll B$ denotes $A\leq c B$ for some small constant $c> 0$. 
We may write  $\les_{\al}$ and $\sim_{\al}$ to 
emphasize the dependence on an external parameter $\al$.
We use $C>0$ to denote various constants, which may vary line by line.

In expressing the dependence of a function $u$
on the time variable, we often use the short-hand notation
$u_t = u(t)$,  which is standard in probability theory and stochastic analysis.

Let $\M = \R$ or $\T$.
We  use $\ft f$ and $\F(f) = \F_\M(f)$
to denote
the  Fourier transform
of a function $f$ on $\M$, 
defined by 
\begin{align*}
\ft  f(n)=\int_{\M}f(x)e^{-i n x} d x_\M.
\end{align*}

\noi
Here, when $\M= \R$, $dx_\R =(2\pi)^{-1} dx$ denotes
the scaled Lebesgue measure on $\R$, 
while, when $\M = \T$, 
$dx_\T$ denotes
 the normalized Lebesgue measure $ dx_\T =  (2\pi)^{-1}dx$.
 With this convention, we have 
\begin{align*}
f(x)
& = \int_{\M^*} \ft f( \xi) e^{ix \xi} d\xi_{\M^*}, \\
\| f\|_{L^2(\M)}^2 
& = \int_{\M}|f(x)|^2 d x_\M
= \int_{\M^*} |\ft f( \xi)|^2 d\xi_{\M^*}, \\
\F_\M(fg)(\xi)
& =
\intt_{\xi = \xi_1 + \xi_2}
 \ft  f(\xi_1)\ft  g(\xi_2) (d\xi_1)_{\M^*}, 
\end{align*}

\noi
where $\M^*$
is the Pontryagin dual of $\M$
(namely, $\M^*=\R$ when $\M = \R$, while  $\M^* = \Z$ when $\M = \T$)
and $d\xi_\R$
 denotes the standard  Lebesgue measure on $\R$, 
while $d\xi_{\Z}$ denotes the counting measure on $\Z$.
On the circle, we impose the mean-zero assumption
for the equations we consider in this paper
and thus we set  $\Z^* = \Z \setminus \{0\}$.

Given $N \in \N$, 
we denote by  $\P_N$
 the Dirichlet projector onto the frequencies 
$\{|\xi| \leq N\}$
defined by 
\begin{align}
 \ft{\P_N f}(\xi)= \ind_{\{|\xi|\leq N\}} \ft f(\xi).
 \label{Diri1}
\end{align}

\noi
We set $\P_N^\perp = \Id - \P_N$.

\subsection{Function spaces}

We use 
$\S(\R)$ to denote the space of Schwartz functions on $\R$
and set
$\S(\T) = C^\infty(\T)$.

Let $\M = \R$ or $\T$.
Given $s \in \R$, 
we define the non-homogeneous and homogeneous Sobolev spaces $H^{s}(\M)$ 
and $\dot H^{s}(\M)$ 
via the norms: 
\begin{align*}
\|  f  \|^2_{H^{s}(\M)}
& =
 \int_{\M^*} \jb{\xi}^{2s} |\ft f( \xi)|^2 d\xi, \\
\|  f  \|^2_{\dot H^{s}(\M)}
&  =
 \int_{\M^*} |\xi|^{2s} |\ft f( \xi)|^2 d\xi.
\end{align*}

Given $k \in \N$, 
let $V, V_1, \dots, V_k$  be separable Hilbert spaces.
We use 
\begin{align*}
\cL_k\Big(\bigotimes_{j = 1}^k V_j; V\Big)
\end{align*}

\noi
 to denote 
the Banach space of bounded $k$-linear operators 
on $\bigotimes_{j = 1}^k V_j$ 
with values in $V$.
When $V_j =  V$ for $j =1, \dots,k$,
we simply set  
 $\cL_k(V)=\cL_k(V^{\otimes k}; V)$.

Let $V$ be a Banach space and $T>0$.
For $n\in\N$, we denote 
\begin{align*}
\Delta_{n, T} = 
\big\{ (t_1, \ldots, t_n) \in [0,T]^n: \ t_i > t_j 
\text{ for } i < j\big\}.
\end{align*}
We denote by $C_{n,T}V$ 
the space of continuous functions 
from $\Delta_{n,T}$ to $V$. When $n=1$,
we may write $C_T V$ for simplicity, 
and equip this space with the supremum norm: 
\begin{align*}
\|f\|_{C_T V} = \|f\|_{L^\infty_T V} = \sup_{0\leq t \leq T} \|f(t)\|_V.
\end{align*}

\noi
We define the coboundary operator 
$\updl: C_{n,T} V  \to C_{n+1,T} V$ 
as follows; 
given  $f\in C_{n,T} V$ 
and $(t_1,\ldots, t_{n+1}) \in \Dl_{{n+1},T}$, 
we set
\begin{align*}
(\updl f)_{t_1,\ldots , t_{n+1}} = 
\sum_{k=1}^{n+1} (-1)^{k} f_{t_1, \ldots, t_{k-1}, t_{k+1}, \ldots, t_{n+1}}.
\end{align*}

\noi
For example, for $f\in C_T V$
and 
$g\in C_{2,T} V$, 
we have 
\begin{align}
\begin{split}
(\updl f )_{t,r} &= f_t - f_r, \\
(\updl g)_{t_1,t_2,t_3} &= 
g_{t_1,t_3} - g_{t_1,t_2} - g_{t_2,t_3}
\end{split}
\label{dl1}
\end{align}

\noi
for $(t,r)\in\Dl_{2,T}$
and $(t_1,t_2,t_3) \in \Dl_{3,T}$.
As noted in \cite{GT10}, 
the sequence 
\begin{align*}
0 \too V \too C_{1, T}V \stackrel{\updl}{\too} C_{2, T}V
\stackrel{\updl}{\too} C_{3, T}V
\stackrel{\updl}{\too} \cdots
\end{align*}

\noi
is exact. 
In particular, we have 
 $\updl\circ\updl =0$ and 
if $f \in C_{n,T} V$ with $\updl f =0$, 
then there exists a $g\in C_{n-1,T}V$ 
such that $f= \updl g$; 
see, for example,  \cite[Lemma 2.1]{GT10}.

Given $0 < \g < 1$, we denote by $C^\g_T V = C^\g([0, T]; V)$ the space of $\g$-H\"older continuous functions taking values in $V$, endowed  with the seminorm:
\begin{align*}
\|f\|_{C^\g_T V} = \sup_{(t,r)\in \Dl_{2,T}} 
\frac{\|(\updl f)_{t,r}\|_V}{|t-r|^\g}.
\end{align*}

\noi
We also define 
 $\CC^\g_T V = \CC^\g([0, T]; V)$ via the norm:
\begin{align}
\| f \|_{\CC^\g_TV} = \| f \|_{L^\infty_TV} + \|f\|_{C^\g_T V}.
\label{Ho2a}
\end{align}

\noi 
We also introduce the spaces 
$C^\g_{n,T}V$, $n=2,3$, 
equipped with the following H\"older-type norms;
 for $g\in C_{2,T}V$ and $h\in C_{3,T}V$, we set
\begin{align}
\begin{split}
\| g\|_{C^\g_{2,T} V} & 
= \sup_{(t,r) \in \Dl_{2,T}} \frac{\|g_{t,r} \|_{V} }{|t-r|^{\g}}, \\
\| h\|_{C^\g_{3,T} V} & 
= \inf\Bigg\{
\sum_{n = 1}^N \sup_{(t_1,t_2,t_3) \in \Dl_{3,T}} 
\frac{ \|h^{(n)}_{t_1,t_2,t_3}\|_V}
{|t_1-t_2|^{\al_n} |t_2-t_3|^{\g-\al_n}}\Bigg\}, 
\end{split}
\label{Ho2}
\end{align}

\noi
where the infimum is taken over all finite representations
$h = \sum_{n= 1}^N h^{(n)}$,   $h^{(n)} \in C_{3,T}V$, of $h$
and over all choices of $\al_n \in (0, \g)$.

\medskip

Let $V$ and $W$ be Banach spaces.
Given $k \in \N$, 
we  use $\Lip_k(V; W)$ 
to denote the  space of 
locally Lipschitz maps $f:V\to W$
with polynomial growth of order $k$
such that 
\begin{align}
\|f\|_{\Lip_k(V; W)} = 
\sup_{x,y\in V} 
\frac{\|f(x)-f(y)\|_W}{\|x-y\|_V \big(1+\|x\|_V+\|y\|_V\big)^{k-1}}
<\infty .
\label{Lip1}
\end{align}

\noi
When $V = W$, 
we simply set $\Lip_k(V)= \Lip_k(V;V)$.
Given an integer $k \ge 2$, 
we say that $f\in\Lip^2_k(V; W)$ 
if 
\smallskip

\begin{itemize}
\item[(i)] $f\in\Lip_k(V; W)$,

\smallskip
\item[(ii)]
 $f$ is Fr\'echet differentiable 
with 
$Df \in \Lip_{k-1}(V;\cL_1(V; W))$.

\end{itemize}

\smallskip

\noi
From \eqref{Ho2} and \eqref{Lip1}, we have 
\begin{align}
\|f \|_{C^\g_{2, T}\Lip_k(V;W)}
= \sup_{(t,r) \in \Dl_{2,T}}
\frac 1{|t-r|^{\g}}
\sup_{x,y\in V} 
\frac{\| f_{t, r}(x)- f_{t, r}(y)\|_W}{\|x-y\|_V \big(1+\|x\|_V+\|y\|_V\big)^{k-1}}.
\label{Ho3}
\end{align}

\noi
In addition, suppose that $V_0 \hookrightarrow V$
is a Banach subspace of $V$.
Given an integer $k \ge 2$, 
we  use $\Lip_k(V, V_0; V_0)\subset  \Lip_k(V)$ 
to denote the space of maps $f\in\Lip_k(V)$
such that $f(V_0)\subset V_0$ and
\begin{align}
\|f\|_{\Lip_k(V, V_0; V_0)}
= 
\sup_{x,y\in V_0} 
\frac{\|f(x)-f(y)\|_{V_0}}
{G_{V, V_0}(x, y)}
<\infty , 
\label{Lip2}
\end{align}

\noi
where $G_{V, V_0}(x, y)$ is given by 
\begin{align*}
G_{V, V_0}(x, y)
& = \|x-y\|_{V_0} \big(1+\|x\|_{V}+\|y\|_{V}\big)^{k-1}\\
& \quad + \|x-y\|_{V} \big(1+\|x\|_{V}+\|y\|_{V}\big)^{k-2}
 \big(1+\|x\|_{V_0}+\|y\|_{V_0}\big).
\end{align*}

\noi
Then, 
from \eqref{Ho2} and \eqref{Lip2}, we have 
\begin{align}
\|f \|_{C^\g_{2, T}\Lip_k(V, V_0; V_0)}
= \sup_{(t,r) \in \Dl_{2,T}}
\frac 1{|t-r|^{\g}}
\sup_{x,y\in V_0} 
\frac{\| f_{t, r}(x)- f_{t, r}(y)\|_{V_0}}{
G_{V, V_0}(x, y)}.
\label{Ho4}
\end{align}

\medskip

Given $s \in \R$, $0 < \g < 1$,  $k \in \N$, and $T > 0$, 
we define the space 
$\cX^{s, \g}_k([0, T]\times\M)$ of drivers 
as follows; for $k \ge 2$, we set
\begin{align}
\begin{split}
& \cX^{s, \g}_k ([0, T]\times \M) \\
& \quad = 
\big\{X \in C^\g_{2, T} \Lip_k^2(H^s(\M)):
X(0) = DX[0]= 0 \text{ and }\updl X = 0\big\}, 
\end{split}
\label{X1}
\end{align}

\noi
endowed with the norm
\begin{align}
\|X\|_{\cX^{s, \g}_k([0, T]\times \M)} 
=  \|X\|_{C^\g_{2, T}\Lip_k(H^s(\M))} + \|DX\|_{C^\g_{2, T}\Lip_{k-1}(H^s(\M);\cL_1(H^s(\M)))}.
\label{X2}
\end{align}

\noi
When $k = 1$, we set
\begin{align}
\begin{split}
 \cX^{s, \g}_1 ([0, T]\times \M) 
& =  \big\{ X \in 
C^\g_{2, T} \L(H^s(\M)):
\updl X = 0\big\}.
\end{split}
\label{X1z}
\end{align}

\noi
In \eqref{X1}, $X(0) = 0$ (and $DX[0]= 0$) means $X_{t, r}(0) = 0$ 
(and $DX_{t, r}[0]= 0$, respectively)
for any $(t, r) \in \Dl_{2, T}$, 
where $DX_{t, r}[0]$
denotes the Fr\'echet derivative of $X_{t, r}$ at $u = 0 \in H^s(\M)$.
Similarly, we define 
$\dot \cX^{s, \g}_k([0, T]\times\M)$ by setting
\begin{align}
\begin{split}
& \dot \cX^{s, \g}_k ([0, T]\times \M) \\
& \quad = 
\big\{X \in C^\g_{2, T} \Lip_k^2(\dot H^s(\M)):
X(0) = DX[0]= 0 \text{ and }\updl X = 0\big\}, 
\end{split}
\label{X1a}
\end{align}

\noi
endowed with the norm
\begin{align*}
\|X\|_{\dot \cX^{s, \g}_k([0, T]\times \M)} 
=  \|X\|_{C^\g_{2, T}\Lip_k(\dot H^s(\M))} + \|DX\|_{C^\g_{2, T}
\Lip_{k-1}(\dot H^s(\M);\cL_1(\dot H^s(\M)))}
\end{align*}

\noi
for  $k \ge 2$.
For  $k = 1$, 
 we set 
$ \dot \cX^{s, \g}_1 ([0, T]\times \M) 
=  C^\g_{2, T} \Lip_1(\dot H^s(\M))$.
We may simply write $\cX^{s, \g}_k ([0, T])$ or $\cX^{s, \g}_k(T)$, 
where there is no confusion about 
an underlying spatial domain.
We write 
$X \in \cX^{s, \g}_k (\R_+)
= \cX^{s, \g}_k (\R_+\times \M)$ 
if  $X \in \cX^{s, \g}_k(T)$ for any $T > 0$
(but $\sup_{T > 0 } \| X\|_{\cX^{s, \g}_k(T)}$ may be infinite).
A similar comment applies to the homogenous space 
$\dot \cX^{s, \g}_k ([0, T]\times \M)$.

More generally, 
given $s, s_0 \in \R$, $0 < \g < 1$,  $ k \in \N$, and $T > 0$
with $s_0 > s$, we define the space
$\cX^{s, s_0, \g}_{k}([0, T]\times\M)$ 
of drivers (for establishing nonlinear smoothing) by 
setting 
\begin{align}
& \cX^{s, s_0, \g}_k ([0, T]\times \M) 
= C^\g_{2, T} \Lip_k(H^s(\M); H^{s_0}(\M)).
\label{X1x}
\end{align}

\noi
Similarly, 
we set 
\begin{align*}
& \dot \cX^{s, s_0, \g}_k ([0, T]\times \M) 
= C^\g_{2, T} \Lip_k(\dot H^s(\M); \dot H^{s_0}(\M)).
\end{align*}

In proving persistence of regularity, we need the following
class of drivers.
Given $s, s_0 \in \R$, $0 < \g < 1$,  
an integer $k \ge 2$,
 and $T > 0$
with $s_0 > s$, we define 
$\cY^{s, s_0, \g}_{k}([0, T]\times\M)$ by 
\begin{align}
\begin{split}
& \cY^{s, s_0, \g}_k ([0, T]\times \M) = 
  C^\g_{2, T} \Lip_k(H^s(\M), H^{s_0}(\M);H^{s_0}(\M)).
\end{split}
\label{X2c}
\end{align}

\noi
Similarly, we set 
\begin{align}
\begin{split}
& \dot \cY^{s, s_0, \g}_k ([0, T]\times \M) = 
  C^\g_{2, T} \Lip_k(\dot H^s(\M), \dot H^{s_0}(\M);\dot H^{s_0}(\M)).
\end{split}
\label{X2d}
\end{align}

\medskip

We often use short-hand notations such as
$C^\g_T H^s_x  = C^\g\big([0, T]; H^s(\M))$, 
 when there is no ambiguity.

\subsection{On continuous and discrete convolutions}

We recall the following basic lemma on continuous and discrete convolutions.

\begin{lemma}\label{LEM:SUM}
Let  $\al, \be \in \R$ satisfy
\begin{align}
 \al \ge \be \ge 0 \qquad \text{and}\qquad  \quad \al+ \be > 1.
 \label{SUM1}
\end{align}

\noi
Then, we have
\begin{align*}
\int_{\xi = \xi_1 + \xi_2}\frac{d\xi_1}{\jb{\xi_1}^\al \jb{\xi_2}^\be}
& \les \frac 1{\jb{\xi}^{ \be - \ld}}, \\
 \sum_{n = n_1 + n_2} \frac{1}{\jb{n_1}^\al \jb{n_2}^\be}
& \les \frac 1{\jb{n}^{ \be - \ld}}
\end{align*}

\noi
for any $\xi \in \R$ and $n \in \Z$, 
where $\ld = 
\max( 1- \al, 0)$ when $\al\ne 1$ and $\ld = \eps$ when $\al = 1$ for any $\eps > 0$.
In particular, if \eqref{SUM1} holds, then we have 
\[
\sup_{\xi  \in \R} \int_{\xi = \xi_1 + \xi_2}\frac{d\xi_1}{\jb{\xi_1}^\al \jb{\xi_2}^\be}
+ 
\sup_{n \in \Z} \sum_{n = n_1 + n_2} \frac{1}{\jb{n_1}^\al \jb{n_2}^\be}
< \infty. \]

\end{lemma}

Lemma \ref{LEM:SUM} follows
from elementary  computations.
See, for example,  
 \cite[Lemma 4.2]{GTV} and \cite[Lemma 4.1]{MWX}.

\section{Young solution theory}
\label{SEC:Young}

Our main goal in this paper is to study well-posedness of 
modulated dispersive equations
in terms of the Duhamel formulation \eqref{mild2} satisfied
by the modulated interaction representation
$\uu(t)=\uw(t)^{-1}u(t)$ defined in \eqref{int1}.
Here,  the main task is to make sense of 
 the integral 
on the right-hand side of \eqref{mild2}.
In this section, we go over the construction of
such an integral as a  nonlinear Young integral
under appropriate assumptions.

Let $p$ be a function on $\R$, smooth on $\R \setminus\{0\}$, 
such that  
\begin{align}
p(0) = 0\qquad 
\text{and}\qquad 
\cj{p(i\xi)} = - p(i\xi), \quad \xi \in \R\setminus \{0\}.
\label{X0a}
\end{align}

\noi
Namely, $p(\dx)$ is formally anti self-adjoint on $L^2(\M)$
(which may be unbounded).
Then, our goal is to give a proper meaning to 
the following integral:
\begin{align}
\label{X0b}
\I^X(\uu)(t) = 
 \int_0^t \uw(t')^{-1}\NN(\uw(t') \uu(t'))dt' , 
\end{align}

\noi
where $U^w(t) = e^{ w(t)p(\dx)}$.
For this purpose, we
introduce a driver $X$, a priori defined on $\S(\M)$,  associated with the integral $\I^X(\uu)$
by setting 
\begin{equation}
X_{t,r}(f) = \int_r^t \uw(t')^{-1} \NN(\uw(t') f) dt', \quad t >  r \ge 0.
\label{X0}
\end{equation}

In the previous work 
\cite{CG1}, 
the first and second authors introduced 
an approach to study modulated dispersive equations, 
where they made sense of the integral in \eqref{X0b}
as a nonlinear Young integral.
While our approach is equivalent
to that presented in \cite{CG1}, 
we construct the integral in \eqref{X0b}
as a nonlinear Young integral, using the sewing lemma
(Lemma \ref{LEM:sew}).
Moreover, we have made our presentation 
in a pedestrian manner, accessible
to non-experts.

In Subsection \ref{SUBSEC:Y1}, 
we provide simple criteria
to verify regularity of a given driver.
In 
Subsection \ref{SUBSEC:Y2}, 
we go over the construction of 
a nonlinear Young integral 
via the sewing lemma in an abstract setting.
In 
Subsection \ref{SUBSEC:Y3}, 
we present a basic local well-posedness
result
along with related properties
for the YDE \eqref{YDE0}
in our functional framework.

\subsection{Regularity  of drivers}
\label{SUBSEC:Y1}

As we see in the next subsection (see Proposition \ref{PROP:young1}), 
in our functional framework, 
we can construct the integral $\I^X(\uu)$ in \eqref{X0b}
with the driver $X$ in~\eqref{X0}
as a nonlinear Young integral, 
provided that 
the driver $X$ belongs to the class
 $\cX^{s, \g}_k([0, T]\times \M)$ 
 defined in~\eqref{X1}.
In this subsection, by restricting our attention to the $k$-linear setting, 
we provide
 simple criteria
for checking regularity of drivers.

In this paper, we consider modulated dispersive equation 
with a nonlinearity 
of the form:\footnote{For the renormalized modulated mKdV \eqref{mkdv2}, the nonlinearity is not $\dx u^3$ but can still be viewed as a trilinear operator.
See \eqref{mK0x} and \eqref{mK0} below.}
\begin{align}
\NN(u) = \dx u^k
\label{non0}
\end{align}

\noi
for  an integer $k \ge 2 $.
With a slight abuse of notation, by setting
\begin{align}
\NN(u_1, \dots, u_k) = \dx \bigg(\prod_{j = 1}^k u_j\bigg), 
\label{non1}
\end{align}

\noi
we can view $\NN$ as a $k$-linear operator, 
a priori bounded on $\S(\M)$.
With this notation,  we can view the driver $X$ defined in \eqref{X0}
as a $k$-linear operator given by 
\begin{equation}
X_{t,r}(f_1, \dots, f_k) = \int_r^t \uw(t')^{-1} \NN(\uw(t') f_1, \dots, \uw(t') f_k) dt', \quad 
t > r \ge 0, 
\label{X3}
\end{equation}

\noi
for functions $f_1, \dots, f_k$ on $\M$.
Note that we have 
\begin{align}
X_{t, r}(f) = X_{t, r}(f, \dots, f)
\quad \text{and}\quad
DX_{t, r}[f] (g) = k X_{t, r}(f, \dots, f, g)
\label{X4}
\end{align}

\noi
for functions $f$ and $g$ on $\M$, 
where $DX_{t, r}$ denotes the  Fr\'echet derivative of $X_{t, r}$.

\medskip


In the current $k$-linear setting, 
the following lemma
provides  simple criteria, 
allowing us to  verify that 
a driver $X$ in \eqref{X0} belongs to 
the class
 $\cX^{s, \g}_k([0, T]\times \M)$
and other related classes of drivers.
Since its proof easily follows
 from 
 \eqref{X4}, 
 the $k$-linearity of $X$ in \eqref{X3}
 (see also \eqref{non1}), 
  \eqref{X1}, \eqref{X2}, \eqref{X1x}, and \eqref{X2c}
with 
\eqref{Ho3} and \eqref{Ho4}, 
 we omit details.

\begin{lemma}\label{LEM:OBS1}

Let $X$ be a driver in \eqref{X0}, 
where $U^w(t) = e^{p(\dx) w(t)}$ with $p$ as in \eqref{X0a}
and $\NN$ is as in \eqref{non0} and \eqref{non1}.
Then, given  $s, s_0 \in \R$, $\frac 12 < \g < 1$, an integer $k \ge 2$,  and $T \ge 1$
such that $s_0 > s$, 
the following statements hold\textup{:}

\smallskip
\begin{itemize}
\item[(i)] 
If we have 
\begin{align*}
\|X_{t,r}\|_{\cL_k(H^s(\M))}
\les 
\|\Phi^w\|_{\W^{\rho,\g}_T}
|t-r|^{\g}
\end{align*}

\noi
for any $0 \le r < t \le T$, 
then
 $X$ belongs to $\cX^{s, \g}_k([0, T]\times \M)$
 defined in 
 \eqref{X1}
with the bound 
\begin{align*}
\|X\|_{\cX^{s, \g}_k(T)} \les \|\Phi^w\|_{\W^{\rho,\g}_T}.
\end{align*}

\smallskip
\item[(ii)] \textup{(persistence of regularity).}
Suppose that $X$ can be written as  $X = \sum_{j = 1}^k X^j$
such that for each $j = 1, \dots, k$, we have 
\begin{align}
\begin{split}
X^j &\in \cL_{k, j}^{s, s_0}\\ :\!&=
\cL_k\Big(\big(\bigotimes_{i  = 1}^{j-1}H^{s}(\M)\big) \otimes
H^{s_0}(\M) 
\otimes \big(\bigotimes_{i  = j+1}^{k}H^{s}(\M)\big); H^{s_0}(\M)\Big)
\end{split}
\label{X5}
\end{align}

\noi
with the bound
\begin{align*}
\|X^j_{t,r}\|_{\cL_{k, j}^{s, s_0}}
\les 
\|\Phi^w\|_{\W^{\rho,\g}_T}
|t-r|^{\g}.
\end{align*}

\noi
Then, 
 $X$ belongs to $ \cY^{s, s_0, \g}_k([0, T]\times \M)$
 defined  in \eqref{X2c}
 with the bound 
\begin{align*}
\|X\|_{\cY^{s,s_0,  \g}_k(T)} \les \|\Phi^w\|_{\W^{\rho,\g}_T}.
\end{align*}

\smallskip
\item[(iii)] \textup{(nonlinear smoothing).}
If we have 
\begin{align*}
\|X_{t,r}\|_{\cL_k(H^s(\M); H^{s_0}(\M))}
\les 
\|\Phi^w\|_{\W^{\rho,\g}_T}
|t-r|^{\g}
\end{align*}

\noi
for any $0 \le r < t \le T$, 
then
 $X$ belongs to $\cX^{s, s_0, \g}_k([0, T]\times \M)$
 defined in 
 \eqref{X1x}
with the bound 
\begin{align*}
\|X\|_{\cX^{s,s_0,  \g}_k(T)} \les \|\Phi^w\|_{\W^{\rho,\g}_T}.
\end{align*}

\smallskip
\item[(iv)] \textup{(convergence).}
Given $N \in \N$, 
define
 the truncated driver $X^N$ by 
\begin{equation*}
X^N_{t,r}(f) = \int_r^t \uw(t')^{-1} \P_N \NN(\P_N \uw(t') f) dt'.
\end{equation*}

\noi
If we have
\begin{align*}
\|X^N_{t, r} - X_{t,r}\|_{\cL_k(H^s(\M))}
\les 
o(1)|t-r|^{\g},
\end{align*}

\noi
as $N\to \infty$, 
uniformly in  $0 \le r < t \le T$, 
then, 
 $X^N$
 converges to $X$ in $\cX^{s, \g}_k([0, T]\times \M)$.

\end{itemize}

\end{lemma}

\begin{remark}\label{REM:OBS2} \rm
(i)
While we stated Part (iv) of Lemma \ref{LEM:OBS1}
for a discrete parameter $N \in \N$, 
it also holds
for a continuous parameter.
In particular, we will use Lemma \ref{LEM:OBS1}\,(iv)
not only for establishing convergence of the Galerkin approximations
but also for establishing deep-water and shallow-water
convergence of the (scaled) modulated ILW~\eqref{ILW1} (and~\eqref{ILW2})
claimed in Theorem \ref{THM:conv}.
See Subsections \ref{SUBSEC:BO3}
and \ref{SUBSEC:BO4}.

\smallskip

\noi
(ii) 
By replacing the non-homogeneous Sobolev spaces
with the corresponding homogeneous ones, 
Lemma \ref{LEM:OBS1}
yields criteria
to check if a driver
belongs to the homogeneous
version of the classes of drivers
such as 
$X \in \dot \cX^{s, \g}_k([0, T]\times \M)$
and   
$ \dot \cY^{s, s_0, \g}_k ([0, T]\times \M)$
  defined in \eqref{X1a}
and 
in \eqref{X2d}, respectively.

\end{remark}

\subsection{Sewing lemma and Young integrals}
\label{SUBSEC:Y2}

In this subsection, 
we go over the construction 
of
 the integral $\I^X(\uu)$ in  \eqref{X0b} 
 as a nonlinear Young integral, 
 using the sewing lemma.
The following version of the sewing lemma
is taken from 
\cite{GT10} 
with slight modifications; 
see  
\cite[Proposition 2.3, Corollary 2.4, Corollary 2.5]{GT10}.
See also \cite[Lemma 4.2]{FH20}.

We set 
$C^{1+}_{n,T}V = \bigcup_{\g>1}C^\g_{n,T}V$.

\begin{lemma}[sewing lemma]
\label{LEM:sew}

Let $V$ be a Banach space and 
fix $T>0$. 
Then,  there exists a unique linear map \textup{(}called the sewing map\textup{)}
$\Lambda:C^{1+}_{3,T} V \cap 
\Ker \updl|_{C_{3,T}V}
\to C^{1+}_{2,T}V$ such that 

\smallskip
\begin{enumerate}
\item[(i)] 
We have 
$\updl \Lambda h = h$
for each  $h\in C^{1+}_{3,T}V\cap \Ker  \updl|_{C_{3,T} V}$.

\smallskip

\item[(ii)]
 For each $\z >1$,
the sewing map $\Lambda$ is continuous
from $C^\z _{3,T}V\cap 
\Ker  \updl|_{C_{3,T} V}$ to 
$C^\z _{2,T}V$ such that 
\begin{align}
\|\Lambda h \|_{C^\z _{2,T}V} 
\le \frac{1}{2^\z - 2}  \| h \|_{C^\z _{3,T}V} 
\label{sew1}
\end{align}

\noi
for any $h\in C^\z _{3,T}V$.

\smallskip
\item[(iii)] 
Given any  $g\in C_{2,T}V$ 
with 
$\updl g\in C^\z _{3,T}V$, 
 there exists  unique
$f\in C([0, T];V)$  \textup{(}modulo an additive  constant\textup{)} 
such that 
$\updl f = (\Id - \Lambda \updl)g$. 
In addition, 
we have 
\begin{align}
(\updl f)_{t,r} = \lim_{|\Pi([r,t])|\to 0} 
\sum_{j=0}^{n-1} g_{t_j,t_{j+1}}
\label{sew2}
\end{align}

\noi
 for any $(t,r)\in \Dl_{2,T}$,
 where 
 the limit is over any partition
 $\Pi([r,t])$  
 of  $[r,t]$\textup{:} 
\[\Pi ([r,t]) = \{r = t_n < \dots < t_1 <  t_0 = t\}\]
whose mesh size 
$|\Pi([r,t])| = \sup_{j} |t_j-t_{j+1}|$ 
tends to $0$.

\end{enumerate}

\end{lemma}

As a consequence of the sewing lemma, 
we have the following lemma
on the construction of nonlinear Young integrals.
See also 
 \cite{HK, G23} (which appeared after the first version of the current paper)
 for the construction of  nonlinear Young integrals
 in a more general setting.

\begin{lemma}[nonlinear Young integral]
\label{LEM:int1}
Let $V$ be a Banach space.
Given $0 < \g < 1$, 
 $T>0$, and $k \in \N$, 
let $X\in C^\g_{2,T}\Lip_k(V)$, where $\Lip_k(V)$ is as in \eqref{Lip1}, such that 
\begin{align}
X_{t,r}(0)=0
\label{Ja0}
\end{align} 

\noi
for any $(t,r) \in  \Dl_{2,T}$, 
and 
\begin{align}
\label{Ja1}
(\updl X)_{t_1,t_2,t_3} = 0
\end{align}

\noi
for any $(t_1,t_2,t_3)\in \Dl_{3,T}$.
Given 
  $f\in C^{\al}([0,T];V)$ for some $0 < \al < 1$
  such that $f(0) \in V$, 
define  $\Theta$ on $\Dl_{2,T}$
by setting
\begin{align}
\Theta_{t,r} = X_{t,r}(f(r)), \quad (t,r)\in \Dl_{2,T}.
\label{Ja1a}
\end{align}

\noi
Suppose that 
 $\z : = \al + \g > 1$, 
Then,  the following statements hold\textup{:}

\smallskip

\begin{enumerate}

\item[(i)] 
We have $\Theta\in C^\g_{2,T}V$
and
$\updl \Theta \in C^\z_{3,T}V$.

\smallskip

\item[(ii)]
There exists 
a unique function 
$\I\in C^\g([0,T]; V)$
such that $\I(0) =0$ and
\begin{align}
\updl \I = (\Id - \Lambda \updl)\Theta, 
\label{Ja1b}
\end{align}

\noi
where $\Lambda$ is the sewing map in Lemma~\ref{LEM:sew}. 
Moreover,
we have 
\begin{align}
\|\updl \I - \Theta\|_{C^\z_{2,T}  V}
& \leq 
C_1(f)
\|X\|_{C^\g_{2,T}\Lip_k(V)},
\label{Ja2}
\\
\|\I\|_{C^\g_T  V} 
& \leq
C_2(f)
\|X\|_{C^\g_{2,T}\Lip_k(V)}, 
\label{Ja3}
\end{align}

\noi
where $C_1(f)$ and $C_2(f)$ are given by 
\begin{align}
\begin{split}
C_1(f)
&= 
\frac{1}{2^\z-2}
\big(1+ 2\|f\|_{ L^\infty_T V}\big)^{k-1} \|f\|_{C^\al_T  V}, \\
C_2 (f)
&= (1 \vee T)^\al C_1(f) + 
\big(1 + \|f\|_{ L^\infty_T V} \big)^{k-1}
\|f\|_{L^\infty_T V}.
\end{split}
\label{Ja3x}
\end{align}

\noi
Here,  $a\vee b := \max(a, b)$.
Finally, for $(t,r)\in \Dl_{2,T}$, 
we have
\begin{align}
(\updl \I)_{t,r} = \lim_{|\Pi([r,t])|\to 0} 
\sum^{n-1}_{j=0} \Theta_{t_j,t_{j+1}},
\label{Ja3a}
\end{align}

\noi
 where 
 the limit is in the sense of Lemma \ref{LEM:sew}\,(iii).
 In particular, we have 
\begin{align*}
\I(t) = \lim_{|\Pi([0,t])|\to 0} 
\sum^{n-1}_{j=0} \Theta_{t_j,t_{j+1}}.
\end{align*}

\noi
In the following, we will refer to  $\I$ 
as the \emph{nonlinear Young integral} 
of $f$ with the driver $X$
and 
 denote it
by $\I^X(f)$.

\smallskip

\item[(iii)]
\textup{(persistence of regularity).}
Let $k \ge 2$ and  $V_0 \hookrightarrow V$
be a Banach subspace of $V$.
In addition, suppose that 
 $X\in C^\g_{2,T}\Lip_k(V, V_0; V_0)$, 
 where $\Lip_k(V, V_0; V_0)$ is as in \eqref{Lip2}
 and that 
   $f\in C^{\al}([0,T];V_0)$
 such that $f(0) \in V_0$.
Then, we have 
\begin{align}
\Theta\in C^\g_{2,T}V_0
\qquad \text{and} 
\qquad 
\updl \Theta \in C^\z_{3,T}V_0.
\label{Jaa1}
\end{align}

\noi
Furthermore, 
we have 
$\I\in C^\g([0,T]; V_0)$, 
satisfying 
\begin{align}
\|\updl \I - \Theta\|_{C^\z_{2,T}  V_0}
& \leq 
C_3(f)
\|X\|_{C^\g_{2,T}\Lip_k(V, V_0; V_0)},
\label{Jaa2}
\\
\|\I\|_{C^\g_T  V_0} 
& \leq
C_4(f)
\|X\|_{C^\g_{2,T}\Lip_k(V, V_0; V_0)}, 
\label{Jaa3}
\end{align}

\noi
where $C_3(f)$ and $C_4(f)$ are given by 
\begin{align}
\begin{split}
C_3(f)
&= 
\frac{1}{2^\z-2}\Big(
\big(1+ 2\|f\|_{L^\infty_T V} \big)^{k-1} \|f\|_{C^\al_T  V_0}\\
& \hphantom{XXXXX} + \big(1+ 2\|f\|_{L^\infty_T V} \big)^{k-2} \|f\|_{C^\al_T  V}
\big(1+ 2\|f\|_{L^\infty_T V_0} \big)\Big),\\
C_4 (f)
&= (1 \vee T)^\al C_3 + 
2\big(1 + \|f\|_{ L^\infty_TV}\big)^{k-1} 
\big(1 + \|f\|_{ L^\infty_T V_0} \big).
\end{split}
\label{Jaa4}
\end{align}

\smallskip

\item[(iv)]
\textup{(convergence).}
Suppose in addition  that a sequence $\{X^N\}_{N \in \N}
\subset C^\g_{2,T}\Lip_k(V)$
with $X^N(0) = 0$ and $\updl X^N = 0$
 converges to $X$ in 
$C^\g_{2,T}\Lip_k(V)$
  as $N \to \infty$.
Then, given
any  $f \in \CC^{\al}([0,T];V)$, 
 the nonlinear Young integral $\I^{X^N}(f)$
 of $f$
 with the driver $X^N$ 
 converges to 
 the nonlinear Young integral $\I^X(f)$
 of $f$
 with the driver $X$ 
 in $\CC^\g([0, T]; V)$
as $N \to \infty$.
Moreover, 
given any $r > 0$, 
the rate of convergence is
uniform in $f \in B_r$,  
where
$B_r$ denotes the ball 
in $\CC^{\al}([0,T];V)$
of radius $r > 0$ centered at the origin.

\end{enumerate}
\end{lemma}

\begin{remark}
\rm

\noi
The nonlinear Young integral $\I^X(f)$
of $f$ with a driver $X$
is  sometimes  
denoted by 
\begin{align*}
\int^t_0 X_{dt'}(f(t')).
\end{align*}

\noi
See, for example, 
\cite{CG1}.
\end{remark}

\begin{remark}\label{REM:nonlin}\rm

Let $X$ be as in Lemma \ref{LEM:int1}.
Suppose in addition that  $X\in C^\g_{2,T}\Lip_k(V; V_0)$,
 where   $V_0 \hookrightarrow V$
is a Banach subspace of $V$.
Then, given 
  $f\in \CC^{\al}([0,T];V)$ for some $0 < \al < 1$
  such that $\z = \al + \g > 1$, 
it follows from a straightforward modification of the proof of Lemma~\ref{LEM:int1}
that 
$\I^X(f)\in C^\g([0,T]; V_0)$
with a bound
\begin{align*}
\|\I^X(f)\|_{C^\g_T  V_0} 
& \leq
C_2(f)
\|X\|_{C^\g_{2,T}\Lip_k(V; V_0)}, 
\end{align*}

\noi
where $C_2(f)$ is as in \eqref{Ja3x}.

\end{remark}

\begin{proof}[Proof of Lemma \ref{LEM:int1}]
(i)
By the definition \eqref{Ja1a} of $\Theta$,
\eqref{Ja0}, 
and \eqref{Ho3}, 
we have 
\begin{equation}
\label{Ja5}
\begin{aligned}
\|\Theta_{t,r}\|_{ V} 
& =\|X_{t,r}(f_r)-X_{t,r}(0)\|_{ V}\\
& \leq \|X\|_{C^\g_{2,T}\Lip_k(V)} |t-r|^\g 
\big(1+\|f_r\|_{V}\big)^{k-1} \|f_r\|_{ V}\\
& \leq A_1(f) \|X\|_{C^\g_{2,T}\Lip_k(V)} |t-r|^\g, 
\end{aligned}
\end{equation}

\noi
for any $(t,r)\in\Dl_{2,T}$, 
where $A_1(f)$ is given by 
\begin{align}
A_1 (f)
= \big(1 + \|f\|_{ L^\infty_T V} \big)^{k-1}
\|f\|_{L^\infty_T V}<\infty. 
\label{Ja5x}
\end{align}

\noi
This proves 
$\Theta\in C^\g_{2,T} V$, 
since $X\in C^\g_{2,T}\Lip_k(V)$.

From \eqref{Ja1a} and \eqref{Ja1}, we have 
\begin{align}
\begin{split}
(\updl \Theta)_{t_1,t_2,t_3}
& = X_{t_1,t_3}(f_{t_3})- X_{t_1,t_2}(f_{t_2}) 
- X_{t_2,t_3}(f_{t_3})\\
& = X_{t_1,t_2}(f_{t_3}) - X_{t_1,t_2}(f_{t_2}),
\end{split}
\label{Ja5b}
\end{align}

\noi
for any $(t_1,t_2,t_3)\in \Dl_{3,T}$.
Then, from \eqref{Ho3}
and the $\al$-H\"older continuity of $f$, 
we have 
\begin{align}
\begin{split}
\| (\updl \Theta)_{t_1,t_2,t_3}\|_{ V}
& \leq \|X\|_{C^\g_{2,T}\Lip_k(V)} |t_1-t_2|^\g 
\big(1+\|f_{t_2}\|_{ V} + \|f_{t_3}\|_{ V}\big)^{k-1}
\\
&\quad \times
\|f_{t_2}-f_{t_3}\|_{ V} \\
& \leq A_2(f) \|X\|_{C^\g_{2,T}\Lip_k(V)} |t_1-t_2|^\g 
|t_2 - t_3|^{\al} , 
\end{split}
\label{Ja5c}
\end{align}

\noi
where $A_2(f)$ is given by 
\begin{align}
A_2(f)
&= 
\big(1+ 2\|f\|_{ L^\infty_T V}\big)^{k-1} \|f\|_{C^\al_T  V}.
\label{Ja5d}
\end{align}

\noi
Consequently, 
from \eqref{Ho2}
and $\z = \al + \g$, 
we have 
\begin{align}
\label{INT2}
\|\updl \Theta\|_{C^\z_{3,T} V}
\leq A_2
\|X\|_{C^\g_{2,T}\Lip_k(V)}
<\infty,
\end{align}

\noi
yielding 
 $\updl \Theta \in C^\z_{3,T}  V$.

\medskip

\noi
(ii)
In view of (i) with $\z = \al + \g > 1$, 
we apply the sewing lemma (Lemma \ref{LEM:sew})
and conclude that 
there exists $\I\in C([0,T]; V)$ 
such that 
$\updl \I = (\Id - \Lambda \updl)\Theta$. 
Moreover, from \eqref{Ja1b} and~\eqref{sew1}, we have 
\begin{align}
\label{Ja6a}
\| \updl \I - \Theta \|_{C^\z _{2,T} V} 
= \|   \Lambda (\updl \Theta) \|_{C^\z _{2,T} V} 
\leq \frac{1}{2^\z -2} 
\|\updl \Theta\|_{C^\z _{3,T} V}.
\end{align}

\noi
Hence, 
by  setting 
\begin{align}
C_1(f)  = \frac{1}{2^\z -2} A_2(f),
\label{Ja6aa} 
\end{align}
the bound \eqref{Ja2}
follows from \eqref{Ja6a}
and 
 \eqref{INT2}.

From 
\eqref{Ja2}, we have 
\begin{align}
\label{Ja7}
\|\I_t - \I_r - \Theta_{t,r}\|_{ V}
\leq C_1(f)
\|X\|_{C^\g_{2,T}\Lip_k(V)}|t-r|^{\z},
\end{align}

\noi
for any  $(t,r)\in \Dl_{2,T}$.
Then, from \eqref{Ja7}, \eqref{Ja5}, 
and the triangle inequality, 
we conclude that 
 $\I\in C^\g_T  V$, satisfying the bound
 \eqref{Ja3}.
Lastly, we note that 
the identity \eqref{Ja3a} follows from~\eqref{sew2}.

\medskip

\noi
(iii)
Proceeding as in \eqref{Ja5}
with \eqref{Ho4}, we have 
\begin{equation}
\begin{aligned}
\|\Theta_{t,r}\|_{V_0} 
& =\|X_{t,r}(f_r)-X_{t,r}(0)\|_{ V_0}\\
& \leq A_3(f) \|X\|_{C^\g_{2,T}\Lip_k(V, V_0; V_0)} |t-r|^\g 
\end{aligned}
\label{Jaa5}
\end{equation}

\noi
for any $(t,r)\in\Dl_{2,T}$, 
where $A_3= A_3(f)$ is given by 
\begin{align*}
A_3(f)
= 2\big(1 + \|f\|_{ L^\infty_TV}\big)^{k-1} 
\big(1 + \|f\|_{ L^\infty_T V_0} \big)
<\infty. 
\end{align*}

\noi
This proves 
$\Theta\in C^\g_{2,T} V_0$, 
since $X\in C^\g_{2,T}\Lip_k(V, V_0; V_0)$. 
Proceeding as in \eqref{Ja5c} with~\eqref{Ja5b} and \eqref{Ho4}, we have 
\begin{align*}
\| (\updl \Theta)_{t_1,t_2,t_3}\|_{ V_0}
& \leq \|X\|_{C^\g_{2,T}\Lip_k(V, V_0;V_0)} |t_1-t_2|^\g \notag \\
&\quad \times
\Big\{\big(1+\|f_{t_2}\|_{ V} + \|f_{t_3}\|_{ V}\big)^{k-1}
\|f_{t_2}-f_{t_3}\|_{ V_0} \notag \\
& \hphantom{XXX}+ \big(1+\|f_{t_2}\|_{ V} + \|f_{t_3}\|_{ V}\big)^{k-2} 
\\
& \hphantom{XXXX} \times \big(1+\|f_{t_2}\|_{ V_0} + \|f_{t_3}\|_{ V_0}\big)
\|f_{t_2}-f_{t_3}\|_{ V} \Big\}\notag \\
& \leq A_4(f) \|X\|_{C^\g_{2,T}\Lip_k(V, V_0; V_0)} |t_1-t_2|^\g 
|t_2 - t_3|^{\al} , \notag 
\end{align*}

\noi
where $A_4 (f)$ is given by 
\begin{align*}
A_4(f)
&= 
\big(1+ 2\|f\|_{L^\infty_T V} \big)^{k-1} \|f\|_{C^\al_T  V_0}\\
& \quad + \big(1+ 2\|f\|_{L^\infty_T V} \big)^{k-2} \|f\|_{C^\al_T  V}
\big(1+ 2\|f\|_{L^\infty_T V_0} \big).
\end{align*}

\noi
Consequently, 
from \eqref{Ho2}
and $\z = \al + \g$, 
we have 
\begin{align}
\label{INT2a}
\|\updl \Theta\|_{C^\z_{3,T} V_0}
\leq A_4(f)
\|X\|_{C^\g_{2,T}\Lip_k(V, V_0; V_0)}
<\infty,
\end{align}

\noi
yielding 
 $\updl \Theta \in C^\z_{3,T}  V_0$.
This proves \eqref{Jaa1}.
Then, from 
 the sewing lemma (Lemma \ref{LEM:sew}), 
 we obtain $\I\in C([0,T]; V_0)$. 
Proceeding as in \eqref{Ja6a}, 
the bound \eqref{Jaa2}
with $C_3(f)  = \frac{1}{2^\z -2} A_4(f)$
follows from \eqref{sew1} and \eqref{INT2a}.
From 
\eqref{Jaa2}, we have 
\begin{align}
\label{Jaa7}
\|\I_t - \I_r - \Theta_{t,r}\|_{ V_0}
\leq C_3(f)
\|X\|_{C^\g_{2,T}\Lip_k(V, V_0;V_0)}|t-r|^{\z},
\end{align}

\noi
for any  $(t,r)\in \Dl_{2,T}$.
Then,
the bound \eqref{Jaa3} follows  from \eqref{Jaa7}, \eqref{Jaa5}, 
and the triangle inequality.

\medskip

\noi
(iv)
Let $\I^N(f) = \I^{X^N}(f)$, $N \in \N$, and $\I(f) = \I^X(f)$
be the nonlinear Young integrals of $f$ 
associated with the drivers
$X^N$ and $X$, respectively, 
constructed in Part (i).
Then, 
from~\eqref{Ja1b} and the linearity of the sewing map $\Ld$, we have 
\begin{align}
\updl \I^N(f) - \updl \I(f)  = \Ta^N - \Ta - \Ld  (\updl \Ta^N - \updl \Ta), 
\label{YE1}
\end{align}

\noi
where $\Ta$ is as in \eqref{Ja1a} and $\Ta^N$ is  given by 
\begin{align}
\Theta_{t,r}^N = X^N_{t,r}(f_r), 
  \quad (t,r)\in \Dl_{2,T}.
\label{YE2}
\end{align}

Proceeding as in \eqref{Ja5}, 
we have
\begin{align}
\| \Ta_{t, r}^N - \Ta_{t, r}\|_{V} 
\le A_1(f) \|X^N - X\|_{C^\g_{2,T}\Lip_k(V)}|t - r |^\g, 
\label{YE3}
\end{align}

\noi
where $A_1(f)$ is as in \eqref{Ja5x}.
On the other hand, 
 from \eqref{sew1} in  Lemma~\ref{LEM:sew}, we have 
\begin{align}
\begin{split}
\|(\Lambda (\updl \Theta^N - \updl \Ta))_{t,r}\|_{V}
& \le 
|t-r|^{\z}
\|\Lambda (\updl \Theta^N - \updl \Ta)\|_{C^\z_{2,T} V}\\
& \le
\frac{1}{2^\z - 2}
|t-r|^{\z}
\| \updl \Theta^N - \updl \Ta\|_{C^\z_{3,T} V}.
\end{split}
\label{YE4}
\end{align}

\noi
Proceeding as in \eqref{Ja5c}
with  \eqref{YE2}
and \eqref{Ja5b} with \eqref{Ja6aa}, 
we have 
\begin{align}
\label{YE5}
\begin{aligned}
 \|  \updl & \Theta^N_{t_1, t_2, t_3} - \updl \Ta_{t_1,t_2,t_3}\|_{V}\\
& = \big\|
(X^N_{t_1,t_2} - X_{t_1,t_2}) (f_{t_3}) - (X^N_{t_1,t_2}- X_{t_1,t_2})(f_{t_2})
\big\|_{V}\\
 & \le A_2(f)
 \|X^N - X\|_{C^\g_{2,T}\Lip_k(V)}|t_1 - t_2|^\g 
 |t_2- t_3|^\al, 
\end{aligned}
\end{align}

\noi
where $A_2(f)$ is as in 
\eqref{Ja5d}.
Thus, from \eqref{YE4} and \eqref{YE5}, we obtain
\begin{align}
\|(\Lambda (\updl \Theta^N - \updl \Ta))_{t,r}\|_{V}
& \le 
C_1(f)
 \|X^N - X\|_{C^\g_{2,T}\Lip_k(V)}
|t-r|^{\z}.
\label{YE6}
\end{align}

Hence, from \eqref{YE1}, \eqref{YE3}, and \eqref{YE6}, we have
\begin{align}
\begin{split}
& \|\updl (\I^N(f))_{t, r} - \updl (\I(f))_{t, r}\|_{V}  \\
& \quad \le 
\Big( A_1(f) |t - r|^\g
+  C_1(f)
|t-r|^{\z}\Big)
\|X^N - X\|_{C^\g_{2,T}\Lip_k(V)}.
\end{split}
\label{YE7}
\end{align}

\noi
from which 
we conclude 
that $\I^N(f)$ converges
to $\I(f)$ in $\CC^\g([0, T]; V)$.

Given $r > 0$, it follows from \eqref{Ja5x}
and 
\eqref{Ja5d}
that 
$A_1(f)$ and $A_2(f)$
are uniformly bounded
for $f \in B_r\subset \CC^{\al}([0,T];V)$.
Then, together with \eqref{YE7}, 
we conclude that  
the rate of convergence is
uniform in $f \in B_r$.  
\end{proof}

In the following proposition, 
we summarize
the properties
of the nonlinear Young integral
$\I^X(\uu)$
for a driver 
 $X\in \cX^{s, \g}_k$
 and $\uu \in  \CC^{\al}([0,T]; H^s(\M))$.

\begin{proposition}\label{PROP:young1}
Given   $s \in \R$, $0 < \g < 1$,  $k \in \N$, and $T >0$, 
let   $X \in \cX^{s, \g}_k([0, T]\times \M)$, 
where $ \cX^{s, \g}_k([0, T]\times \M)$ is  defined  in \eqref{X1}
and \eqref{X1z}.
Let
$0 < \al < 1$
such that 
 $\z  = \al + \g > 1$.

\smallskip

\noi
\textup{(i)}
Let  
  $\uu \in \CC^{\al}([0,T];H^s(\M))$.
 Then, the nonlinear Young integral $\I^X(\uu)$
 with the driver $X$ exists
 as the unique function  
$\I^X(\uu)\in C^\g([0,T]; H^s(\M))$
 with $\I^X(\uu)(0) = 0$
 such that 
 \begin{align}
\updl \I^X(\uu) = (\Id - \Lambda \updl)\Theta, 
\label{Y1}
\end{align}

\noi
where $\Ta$ is given by   
\begin{align*}
\Theta_{t,r} = X_{t,r}(\uu(r)), \quad (t,r)\in \Dl_{2,T}.
\end{align*}

\noi
Moreover, we have
\begin{align}
\label{Y2}
\|\updl \I^X(\uu) - \Theta\|_{C^\z_{2,T} H^s_x}
\le
C_1(\uu)\|X\|_{\cX^{s, \g}_k(T)},\\
\label{Y3}
\|\I^X(\uu)\|_{C^\g_T  H^s_x} 
\le
C_2(\uu)\|X\|_{\cX^{s, \g}_k(T)}, 
\end{align}

\noi
where $C_1(\uu)$ and $C_2(\uu)$ are as in \eqref{Ja3x}
with $V = H^s(\M)$, 
and 
\begin{align*}
(\updl \I^X(\uu))_{t,r} = \lim_{|\Pi([r,t])|\to 0} 
\sum^{n-1}_{j=0} \Theta_{t_j,t_{j+1}}
\end{align*}

\noi
for any  $(t,r)\in \Dl_{2,T}$.

\smallskip

\noi
\textup{(ii)}
Suppose that a sequence $\{X^N\}_{N \in \N}\subset 
 \cX^{s, \g}_k([0, T]\times \M)$
 converges to $X$ in 
 $ \cX^{s, \g}_k([0, T]\times \M)$ as $N \to \infty$.
Then, given
any  $\uu \in \CC^{\al}([0,T];H^s(\M))$, 
 the nonlinear Young integral $\I^{X^N}(\uu)$
 of~$\uu$ 
 with the driver $X^N$ 
 converges to 
 the nonlinear Young integral $\I^X(\uu)$
 of~$\uu$ 
 with the driver $X$ 
 in $\CC^\g([0, T]; H^s(\M))$
as $N \to \infty$.
Moreover, 
given any $r > 0$, 
the rate of convergence is
uniform in $\uu \in B_r$ 
where
$B_r$ denotes the ball 
in $\CC^{\al}([0,T];H^s(\M))$
of radius $r > 0$ centered at the origin.

\end{proposition}

\begin{proof}
 By noting from \eqref{X1} and \eqref{X2} 
 (see also \eqref{X1z}) that
$\cX^{s, \g}_k([0, T]\times \M) \subset  C^\g_{2, T}\Lip_k(H^s(\M))$, 
the claims follow from Lemma \ref{LEM:int1} with $V = H^s(\M)$.
\end{proof}

\subsection{Young differential equation}
\label{SUBSEC:Y3}

Our goal in this paper is to 
 study the 
following YDE:
\begin{equation}
\label{YDE1}
\uu(t) = u_0 + \I^X(\uu)(t), 
\end{equation}

\noi
where $\I^X(\uu)$ is the nonlinear Young integral of 
$\uu$ with a driver $X$  of the form \eqref{X0}.
In the following proposition, 
we establish local well-posedness
and related properties
of 
the YDE~\eqref{YDE1} in the abstract setting  with a driver $X\in \cX^{s, \g}_k([0, T] \times \M)$. 
We refer interested readers
 to \cite{G23} for a recent review on
nonlinear YDEs.
In Sections~\ref{SEC:LWP}
and~\ref{SEC:BO}, 
by verifying the hypotheses
in Lemma \ref{LEM:OBS1}, 
we apply Proposition \ref{PROP:main}
to prove local well-posedness
and other related properties
of modulated dispersive equations.

\begin{proposition}
\label{PROP:main}

\textup{(i) (local well-posedness).}
Given   $s \in \R$, $\frac 12 < \g < 1$,  $k \in \N$, and $T \ge 1$, 
let     $X \in \cX^{s, \g}_k([0, T]\times \M)$, 
where $ \cX^{s, \g}_k([0, T]\times \M)$ is  defined  in \eqref{X1}
and \eqref{X1z}.
Then, 
 the Young differential equation \eqref{YDE1} with 
 the driver $X$ 
is locally well-posed in $H^s(\M)$.
More precisely, given  $u_0\in H^s(\M)$, 
there exist $C_0 >0$ and   $\ta>0$, independent of $u_0$ and $X$, 
and 
 a unique solution $\uu \in \cC^\g([0, \tau]; H^s(\M))$
   to 
 \eqref{YDE1}  
with $\uu|_{t= 0} = u_0$, where the local existence time $\tau = 
\tau\big(\|u_0\|_{H^s(\M)}, \|X\|_{\cX^{s, \g}_k(T)}\big) \in (0, 1]$
satisfies 
\begin{align}
\tau \ge  C_0 
\Big(\big(1+ \|X\|_{\cX^{s,\g}_k(T)} \big) \big(1+ \|u_0\|_{H^s(\M)}\big)^{k-1}\Big)^{-\ta}.
\label{YD1}
\end{align}

\noi
Moreover, given $0 < \al < \g$
with $\al + \g > 1$, 
there exists $C_\al > 0$ 
such that 
\begin{align}
\| \uu \|_{\CC^\al_\tau H^s_x}\le C_\al \|u_0\|_{H^s}, 
\label{YD1a}
\end{align}

\noi
while we have
\begin{align}
\| \uu \|_{C^\g_\tau H^s_x}\les \|X\|_{\cX^{s,\g}_k(T)}  (1+\|u_0\|_{H^s})^k.
\label{YD1b}
\end{align}

\smallskip

\noi
\textup{(ii) (persistence of regularity).}
In addition, suppose that $u_0 \in H^{s_0}(\M)$ for some $s_0 > s$
and that  $X \in \cY^{s, s_0, \g}_k([0, T]\times \M)$, 
where $\cY^{s, s_0, \g}_k([0, T]\times \M)$  is  as in \eqref{X2c}.
Then, 
by possibly making $\tau > 0$ smaller by a multiplicative constant 
\textup{(}still satisfying \eqref{YD1}\textup{;} in particular $\tau > 0$
depends on the $H^s$-norm of the initial data $u_0$
but not on its $H^{s_0}$-norm\textup{)}, 
we have $\uu \in \cC^\g([0, \tau]; H^{s_0}(\M))$.

\smallskip

\noi
\textup{(iii) (nonlinear smoothing).}
In addition to the hypotheses in Part (i), 
suppose  that  $X \in \cX^{s, s_0, \g}_k([0, T]\times \M)$, 
where $\cX^{s, s_0, \g}_k([0, T]\times \M)$  is  as in \eqref{X1x}.
Then, we have
 $\I^X(\uu) \in \cC^\g([0, \tau]; H^{s_0}(\M))$, 
 where $\uu$ is the solution to \eqref{YDE1}
 constructed in Part (i).

\smallskip

\noi
\textup{(iv) (convergence).}
In addition to the hypotheses in Part (i), 
suppose  that a sequence $\{X^N\}_{N \in \N}\subset 
 \cX^{s, \g}_k([0, T]\times \M)$
 converges to $X$ in 
 $ \cX^{s, \g}_k([0, T]\times \M)$ as $N \to \infty$.
Then, 
by possibly making $\tau > 0$ smaller by a multiplicative constant 
\textup{(}still satisfying~\eqref{YD1}\textup{)},  the solution 
$\uu^N$  
to the following Young differential equation\textup{:}
\begin{equation}
\label{YDE2}
\uu^N(t) = u_0 + \I^{X^N}(\uu^N)(t)
\end{equation}

\noi
converges to the solution $\uu$ to \eqref{YDE1}
constructed in Part (i)  
 in $\CC^\g([0, \tau]; H^s(\M))$
as $N \to \infty$.
Moreover, 
given any $r > 0$, 
the rate of convergence of $\uu^N$ to $\uu$ is
uniform in $u_0 \in B_r$ 
where
$B_r$ denotes the ball 
in $H^s(\M)$
of radius $r > 0$ centered at the origin.

\end{proposition}

In Proposition \ref{PROP:main}, 
we assumed $T \ge 1$
for simplicity.
With the same proof, 
we have the same results even when $0 < T < 1$, 
in which case, 
the local existence time $\tau > 0$ 
needs to satisfy an additional constraint $\tau \le T$.

\begin{remark}\label{REM:main2}\rm
(i)
Suppose that $X \in \cX^{s, \g}_k(\R_+\times \M)$.
Then, from \eqref{YD1}, we obtain the following blowup alternative;
let $T_* \in (0, \infty]$ be the maximal time of existence
of the solution~$\uu$ to the YDE~\eqref{YDE1}.
Then, we have either 
\begin{align}
T_* = \infty \qquad \text{or} \qquad \lim_{t \nearrow T_*} \|\uu(t)\|_{H^s} = \infty.
\label{BA1}
\end{align}

\smallskip

\noi
(ii)
In Proposition \ref{PROP:main}, we presented 
a local well-posedness result
for the YDE~\eqref{YDE1}
in the non-homogeneous Sobolev space $H^s(\M)$.
Given 
  $X \in \dot \cX^{s, \g}_k([0, T]\times \M)$
  defined in \eqref{X1a}
  and 
  $u_0 \in \dot H^s(\M)$, 
the same proof yields
the corresponding local well-posedness
result in the 
 homogeneous Sobolev spaces.

\smallskip

\noi
(iii) 
Given $\rho >0$,  $\frac12< \g < 1$, and $T> 0$, 
let  $w$ and $w_n$, $n \in \N$,  be $(\rho,\g)$-irregular on $[0, T]$ in the sense of Definition~\ref{DEF:ir}.
Let $X = X^w$ and $X^n = X^{w_n}$ be 
the associated  drivers defined in~\eqref{X0}
with the nonlinearity of the form~\eqref{non0}.
 Suppose that convergence of $w_n$ to $w$ in 
$\W^{\rho,\g}_T$  implies
convergence of 
$X^n$ to $X$ in 
$ \cX^{s, \g}_k([0, T]\times \M)$.
Then, 
from  (the proof of) Proposition~\ref{PROP:main}, 
we see that
the solution map: 
 \[(u_0, w) \in H^s(\M) 
\times  \W^{\rho,\g}_T \longmapsto \uu \in \cC^\g([0, \tau]; H^s(\M))\]

\noi
is continuous.
In our application, 
 convergence of $w_n$ to $w$ in 
$\W^{\rho,\g}_T$ indeed  implies
convergence of the corresponding driver $X^n$ to $X$
 in 
$ \cX^{s, \g}_k([0, T]\times \M)$.
See
Remark \ref{REM:conv1}.

\end{remark}


\begin{proof}[Proof of Proposition \ref{PROP:main}]
(i) 
Let $0 < \al < \g$ such that 
 $\al + \g>1$. 
Given $0 < \tau \le 1$ (to be chosen later), 
 define 
 the map
$\G = \G_{X, u_0}$ on $\cC^\al([0, \tau]; H^s(\M))$
by setting
\begin{align}
\G(\uu)(t) = u_0 + \I^X(\uu)(t), 
\label{YD3}
\end{align}

\noi

\noi
where $\I^X(\uu)$ is the nonlinear Young integral of 
$\uu$ with a driver 
$X \in \cX^{s, \g}_k([0, T]\times \M)$, 
constructed in Proposition~\ref{PROP:young1}.
We show that $\G$ is a contraction on 
the ball $B_R\subset 
\cC^\al([0, \tau]; H^s(\M))$
of radius $R > 0$ 
(to be chosen later)
centered at the origin.

Recalling that $\al + \g > 1$, 
it follows 
from  \eqref{Y3} in 
Proposition \ref{PROP:young1}
with \eqref{Ja3x} that 
\begin{align*}
\|\I^X(\uu)(t) - \I^X(\uu)(r)\|_{H^s}
& \les 
\|X\|_{\cX^{s, \g}_k(T)}
|t-r|^{\g}
\big(1+ 
\|\uu\|_{\CC^\al_{\tau } H^s_x}\big)^{k-1}
\|\uu\|_{\CC^\al_{\tau } H^s_x}
\end{align*}

\noi
for any $0 \le r < t \le \tau$.
Thus, we obtain 
\begin{align}
\label{YD5}
\|\I^X(\uu)\|_{C^\al_\tau H^s_x}
\les \tau^{\g - \al}
\|X\|_{\cX^{s, \g}_k(T)}
\big(1+ 
\|\uu\|_{\CC^\al_{\tau } H^s_x}\big)^{k-1}
\|\uu\|_{\CC^\al_{\tau } H^s_x}.
\end{align}

\noi
Hence, from \eqref{YD3}
and \eqref{YD5} with 
\eqref{Ho2a} and 
$\I^X(\uu)(0) = 0$, 
we have 
\begin{align}
\label{YD5a}
\|\G(\uu)\|_{\CC^\al_\tau H^s_x}
\le \|u_0\|_{H^s} +  C\tau^{\g - \al}
\|X\|_{\cX^{s, \g}_k(T)}
\big(1+ 
\|\uu\|_{\CC^\al_{\tau } H^s_x}\big)^{k-1}
\|\uu\|_{\CC^\al_{\tau } H^s_x}
\end{align}

\noi
for any $0 < \tau \le 1$.
Then, by setting\footnote{In the following, we assume that $u_0 \ne 0$.
If $u_0 = 0$, we can simply set $R = 1$ and carry out the contraction argument
to show that  $\uu \equiv 0$ is a unique solution to \eqref{YDE1}.} 
\begin{align}
R = 2 \|u_0\|_{H^s} 
\label{YD5b}
\end{align}

\noi
and choosing 
$\tau = 
\tau\big(R, \|X\|_{\cX^{s, \g}_k(T)}\big) 
= 
\tau\big(\|u_0\|_{H^s}, \|X\|_{\cX^{s, \g}_k(T)}\big) >0$
sufficiently small, 
we see that 
$\G$ maps the ball $B_R\subset 
\cC^\al([0, \tau]; H^s(\M))$
of radius $R > 0$ 
centered at the origin into itself.

Next, we show that $\G$ is a contraction. 
Given $\uu^1, \uu^2 \in B_R$, let 
\begin{align}
\Theta_{t,r}^1 = X_{t,r}(\uu^1_r)
\qquad \text{and}\qquad 
\Ta_{t,r}^2 = X_{t,r}(\uu^2_r).
\label{YD6}
\end{align}

\noi
Then, 
from \eqref{YD3}, we have 
\begin{align}
\|\G(\uu^1)- \G(\uu^2)\|_{\CC^\al_{\tau }H^s_x}
=
\|\I^X(\uu^1)- \I^X(\uu^2)\|_{\CC^\al_{\tau}H^s_x}
\label{YD7}
\end{align}

\noi
for any $0 < \tau \le 1$.
From \eqref{Y1} and the linearity of $\Ld$, we have 
\begin{align}
\begin{split}
& \|\updl(\I^X(\uu^1))_{t, r}- 
\updl(\I^X(\uu^2))_{t, r}\|_{H^s}\\
& \quad \le
\|\Theta^1_{t,r} - \Ta^2_{t,r} \|_{H^s} 
+
\|(\Lambda (\updl \Theta^1 - \updl \Ta^2))_{t,r}\|_{H^s}
\end{split}
\label{YD8}
\end{align}

\noi
for any $0 \le r < t \le \tau$.
From \eqref{YD6} and  \eqref{Ho3}, we have 
\begin{equation}
\label{YD9}
\begin{aligned}
\|\Theta^1_{t,r} - \Ta^2_{t,r} \|_{H^s} 
& \le 
 \|X\|_{\cX^{s,\g}_k(T)}|t-r|^{\g}\\
& \quad \times 
\big(1+\|\uu^1_r\|_{H^s} + \|\uu^2 _r\|_{H^s}\big)^{k-1}
\|\uu^1_r - \uu^2_r\|_{H^s}
\\
& \les 
\|X\|_{\cX^{s,\g}_k(T)}|t-r|^{\g}
(1+R)^{k-1}
\|\uu^1 - \uu^2\|_{\CC^\al_\tau H^s_x}
\end{aligned}
\end{equation}

\noi
for any $0 \le r < t \le \tau$
and $\uu^1, \uu^2 \in B_R$.

Now, let us estimate the second term on the right-hand side
of \eqref{YD8}.
Let  $\z = \al + \g>1$.
Then, from \eqref{sew1} in  Lemma~\ref{LEM:sew}, we have 
\begin{align}
\label{YD10}
\begin{split}
\|(\Lambda (\updl \Theta^1 - \updl \Ta^2))_{t,r}\|_{H^s}
& \le 
|t-r|^{\z}
\|\Lambda (\updl \Theta^1 - \updl \Ta^2)\|_{C^\z_{2,\tau} H^s_x}\\
& \les
|t-r|^{\z}
\| \updl \Theta^1 - \updl \Ta^2\|_{C^\z_{3,\tau} H^s_x}.
\end{split}
\end{align}

\noi
In the following, we first consider the case $k \ge 2$.
From \eqref{YD6}, 
\eqref{Ja5b}, 
and the fundamental theorem of calculus, 
we have 
\begin{align}
\label{YD11}
\begin{aligned}
 \|  \updl & \Theta^1_{t_1, t_2, t_3} - \updl \Ta^2_{t_1,t_2,t_3}\|_{H^s}\\
& = 
\big\|
X_{t_1,t_2}(\uu^1_{t_3}) - X_{t_1,t_2}(\uu^2 _{t_3})
-\big(X_{t_1,t_2}(\uu^1_{t_2}) - X_{t_1,t_2}(\uu^2_{t_2})\big)\big\|_{H^s}
\\
& =
\bigg\| 
\int^1_0 DX_{t_1,t_2}[\ld \uu^1_{t_3} + (1-\ld)\uu^2 _{t_3}]
(\uu^1_{t_3} - \uu^2_{t_3}) d \ld 
\\
& \hphantom{XX}
- \int^1_0 DX_{t_1,t_2}[\ld \uu^1_{t_2} + (1-\ld)\uu^2 _{t_2}]
 (\uu^1_{t_2} - \uu^2 _{t_2}) d \ld
\bigg\|_{H^s}
\end{aligned}
\end{align}

\noi
for any  $(t_1,t_2,t_3)\in \Dl_{3,\tau}$, where
$D X_{t_1, t_2}$ denotes the 
Fr\'echet derivative of $X_{t_1, t_2}$.
Then, from~\eqref{X2} with \eqref{Lip1}, 
we bound the right-hand side of \eqref{YD11}
by 
\begin{align}
\bigg\| & 
\int^1_0 DX_{t_1,t_2}[\ld \uu^1_{t_3} + (1-\ld)\uu^2 _{t_3}]
\big((\updl \uu^1)_{t_2,t_3} - (\updl \uu^2 )_{t_2,t_3}\big) d \ld \bigg\|_{H^s}
\notag \\
& \quad
+ \bigg\|\int^1_0 
\Big(DX_{t_1,t_2}[\ld \uu^1_{t_3} + (1-\ld)\uu^2 _{t_3}]\notag  \\
& \hphantom{XXXXX}
-  DX_{t_1,t_2}[\ld \uu^1_{t_2} + (1-\ld)\uu^2 _{t_2}]\Big)
 (\uu^1_{t_2} - \uu^2 _{t_2}) d \ld\bigg\|_{H^s}
\label{YD12} \\
& \les
\|X\|_{\cX^{s, \g}_k(T)} |t_1-t_2|^\g 
(1+R)^{k-1}
\|(\updl \uu^1)_{t_2,t_3} - (\updl \uu^2 )_{t_2,t_3}\|_{H^s}
\notag \\
& \quad +
\|X\|_{\cX^{s, \g}_k(T)} |t_1-t_2|^\g 
(1+R)^{k-2}
\big(\|\updl \uu^1_{t_2, t_3}\|_{H^s}+ \|\updl \uu^2_{t_2, t_3}\|_{H^s}\big)
\|\uu^1_{t_2} - \uu^2 _{t_2}\|_{H^s} \notag \\
& \les 
\|X\|_{\cX^{s, \g}_k(T)} |t_1-t_2|^\g 
|t_2 - t_3|^{\al} (1+R)^{k-1}\|\uu^1- \uu^2 \|_{\CC^\al _\tau H^s_x}
\notag 
\end{align}

\noi
for 
 $\uu^1, \uu^2 \in B_R\subset \cC^\al([0, \tau]; H^s(\M))$.
Thus, from \eqref{YD11} and \eqref{YD12} with $\z = \al + \g > 1$, we obtain
\begin{align}
\label{YD13}
\| \updl (\Theta^1 - \Ta^2)\|_{C^\z_{3,\tau} H^s_x}
\les 
\|X\|_{\cX^{s, \g}_k(T)}  (1+R)^{k-1}\|\uu^1- \uu^2 \|_{\CC^\al _\tau H^s_x}.
\end{align}

\noi
When $k = 1$, the conclusion \eqref{YD13}
follows from the second line of \eqref{YD11} and the linearity of $X_{t_1, t_2}$.

Hence, putting
\eqref{YD7}, \eqref{YD8}, 
\eqref{YD9}, \eqref{YD10}, and \eqref{YD13}
together, 
we obtain
\begin{align}
\begin{split}
\|\G(\uu^1)- \G(\uu^2)\|_{\cC^\al _{\tau} H^s_x}
& = 
\|\I^X(\uu^1)- \I^X(\uu^2)\|_{\CC^\al_{\tau}H^s_x}\\
&  \les 
\tau^{\g - \al}\|X\|_{\cX^{s, \g}_k(T)}  (1+R)^{k-1}\|\uu^1- \uu^2 \|_{\CC^\al _\tau H^s_x}.
\end{split}
\label{YD14}
\end{align}

\noi
Then, by 
 choosing 
$\tau = 
\tau\big(R, \|X\|_{\cX^{s, \g}_k(T)}\big) 
= 
\tau\big(\|u_0\|_{H^s}, \|X\|_{\cX^{s, \g}_k(T)}\big) >0$
sufficiently small, 
we conclude that $\G$ is a contraction on $B_R$.
The bound \eqref{YD1} follows from \eqref{YD5a} and \eqref{YD14}.

By the Banach fixed point theorem, 
there exists a unique fixed point $\uu \in B_R$
such that 
\begin{align}
\uu (t)= \G(\uu)(t) = u_0 + \I^X(\uu)(t).
\label{YD15}
\end{align}

\noi
Namely, 
$\uu$ is a solution to the YDE \eqref{YDE1}.
Here, the uniqueness a priori holds only in the ball $B_R$
but, by a standard continuity argument, 
we can extend uniqueness to the entire space $\CC^\al([0, \tau]; H^s(\M))$
(by possibly making $\tau$ smaller by a multiplicative constant).
Moreover,  from Proposition \ref{PROP:young1}, 
we know that the right-hand side of~\eqref{YD15}
belongs to $\CC^\g([0, \tau]; H^s(\M))$, 
from which we conclude that $\uu \in \CC^\g([0, \tau]; H^s(\M))$.

Finally the bound \eqref{YD1a}
follows from the fact that $\uu \in B_R$ with \eqref{YD5b}, 
while 
 the bound~\eqref{YD1b}
follows from 
\eqref{YD15}, 
\eqref{Y3}, and  \eqref{Ja3x} 
with the fact that $\uu \in B_R$ and \eqref{YD5b}.

\medskip

\noi
(ii)
In Part (i), we constructed a solution to \eqref{YDE1}
on the time interval $[0, \tau]$ with the bound 
\begin{align*}
\| \uu\|_{\CC^\al_\tau H^s_x} \le R.
\end{align*}

\noi
Proceeding as in \eqref{YD5}  and \eqref{YD5a} with 
\eqref{YD15}, 
\eqref{Jaa3},  \eqref{Jaa4}, 
and \eqref{X2c}, 
we have
\begin{align*}
\| \uu\|_{\CC^\al_\tau H^{s_0}_x} \le  \|u_0\|_{H^{s_0}} 
+ C \tau^{\g- \al}
\|X\|_{\cY^{s, s_0, \g}_k(T)}
(1+ R)^{k-1} \| \uu\|_{C^\al_\tau H^{s_0}_x}.
\end{align*}

\noi
Then, by possibly making $\tau > 0$ smaller (by a multiplicative constant, and hence \eqref{YD1}
still holds), 
we conclude that 
$\uu \in  \CC^\al([0, \tau]; H^{s_0}(\M))$.
Then, from 
Lemma \ref{LEM:int1}\,(iii), 
we obtain 
 $\I^X(\uu) \in 
 \CC^\g([0, \tau]; H^{s_0}(\M))$
 and thus
  $\uu \in 
 \CC^\g([0, \tau]; H^{s_0}(\M))$
 in view of \eqref{YD15}.


\medskip

\noi
(iii)
Let   $\uu \in \cC^\g([0, \tau]; H^s(\M))$ be the solution constructed in Part (i).
Then, since $2\g > 1$, 
it follows from 
Remark \ref{REM:nonlin}
that 
 $\I^X(\uu) \in \cC^\g([0, \tau]; H^{s_0}(\M))$.

\medskip

\noi
(iv)
For simplicity, we assume that $X\ne 0$.
We first note that the convergence of $X^N$ to $X$ in 
$ \cX^{s, \g}_k([0, T]\times \M)$ implies
that there exists $N_0 \in \N$ such that 
\begin{align}\|X^N\|_{\cX^{s, \g}_k(T)} \le
2 \|X\|_{\cX^{s, \g}_k(T)} 
\label{YDE3a}
\end{align}

\noi
for any $N \ge N_0$.
Thus, by making $\tau$ smaller by a multiplicative constant, 
Part (i) guarantees existence of 
a unique solution $\uu^N \in \CC^\g([0, \tau]; H^s(\M))$
to \eqref{YDE2} for each $N \ge N_0$.
Moreover, we have 
\begin{align}
\|\uu^N \|_{\CC^\al_\tau H^s_x} , \|\uu \|_{\CC^\al_\tau H^s_x} \le R
\label{YDE3b}
\end{align}

\noi
for any $N \ge N_0$, where $R$ is as in Part (i).

From \eqref{YDE1} and \eqref{YDE2}, we have
\begin{align}
\label{YDE3}
\begin{split}
\uu(t) - \uu^N(t) 
& = 
 \I^{X}(\uu)(t)-  
 \I^{X^N}(\uu^N)(t)\\
& = 
\Big( \I^{X}(\uu)(t)-  
 \I^{X^N}(\uu)(t)\Big)
 + \Big(\I^{X^N}(\uu)(t)- 
 \I^{X^N}(\uu^N)(t)\Big).
\end{split}
\end{align}

\noi
From 
Proposition \ref{PROP:young1}\,(ii), 
we have 
\begin{align}
\lim_{N\to \infty }\| \I^{X}(\uu)-  
 \I^{X^N}(\uu)\|_{\CC^\g_\tau H^s_x} = 0.
\label{YDE4}
\end{align}

\noi
By applying
 \eqref{YD8}, 
\eqref{YD9}, \eqref{YD10}, and \eqref{YD13}
with \eqref{YDE3a} and \eqref{YDE3b}, 
we obtain
\begin{align}
\| \I^{X^N}(\uu)- 
 \I^{X^N}(\uu^N)\|_{\cC^\al _{\tau} H^s_x}
\les 
\tau^{\g - \al}\|X\|_{\cX^{s, \g}_k(T)}  (1+R)^{k-1}\|\uu- \uu^N \|_{\CC^\al _\tau H^s_x}
\label{YDE5}
\end{align}

\noi
for any $N \ge N_0$.
Hence, 
 by possibly taking $\tau > 0$ smaller (by a multiplicative constant), 
it follows 
from \eqref{YDE3}, 
\eqref{YDE4}, 
and \eqref{YDE5} that 
\begin{align}
\|\uu- \uu^N \|_{\CC^\al _\tau H^s_x}
\les 
\| \I^{X}(\uu)-  
 \I^{X^N}(\uu)\|_{\CC^\g_\tau H^s_x} \too 0, 
\label{YDE6}
\end{align}

\noi
as $N \to \infty$.
This proves convergence of $\uu^N$ to $\uu$ in 
 $\CC^\al([0, \tau]; H^s(\M))$.

We now upgrade this convergence to that 
 in $\CC^\g([0, \tau]; H^s(\M))$.
By applying  \eqref{YD8}, 
\eqref{YD9}, \eqref{YD10}, and \eqref{YD13}
with \eqref{YDE3a} and \eqref{YDE3b}
 once again, 
we have
\begin{align}
\| \I^{X^N}(\uu)- 
 \I^{X^N}(\uu^N)\|_{\cC^\g _{\tau} H^s_x}
\les 
\|X\|_{\cX^{s, \g}_k(T)}  (1+R)^{k-1}\|\uu- \uu^N \|_{\CC^\al _\tau H^s_x}.
\label{YDE7}
\end{align}

\noi
Then, the claimed convergence of $\uu^N$ to $\uu$ 
 in $\CC^\g([0, \tau]; H^s(\M))$
 follows from \eqref{YDE3}, 
 \eqref{YDE4}, and \eqref{YDE7} with \eqref{YDE6}. 
Moreover, given $r > 0$, 
the uniform convergence rate for $\|u_0\|_{H^s}\le r$
follows
from the fact that the 
rate of the convergence \eqref{YDE4}
is uniform in 
$\|u_0\|_{H^s}\le r$ 
(Proposition \ref{PROP:young1}\,(ii))
and the bound \eqref{YDE3b}
which is uniform with respect to the initial data
$u_0$ 
with $\|u_0\|_{H^s}\le r$.

This completes the proof of Proposition \ref{PROP:main}.
\end{proof}

\section{Local well-posedness
of the modulated KdV and mKdV}
\label{SEC:LWP}

In this section, we prove local well-posedness
and related properties
(persistence of regularity, nonlinear smoothing,  and convergence of the Galerkin approximation)
of the modulated KdV-type equations:

\smallskip

\begin{itemize}
\item
Theorem~\ref{THM:1}\,(i) and Theorem \ref{THM:G1} 
for the modulated 
KdV \eqref{kdv1}
on the circle,

\smallskip

\item 
Theorems~\ref{THM:3} and \ref{THM:G2}
for the modulated 
mKdV \eqref{mkdv1}
on the circle,

\smallskip

\item 
Theorem~\ref{THM:2}\,(i)
for  the modulated KdV \eqref{kdv1}
on the real line.

\end{itemize}

\smallskip

\noi
These results follow
from 
Proposition \ref{PROP:main}, 
once we verify the conditions
in 
Lemma \ref{LEM:OBS1}
for each of the equations.

\subsection{Modulated KdV equation on the circle}
\label{SUBSEC:K1}
We consider the  modulated KdV~\eqref{kdv1} on~$\T$.
More precisely, as we discussed in Section \ref{SEC:1},
we consider the following equation satisfied by 
the modulated interaction representation $\uu$
of $u$ defined in~\eqref{int1}:
\begin{equation}
\dt \uu = \uw(t)^{-1} \dx\big( (\uw(t) \uu)^2\big), 
\label{kdv1a}
\end{equation}

\noi
where 
$\uw (t)=e^{-   w(t)\dx^3}  $
denotes the modulated linear propagator for \eqref{kdv1}.
By writing \eqref{kdv1a} (with initial data $u_0$) in the integral form, we have
\begin{equation}
\uu(t) = u_0 + \int_0 ^t \uw(t')^{-1}
\dx\big( (\uw(t') \uu(t'))^2\big)dt'.
\label{mild3x}
\end{equation}

\noi
Then, 
the bilinear driver $X= X^\KDV$ associated with the modulated KdV  is 
given by 
 \begin{align}
X_{t,r}(f_1,f_2)
= X^\KDV_{t,r}(f_1,f_2)
=\int_r^t  \uw(t')^{-1}  \partial_x
\big( (\uw(t')  f_1)( \uw(t') f_2 )\big) dt'
\label{K1}
\end{align}

\noi
for functions $f_1$ and $f_2$ on $\T$.
By taking the Fourier transform, we have
\begin{align}
\F\big(X_{t,r} (f_1,f_2)\big) (n)
= in \sum_{ \substack{n_1, n_2 \in \Z^*\\n = n_1+n_2}}
 \Phi^w_{t,r}(\Xi_\KDV (\bar n))
\ft f_1(n_1)  \ft f_2(n_2),
\label{K2}
\end{align}

\noi
where $\Phi^w_{t,r}$ is as in \eqref{rho2}
and $\Xi_\KDV (\bar n)$ denotes the {\it resonance function}\footnote{Here, we follow
the terminology in \cite{Tao2}.  
We point out that $\Xi_\KDV (\bar n)$ is also called the modulation function (see \cite{KOY})
but we do not use this latter terminology to avoid
confusion with a modulation function $w(t)$.
}
 for KdV given by 
\begin{align}
\Xi_\KDV (\bar n) &  = \Xi_\KDV (n,n_1,n_2) 
= - n^3+ n_1^3 +n_2^3.
\label{K3}
\end{align}

\noi
Recall that, under $n = n_1+ n_2$, we have
\begin{align}
\Xi_\KDV (\bar n) = - 3n n_1 n_2.
\label{K3a}
\end{align}

The following proposition establishes
basic mapping properties
of the bilinear driver $X$ defined in \eqref{K1}.

\begin{proposition}
\label{PROP:kdv1}

Given $\rho \ge\frac 12$, $\frac 12 < \g < 1$, and $T> 0$, 
let  $w$ be $(\rho,\g)$-irregular on $[0, T]$ in the sense of Definition~\ref{DEF:ir}.

\smallskip

\begin{itemize}
\item[(i)]
Suppose that $\rho \ge \frac 12$ and $s \in \R$
satisfy 
\eqref{reg1}.
Then, 
the driver $X$ defined in \eqref{K1}
belongs to $ \cX^{s, \g}_2([0, T]\times \T)$
defined in~\eqref{X1}.

\smallskip
\item[(ii)] \textup{(persistence of regularity).}
Suppose that $\rho \ge \frac 12$ and $s \in \R$
satisfy 
\eqref{reg1}.
Then, 
for any $s_0 > s$, 
the driver $X$ 
belongs to 
$\cY^{s, s_0, \g}_{2}([0, T]\times \T)$
defined in~\eqref{X2c}.

\smallskip
\item[(iii)]
\textup{(nonlinear smoothing).}
In addition, suppose that $s_0 > s$ satisfies
\eqref{reg3}.
Then, 
the driver $X$ 
belongs to 
$\cX^{s, s_0, \g}_{2}([0, T]\times \T)$
defined in~\eqref{X1x}.

\smallskip
\item[(iv)]
\textup{(Galerkin approximation).}
Suppose that $\rho \ge \frac 12$ and $s \in \R$
satisfy 
\eqref{reg3a}.
Then, 
the truncated driver $X^N = X^{\KDV,  N}$, defined by 
\begin{align}
X_{t,r}^N(f_1,f_2)=\int_r^t \uw(t')^{-1} \P_N \partial_x\big(
( \P_N \uw(t')  f_1 )(\P_N \uw(t') f_2 )\big) dt', 
\label{XN}
\end{align}

\noi
 converges
to  $X$ 
in  $ \cX^{s, \g}_2([0, T]\times \T)$ as $N \to \infty$.

\end{itemize}

 \end{proposition}

Once we prove Proposition \ref{PROP:kdv1}, 
local well-posedness of the modulated KdV \eqref{kdv1}
on $\T$
and other related properties
(Theorem \ref{THM:1}\,(i) and Theorem \ref{THM:G1}) follow
from Proposition~\ref{PROP:main}
with Proposition~\ref{PROP:kdv1}
applied to \eqref{mild3x}.
We omit details.

\begin{proof}[Proof of Proposition \ref{PROP:kdv1}]

In the following, we verify the hypotheses
in Lemma \ref{LEM:OBS1} 
with $k = 2$.
Then, the claimed results follow from 
Lemma \ref{LEM:OBS1}.


\medskip

\noi
(i)
From \eqref{K2}
and H\"older's inequality, 
we have
\begin{align}
\begin{aligned}
\| X_{t,r} (f_1,f_2)\|_{H^s}^2 
&\le \sum_{n \in \Z^*} \jb{n}^{2s+2} \bigg |  \sum_{ \substack{n_1, n_2 \in \Z^*\\n = n_1+n_2}}
\Phi^{w}_{t,r}(\Xi_\KDV(\bar n) ) 
\ft f_1(n_1)  \ft f_2 (n_2)  \bigg|^2 \\
&\le \sum_{n \in \Z^*} \jb{n}^{2s+2}   \sum_{ \substack{n_1, n_2 \in \Z^*\\n = n_1+n_2}}
\frac{|\Phi^{w}_{t,r}(\Xi_\KDV(\bar n) )|^2 }{\jb{n_2}^{2s}}
|\ft f_1(n_1)|^2
\|f_2\|_{H^s(\T)}^2
\\
&\leq  
\bigg(\sup_{n_1 \in \Z^*}
\sum_{\substack{n, n_2 \in \Z^*\\n_1 = n - n_2}} \jb{n}^{2s+2}
\frac{|\Phi^{w}_{t,r} (\Xi_\KDV(\bar n)) |^2}{\jb{n_1}^{2s}\jb{n_2}^{2s}} \bigg)
\|f_1\|_{H^s}^2
\|f_2\|_{H^s}^2
\end{aligned}
\label{K4}
\end{align}

\noi
for any $0 \le r <  t \le T$.
Then, from 
\eqref{K4}, 
 \eqref{rho1}, and
 \eqref{K3a},
 we have 
\begin{align}
\begin{aligned}
\|X_{t,r} \|_{\cL_2 (H^s)}
&\les
 \| \Phi^w \|_{\W_T^{\rho,\g}} 
  |t-r|^{\g} \\
  & \quad
  \times 
\bigg(\sup_{n_1 \in \Z^*}
\sum_{\substack{n, n_2 \in \Z^*\\n_1 = n - n_2}} 
 \jb{n}^{2- 4\rho}
 \frac{\jb{n}^{2s+2\rho}}{\jb{n_1}^{2s+2\rho} \jb{n_2}^{2s +2\rho}}\bigg)^\frac 12
\end{aligned}
\label{K5}
\end{align}

\noi
for any $0 \le r <  t \le T$.
Without loss of generality, 
assume that $|n_1|\ges |n_2|$ in \eqref{K5}.
Then, we have $|n_1|\ges |n|$ under $n = n_1 + n_2$.
Thus, we have 
\begin{align}
\sup_{n_1 \in \Z^*}
\sum_{\substack{n, n_2 \in \Z^*\\n_1 = n - n_2}} 
 \jb{n}^{2- 4\rho}
 \frac{\jb{n}^{2s+2\rho}}{\jb{n_1}^{2s+2\rho} \jb{n_2}^{2s +2\rho}}
\les 
\sup_{n_1 \in \Z^*}
\sum_{n\in \Z^*}
 \frac{1}{\jb{n}^{4\rho- 2}\jb{n - n_1}^{2s +2\rho}}, 
\label{K5a}
\end{align}

\noi
provided that 
 $s + \rho \ge 0$.
 In view of Lemma \ref{LEM:SUM}, 
we see that the right-hand side of \eqref{K5a}
is finite if 
 \begin{align*}
4\rho- 2 \ge 0, \quad 2s + 2\rho  \ge 0, 
\quad \text{and}\quad 
 4\rho - 2 + 2s  +2\rho > 1, 
 \end{align*}

\noi
Namely, 
 \begin{align}
\rho \ge \frac 12, \quad s \ge - \rho,
\quad \text{and}\quad 
s > \frac 32 - 3 \rho, 
\label{K5b}
 \end{align}

\noi
which is satisfied under  \eqref{reg1}.
Hence,
from  Lemma \ref{LEM:OBS1}\,(i),  
we conclude that 
 $X\in  \cX^{s, \g}_2([0, T]\times \T)$ under~\eqref{reg1}.

\medskip

\noi
(ii)  Write $X = X^1 + X^2$, 
where $X^1$ is the contribution from $|n_1| \ge |n_2|$ in \eqref{K2}
and $X^2 = X - X^1$.
We only consider $X^1$ since a similar bound holds for $X^2$ by symmetry.
Proceeding as in \eqref{K4}, \eqref{K5}, and \eqref{K5a}
with $|n_1|\ges |n|$
and \eqref{X5}, 
we have 
\begin{align*}
\|X^1_{t,r} \|_{\cL_{2, 1}^{s, s_0}}
&\les
 \| \Phi^w \|_{\W_T^{\rho,\g}} 
  |t-r|^{\g} \\
  & \quad
  \times 
\bigg(
\sup_{n_1 \in \Z^*}
\sum_{n\in \Z^*}
 \frac{1}{\jb{n}^{4\rho- 2}\jb{n - n_1}^{2s +2\rho}}
\bigg)^\frac 12 \\
 &\les
 \| \Phi^w \|_{\W_T^{\rho,\g}} 
  |t-r|^{\g}
\end{align*}

\noi
for any $0 \le r <  t \le T$, 
provided that \eqref{reg1} holds.
Hence,
from  Lemma \ref{LEM:OBS1}\,(ii),  
we conclude that, if \eqref{reg1} holds, 
then
 $X\in  \cY^{s, s_0, \g}_2([0, T]\times \T)$
 for any $s_0 > s$.

\medskip

\noi
(iii) Without loss of generality, 
let us assume $|n_1|\ges |n|$ in \eqref{K2}.
Then, proceeding as in Part~(i), 
we have 
\begin{align}
\begin{aligned}
\|X_{t,r} \|_{\cL_2 (H^s; H^{s_0})}
&\les
 \| \Phi^w \|_{\W_T^{\rho,\g}} 
  |t-r|^{\g} \\
  & \quad
  \times 
\bigg(
\sup_{n_1 \in \Z^*}
\sum_{n\in \Z^*}
 \frac{1}{\jb{n}^{4\rho- 2s_0 + 2s- 2}\jb{n - n_1}^{2s +2\rho}}
\bigg)^\frac 12
\end{aligned}
\label{K6}
\end{align}

\noi
for any $0 \le r <  t \le T$.
 In view of Lemma \ref{LEM:SUM}, 
we see that the right-hand side of \eqref{K6}
is finite if 
 \begin{align*}
s_0 \le s + 2\rho - 1, \quad s \ge - \rho, 
\quad \text{and}\quad 
s_0 < 2s + 3 \rho - \frac 32, 
 \end{align*}

\noi
which is satisfied under  \eqref{reg3}.
Hence,
from  Lemma \ref{LEM:OBS1}\,(iii),  
we conclude that 
 $X\in  \cX^{s,s_0,  \g}_2([0, T]\times \T)$ under~\eqref{reg3}.

\medskip

\noi
(iv)
From  \eqref{K2} with \eqref{XN}, we have 
\begin{align*}
& \F\big( X_{t,r}(f_1,f_2)\big)(n) - \F\big( X^N_{t,r}(f_1,f_2)\big)(n)\\
& \quad =in \sum_{ \substack{n_1, n_2 \in \Z^*\\n = n_1+n_2}}
 \ind_{\{\max(|n|, |n_1|, |n_2| )>  N\}}
\Phi^w_{t,r}(\Xi_\KDV(\bar n))
\ft f_1(n_1) \ft f_2(n_2).
\end{align*}

\noi
Note that 
\begin{equation}
 \ind_{\{\max(|n|, |n_1|, |n_2|) >  N\}} 
 \les N^{-\eps} 
\max(|n|, |n_1|, |n_2| )^{\eps}.
\label{K8}
\end{equation}

\noi
Assume
 that $|n_1|\ges |n_2|$ as in Part (i)
(which implies $|n_1| \sim \max(|n|, |n_1|, |n_2| )$), 
Then, by proceeding as in \eqref{K5} and \eqref{K5a}
\begin{align}
\begin{aligned}
& \|X_{t,r} - X^N_{t, r} \|_{\cL_2 (H^s)}\\
&\quad \les N^{-\eps}
 \| \Phi^w \|_{\W_T^{\rho,\g}} 
  |t-r|^{\g} 
\bigg(\sup_{n_1 \in \Z^*}
\sum_{\substack{n, n_2 \in \Z^*\\n_1 = n - n_2}} 
 \jb{n}^{2- 4\rho+2\eps}
 \frac{\jb{n}^{2s+2\rho-2\eps}}{\jb{n_1}^{2s+2\rho-2\eps} \jb{n_2}^{2s +2\rho}}\bigg)^\frac 12\\
&\quad \les N^{-\eps}
 \| \Phi^w \|_{\W_T^{\rho,\g}} 
  |t-r|^{\g} 
\bigg(\sup_{n_1 \in \Z^*}
\sum_{n\in \Z^*}
 \frac{1}{\jb{n}^{4\rho- 2-2\eps}\jb{n - n_1}^{2s +2\rho}}\bigg)^\frac 12.
 \label{K8a}
\end{aligned}
\end{align}

\noi
for any $0 \le r <  t \le T$, 
provided that 
 $s + \rho \ge \eps$.
Moreover, 
in view of Lemma \ref{LEM:SUM}, 
we see that the last factor on the right-hand side of \eqref{K8a}
is finite if 
 \begin{align*}
\rho \ge \frac 12 + \frac 12\eps, \quad s \ge - \rho,
\quad \text{and}\quad 
s > \frac 32 - 3 \rho, 
 \end{align*}

\noi
which is satisfied under  \eqref{reg3a}.
Hence,
from  Lemma \ref{LEM:OBS1}\,(iv),  
we conclude that 
$X^N$ converges to $X$ in 
 $  \cX^{s, \g}_2([0, T]\times \T)$ under~\eqref{reg3a}
 by choosing $\eps > 0$ sufficiently small.

This concludes the proof of Proposition \ref{PROP:kdv1}.
\end{proof}

\begin{remark}\label{REM:conv1} \rm

Given $\rho \ge\frac 12$,  $\frac12< \g < 1$, and $T> 0$, 
let  $w$ and $w_n$, $n \in \N$,  be $(\rho,\g)$-irregular on $[0, T]$ in the sense of Definition~\ref{DEF:ir}.
Let $X = X^w$ and $X^n = X^{w_n}$ be 
the   bilinear drivers 
for the modulated KdV on the circle
defined in \eqref{K1}
associated with the modulations $w$ and $w_n$, respectively.
Let  $s \in \R$
satisfy 
\eqref{reg1}
such that Proposition \ref{PROP:kdv1}\,(i) guarantees
$X, X^n \in  \cX^{s, \g}_2([0, T]\times \T)$.
Now,  suppose that $w_n$ converges to $w$ in 
$\W^{\rho,\g}_T$ as $n \to \infty$.
Then, a straightforward  modification of 
\eqref{K4} and \eqref{K5} with the finiteness of the sum
on the right-hand side of~\eqref{K5}
under \eqref{reg1} yields
\begin{align}
\|X^n_{t,r} - X_{t,r} \|_{\cL_2 (H^s)}
&\les
 \| \Phi^{w_n} - \Phi^w \|_{\W_T^{\rho,\g}} 
  |t-r|^{\g}
  = o(1)   |t-r|^{\g}
\label{YE8}
\end{align}

\noi
as $n \to \infty$.
Hence, we conclude from Lemma \ref{LEM:OBS1}\,(iv)
that 
 $X^n$ converges to $X$ 
in $ \cX^{s, \g}_2([0, T] \times \T)$
with the Lipschitz bound \eqref{YE8}.
Then, from Remark \ref{REM:main2}\,(iii)
with \eqref{YE7} and\eqref{YE8}, 
we see that the solution map:
 \[(u_0, w) \in H^s(\T) 
\times  \W^{\rho,\g}_T \longmapsto \uu \in \cC^\g([0, \tau]; H^s(\T))\]

\noi
is in fact locally Lipschitz continuous.
 A similar comment applies
to all the modulated dispersion equations
we consider in this paper.

\end{remark}

\subsection{Modulated mKdV equation on the circle}
\label{SUBSEC:K2}

In this subsection,
we consider the modulated (renormalized) mKdV  \eqref{mkdv2} on $\T$.
The renormalized nonlinearity in \eqref{mkdv2} is given by 
\begin{align}
\NN(f)  = 3   \big(f^2- \|f\|_{L^2}^2\big)\dx f.
\label{mK0x}
\end{align}

\noi
By taking the Fourier transform and symmetrization, we have 
\begin{align}
\ft{\NN(f)}(n) 
& = in \sum_{A_n}
\ft f(n_1)\ft f(n_2)\ft f(n_3)
- 3 in |\ft f(n)|^2 \ft f(n), 
\label{mK0}
\end{align}

\noi
where 
\begin{align}
A_n =\big\{(n_1, n_2,n_3)\in (\Z^*)^3 : 
n = n_1+n_2+n_3, \ n\ne n_1, n_2, n_3\big\}.
\label{mK2b}
\end{align}

\noi
See
\cite[Lemma 8.16]{BO93}
and \cite[(1.9)]{TT}.
With a slight abuse of notation, 
define 
the trilinear version of $\NN(f)$ by 
setting
\begin{align}
\begin{split}
\F\big(\NN(f_1, f_2, f_3)\big)(n) 
& = in \sum_{A_n}
\ft f_1(n_1)\ft f_2(n_2)\ft f_3(n_3)\\
& \quad - 3 in \ft f_1(n)\cj{\ft f_2(n)}\ft f_3(n).
\end{split}
\label{mK0a}
\end{align}

\noi
Then, 
the trilinear driver associated with the modulated mKdV \eqref{mkdv2} is given by 
 \begin{align}
X_{t,r}(f_1,f_2, f_3)=\int_r^t  \uw(t')^{-1}
\NN( \uw(t')  f_1 , \uw(t') f_2,  \uw(t') f_3) dt',
\label{mK1}
\end{align}

\noi
\noi
where 
$\uw (t)=e^{- w(t)\dx^3 }  $
and $\NN$ is as in \eqref{mK0a}.
By taking the Fourier transform, we have 
\begin{align}
\begin{split}
\F\big(X_{t,r} (f_1,f_2, f_3) \big)(n)
& = in  \sum_{ A_n }
 \Phi^w_{t,r}(\Xi_\MKDV )
\ft f_1(n_1)  \ft f_2(n_2) \ft f_3(n_3)\\
& \quad - 3 i (t-r)n \ft f_1(n)\cj{\ft f_2(n)}\ft f_3(n)\\
& =:
\F\big(X_{t,r}^{(1)} (f_1,f_2, f_3) \big)(n)
+ \F\big(X_{t,r}^{(2)} (f_1,f_2, f_3) \big)(n), 
\end{split}
\label{mK2}
\end{align}

\noi

\noi
where $ \Phi^w_{t,r} (\cdot)$ is as in \eqref{rho2} and 
  $X^{(1)}$ corresponds to the non-resonant part of $X$ (namely 
$\Xi_\MKDV (\bar n) \ne 0$), while 
 $X^{(2)}$ corresponds to the resonant part of $X$ (namely 
$\Xi_\MKDV (\bar n) = 0$).
Here, $\Xi_\MKDV (\bar n)$ denotes the resonance function for mKdV given by 
\begin{align*}
\Xi_\MKDV (\bar n)
= \Xi_\MKDV (n,n_1,n_2, n_3) 
&= - n^3+ n_1^3 +n_2^3 + n_3^3.
\end{align*}

\noi
Recall that, under $n = n_1+ n_2+n_3$, we have
\begin{align}
\Xi_\MKDV (\bar n)
& = - 3 (n-n_1)(n-n_2)(n-n_3).
\label{mK3a}
\end{align}

As in the case of the modulated KdV studied in the previous subsection, 
once we prove the following proposition 
on basic mapping properties
of the trilinear driver $X$ defined in \eqref{mK1}, 
local well-posedness of the modulated mKdV 
and other related properties
(Theorem \ref{THM:3} and Theorem \ref{THM:G2}) follow
from  Proposition~\ref{PROP:main}.

\begin{proposition}\label{PROP:mkdv1}
Given $\rho \ge\frac 12$,  $\frac 12 < \g < 1$, and $T> 0$, 
let  $w$ be $(\rho,\g)$-irregular on $[0, T]$ in the sense of Definition~\ref{DEF:ir}.

\smallskip

\begin{itemize}
\item[(i)]
Suppose that $s \ge \frac 12$. 
Then, 
the driver $X$ defined in \eqref{mK1}
belongs to $ \cX^{s, \g}_3([0, T]\times \T)$
defined in~\eqref{X1}.

\smallskip
\item[(ii)] \textup{(persistence of regularity).}
Suppose that $s \ge \frac 12$. 
Then, 
for any $s_0 > s$, 
the driver $X$ 
belongs to 
$\cY^{s, s_0, \g}_{3}([0, T]\times \T)$
defined in~\eqref{X2c}.

\smallskip
\item[(iii)]
\textup{(nonlinear smoothing).}
In addition, suppose that $s_0 > s \ge \frac 12$ satisfies
\eqref{regM2}.
Then, 
the driver $X$ 
belongs to 
$\cX^{s, s_0, \g}_{3}([0, T]\times \T)$
defined in~\eqref{X1x}.

\smallskip
\item[(iv)]
\textup{(Galerkin approximation).}
Suppose that $\rho > \frac12$ and $s > \frac 12$.
Then, 
the truncated driver $X^N$, defined by 
 \begin{align}
X_{t,r}^N(f_1,f_2, f_3)=\int_r^t  \uw(t')^{-1}
\P_N \NN( \P_N \uw(t')  f_1 ,\P_N  \uw(t') f_2, \P_N  \uw(t') f_3) dt',
\label{XN2}
\end{align}

\noi
 converges
to  $X$ 
in  $ \cX^{s, \g}_3([0, T]\times \T)$ as $N \to \infty$, 
where $\NN$ is as in \eqref{mK0a}.

\end{itemize}

\end{proposition}

\begin{proof}
Let $\rho, s \ge \frac 12$.
As in the proof of Proposition \ref{PROP:kdv1}, 
we verify  the hypotheses
in Lemma~\ref{LEM:OBS1} 
with $k = 3$.

\medskip

\noi
(i) 
We first consider  the resonant part $X^{(2)}$.
From \eqref{mK2}, 
we have
\begin{align*}
\|X_{t,r}^{(2)} \|_{\cL_3 (H^s)}
\les 
  |t-r| 
\end{align*}

\noi
for any $0 \le r < t \le T$, 
provided that 
 $s \ge \frac 12$, just as in the usual (unmodulated) mKdV case.

Next, we consider the non-resonant part  $X^{(1)}$.
From \eqref{mK2} and H\"older's inequality, we have 
\begin{equation}
\begin{split}
& \|X_{t,r}^{(1)}(f_1, f_2, f_3)\|^2_{H^{s} }
\le\sum_{n\in \Z^*} \jb{n}^{2s+2}
\bigg| \sum_{A_n} 
\Phi^{w}_{t,r} (\Xi_\MKDV(\bar n ))
\prod_{j = 1}^3 \ft f_j(n_j)
 \bigg|^2\\
&\quad \les \bigg(
\sup_{n\in \Z^*} 
\sum_{A_n} 
\frac{\jb{n}^{2s+2} }{\jb{n_1}^{2s} \jb{n_2}^{2s}\jb{n_3}^{2s}}
|\Phi^{w}_{t,r}(\Xi_\MKDV(\bar n))|^2  \bigg) 
\prod_{j=1}^3  \| f_j\|^2_{H^{s}  }   .
\end{split}
\label{mK4x}
\end{equation}

\noi
By symmetry, 
assume that $|n_1| \ges \max (|n_1|, |n_2|,  |n_3|)$, 
which implies $|n_1|\ges |n|$.
Then, from~\eqref{rho1}
and   \eqref{mK3a} with $s \ge 0$ and \eqref{mK2b}, we have 
\begin{align}
\begin{aligned}
\|X_{t,r}^{(1)}\|_{\cL_3 (H^{s})} 
&\leq \bigg( \sup_{n\in \Z^*} 
\sum_{A_n}
\frac{\jb{n}^{2} }{ \jb{n_2}^{2s}\jb{n_3}^{2s}}
|\Phi^{w}_{t,r}( \Xi_\MKDV (\bar n)) |^2\bigg)^\frac 12\\
&\leq A_\MKDV \| \Phi^w \|_{\W_T^{\rho,\g}}|t-r|^{\g}
\end{aligned}
\label{mK4a}
\end{align}

\noi
for any $0 \le r <  t \le T$, 
where $A_\MKDV$ is given by 
\begin{align*}
A_\MKDV^2 =  \sup_{n\in \Z^*} 
\sum_{A_n}
\frac{\jb{n}^{2} }{\jb{n_2 + n_3}^{2\rho}\prod_{j = 2}^3\jb{n_j}^{2s} \jb{n - n_j}^{2\rho}}.
\end{align*}

\noi
Fix small $\eps > 0$.
From H\"older's inequality, Young's inequality (with $\rho \ge \frac 12$),
and Lemma \ref{LEM:SUM}, 
 we have
\begin{align}
\begin{split}
 A_\MKDV^2 
& \le  \sup_{n\in \Z^*} \jb{n}^2
\bigg\| \frac1{\jb{n_2}^{2s}\jb{n - n_2}^{2\rho}}\bigg\|_{\l^\frac{2(1+\eps)}{1+2\eps}_{n_2}}^2
\bigg\| \frac1{\jb{n_2}^{2\rho}}\bigg\|_{\l^{1+\eps}_{n_2}}\\
& \les  \sup_{n\in \Z^*} \jb{n}^2
\bigg\| \frac1{\jb{n_2}^{2s}\jb{n - n_2}^{2\rho}}\bigg\|_{\l^\frac{2(1+\eps)}{1+2\eps}_{n_2}}^2\\
& \les 
\sup_{n\in \Z^*} 
\jb{n}^{2 - 4\min(s, \rho) } < \infty, 
\end{split}
\label{mK8}
\end{align}

\noi
\noi
provided that 
$\min(s, \rho)\ge \frac 12$.

\medskip

\noi
(ii) 
As for   the resonant part $X^{(2)}$, it suffices
to note that 
\begin{align*}
\frac{\jb{n}^{2s_0+2}}{\jb{n}^{2s_0 + 4s}}
\les 1, 
\end{align*}

\noi
uniformly in $n \in \Z^*$, 
provided that 
 $s \ge \frac 12$.
As for the non-resonant part $X^{(1)}$, 
without loss of generality,
assume that $|n_1| \ges \max (|n_1|, |n_2|,  |n_3|)$, 
which implies $|n_1|\ges |n|$.
%
%
%
Then, proceeding as in \eqref{mK4x} and~\eqref{mK4a}
with $s\ge 0$ (see also \eqref{X5}), 
we have 
\begin{align*}
\|X_{t,r}^{(1)}\|_{\cL_{3, 1}^{s, s_0}}
&\leq \bigg( \sup_{n\in \Z^*} 
\sum_{\substack{A_n\\|n_1|\ges |n|}}
\frac{\jb{n}^{2s_0+2} }{\jb{n_1}^{2s_0} \jb{n_2}^{2s}\jb{n_3}^{2s}}
|\Phi^{w}_{t,r}( \Xi_\MKDV (\bar n)) |^2\bigg)^\frac 12\\
&\leq A_\MKDV \| \Phi^w \|_{\W_T^{\rho,\g}}|t-r|^{\g}\\
&\les \| \Phi^w \|_{\W_T^{\rho,\g}}|t-r|^{\g}
\end{align*}

\noi
for any $0 \le r <  t \le T$, 
provided that 
$\min(\rho, s) \ge  \frac 12$.

\medskip

\noi
(iii) 
As for   the resonant part $X^{(2)}$, it suffices
to note that 
\begin{align*}
\frac{\jb{n}^{2s_0+2}}{\jb{n}^{6s}}
\les 1, 
\end{align*}

\noi
\noi
uniformly in $n \in \Z^*$, 
provided that 
 $s_0 \le 3s - 1$.
As for the non-resonant part $X^{(1)}$, 
without loss of generality,
assume that $|n_1| \ges \max (|n_1|, |n_2|,  |n_3|)$, 
which implies $|n_1|\ges |n|$.
Then, 
proceeding as in~\eqref{mK4x} and \eqref{mK4a}
with $s_0 > s\ge 0$, 
we have 
\begin{align*}
\|X_{t,r}^{(1)}\|_{\cL_3(H^s; H^{s_0})}
&\leq A_\MKDV^{(1)} \| \Phi^w \|_{\W_T^{\rho,\g}}|t-r|^{\g}
\end{align*}

\noi
for any $0 \le r <  t \le T$, 
where $ A_\MKDV^{(1)}$ is given by 
\begin{align*}
( A_\MKDV^{(1)})^2
 =
 \sup_{n\in \Z^*} 
\sum_{A_n}
\frac{\jb{n}^{2s_0 - 2s +2} }{\jb{n_2 + n_3}^{2\rho} \prod_{j = 2}^3\jb{n_j}^{2s} \jb{n - n_j}^{2\rho}}.
\end{align*}

\noi
Here, by a crude estimate, 
we have dropped 
the factor $\jb{n - n_1}^{2\rho}$ (see Remark \ref{REM:mkdv2} below).
By proceeding as in \eqref{mK8} with 
Lemma \ref{LEM:SUM}, 
we have 
\begin{align*}
( A_\MKDV^{(1)})^2 \les 
\sup_{n\in \Z^*} 
\jb{n}^{2 s_0 -2s + 2 - 4\min(\rho, s)} < \infty, 
\end{align*}

\noi
provided that \eqref{regM2} holds.

\medskip

\noi
(iv) 
From \eqref{mK2} with \eqref{XN2}, we have 
\begin{align*}
 \F & \big(X_{t,r} (f_1,f_2, f_3) \big)(n)
- \F\big(X^N_{t,r} (f_1,f_2, f_3) \big)(n)\\
& = in  \sum_{ \wt A_n }
 \Phi^w_{t,r}(\Xi_\MKDV )
\ft f_1(n_1)  \ft f_2(n_2) \ft f_3(n_3)\\
& \quad - \ind_{\{|n| > N\}}\cdot 3 i (t-r)n \ft f_1(n)\cj{\ft f_2(n)}\ft f_3(n).
\end{align*}

\noi
Here,  $\wt A_n$ is given by 
\begin{align*}
\wt A_n = 
\big\{(n_1, n_2,n_3)\in A_n : 
\max(|n|, |n_1|, |n_2|, |n_3|) > N\big\}, 
\end{align*}

\noi
where  $A_n$ is as in \eqref{mK2b}.
Let $X^{N, (1)}$ and $X^{N, (2)}$ denote the non-resonant and resonant 
parts of $X^N$, respectively.
Then, we have
\begin{align*}
\|X_{t,r}^{(2)}- X_{t,r}^{N, (2)} \|_{\cL_3 (H^s)}
\les N^{-\eps}
  |t-r| 
\end{align*}

\noi
for any $0 \le r <  t \le T$, 
provided that $s \ge \frac 12 + \frac 1 2\eps$.

Let us consider the non-resonant part.
From~\eqref{mK3a}, we have 
\[
\max(|n|,|n_1|,|n_2|,|n_3|)
\les \jb{\Xi_\MKDV(\bar n)}
\]

\noi
on $A_n$ defined in \eqref{mK2b}.
Consequently, we have 
\begin{align}
\ind_{\wt A_n}
\les N^{-\eps_0}\jb{\Xi_\MKDV(\bar n)}^{\eps_0}
\label{mKX1}
\end{align}

\noi
for $\eps_0 > 0$.
By repeating 
the computations in \eqref{mK4a} and \eqref{mK8}
with $\min(s, \rho- \eps_0) \ge \frac 12$, we obtain
\begin{align*}
\|X_{t,r}^{(1)}- X_{t,r}^{N, (1)} \|_{\cL_3 (H^s)}
\les N^{-\eps_0}
\| \Phi^w \|_{\W_T^{\rho,\g}}|t-r|^{\g}
\end{align*}

\noi
for any $0 \le r <  t \le T$.
\end{proof}

\subsection{Modulated KdV equation  on the real line}
\label{SUBSEC:K3}

We conclude this section by 
considering
 the  modulated KdV  \eqref{kdv1} on the real line.
The bilinear driver associated with the modulated KdV  is 
given by $X$ in \eqref{K1}, 
where we view all the operators and functions 
 as those defined on $\R$.
By taking the Fourier transform, we have
\begin{align}
\F_\R \big(X_{t,r} (f_1,f_2)\big) (\xi)
= i\xi \intt_{\xi = \xi_1 + \xi_2}
 \Phi^w_{t,r}(\Xi_\KDV (\bar \xi))
\ft f_1(\xi_1)  \ft f_2(\xi_2)d\xi_1,
\label{KR1}
\end{align}

\noi
where $\Phi^w_{t,r}$ is as in \eqref{rho2}
and 
\begin{align*}
\Xi_\KDV (\bar \xi) &  = \Xi_\KDV (\xi,\xi_1,\xi_2) 
= - \xi^3+ \xi_1^3 +\xi_2^3.
\end{align*}

\noi
As in \eqref{K3a}, 
 we have
\begin{align}
\Xi_\KDV (\bar \xi) = - 3\xi \xi_1 \xi_2
\label{KR1b}
\end{align}

\noi
 under $\xi = \xi_1+ \xi_2$.
Then, as in the previous subsections, 
local well-posedness of the modulated KdV \eqref{kdv1}
on $\R$
(Theorem \ref{THM:2}\,(i)) follows
from Proposition~\ref{PROP:main}
once we prove
the following proposition.

\begin{proposition}
\label{PROP:KR1}

Given $\rho > \frac 12$,  $\frac 12 < \g < 1$, and $T> 0$, 
let  $w$ be $(\rho,\g)$-irregular on $[0, T]$ in the sense of Definition~\ref{DEF:ir}.

\smallskip

\begin{itemize}
\item[(i)]
Suppose that $\rho >  \frac 12$ and $s \in \R$
satisfy 
\eqref{regR1}.
Then, 
the driver $X$ defined in \eqref{K1}
belongs to $ \cX^{s, \g}_2([0, T]\times \R)$
defined in~\eqref{X1}.

\smallskip
\item[(ii)] \textup{(persistence of regularity).}
Suppose that $\rho >  \frac 12$ and $s \in \R$
satisfy 
\eqref{regR1}.
Then, 
for any $s_0 > s$, 
the driver $X$ 
belongs to 
$\cY^{s, s_0, \g}_{2}([0, T]\times \R)$
defined in~\eqref{X2c}.

\end{itemize}

\end{proposition}

The regularity restriction $s > - \frac 32$ in \eqref{regR1}
(for $\rho \ge \frac 32$) is (essentially) sharp.
See Remark~\ref{REM:low1} below.

\begin{proof}[Proof of Proposition \ref{PROP:KR1}]
While we proceed as in the proof of Proposition \ref{PROP:kdv1}
for the modulated KdV on $\T$, 
we need to pay particular attention to low frequencies 
in the current real line case.

\medskip

\noi
(i) 
Proceeding as in \eqref{K4}, we have 
\begin{align}
\| X_{t,r} \|_{\cL_2(H^s(\R))}^2
&\leq  
\sup_{\xi_1 \in \R}
\int_\R \jb{\xi}^{2s}
\xi^2
\frac{|\Phi^{w}_{t,r} (\Xi_\KDV(\bar \xi)) |^2}{\jb{\xi_1}^{2s}\jb{\xi - \xi_1}^{2s}} d\xi.
\label{KR0a}
\end{align}

\noi
On the other hand, by H\"older's inequality, we have 
\begin{align}
\begin{split}
\| X_{t,r} (f_1,f_2)\|_{H^s(\R)}^2
&\leq
 \int_\R \jb{\xi}^{2s}
 \xi^2
 \bigg(  \intt_{\xi = \xi_1 + \xi_2}
\frac{ |\Phi^w_{t,r} ( \Xi_\KDV(\bar \xi) )|^2 }{\jb{\xi_1}^{2 s}  \jb{\xi_2}^{2s} } d\xi_1  \bigg) \\
& \quad\times 
 \Big(\intt_{\xi = \xi_1 + \xi_2}   \jb{\xi_1}^{2 s} \jb{\xi_2}^{2s}
|\ft f_1(\xi_1)  |^2  |\ft f_2(\xi_2)  |^2  d\xi_1 \Big)  d\xi  \\
&\leq
\bigg(\sup_{\xi\in\R} 
  \int_\R 
  \jb{\xi}^{ 2s }
  \xi^2
   \frac{|\Phi^{w}_{t,r} (\Xi_\KDV(\bar \xi)) |^2}{ \jb{\xi_1}^{2s}  \jb{\xi - \xi_1}^{2s}}  d\xi_1\bigg)
   \|f_1\|_{H^s(\R)}^2
\|f_2\|_{H^s(\R)}^2.
\end{split}
\label{KR0x}
\end{align}

\noi
Thus, we have 
\begin{align}
\| X_{t,r} \|_{\cL_2(H^s(\R))}^2
&\leq  
\sup_{\xi\in\R} 
  \int_\R 
  \jb{\xi}^{ 2s}
  \xi^2
   \frac{|\Phi^{w}_{t,r} (\Xi_\KDV(\bar \xi)) |^2}{ \jb{\xi_1}^{2s}  \jb{\xi - \xi_1}^{2s}}  d\xi_1.
\label{KR0b}
\end{align}

\noi
Then, 
by dividing the frequency space $\R^2_{\xi, \xi_1}$ into the following four regions:
\begin{align}
\begin{aligned}
D_1 &=  \{| \xi |< 1, |\xi_1|< 2\} , \\
D_2 &= \{| \xi |< 1, |\xi_1|\ge 2\}, \\
D_3 &= \Big\{ | \xi |\ge 1, |\xi_1|\ge \frac12, |\xi_2|\ge \frac12  \Big\},   \\
D_4 &= 
 \Big\{   |  \xi  |\ge 1, |\xi_1|\ge \frac12, |\xi_2|< \frac12  \Big\} 
 \cup
\Big\{| \xi  |\ge 1, | \xi_1 |< \frac12, |\xi_2| \ge \frac12  \Big\}
\end{aligned}
\label{KR2b}
\end{align}

\noi
and applying \eqref{KR0a} and \eqref{KR0b}
(with symmetry $\xi_1\leftrightarrow\xi_2$ on $D_4$), 
we obtain
\begin{align}
\|X_{t , r }\|^2_{\cL _2 (H^s (\R))}
 &
\le
\sup_{|\xi|< 1} 
  \intt_{|\xi_1|< 2} 
  \jb{\xi}^{ 2s }\xi^2
   \frac{|\Phi^{w}_{t,r} (\Xi_\KDV(\bar \xi)) |^2}{ \jb{\xi_1}^{2s}  \jb{\xi - \xi_1}^{2s}}  d\xi_1
\notag \\
& \quad + 
\sup_{|\xi_1|\ge 2}
\intt_{|\xi|< 1} \jb{\xi}^{2s}
\xi^2
\frac{|\Phi^{w}_{t,r} (\Xi_\KDV(\bar \xi)) |^2}{\jb{\xi_1}^{2s}\jb{\xi - \xi_1}^{2s}} d\xi
\notag \\ 
& \quad +
\sup_{|\xi_1|\ge \frac 12}
\intt_{\substack{|\xi|\ge 1\\|\xi-\xi_1|\ge \frac 12 }} 
\jb{\xi}^{2s}\xi^2
\frac{|\Phi^{w}_{t,r} (\Xi_\KDV(\bar \xi)) |^2}{\jb{\xi_1}^{2s}\jb{\xi - \xi_1}^{2s}} d\xi
\label{KR2} \\
& \quad +
2\sup_{|\xi|\ge 1} 
  \intt_{|\xi - \xi_1|< \frac 12} 
  \jb{\xi}^{ 2s }\xi^2
   \frac{|\Phi^{w}_{t,r} (\Xi_\KDV(\bar \xi)) |^2}{ \jb{\xi_1}^{2s}  \jb{\xi - \xi_1}^{2s}}  d\xi_1
\notag \\
& =:
A_1+ A_2 + A_3 + A_4.
\notag
\end{align}

From \eqref{rho1} with \eqref{KR1b}, we have 
\begin{align}
A_1 
\les  \|\Phi^w\|^2_{\W_T^{\rho,\g}} |t-r |^{2\g}
\label{KR3}
\end{align}

\noi
\noi
for any $0 \le r <  t \le T$
without any restriction on  $\rho > 0$ and  $s \in \R$.

Next, we consider $A_2$.
When $s \ge 0$, it follows from the triangle inequality and \eqref{rho1}
that 
\begin{align}
A_2
\les  \|\Phi^w\|^2_{\W_T^{\rho,\g}} |t-r |^{2\g}
\label{KR3a}
\end{align}
 
\noi
for any $0 \le r <  t \le T$
without any restriction on  $\rho > 0$ (and  $s\ge 0$).
Now, suppose that $ s < 0$. Note from \eqref{rho1} that
\begin{align}
|\Phi^w_{t,r}(a)|
\le
\|\Phi^w\|_{\W^{\rho,\g}_T}
\frac{|t-r|^\g}{ \jb{a}^{\rho}}
\le
\|\Phi^w\|_{\W^{\rho,\g}_T}
\frac{|t-r|^\g}{ \jb{a}^{\rho'}}
\label{KR3b}
\end{align}

\noi
 for any $a\in \R$
and 
 $0<\rho'\le\rho$.
Under $|\xi| < 1$ and $|\xi_1|\ge 2$, 
we have $|\xi - \xi_1|> 1$.
In particular, we have $|\xi_1|\sim \jb{\xi_1}$
and $|\xi - \xi_1|\sim \jb{\xi- \xi_1}$.
Moreover, 
we have 
$\jb{\xi\xi_1(\xi-\xi_1)}\ge
 |\xi \xi_1(\xi- \xi_1)|$.
Then, with \eqref{KR3b}, we have
\begin{align}
\begin{split}
A_2
& \les
\|\Phi^w\|_{\W^{\rho,\g}_T}^2
|t-r|^{2\g}
\intt_{|\xi|<1}  \xi^2 
\sup_{\substack{|\xi_1|\ge 2\\|\xi - \xi_1|> 1}}
\frac{1}{\jb {\xi \xi_1 (\xi-\xi_1)}^{2\rho'}\jb {\xi_1}^{2s} \jb {\xi-\xi_1}^{2s}} d\xi
\\
& \les 
\|\Phi^w\|_{\W^{\rho,\g}_T}^2
|t-r|^{2\g}
\intt_{|\xi|<1} |\xi|^{2-2\rho' } 
\sup_{\substack{|\xi_1|\ge 2\\|\xi - \xi_1|> 1}}
\frac{1}{| \xi_1 (\xi-\xi_1)|^{2s+2\rho'} } d\xi\\
& \les 
\|\Phi^w\|_{\W^{\rho,\g}_T}^2
|t-r|^{2\g} 
\end{split}
\label{KR3c}
\end{align}

\noi
for any $0 \le r <  t \le T$, 
provided that $2-2\rho'>-1$ and $ s+\rho' \ge 0$, 
which can be satisfied by choosing 
suitable $0 < \rho' \le \rho$
if 
\begin{align}
\rho > 0, \quad s > -\frac 32,
\quad \text{and} \quad  
s\ge -\rho.
\label{KR3d}
\end{align}

As for $A_3$, 
there is no low frequency issue.
Thus, by 
proceeding as in \eqref{K5} and \eqref{K5a}, we have
\begin{align}
\begin{aligned}
A_3
& \les
\|\Phi^w\|^2_{\W_T^{\rho,\g}}|t-r |^{2\g}
\sup_{|\xi_1|\ge \frac 12}
\intt_{\substack{|\xi|\ge 1\\|\xi-\xi_1|\ge \frac 12 }}
\frac{ \jb{\xi}^{2s+2-2\rho}}{\jb{\xi_1}^{2s+2\rho}\jb{\xi - \xi_1}^{2s+2\rho}} d\xi\\
& \les
\|\Phi^w\|^2_{\W_T^{\rho,\g}}|t-r |^{2\g}
\end{aligned}
\label{KR4}
\end{align}

\noi
for any $0 \le r <  t \le T$, 
provided that
\eqref{K5b} holds.

As for $A_4$, 
we first note that 
$ |\xi_1|\sim |\xi| \ge 1$
and $\jb{\xi - \xi_1}\sim 1$
under 
 $|\xi|\ge1$ and $|\xi - \xi_1|< \frac12 $.
 Then, from~\eqref{rho1}, \eqref{KR1b},  and 
 a change of variables $\zeta_2 = \xi^2 (\xi - \xi_1)$, we have 
 \begin{align}
\begin{aligned}
A_4
&\les
\|\Phi^w\|^2_{\W_T^{\rho,\g}}
|t-r |^{2\g }    \sup_{|\xi|\ge 1}  \xi^2
  \intt_{|\xi - \xi_1|< \frac 12} 
\frac{ 1}{  \jb{\xi^2(\xi - \xi_1)}^{2\rho}} d \xi_1\\
&= 
\|\Phi^w\|^2_{\W_T^{\rho,\g}}
|t-r |^{2\g }  
   \sup_{|\xi|\ge 1} 
  \intt_{|\z_2|< \frac {\xi^2}2} 
 \frac{1}{\jb{\z_2}^{2\rho} }d \z_2
 \les \|\Phi^w\|^2_{\W_T^{\rho,\g}}
|t-r |^{2\g } 
\end{aligned}
\label{KR5}
\end{align}

\noi
for any $0 \le r <  t \le T$, provided that $\rho > \frac 12$
(and any $s \in \R$).

Therefore, putting \eqref{KR2}, 
\eqref{KR3}, 
\eqref{KR3a}, 
\eqref{KR3c}, 
\eqref{KR4}, and \eqref{KR5}
together, we conclude 
from  Lemma \ref{LEM:OBS1}\,(i)
that 
 $X\in  \cX^{s, \g}_2([0, T]\times \R)$ under~\eqref{regR1}.

\medskip

\noi
(ii)  Write $X = X^1 + X^2$, 
where $X^1$ is the contribution from $|\xi_1| \ge |\xi_2|$ in \eqref{KR1}
and $X^2 = X - X^1$.
We only consider $X^1$ since $X^2$ can be handled 
in an analogous manner.
From \eqref{KR2b}
with \eqref{X5}, we have 
\begin{align*}
\|X^1_{t,r} \|_{\cL_{2, 1}^{s, s_0}}^2
 \le \wt A_1+ \wt A_2 + \wt A_3 +\wt  A_4, 
\end{align*}

\noi
where $\wt A_j$, $j = 1, \dots, 4$, 
is obtained from $A_j$ in \eqref{KR2}
by replacing 
$\jb{\xi}^{2s}$
and $\jb{\xi_1}^{2s}$
with $\jb{\xi}^{2s_0}$
and $\jb{\xi_1}^{2s_0}$, respectively.

Suppose that  \eqref{regR1} holds and that $s_0 > s$.
From 
\eqref{KR3}, 
a slight modification of \eqref{KR4}
(as in the proof of Proposition \ref{PROP:kdv1}\,(ii)), and 
\eqref{KR5}, 
we have 
\begin{align*}
\wt A_1 + 
\wt A_3 + 
\wt A_4
& \les
\|\Phi^w\|^2_{\W_T^{\rho,\g}}|t-r |^{2\g}
\end{align*}

\noi
for any $0 \le r <  t \le T$.
When $s \ge 0$, we also have 
\begin{align*}
\wt A_2
& \les
\|\Phi^w\|^2_{\W_T^{\rho,\g}}|t-r |^{2\g}
\end{align*}

\noi
for any $0 \le r <  t \le T$.
When $s < 0$, 
a slight modification of~\eqref{KR3c} with $|\xi_1|\sim |\xi - \xi_1|$ yields
\begin{align*}
\wt A_2
& \les
\|\Phi^w\|_{\W^{\rho,\g}_T}^2
|t-r|^{2\g}
\intt_{|\xi|<1}  \xi^2 
\sup_{\substack{|\xi_1|\ge 2\\|\xi - \xi_1|> 1}}
\frac{1}{\jb {\xi \xi_1 (\xi-\xi_1)}^{2\rho'}\jb {\xi_1}^{2s_0} \jb {\xi-\xi_1}^{2s}} d\xi
\\
& \les 
\|\Phi^w\|_{\W^{\rho,\g}_T}^2
|t-r|^{2\g}
\intt_{|\xi|<1} |\xi|^{2-2\rho' } 
\sup_{|\xi_1|\ge 2}
\frac{1}{| \xi_1|^{2s_0+2s + 4\rho'} } d\xi\\
& \les 
\|\Phi^w\|_{\W^{\rho,\g}_T}^2
|t-r|^{2\g} 
\end{align*}

\noi
for any $0 \le r <  t \le T$, 
provided that $2-2\rho'>-1$ and $ s_0 + s+2\rho' \ge 0$, 
which can be satisfied by choosing 
suitable $0 < \rho' \le \rho$
under  \eqref{regR1} and  $s_0 > s$.
Hence,
from  Lemma \ref{LEM:OBS1}\,(ii),  
we conclude that, if \eqref{regR1} holds, 
then
 $X\in  \cY^{s, s_0, \g}_2([0, T]\times \R)$
 for any $s_0 > s$.
\end{proof}

\begin{remark} \label{REM:low1} \rm

The regularity restriction $s > - \frac 32$
appears in the analysis of $A_2$ in \eqref{KR3c}; 
see~\eqref{KR3d}.
In this remark, we show that 
the regularity restriction $s > - \frac 32$ is essentially necessary
(modulo the endpoint $s = -\frac 32$), 
regardless of the value of $\rho$.

Given $N \gg 1$, define a function $f$ on $\R$ by setting
\begin{align}
\ft f(\xi) = N^{-s-\frac 12} \cdot \ind_{|\xi|\sim N}
\label{AZ0}
\end{align}
such that 
\begin{align}
\| f \|_{H^s} \sim 1.
\label{AZ1}
\end{align}
We now estimate the $H^s$-norm
of $X_{t,r} (f,f)$ from below
by considering the contribution coming from $|\xi| \le c_0 N^{-2}$
for some small constant $c_0 > 0$.
Namely, 
\begin{align}
|\Xi_\KDV (\bar \xi)|
= 3|\xi \xi_1\xi_2|\les c_0 \ll 1.
\label{AZ1a}
\end{align}

\noi
In this case, it follows from~\eqref{rho2}
and the continuity of the modulation function $w$
that 
\begin{align}
\Re \Phi^w_{t,r}(\Xi_\KDV(\bar \xi))
\ges_w  |t-r|,
\label{AZ2}
\end{align}

\noi
provided that $t$ is sufficiently close to (given) $r$
(and by choosing $c_0 > 0$ sufficiently small, depending on the value of $w(r)$).
Then, from \eqref{KR1} with \eqref{AZ0} and \eqref{AZ2}, we have 
\begin{align}
\begin{split}
\| X_{t,r} (f,f)\|_{H^s(\R)}^2 
& \ges_{w, t, r}
N^{-4s-2} 
 \intt_{|\xi| \le c_0 N^{-2}}
|\xi|^2 
\bigg(\int
\ind_{|\xi_1|\sim N} \cdot \ind_{|\xi-\xi_1|\sim N}
d\xi_1\bigg)^2 d\xi \\
& \ges
N^{-4s} 
 \intt_{|\xi| \le c_0 N^{-2}}
|\xi|^2 
 d\xi 
\sim  
N^{-4s-6}. 
\end{split}
\label{AZ3}
\end{align}

\noi
Hence, from \eqref{AZ1} and \eqref{AZ3}
and taking $N \to \infty$, 
we conclude that if $X_{t, r} \in \L_2(H^s(\R))$, we must have $s \ge - \frac 32$, 
showing the essential sharpness 
of the condition $s > - \frac 32$ for $\rho \ge \frac 32$ in~\eqref{regR1}
(modulo the endpoint $s = -\frac 32$).
We point out that by following
\cite[Section~6]{BO97}, 
the argument above (with $r = 0$) yields the failure of $C^2$-smoothness
of the solution map for the modulated KdV \eqref{kdv1} on $H^s(\R)$
for $s < -\frac 32$.

In the argument presented above, 
we considered the case where 
the resonance function 
$\Xi_\KDV (\bar \xi)$
is small; see \eqref{AZ1a}.
In this case, we can not make use
of the irregularity of the modulation function $w$, 
thus reducing the problem to the (unmodulated) KdV case.
This is the reason why we encountered the scaling-critical regularity $s = - \frac 32$
for  the (unmodulated) KdV.

The argument presented above exploits the high-to-low energy transfer, 
leading to the divergence of the Picard second iterate for $s < - \frac 32$.
Such an argument also appears in showing norm inflation
for dispersive equations in negative regularity; see, for example, \cite{OH3, COW}.
By considering
a (suitably weighted) lacunary series of initial data of the form \eqref{AZ0}, 
it may be possible to 
construct a counterexample for  the endpoint case $s = - \frac 32$.
See \cite[Section~4]{COW}
for such an argument.
For simplicity of the presentation, however, we do not pursue this issue here.

\end{remark}

\begin{remark}\label{REM:R1}\rm
In order to prove analogues of Proposition \ref{PROP:kdv1}\,(iii) and (iv)
as a modification of the proof of Proposition \ref{PROP:KR1}, 
we would need to 
(at least) bound
\begin{align*}
A_4' = \sup_{|\xi|\ge 1} 
  \intt_{|\xi - \xi_1|< \frac 12} 
  \jb{\xi}^{ 2(s+ \eps) }\xi^2
   \frac{|\Phi^{w}_{t,r} (\Xi_\KDV(\bar \xi)) |^2}{ \jb{\xi_1}^{2s}  \jb{\xi - \xi_1}^{2s}}  d\xi_1
\end{align*}

\noi
for some $\eps > 0$.
By proceeding as in \eqref{KR5}, we have
 \begin{align*}
A_4'
&\les
\|\Phi^w\|^2_{\W_T^{\rho,\g}}
|t-r |^{2\g }    \sup_{|\xi|\ge 1}  |\xi|^{2 + 2\eps}
  \intt_{|\xi - \xi_1|< \frac 12} 
\frac{ 1}{  \jb{\xi^2(\xi - \xi_1)}^{2\rho}} d \xi_1\\
& = 
\|\Phi^w\|^2_{\W_T^{\rho,\g}}
|t-r |^{2\g }  
    \sup_{|\xi|\ge 1}  |\xi|^{2\eps}
  \intt_{|\xi_2|< \frac {\xi^2}2} 
 \frac{1}{\jb{\xi_2}^{2\rho} }d \xi_2
 = \infty.
\end{align*}

\noi
This shows that 
we can not  establish nonlinear smoothing or convergence
of the Galerkin approximation
as a straightforward modification of the proof of Proposition \ref{PROP:KR1}.

\end{remark}

\section{Modulated BO and ILW equations on the circle}
\label{SEC:BO}

In this section, we consider the
modulated BO  and (scaled) ILW  equations
on the circle.
By establishing regularity properties
of the associated drivers, we first prove 
their local well-posedness and related properties
(Theorems \ref{THM:4} and \ref{THM:G3})
in Subsections \ref{SUBSEC:BO1} and
 \ref{SUBSEC:BO2}, 
 where the claims for the modulated (scaled) ILW follow
 from a perturbation of the BO case.
This perturbative viewpoint
allows us to 
easily prove
deep-water convergence (Theorem \ref{THM:conv}\,(i)) in Subsection~\ref{SUBSEC:BO3}.
On the other hand, 
 careful analysis
on the resonance function $\wt \Xi_\dl(\bar n)$ 
for the scaled ILW defined in \eqref{L11}
is required to establish
shallow-water convergence
(Theorem~\ref{THM:conv}\,(ii)); 
see  Subsection  \ref{SUBSEC:BO4}.

\subsection{Modulated BO equation}
\label{SUBSEC:BO1}

The bilinear driver associated with 
the modulated BO~\eqref{BO} on $\T$ is given by 
\begin{align}
X^\BO_{t,r}( f _1, f _2)
=\int_r^t  \uw  (t')^{-1}
\dx \big( (\uw (t')   f _1)( \uw (t')  f_2)\big) dt',
\label{B1}
\end{align}

\noi
where 
\begin{align}
\uw (t)=e^{w(t) \H\dx^2 }
\label{B1a}
\end{align}

\noi
denotes 
the modulated linear propagator
for \eqref{BO}.
By taking the Fourier transform, we have 
\begin{align}
\F\big( X^\BO_{t,r}  ( f _1, f _2) \big)(n)  =
in  \sum_{ \sub {n_1, n_2 \in \Z^*\\ n = n_1+n_2}}
 \Phi_{t,r}^{w}(\Xi_\BO (\bar n))\ft  f_1 ( n_1 )  \ft   f_2 (n_2)
\label{B2}
\end{align}

\noi
where $\Phi_{t,r}^w$ is as in  \eqref{rho2}
and 
$\Xi_\BO (\bar n)$ denotes the resonance function for BO given by 
\begin{align}
\Xi_\BO (\bar n)
= \Xi_\BO(n,n_1,n_2) = - |n| n +  | n_1 |n_1 +  |n_2| n_2.
\label{B3}
\end{align}

As in the previous section, 
Theorem \ref{THM:4}\,(i) and Theorem \ref{THM:G3}\,(i) and (ii)  follow
from  Proposition~\ref{PROP:main}, 
once we prove the following proposition.

\begin{proposition}

\label{PROP:BO}

Given $\rho \ge1$,  $\frac 12 < \g < 1$, and $T> 0$, 
let  $w$ be $(\rho,\g)$-irregular on $[0, T]$ in the sense of Definition~\ref{DEF:ir}.

\smallskip

\begin{itemize}
\item[(i)]
Suppose that $\rho \ge 1$ and $s \in \R$
satisfy 
\eqref{regBO1}.
Then, 
the driver $X^\BO$ defined in \eqref{B1}
belongs to $ \cX^{s, \g}_2([0, T]\times \T)$
defined in~\eqref{X1}.

\smallskip
\item[(ii)] \textup{(persistence of regularity).}
Suppose that $\rho \ge 1$ and $s \in \R$
satisfy 
\eqref{regBO1}.
Then, 
for any $s_0 > s$, 
the driver $X^\BO$ 
belongs to 
$\cY^{s, s_0, \g}_{2}([0, T]\times \T)$
defined in~\eqref{X2c}.

\smallskip
\item[(iii)]
\textup{(nonlinear smoothing).}
In addition, suppose that $s_0 > s$ satisfies
\eqref{regBO2}.
Then, 
the driver $X^\BO$ 
belongs to 
$\cX^{s, s_0, \g}_{2}([0, T]\times \T)$
defined in~\eqref{X1x}.

\smallskip
\item[(iv)]
\textup{(Galerkin approximation).}
Define the truncated driver $X^{\BO, N}$
as in \eqref{XN}, 
where we replace
 $\uw(t)$ by the modulated linear propagator
for BO defined in \eqref{B1a}.
Suppose that $\rho > 1$ and $s > -\frac 12\rho$.
Then, $X^{\BO, N}$
 converges
to  $X^\BO$ 
in  $ \cX^{s, \g}_2([0, T]\times \T)$ as $N \to \infty$.

\end{itemize}

\end{proposition}

\begin{proof}
As in the proof of Proposition \ref{PROP:kdv1}, 
we verify  the hypotheses
in Lemma \ref{LEM:OBS1} 
with $k = 2$.

\medskip

\noi
(i) 
Proceeding as in \eqref{K4} (see also  \eqref{K5}), 
we have 
\begin{align*}
\|X^\BO_{t,r} \|_{\cL_2 (H^s)}
&\les
A_\BO \| \Phi^w \|_{\W_T^{\rho,\g}} 
  |t-r|^{\g}
\end{align*}

\noi
\noi
for any $0 \le r <  t \le T$, 
where $A_\BO$ is given by
\begin{align*}
A_\BO = \bigg(\sup_{n_1 \in \Z^*}
\sum_{\substack{n, n_2 \in \Z^*\\n_1 = n - n_2}} 
 \jb{n}^{2s + 2}
 \frac{1}{\jb{n_1}^{2s} \jb{n_2}^{2s} |  \Xi_\BO (\bar n)|^{2\rho}}\bigg)^\frac 12.
\end{align*}

\noi
From \eqref{B3}, we have, under $n = n_1 + n_2$,  

\smallskip

\begin{itemize}
\item[(a)]
If $n_1n_2>0$ (i.e.~$n$, $n_1$, and $n_2$ all have the same sign), then  $|\Xi_\BO(\bar n)|=2|n_1n_2|$.

\smallskip

\item[(b)]
If $n, n_1>0>n_2$,  then  $|\Xi_\BO(\bar n)|=2|n n_2|$.

\smallskip

\item[(c)]

If $n, n_2>0>n_1$,  then  $|\Xi_\BO(\bar n)|=2|n n_1|$.

\smallskip

\item[(d)]
If $n_1>0>n, n_2$,  then  $|\Xi_\BO(\bar n)|=2|n n_1|$.

\smallskip

\item[(e)]
If $n_2>0>n, n_1$,  then  $|\Xi_\BO(\bar n)|=2|n n_2|$.

\end{itemize}

\smallskip

\noi
By symmetry, 
we only consider Cases (a) and (b) in the following.

\medskip

\noi
$\bul$ {\bf Case (i.a):} $n_1 n_2 >0$.\\
\indent
Without loss of generality, assume that $|n_1|\ges |n_2|$.
Then, we have $|n_1|\ges |n|$ under $n = n_1 + n_2$.
In view of Lemma \ref{LEM:SUM}, we have 
\begin{align}
\begin{split}
A_\BO^2
& \les 
 \sup_{n_1\in \Z^*}
 \sum_{\substack{n\in \Z^*\\n_1 = n - n_2}}  |n|^{2-2\rho}
\frac{  |n|^{2s+2\rho}  }{  |n_1 n_2|^{2s+2\rho} }\\
& \les \sup_{n_1\in \Z^*}
\sum_{n\in \Z^*} 
\frac{1}{ |n|^{2\rho-2}
|n - n_1|^{2s+2\rho} }
< \infty, 
\end{split}
\label{B6}
\end{align}

\noi
provided that 
\begin{align}
\rho \ge 1, \quad s \ge - \rho, \quad
\text{and}\quad 
s > \frac 32 - 2 \rho.
\label{B5}
\end{align}

\medskip

\noi
$\bul$ {\bf Case (i.b):} $n, n_1>0>n_2$.\\
\indent
In this case, we have $|n_1| > |n_2|$ and thus $|n_1| \ges |n|$.
When  $s \ge 0$, it follows from Lemma~\ref{LEM:SUM} that 
\begin{align}
A_\BO^2\les  
\sup_{n_1\in \Z^*}
\sum_{\substack{n\in \Z^*\\ n_1 = n- n_2}}     \frac{ 1 }{  |n|^{ 2\rho - 2}  | n_2|^{2s+2\rho} }
< \infty, 
\label{B7}
\end{align}

\noi
provided that \eqref{B5} holds.

Next, 
we consider the case $s <  0$.
We separately consider the cases $|n_1| \sim |n_2|$ and $|n_1|\sim |n|$ ($\gg |n_2|$).
When $ |n_1| \sim |n_2|$, 
it follows from Lemma \ref{LEM:SUM} that 
\begin{align}
A_\BO^2\les  
\sup_{n_1\in \Z^*}
\sum_{\substack{n\in \Z^*\\ n_1 = n- n_2}}    
 \frac{ 1 }{  |n|^{ 2\rho -2s - 2}  | n_2|^{4s+2\rho} }
< \infty, 
\label{B7a}
\end{align}

\noi
provided that 
\begin{align}
\rho \ge s + 1, \quad s \ge -\frac 12  \rho, \quad
\text{and}\quad 
s > \frac 32 - 2 \rho.
\label{B7b}
\end{align}

\noi
When $ |n_1| \sim |n|$, 
\eqref{B6} holds under 
 \eqref{B5}.

A simple computation shows that \eqref{B5}
and \eqref{B7b} (with $s < 0$ for the latter)
are satisfied under \eqref{regBO1}.
Hence,
from  Lemma \ref{LEM:OBS1}\,(i),  
we conclude that 
 $X^\BO\in  \cX^{s, \g}_2([0, T]\times \T)$ under~\eqref{regBO1}.

\medskip

\noi
(ii) 
 Write $X^\BO = X^{\BO, 1} + X^{\BO, 2}$, 
where $X^{\BO, 1}$ is the contribution from $|n_1| \ge |n_2|$ in \eqref{B2}
and $X^{\BO, 2} = X^\BO - X^{\BO, 1}$.
We only consider $X^{\BO, 1}$ since a similar bound holds for $X^{\BO, 2}$ by symmetry.
Since $|n_1| \ge |n_2|$, 
we only need to consider the cases (a), (b), and (e) in Part~(i).
Moreover, by symmetry, 
it suffices to consider the cases (a) and (b).

As before, we have 
\begin{align*}
\|X^{\BO, 1}_{t,r} \|_{\cL_{2, 1}^{s, s_0}}
&\les
A^{(1)}_\BO \| \Phi^w \|_{\W_T^{\rho,\g}} 
  |t-r|^{\g},  
\end{align*}

\noi
where $A^{(1)}_\BO$ is given by
\begin{align}
A^{(1)}_\BO = \bigg(\sup_{n_1 \in \Z^*}
\sum_{\substack{n, n_2 \in \Z^*\\n_1 = n - n_2\\|n_1|\ge |n_2|}} 
 \jb{n}^{2{s_0} + 2}
 \frac{1}{\jb{n_1}^{2{s_0}} \jb{n_2}^{2s} |  \Xi_\BO (\bar n)|^{2\rho}}\bigg)^\frac 12.
\label{B8x}
\end{align}

\medskip

\noi
$\bul$ {\bf Case (ii.a):} $n_1 n_2 >0$.\\
\indent
In this case, we can proceed as in Case (i.a) 
and obtain 
$A^{(1)}_\BO < \infty$ under  \eqref{regBO1}.

\medskip

\noi
$\bul$ {\bf Case (ii.b):} $n, n_1>0>n_2$.\\
\indent
When $s_0 \ge 0$ or $|n_1| \sim |n|$, 
we  can proceed as in Case (i.b) 
and obtain 
$A^{(1)}_\BO < \infty$ under~\eqref{B5}.
Thus, it suffices to consider the case 
$s_0 < 0$ and $ |n_1| \sim |n_2|$.
From \eqref{B8x} and Lemma \ref{LEM:SUM}, we have 
\begin{align*}
(A_\BO^{(1)})^2\les  
\sup_{n_1\in \Z^*}
\sum_{\substack{n\in \Z^*\\ n_1 = n- n_2}}    
 \frac{ 1 }{  |n|^{ 2\rho -2s_0 - 2}  | n_2|^{2s_0 + 2s +2\rho} }
< \infty, 
\end{align*}

\noi
provided that 
\begin{align*}
\rho \ge s_0 + 1, \quad s_0 + s  \ge -  \rho, \quad
\text{and}\quad 
s > \frac 32 - 2 \rho
\end{align*}

\noi
which is satisfied under \eqref{regBO1}.
Hence,
from  Lemma \ref{LEM:OBS1}\,(ii),  
we conclude that, if \eqref{regBO1} holds, 
then
 $X^\BO\in  \cY^{s, s_0, \g}_2([0, T]\times \T)$
 for any $s_0 > s$.

\medskip

\noi
(iii) As in Part (i), 
by symmetry, 
we only consider Cases (a) and (b). 
As before, we have 
\begin{align*}
\|X^\BO_{t,r} \|_{\cL_2 ( H^s; H^{s_0})}
&\les
A^{(2)}_\BO \| \Phi^w \|_{\W_T^{\rho,\g}} 
  |t-r|^{\g},  
\end{align*}

\noi
where $A^{(2)}_\BO$ is given by
\begin{align}
A^{(2)}_\BO = \bigg(\sup_{n_1 \in \Z^*}
\sum_{\substack{n, n_2 \in \Z^*\\n_1 = n - n_2}}
 \jb{n}^{2{s_0} + 2}
 \frac{1}{\jb{n_1}^{2{s}} \jb{n_2}^{2s} |  \Xi_\BO (\bar n)|^{2\rho}}\bigg)^\frac 12.
\label{B8a}
\end{align}

\medskip

\noi
$\bul$ {\bf Case (iii.a):} $n_1 n_2 >0$.\\
\indent
Without loss of generality, assume that $|n_1|\ge |n_2|$.
Then, we have $|n_1|\ges |n|$ under $n = n_1 + n_2$.
In view of Lemma \ref{LEM:SUM}, we have 
\begin{align*}
(A_\BO^{(2)})^2
& \les 
 \sup_{n_1\in \Z^*}
 \sum_{\substack{n\in \Z^*\\n_1 = n - n_2}}  |n|^{2s_0 - 2s + 2-2\rho}
\frac{  |n|^{2s+2\rho}  }{  |n_1 n_2|^{2s+2\rho} }\\
& \les \sup_{n_1\in \Z^*}
\sum_{n\in \Z^*} 
\frac{1}{ |n|^{2\rho-2s_0 + 2s-2}
|n - n_1|^{2s+2\rho} }
< \infty, 
\end{align*}

\noi
provided that 
\begin{align}
s_0 \le s + \rho - 1, \quad  s \ge - \rho, \quad
\text{and}\quad 
s_0< 2s + 2\rho -\frac 32.
\label{B9}
\end{align}

\medskip

\noi
$\bul$ {\bf Case (iii.b):} $n, n_1>0>n_2$.\\
\indent
When $s_0 \ge 0$ or $|n_1| \sim |n|$, 
we  can proceed as in Case (iii.a) 
and obtain 
$A^{(2)}_\BO < \infty$ under~\eqref{B5}.
Thus, it suffices to consider the case 
$s_0 < 0$ and $ |n_1| \sim |n_2|$.
From \eqref{B8a} and Lemma \ref{LEM:SUM}, we have 
\begin{align*}
(A_\BO^{(2)})^2\les  
\sup_{n_1\in \Z^*}
\sum_{\substack{n\in \Z^*\\ n_1 = n- n_2}}    
 \frac{ 1 }{  |n|^{ 2\rho -2s_0 - 2}  | n_2|^{4s +2\rho} }
< \infty, 
\end{align*}

\noi
provided that 
\begin{align}
s_0 \le \rho - 1, \quad  s  \ge - \frac 12  \rho, \quad
\text{and}\quad 
s_0< 2s + 2\rho -\frac 32.
\label{B10}
\end{align}

By separately considering the cases $s \ge 0$ and $s < 0$, 
we obtain \eqref{regBO2}
from \eqref{B9} and~\eqref{B10}.
Hence,
from  Lemma \ref{LEM:OBS1}\,(iii),  
we conclude that 
 $X\in  \cX^{s,s_0,  \g}_2([0, T]\times \T)$ under~\eqref{regBO2}.

\medskip

\noi
(iv) 
By proceeding as in the proof of Proposition \ref{PROP:kdv1}\,(iv)
with \eqref{K8}, 
we see that the claim follows
from a slight modification of Part (i), 
provided that 
\begin{align*}
\rho >  1, \quad s > - \rho, \quad
\text{and}\quad 
s > \frac 32 - 2 \rho.
\end{align*}

\noi
and 
\begin{align*}
\rho \ge s + 1, \quad s >  -\frac 12  \rho, \quad
\text{and}\quad 
s > \frac 32 - 2 \rho.
\end{align*}

\noi
where the latter set of the conditions is for $s < 0$.
Note that these conditions are  satisfied 
if  $\rho > 1$ and $s > -\frac 12\rho$.

This concludes the proof of Proposition \ref{PROP:BO}.
\end{proof}

\subsection{Modulated ILW equation}
\label{SUBSEC:BO2}

In this subsection, we consider 
the well-posedness and related issues for 
the modulated ILW~\eqref{ILW1}
and the modulated scaled ILW \eqref{ILW2}
for fixed $0 < \dl < \infty$.
Since these two equations are related by 
the scaling transform \eqref{scaling0}, 
we only consider the modulated ILW \eqref{ILW1}.

Fix $0 < \dl < \infty$. The bilinear driver associated with
the modulated ILW \eqref{ILW1} is given by 
\begin{align}
X^\dl_{t,r}( f _1, f _2)
=\int_r^t  \uw  (t')^{-1}
\dx \big( (\uw (t')   f _1 )(\uw (t')  f_2) \big)dt',
\label{IL1}
\end{align}

\noi
where 
\begin{align}
\uw (t)=e^{ w(t)\Gdl \dx^2}
\label{IL2}
\end{align}

\noi
denotes 
the modulated linear propagator
for \eqref{ILW1}.
By taking the Fourier transform, we have 
\begin{align}
\F\big( X^\dl_{t,r}  ( f _1, f _2) \big)(n)  =
in  \sum_{ \sub {n_1, n_2 \in \Z^*\\ n = n_1+n_2}}
 \Phi_{t,r}^{w}(\Xi_\dl (\bar n))\ft  f_1 ( n_1 )  \ft   f_2 (n_2)
\label{IL3}
\end{align}

\noi
where $\Phi_{t,r}^w$ is as in  \eqref{rho2}
and 
$\Xi_\dl (\bar n)$ denotes the resonance function for ILW given by 
\begin{align}
\Xi_\dl (\bar n)
= \Xi_\dl(n,n_1,n_2) = - p_\dl(n)  +   p_\dl(n_1) +  p_\dl(n_2).
\label{IL4}
\end{align}

\noi
Here, $p_\dl(n)$ denotes the multiplier for $-i \Gdl \dx^2$ given by 
\begin{align}
p_\dl(n)=
n^2 \coth(\dl n) -\frac{n}{\dl} .
\label{IL5}
\end{align}

\noi
See \eqref{GG1}.
Then, under $n = n_1 + n_2$, we have
\begin{align}
\Xi_\dl (\bar n)
= 
- n^2 \coth(\dl n)
+ n_1^2 \coth(\dl n_1)
+ n_2^2 \coth(\dl n_2).
\label{IL6}
\end{align}

\noi
By applying 
\cite[(2.7)]{CLOP}
and 
the fact that  $x^{\al +1} \le C_\al( e^{2x} - 1)$ for any $x \ge  0$
and $\al \ge 0$, 
we have 
\begin{align}
n^2 \coth (\dl n) - |n|n = \sgn(n) \cdot \frac {2n^2} {e^{2|\dl n|} - 1}
= O( \dl^{-2})
\label{IL7}
\end{align}

\noi
for any $n \in \Z^*$.
Hence, from \eqref{B3}, \eqref{IL6} and \eqref{IL7}, 
we have
\begin{align}
\Xi_\dl(\bar n) =  \Xi_\BO(\bar n) + O(\dl^{-2}), 
\label{IL8}
\end{align}

\noi
uniformly in $\bar n = (n, n_1, n_2) \in (\Z^*)^3$.
Then, a slight modification of the proof of Proposition~\ref{PROP:BO}
for the modulated BO
yields the following proposition 
on the driver $X^\dl$ for the modulated ILW, 
from which 
Theorem \ref{THM:4}\,(ii) and Theorem \ref{THM:G3}\,(iii) follow.

\begin{proposition}

\label{PROP:ILW1}

Given $\rho \ge1$,  $\frac 12 < \g < 1$, and $T> 0$, 
let  $w$ be $(\rho,\g)$-irregular on $[0, T]$ in the sense of Definition~\ref{DEF:ir}.
Fix $0 < \dl < \infty$.

\smallskip

\begin{itemize}
\item[(i)]
Suppose that $\rho \ge 1$ and $s \in \R$
satisfy 
\eqref{regBO1}.
Then, 
the driver $X^\dl$ defined in \eqref{IL1}
belongs to $ \cX^{s, \g}_2([0, T]\times \T)$
defined in~\eqref{X1}.

\smallskip
\item[(ii)] \textup{(persistence of regularity).}
Suppose that $\rho \ge 1$ and $s \in \R$
satisfy 
\eqref{regBO1}.
Then, 
for any $s_0 > s$, 
the driver $X^\dl$ 
belongs to 
$\cY^{s, s_0, \g}_{2}([0, T]\times \T)$
defined in~\eqref{X2c}.

\smallskip
\item[(iii)]
\textup{(nonlinear smoothing).}
In addition, suppose that $s_0 > s$ satisfies
\eqref{regBO2}.
Then, 
the driver $X^\dl$ 
belongs to 
$\cX^{s, s_0, \g}_{2}([0, T]\times \T)$
defined in~\eqref{X1x}.

\smallskip
\item[(iv)]
\textup{(Galerkin approximation).}
Define the truncated driver $X^{\dl, N}$
as in \eqref{XN}, 
where we replace
 $\uw(t)$ by the modulated linear propagator
for ILW defined in \eqref{IL2}.
Suppose that $\rho > 1$ and $s > -\frac 12\rho$.
Then, 
the truncated driver $X^{\dl, N}$
 converges
to  $X^\dl$ 
in  $ \cX^{s, \g}_2([0, T]\times \T)$ as $N \to \infty$.

\end{itemize}

\end{proposition}

 \begin{proof}
 We only show how the claim in  (i) 
 follows as a modification of the proof of 
 Proposition~\ref{PROP:BO}\,(i).
 The claims in (ii), (iii), and (iv) follow from 
 similar modifications. 

 Proceeding as in \eqref{K4} (see also  \eqref{K5}), 
we have 
\begin{align}
\|X^\dl_{t,r} \|_{\cL_2 (H^s)}
&\les
A_\dl \| \Phi^w \|_{\W_T^{\rho,\g}} 
  |t-r|^{\g}
\label{IL9}
\end{align}

\noi
\noi
for any $0 \le r <  t \le T$, 
where $A_\dl$ is given by
\begin{align}
A_\dl = \bigg(\sup_{n_1 \in \Z^*}
\sum_{\substack{n, n_2 \in \Z^*\\n_1 = n - n_2}} 
 \jb{n}^{2s + 2}
 \frac{1}{\jb{n_1}^{2s} \jb{n_2}^{2s} \jb{  \Xi_\dl (\bar n)}^{2\rho}}\bigg)^\frac 12.
\label{IL9a}
\end{align}

 From (a)-(e) in the proof of Proposition \ref{PROP:BO}\,(i)
 we have
 \begin{align}
 |\Xi_\BO(\bar n)| \ges \max (|n|, |n_1|, |n_2|)
\label{IL9b}
 \end{align} 
 
 \noi
 for any $\bar n = (n, n_1, n_2) \in (\Z^*)^3$
 with $n = n_1 + n_2$.
 In particular, from \eqref{IL8} and \eqref{IL9b}, 
 we obtain
 \begin{align}
|\Xi_\dl(\bar n)| \ge   |\Xi_\BO(\bar n)| -  O(\dl^{-2})
\sim  |\Xi_\BO(\bar n)|
\label{IL9c}
\end{align}

\noi
 when 
$ \max (|n|, |n_1|, |n_2|) \gg \dl^{-2}$.
Now, write 
\begin{align}
A_\dl = A_\dl^1 + A_\dl^2,
\label{IL9x} 
\end{align}
where
$A_\dl^1$ is the restriction 
of $A_\dl$ to the case 
$ \max (|n|, |n_1|, |n_2|) \gg \dl^{-2}$
and $A_\dl^2 = A_\dl - A_\dl^1$.
 Then, from \eqref{IL9a}, \eqref{IL9c}, 
 and
 the proof of Proposition \ref{PROP:BO}\,(i), 
 we have 
\begin{align}
A_\dl^1 < \infty, 
\label{IL9y}
\end{align}

\noi
 provided that \eqref{regBO1} holds.
 On the other hand, by a crude estimate, 
 we have 
 \begin{align}
A_\dl^2 
& \le  \bigg(\sup_{|n_1|\les \dl^{-2}}
\sum_{|n|\les \dl^{-2} }
 \frac{\jb{n}^{2s + 2}}{\jb{n_1}^{2s} \jb{n - n_1}^{2s} }\bigg)^\frac 12
 \les 1 + \dl^{-6|s| - 3} < \infty.
\label{IL9z}
\end{align}

\noi
Hence,
from  Lemma \ref{LEM:OBS1}\,(i)
with \eqref{IL9}, 
 \eqref{IL9x},  \eqref{IL9y}, and  \eqref{IL9z}, 
we conclude that 
 $X\in  \cX^{s, \g}_2([0, T]\times \T)$ under~\eqref{regBO1}.
 \end{proof}

\begin{remark}\label{REM:shallow1}\rm
It is clear from the proof 
that 
the claims in 
Proposition \ref{PROP:ILW1}
(such as the size of the $\cX^{s, \g}_2(T)$-norm, convergence speed, etc.)
are uniform in $\dl \gg 1$.
However, there is no uniform control
for  $0 < \dl \ll1$, which
presents an additional difficulty in establishing
shallow-water convergence
(as compared to  deep-water convergence);
see Subsection \ref{SUBSEC:BO4}.

\end{remark}

\subsection{Deep-water convergence of the modulated ILW}
\label{SUBSEC:BO3}

In this subsection, we establish
deep-water convergence ($\dl \to \infty$) of the modulated ILW \eqref{ILW1}
to the modulated BO~\eqref{BO}
(Theorem~\ref{THM:conv}\,(i)).

Given $u_0 \in H^s (\T)$,
let $u^\dl$, $0<   \dl < \infty$,  and $u$ denote the local-in-time solutions
to the modulated ILW \eqref{ILW1}
and the modulated BO \eqref{BO}
with $u^\dl|_{t = 0} = u|_{t = 0} = u_0$, respectively.
Let 
$\uu^\dl(t) = e^{- w(t) \Gdl \dx^2} u^\dl(t)$
and 
 $\uu(t) = e^{-  w(t) \H \dx^2} u(t)$
be the modulated interaction representations
of $u^\dl$ and $u$ (under the modulated ILW and the modulated BO, respectively).
Then, by the unitarity of $e^{- w(t) \Gdl \dx^2}$ on $H^s(\T)$, 
we have 
\begin{align}
\| u^\dl - u\|_{C_\tau H^s_x}
\le \| \uu^\dl - \uu\|_{C_\tau H^s_x}
+ \| ( e^{ - w(t) \Gdl \dx^2} -  e^{ -w(t) \H \dx^2}) u \|_{C_\tau H^s_x}
\label{ILD1a}
\end{align}

\noi
From the mean value theorem and \cite[Lemma 2.1]{LOZ}, we have,
for each $n \in \Z^*$, 
\begin{align*}
\big| \F\big( (e^{ - w(t) \Gdl \dx^2} -  e^{-  w(t) \H \dx^2}) f \big)(n)\big|
& \le
 |w(t)|\big| n^2 \ft \Gdl(n) - ( - i |n|n) \big|
|\ft f (n)|\\
& \too0,
\end{align*}

\noi
as $\dl \to \infty$.
Then, from the dominated convergence theorem, 
we see that the second term on the right-hand side
of \eqref{ILD1a} tends to $0$ 
as $\dl \to \infty$.

Hence, 
Theorem~\ref{THM:conv}\,(i)
follows once we show
\begin{align}
\lim_{\dl \to \infty} \| \uu^\dl - \uu\|_{C_\tau H^s_x} = 0.
\label{ILD1b}
\end{align}

\noi
In view of 
Lemma \ref{LEM:OBS1}\,(iv), Remark \ref{REM:OBS2}\,(i), and 
Proposition \ref{PROP:main}\,(iv), 
it then suffices to prove 
\begin{align}
\|X^\dl_{t, r} - X^\BO_{t,r}\|_{\cL_2(H^s)}
\les 
o(1)|t-r|^{\g},
\label{ILD1}
\end{align}

\noi
as $\dl\to \infty$, 
uniformly in  $0 \le r < t \le T$, 
where $X^\dl$ and $X^\BO$
are as in 
\eqref{IL1} and 
\eqref{B1}, respectively.

Proceeding as in \eqref{K4} with \eqref{IL3} and \eqref{B2}
we have 
\begin{align}
\| X_{t, r}^\dl - X^\BO_{t,r} \|_{\cL_2 (H^s)}
&\leq  
\bigg(\sup_{n_1 \in \Z^*}
\sum_{\substack{n, n_2 \in \Z^*\\n_1 = n - n_2}} \jb{n}^{2s+2}
\frac{|\Phi^{w}_{t,r} (\Xi_\dl(\bar n))- \Phi^{w}_{t,r} (\Xi_\BO(\bar n)) |^2}{\jb{n_1}^{2s}\jb{n_2}^{2s}} \bigg)^\frac 12 .
\label{ILD2}
\end{align}

\noi
From 
 \eqref{rho1}  and \eqref{IL8}, we have 
\begin{align}
\begin{split}
|\Phi_{t,r}^w(\Xi_\dl (\bar n)) - \Phi_{t,r}^w(\Xi_\BO(\bar n) )|
& \le |\Phi_{t,r}^w(\Xi_\dl (\bar n))| +| \Phi_{t,r}^w(\Xi_\BO(\bar n) )|\\
& \les \|\Phi^w\|_{  \W^{\rho,\g}_T} 
|t - r|^\g
\frac 1{\jb{\Xi_\BO (\bar n)}^\rho}
\end{split}
\label{ILD3}
\end{align}

\noi
for any $\dl \gg1 $ and 
 $0 \le r < t \le T$.
On the other hand, 
from  \eqref{rho2},  \eqref{IL8}, 
and the mean value theorem, we have 
\begin{align}
\begin{aligned}
|\Phi_{t,r}^w(\Xi_\dl (\bar n)) - \Phi_{t,r}^w(\Xi_\BO(\bar n) )|
&=\bigg|
\int_r^t   e^{i \Xi_{\rm BO}  (\bar n)   w(t') } 
\big(e^{i  O(\dl^{-2})  w(t') }-1 \big)
 d t' \bigg|\\
& \les \dl^{-2}  |t- r|.
\end{aligned}
\label{ILD4}
\end{align}

\noi
By interpolating \eqref{ILD3} and \eqref{ILD4}, 
we have 
\begin{align}
|\Phi_{t,r}^w(\Xi_\dl (\bar n)) - \Phi_{t,r}^w(\Xi_\BO(\bar n) )|
\les_T \dl^{-2\eps} \|\Phi^w\|_{  \W^{\rho,\g}_T}^{1-\eps} 
|t - r|^\g
\frac 1{\jb{\Xi_\BO (\bar n)}^{\rho(1-\eps)}}
\label{ILD5}
\end{align}

\noi
for any $0 \le \eps \le 1$.
Hence, 
by a slight modification of \eqref{B6}, \eqref{B7}, and \eqref{B7a}
(where we replace $\rho$ by $\rho (1-\eps)$ for small $\eps > 0$)
with \eqref{ILD2} and \eqref{ILD5}, 
we obtain 
\begin{align*}
\| X_{t, r}^\dl - X^\BO_{t,r} \|_{\cL_2 (H^s)}
&\les_T \dl^{-2\eps}  
\|\Phi^w\|_{  \W^{\rho,\g}_T}^{1-\eps} 
|t - r|^\g
= o(1) |t - r|^\g, 
\end{align*}

\noi
as $\dl \to \infty$, 
uniformly in  $0 \le r < t \le T$, 
provided that 
\begin{align*}
\rho > 1\quad \text{and}\quad
s > - \frac 12 \rho.
\end{align*}

\noi
This proves \eqref{ILD1} (and \eqref{ILD1b})
and thus concludes the proof of Theorem \ref{THM:conv}\,(i).

\subsection{Shallow-water convergence of  the modulated scaled ILW}
\label{SUBSEC:BO4}

In this subsection,
we establish shallow-water 
convergence ($\dl \to 0$) of the modulated scaled ILW
\eqref{ILW2} to the modulated KdV \eqref{kdv1} 
(Theorem~\ref{THM:conv}\,(ii)).
In the previous subsection, 
we proved the deep-water convergence
 by viewing the modulated ILW
as a perturbation of the modulated BO.
In the shallow-water regime, 
we can not employ this simple perturbative viewpoint
and thus, we need to carry out careful analysis
on the resonance function $\wt \Xi_\dl(\bar n)$ 
for the scaled ILW (see  \eqref{L11}).
In fact, in the low frequency regime, 
we view the modulated scaled ILW
as a perturbation of the modulated KdV, 
while 
we view it 
as a perturbation of the modulated BO
in the high frequency regime.
When there is a mixture of high and low frequencies, 
we proceed with a direct computation, 
exhibiting an extra smoothing property.
In the remaining part of this section, 
we assume that  $0<\dl\ll1$.

Let $\Gd$ be as in \eqref{GG2} (see also \eqref{GG1}).
Then, 
 the Fourier multiplier $\wt p_\dl(n)$ of $-i\Gd \dx^2$ is given by 
 \begin{align}
\wt p_\dl(n)
&  =   \frac 3 \dl p_\dl(n) =  n^3 \big(1- h(n,\dl)  \big), 
\label{L9}
\end{align}

\noi
where $p_\dl(n)$ is as in \eqref{IL5} and $h(n,\dl)$ is given by\footnote{Recall that  
$ \sum_{k=1}^{\infty}   \frac1 {k^2} = \frac{\pi^2}6$.}
\begin{align}
h(n,\dl) = \sum_{k=1}^{\infty}   \frac{6\dl^2n^2}{\pi^2k^2 (\pi^2k^2 +\dl^2n^2)}
\in (0, 1).
\label{L9b}
\end{align}

\noi
See \cite[Lemma 8.2.1]{ABFS}, 
\cite[(the proof of) Lemma 2.3]{LOZ}, 
and \cite[Lemma 2.3]{CLO}.
Fix $N \in \N$.
Then, 
from \eqref{L9b}, we have 
\begin{align}
h(n,\dl)  \les  \dl^2 N^2
\label{L9a}
\end{align}

\noi
for any $|n| \le N$.

Let $\wt X^\dl$ be the bilinear driver associated with
the modulated scaled ILW \eqref{ILW2},  given by 
\begin{align}
\wt X^\dl_{t,r}( f _1, f _2)
=\int_r^t  \uw  (t')^{-1}
\dx \big( (\uw (t')   f _1 )(\uw (t')  f_2) \big)dt',
\label{L9x}
\end{align}

\noi
where  $\uw (t)=e^{ w(t)\Gd \dx^2}$.
By taking the Fourier transform, we have 
\begin{align*}
\F\big( \wt X^\dl_{t,r}  ( f _1, f _2) \big)(n)  =
in  \sum_{ \sub {n_1, n_2 \in \Z^*\\ n = n_1+n_2}}
 \Phi_{t,r}^{w}(\wt \Xi_\dl (\bar n))\ft  f_1 ( n_1 )  \ft   f_2 (n_2)
\end{align*}

\noi
where $\Phi_{t,r}^w$ is as in  \eqref{rho2}
and 
$\wt \Xi_\dl (\bar n)$ denotes the resonance function for the scaled ILW given by 
\begin{align}
\begin{split}
\wt \Xi_\dl (\bar n)
& = \wt \Xi_\dl(n,n_1,n_2) = - \wt p_\dl(n)  +  \wt  p_\dl(n_1) +  \wt p_\dl(n_2)\\
& = - n^3 \big(1- h(n,\dl)  \big)
+ n_1^3 \big(1- h(n_1,\dl)  \big)
+ n_2^3 \big(1- h(n_2,\dl)  \big).
\end{split}
\label{L11}
\end{align}

Given $u_0 \in H^s (\T)$,
let $v^\dl$, $0<   \dl \ll 1$,  and $u$ denote the local-in-time solutions
to the modulated scaled ILW \eqref{ILW2}
and the modulated KdV \eqref{kdv1}
with $v^\dl|_{t = 0} = u|_{t = 0} = u_0$, respectively.
Let 
$\vv^\dl(t) = e^{- w(t) \Gd \dx^2} v^\dl(t)$
and 
 $\uu(t) = e^{ w(t) \dx^3} u(t)$
be the modulated interaction representations
of $v^\dl$ and $u$ (under the modulated scaled ILW and the modulated KdV, respectively).
Then, by the unitarity of $e^{- w(t) \Gd \dx^2}$ on $H^s(\T)$, we have
\begin{align}
\| v^\dl - u\|_{C_\tau H^s_x}
\le \| \vv^\dl - \uu\|_{C_\tau H^s_x}
+ \| ( e^{ - w(t) \Gd \dx^2} -  e^{ w(t)  \dx^3}) u \|_{C_\tau H^s_x}
\label{ILS1}
\end{align}

\noi
From the mean value theorem and \cite[Lemma 2.3]{LOZ}, we have,
for each $n \in \Z^*$, 
\begin{align*}
\big| \F\big( (e^{ - w(t) \Gd \dx^2} -  e^{  w(t)  \dx^3}) f \big)(n)\big|
& \le
 |w(t)|\big| n^2 \ft \Gd(n) - (-i n^3) \big|
|\ft f (n)|\\
& \too0,
\end{align*}

\noi
as $\dl \to 0$.
Then, from the dominated convergence theorem, 
we see that the second term on the right-hand side
of \eqref{ILS1} tends to $0$ 
as $\dl \to 0$.
Hence, 
Theorem~\ref{THM:conv}\,(ii)
follows once we show
\begin{align}
\lim_{\dl \to 0} \| \vv^\dl - \uu\|_{C_\tau H^s_x} = 0.
\label{ILS2}
\end{align}

As in  the case of the deep-water convergence, 
in view of Lemma \ref{LEM:OBS1}\,(iv), Remark \ref{REM:OBS2}, and 
Proposition \ref{PROP:main}\,(iv),
\eqref{ILS2} follows once we  prove
\begin{align}
\|\wt X^{\dl}_{t,r} - X^{\KDV}_{t,r}\|_{\cL_2(H^s)}
\les 
o(1)|t-r|^{\g},
\label{L11a}
\end{align}

\noi
as $\dl\to 0$, 
uniformly in  $0 \le r \le  t \le T$, 
where $X^{\KDV}$
is  the driver for the modulated KdV defined  in \eqref{K1}.

From  \eqref{K3} and \eqref{L11},  we have
\begin{align}
\Dl \wt \Xi_\dl (\bar n): = \wt \Xi_\dl (\bar n) - \Xi_\KDV (\bar n)
 =
n^3 h(n,\dl)  - n_1^3 h(n_1,\dl)   - n_2^3 h(n_2,\dl) .
\label{L11b}
\end{align}

\noi
Then,  from \eqref{L9a}, 
and  \eqref{L11b}, we have 
\begin{align}
|\Dl \wt \Xi_\dl (\bar n)|
\les  \dl^2N^5, 
 \label{L14a}
\end{align}

\noi
uniformly in   $|n|, |n_1|, |n_2|\le N$.
Thus, we have 
\begin{align}
\jb{\wt \Xi_\dl (\bar n)}
\sim \jb{ \Xi_\KDV (\bar n)}, 
 \label{L14aa}
\end{align}

\noi
uniformly in  $|n|, |n_1|, |n_2|\le N$
and $\dl \ll N^{-\frac 52}$.
Then, 
from \eqref{rho1}  and \eqref{L14aa}, we have 
\begin{align}
\begin{split}
|\Phi_{t,r}^{w}( \wt \Xi_\dl (\bar n) )
- \Phi_{t,r}^{w} (\Xi_\KDV (\bar n) )|
& \le 
 |\Phi_{t,r}^{w}( \wt \Xi_\dl (\bar n) )|
+ 
|\Phi_{t,r}^{w} (\Xi_\KDV (\bar n) )|\\
& 
\les \|\Phi^w\|_{  \W^{\rho,\g}_T} |t-r|^\g
\frac{1}{\jb{\Xi_\KDV (\bar n )}^\rho}
\end{split}
\label{L14b}
\end{align}

\noi
\noi
uniformly in  $|n|, |n_1|, |n_2|\le N$, 
 $\dl \ll N^{-\frac 52}$, 
 and 
 $0 \leq  r < t \leq T$.
Moreover, 
from  \eqref{rho2},  \eqref{L14a}, and the mean value theorem, 
we have 
\begin{align}
\begin{aligned}
|\Phi_{t,r}^{w}( \wt \Xi_\dl (\bar n) )
- \Phi_{t,r}^{w} (\Xi_\KDV (\bar n) )|
&=
\bigg| \int_r^t  
e^{i  \Xi_\KDV(\bar n)    w(t') } \big( e^{i   \Dl \wt \Xi_\dl (\cj n)   w(t') } - 1 \big)
 d t' 
 \bigg|
 \\
 &\les 
\dl^2  N^5  |t-r|, 
 \end{aligned}
\label{L12}
\end{align}

\noi
uniformly in   $|n|, |n_1|, |n_2|\le N$.
By interpolating \eqref{L14b} and \eqref{L12}, we have
\begin{align}
\begin{aligned}
& |\Phi_{t,r}^{w}( \wt \Xi_\dl (\bar n) ) - 
\Phi_{t,r}^{w} (\Xi_\KDV (\bar n) )| \cdot
\ind_{ \{ |n|, |n_1|, |n_2|\le N \} } \\
&\quad \les_{ T}
\dl^ {2\eps}  N^{5\eps} \|\Phi^w\|_{  \W^{\rho,\g}_T}^{1-\eps}
 |t-r|^\g \frac{\ind_{ \{ |n|, |n_1|, |n_2|\le N \} }}{\jb{\Xi_\KDV (\bar n )}^{\rho(1-\eps)}}.
 \end{aligned}
 \label{L12a}
\end{align}

\noi
for any $0 \le \eps \le 1$.

Given $N \in \N$, 
let $X^{\KDV,N}_{t,r} $
be  the frequency-truncated driver
defined in \eqref{XN}.
Similarly, 
let $\wt X^{\dl,N}_{t,r} $ 
be 
 the frequency-truncated version of $\wt X^\dl$  in \eqref{L9x}.
Then, proceeding as in the proof of Proposition \ref{PROP:kdv1}\,(i)
(where we replace $\rho$ by $\rho(1-\eps)$
for small $\eps > 0$)
with \eqref{L12a}, 
\begin{align}
\|   \wt X^{\dl,N}_{t,r} -X^{\KDV,N}_{t,r}\|_{\cL_2({H^s})} 
&
\les_{T}
\dl^ {2\eps} 
 N^{5\eps}
\|\Phi^w\|_{  \W^{\rho,\g}_T}^{1-\eps}
 |t-r|^\g = o(1) |t-r|^\g, 
\label{L7}
\end{align}

\noi
as $\dl \to 0$, 
uniformly in  $0 \le r <  t \le T$, 
provided that  $\rho > \frac 12$ and $s \in \R$
satisfy 
\eqref{reg3a}.
From the proof of 
Proposition \ref{PROP:kdv1}\,(iv), we  have 
\begin{align}
\|X^{\KDV, N}_{t,r} -X^\KDV_{t,r}\|_{\cL_2 (H^s)} 
&
\les
 o(1) |t-r|^\g, 
\label{L7x}
\end{align}

\noi
as $N \to \infty$, 
uniformly in  $0 \le r <  t \le T$.
Hence, by the triangle inequality with \eqref{L7}
and~\eqref{L7x}, 
the bound \eqref{L11a}
follows once we prove the following uniform bound:
\begin{align}
\sup_{0 < \dl \ll 1} \|  \wt X^\dl_{t,r} - \wt X^{\dl,N}_{t,r} \|_{\cL_2 (H^s)}
\les
 o(1) |t-r|^\g, 
\label{L7y}
\end{align}

\noi
as $N \to \infty$.

Let us first take a close look at the resonance function 
$\wt \Xi_\dl(\bar n)$ defined in \eqref{L11}.
A direct computation  with \eqref{L9} and \eqref{L9b}
yields
\begin{align}
\wt \Xi_\dl(\bar n)
= -6n n_1 n_2
 \sum_{k=1}^{\infty}   
 \frac{ \pi^2 k^2 \big(3 \pi^2 k^2 + \dl^2(n_1^2 + n_1 n_2 + n_2^2)\big)}
 {\prod_{j = 0}^2(\pi^2 k^2 +\dl^2n_j^2)}
\label{LL1}
\end{align}

\noi
under  $n = n_1 + n_2$, where $n_0 := n$.
Write  $|n|$, $|n_1|$, and $|n_2|$
in the non-increasing order and denote them by 
 $n_{\max}$, 
$n_{\med}$, and $n_{\min}$, respectively.

\medskip

\noi
$\bul$ {\bf Case 1:} $n_{\max} \les \dl^{-1}$.\\
\indent
Denote by 
 $\wt X^{\dl, (1)}$ and $ \wt X^{\dl,N, (1)}$
the contributions to 
$\wt X^\dl$ and $\wt X^{\dl,N}$, respectively, 
from the frequency region $n_{\max} \les \dl^{-1}$. 
In this case, from \eqref{LL1}, we have
\begin{align}
|\wt \Xi_\dl(\bar n)|
\sim  |n n_1 n_2|
 \sum_{k=1}^{\infty}   
 \frac1{\pi^2 k^2}
\sim  |n n_1 n_2|
\sim |\Xi_\KDV(\bar n)|, 
\label{LL2}
\end{align}

\noi
uniformly in $0 < \dl \ll 1$.
Hence, by repeating the proof of Proposition 
\ref{PROP:kdv1}\,(iv) with \eqref{LL2}, we have 
\begin{align}
\sup_{0 < \dl \ll 1} \|  \wt X^{\dl, (1)}_{t,r} - \wt X^{\dl,N, (1)}_{t,r} \|_{\cL_2 (H^s)}
\les
 o(1) |t-r|^\g, 
\label{LL3}
\end{align}

\noi
as $N \to \infty$, 
provided that 
 $\rho > \frac 12$ and $s \in \R$
satisfy 
\eqref{reg3a}.

\medskip

\noi
$\bul$ {\bf Case 2:} $ n_{\min} \gg \dl^{-1}$.\\
\indent
Denote by 
 $\wt X^{\dl, (2)}$ and $ \wt X^{\dl,N, (2)}$
the contributions to 
$\wt X^\dl$ and $\wt X^{\dl,N}$, respectively, 
from the frequency region 
$ n_{\min} \gg \dl^{-1}$.
 From (a)-(e) in the proof of Proposition \ref{PROP:BO}\,(i)
 we have
\begin{align}
|\Xi_\BO(\bar n)|\ges n_{\max}n_{\min} \gg \dl^{-2}.
\label{LL4}
\end{align}

\noi
Then, 
from 
  \eqref{IL8} and  \eqref{LL4}, 
we obtain
\begin{align*}
|\Xi_\dl(\bar n) - \Xi_\BO(\bar n)| = O(\dl^{-2}) \ll 
 |\Xi_\BO(\bar n)|.
\end{align*}

\noi
from which we conclude that 
\begin{align}
|\Xi_\dl(\bar n)|\sim | \Xi_\BO(\bar n)|
\label{LL5a}
\end{align}

\noi
uniformly in $0 < \dl \ll 1$.
Thus, from 
\eqref{LL5a} with 
\eqref{L11},  \eqref{L9}, 
and  \eqref{IL4},
we obtain
\begin{align}
|\wt \Xi_\dl(\bar n)|\sim \dl^{-1} | \Xi_\BO(\bar n)|
\gg  | \Xi_\BO(\bar n)|, 
\label{LL5}
\end{align}

\noi
uniformly in $0 < \dl \ll 1$.
Hence, by repeating the proof of Proposition 
\ref{PROP:BO}\,(iv) with \eqref{LL5}, we have 
\begin{align}
\sup_{0 < \dl \ll 1} \|  \wt X^{\dl, (2)}_{t,r} - \wt X^{\dl,N, (2)}_{t,r} \|_{\cL_2 (H^s)}
\les
 o(1) |t-r|^\g, 
\label{LL6}
\end{align}

\noi
as $N \to \infty$, 
provided that 
 $\rho > 1$ and $s > -\frac 12\rho$.

\medskip

\noi
$\bul$ {\bf Case 3:} $n_{\max} \gg    \dl^{-1} \ges  n_{\min}$.\\
\indent
Denote by 
 $\wt X^{\dl, (3)}$ and $ \wt X^{\dl,N, (3)}$
the contributions to 
$\wt X^\dl$ and $\wt X^{\dl,N}$, respectively, in this case.
Assume that $|n_1|\le |n_2|$.
Since $n_{\med} \sim n_{\max}$ under $n = n_1 + n_2$, 
it follows from 
 \eqref{LL1} with $n_1 n_2 \le \frac 12 n_1^2 + \frac 12 n_2^2$ that 
\begin{align*}
|\wt \Xi_\dl(\bar n)|
\sim |n n_1 n_2|
 \sum_{k=1}^{\infty}   
 \frac{ \pi^2 k^2 }
 {\prod_{j = 0}^1(\pi^2k^2  +\dl^2n_j^2)}
\end{align*}

\noi
Without loss of generality, assume that $|n| = |n_0| \le |n_1|$.
Then, we have $|n_1| \sim n_{\max}
\gg \dl^{-1} \ges n_{\min} = |n|$. 
In particular, we have 
$\pi^2k^2  +\dl^2n^2 \sim \pi^2 k^2$.
Then, by  estimating from below by an integral, we have 
\begin{align}
\begin{split}
|\wt \Xi_\dl(\bar n)|
& \sim \dl^{-1} |n  n_2|
 \sum_{k=1}^{\infty}   
 \frac{ 1}
 {\pi^2(\frac k{\dl n_1})^2  +1} \frac{1}{\dl |n_1|}\\
 & \ges  n_{\min}^2n_{\max}
\int_{0}^\infty
 \frac{ 1}
 {\pi^2x^2 +1}dx\\
& \sim   n_{\min}^2n_{\max}, 
\end{split}
\label{BOX1}
\end{align}

\noi
uniformly in $0 < \dl \ll 1$.
Hence, by repeating the proof of Proposition 
\ref{PROP:BO}\,(iv) with \eqref{BOX1}, we have 
\begin{align}
\sup_{0 < \dl \ll 1} \|  \wt X^{\dl, (3)}_{t,r} - \wt X^{\dl,N, (3)}_{t,r} \|_{\cL_2 (H^s)}
\les
 o(1) |t-r|^\g, 
\label{LL9}
\end{align}

\noi
as $N \to \infty$, 
provided that $\rho > 1$ and $s > -\frac12\rho$.

\medskip

Putting 
\eqref{LL3}, 
\eqref{LL6},  and 
\eqref{LL9} together, 
we obtain \eqref{L7y}, 
provided that 
 $\rho > 1$ and $s > -\frac 12\rho$.
Therefore, under the same condition, 
we obtain
\eqref{L11a}
from the triangle inequality,~\eqref{L7}, \eqref{L7x}, and \eqref{L7y}, 
where we first take $\dl \to 0$ and then $N \to \infty$.
This  concludes the proof of Theorem \ref{THM:conv}\,(ii).

\section{Global well-posedness in 
$L^2$}
\label{SEC:GWP1}

In this section, 
our main goal is to  establish the following proposition on the conservation of the $L^2$-norm
of a solution to the modulated dispersive equations
we considered in Sections~\ref{SEC:LWP} and~\ref{SEC:BO}.

\begin{proposition}
\label{PROP:L2}

Given $\rho > 0$,  $\frac 12 < \g < 1$, and $T> 0$, 
let  $w$ be $(\rho,\g)$-irregular on $[0, T]$ in the sense of Definition~\ref{DEF:ir}.
Given $s \ge 0$, suppose that 
the driver $X$  defined in 
\eqref{K1} belongs to 
   $ \cX^{s, \g}_2([0, T]\times \T)$
   and that 
$u $ is a solution to the modulated KdV equation \eqref{kdv1} on $\T$,  whose
modulated interaction representation $\uu$ defined in \eqref{int1}
belongs to  $\CC^\g([0, T];H^s(\T)) $, 
satisfying 
\begin{align}
\uu  = u_0 + \I^X(\uu).
\label{GK0a}
\end{align}

\noi
Here,  $\I^X(\uu)$ denotes the nonlinear Young integral of $\uu$
with the driver $X$.
Then, we have
\begin{align}
\| u(t) \|_{L^2} = \| u(0) \|_{L^2}
\label{GK0b}
\end{align}

\noi
for any $0 \le t \le T$.

With appropriate minor modifications
on the driver space 
   $X \in \cX^{s, \g}_k([0, T]\times \M)$, 
the  $L^2$-conservation \eqref{GK0b} also  holds
for the modulated KdV equation on $\R$, 
the modulated mKdV, BO, ILW, and scaled ILW equations on $\T$.

\end{proposition}

Once we prove 
 Proposition \ref{PROP:L2}, 
 global well-posedness in $H^s(\M)$, $s \ge 0$, 
holds
 for 
 the following cases:
 
  \smallskip
 \begin{itemize}
 \item[(i)]
 the modulated KdV \eqref{kdv1} on $\M = \T$
for   $\rho >  \frac 12$,

 \smallskip
  \item[(ii)]
 the modulated KdV \eqref{kdv1} on $\M = \R$
for   $\rho >  \frac 12$,

 \smallskip  
   \item[(iii)]
the modulated BO~\eqref{BO}, 
the modulated ILW~\eqref{ILW1}, 
and the modulated scaled ILW~\eqref{ILW2} on $\M = \T$
for  $\rho \ge 1 $
(Theorem \ref{THM:4}\,(ii)).

 \end{itemize}
 
  
  \noi
These global well-posedness results
follow from Proposition \ref{PROP:L2}
together with 
 the persistence-of-regularity results (Proposition~\ref{PROP:main}\,(ii)
with 
Propositions \ref{PROP:kdv1}\,(ii), 
 \ref{PROP:KR1}\,(ii), 
 \ref{PROP:BO}\,(ii), 
and  \ref{PROP:ILW1}\,(ii))
which state that the local existence
time of an $H^s$-solution depends
only on the $L^2$-norm of initial data
(besides the $\cX^{s, \g}_k(T)$-norm of the associated driver $X$).
Since the argument is standard, we omit details.

\begin{proof}[Proof of Proposition \ref{PROP:L2}]

We only consider the periodic modulated KdV case.
The claim for the other equations follows 
from a similar computation.

Let $X$ be as in 
\eqref{K1}.
For simplicity, we set $X_{t, r}(f) = X_{t, r}(f, f)$.
Then, 
from the unitarity of $\uw(t)$ on $L^2(\T)$, 
we have 
\begin{align}
\begin{split}
\jb{f, X_{t,r}(f)  }_{L^2}
&=\int_{r}^{t} 
\int_{\T}  f\cdot 
 \uw(t')^{-1}  
\dx \big((\uw(t') f)^2 \big) 
dx  dt'\\
&=\int_{r}^{t}
 \int_{\T}   \uw(t')  f \cdot 
  \dx \big((\uw(t') f)^2\big)   dxdt'
\\&=\frac 23 \int_{r}^{t} \int_{\T}
  \dx \big((\uw(t') f)^3\big)
   dx dt'=0
\end{split}
\label{GK0}
\end{align}

\noi
for any  $0\leq r<t\leq T$ and $f \in L^2(\T)$.
Consequently,
we obtain 
\begin{align}
\label{GK1}
\begin{aligned}
\| \uu(t) \|^2_{L^2}
&=\|  \uu(r)  \|_{L^2}^2
+2\big\langle \uu(r),X_{t,r}(\uu(r)) 
\big\rangle_{L^2}
 +  \|  \uu(t) -\uu(r)  \|^2_{L^2}+R_{t,r}\\
&=\|  \uu(r)  \|^2_{L^2} +
\|   \uu(t)-\uu(r)   \|^2_{L^2}
+R_{t,r},
\end{aligned}
\end{align}

\noi
where 
\begin{align*}
R_{t,r}=2\big\langle 
\uu(r) ,  \uu(t)-  \uu(r) -X_{t,r}(\uu(r))\big\rangle_{L^2}. 
\end{align*}

\noi
By assumption, we have 
 $\uu\in \CC^\g([0, T];  L^2(\T))$
with $\g>\frac 12$. 
Then, 
from \eqref{GK0a} and~\eqref{Y2} in Proposition \ref{PROP:young1}
(see also \eqref{Ja3x}), 
we have 
\begin{align}
\begin{split}
|R_{t,r}| 
&\leq 2 
\|\uu(r)\|_{L^2}
\|\uu(t)- \uu(r) -X_{t,r}(\uu(r))\|_{L^2}
\\
&= 2
\|\uu(r)\|_{L^2}
\|(\updl\I^X(\uu))_{t, r}  -X_{t,r}(\uu(r))\|_{L^2}
\\
&\le
C\big(\|\uu\|_{\CC^\g_{T} L^2_x}\big)  \|X\|_{\cX^{s, \g}_k(T)} |t-r|^{2\g}.
\end{split}
\label{GK2}
\end{align}

\noi 
In addition, we have 
\begin{align}
\|\uu(t)-\uu(r)\|_{L^2}^2 
\leq 
\|\uu\|_{C^\g_T L^2_x}^2 
|t-r|^{2\g}
\label{GK3}
\end{align}

\noi
for any  $0\leq r<t\leq T$.
Hence, from \eqref{GK1}, \eqref{GK2}, and \eqref{GK3}, 
we obtain 
\begin{align}
\big| 
 \|  \uu(t)  \|^2_{L^2}- 
\| \uu(r)  \|^2_{L^2}
 \big|
 \lesssim
 C\big(\|\uu\|_{\CC^\g_{T} L^2_x},   \|X\|_{\cX^{s, \g}_k(T)}\big)
 |t-r |^{2\g}.
\label{GK4}
\end{align}

\noi
By recalling $2\g > 1$, 
the $L^2$-conservation \eqref{GK0b}
then follows from \eqref{GK4}.
\end{proof}

\section{Global well-posedness of the modulated KdV in negative Sobolev spaces}
\label{SEC:I}

In this section, we prove global well-posedness of 
the periodic modulated KdV \eqref{kdv1}  in negative Sobolev spaces 
(Theorem \ref{THM:1}\,(ii)).
Our strategy is to adapt 
the so-called $I$-method ($=$ the method of almost conservation laws), 
introduced by Colliander, Keel, Staffilani, Takaoka, and Tao
\cite{CKSTT03}, to the current modulated setting.

Fix   $s<0$
in the remaining part of this section.
 Given $N \geq 1$,
 we define a smooth, even, 
 non-increasing (in $|\xi|$) 
 function 
 $m_{s, N} \in C^\infty(\R; [0, 1])$ 
 by setting
\begin{align}
m_{s, N}(\xi)=
\begin{cases}
1, 
&\text{for }
|\xi|\leq N,   \\
\frac{ |\xi|^s }{ N^s }, 
&\text{for }
|\xi|\geq  2N.
\end{cases}
\label{Iop1}
\end{align}

\noi
Then, we define the so-called $I$-operator $I = I_{s, N}$
to be the Fourier multiplier operator with the multiplier $m_{s, N}$:
\begin{align}
\ft{I_{s, N}f}(\xi)=m_{s, N}(\xi)\ft{f}(\xi).
\label{Ix2}
\end{align}

\noi
Note that 
$I_{s, N}$ acts as the {\bi i}\hspace{0.25mm}dentity operator on low frequencies $\{ |\xi| \leq N\}$, 
while it acts as a (fractional)  
{{\bi i}\hspace{0.25mm}\nolinebreak ntegration}   operator 
(recall  $s < 0$)
 on high frequencies $\{ |\xi| \geq 2N\}$
 with the following bound:
\begin{align}
\| f \|_{H^s (\T) } \leq\| I f \|_{L^2(\T)}
\les N^{-s} \| f \|_{H^s(\T)}, 
\label{I1}
\end{align}

\noi
 and hence the name $\pmb{I}$-method.

We proved local well-posedness
of 
the modulated KdV \eqref{kdv1} on $\T$
by applying the abstract subcritical  local well-posedness
result (Proposition \ref{PROP:main}).
Thus, it follows from the blowup alternative \eqref{BA1}
in Remark \ref{REM:main2}\,(i)
that 
global well-posedness of the periodic modulated KdV follows
once we show that the $H^s$-norm 
of a solution $u(t)$ remains finite 
for each finite time $t > 0$.
In view of \eqref{I1}, 
it then suffices to control the so-called
modified energy $\| Iu(t)\|_{L^2}^2$.\footnote{As we see below, given a target time $T> 0$, 
we in fact apply a scaling first and study the scaled modified energy in \eqref{ME1}.}

In Subsection \ref{SUBSEC:I1}, 
we go over basic scaling properties
of the periodic modulated KdV.
In Subsection~\ref{SUBSEC:I2}, 
we then study  local well-posedness
of the so-called (scaled) modulated $I$-KdV equation~\eqref{kdv4}.
After establishing a crucial commutator estimate
(Proposition \ref{PROP:com}) in Subsection~\ref{SUBSEC:I3}, 
we establish global well-posedness
by combining the commutator estimate (Proposition \ref{PROP:com})
with the sewing lemma (Lemma \ref{LEM:sew})
in Subsection~\ref{SUBSEC:I4}.
It is worthwhile to note that we need to make use
of the sewing lemma even in the globalization argument via the $I$-method
which is a new feature in the current modulated setting.
We also discuss the real line case in Subsection~\ref{SUBSEC:I5}.

\subsection{Scaling properties of the modulated KdV}
\label{SUBSEC:I1}

In this subsection, we discuss
the scaling properties
of the modulated KdV \eqref{kdv1}
on $\T$
which will be used to make the problem at hand
into a small data problem; see \eqref{scaling3},  \eqref{IT1}, and \eqref{IT2a}.
We first recall the scaling invariance
of the (unmodulated) KdV \eqref{kdv2};
a function  $u$ is a solution to KdV \eqref{kdv2} on  $[0, T]\times \T$
 if and only if  
  the scaled function 
\begin{align}
 u^{\ld}(t,x)=\ld^{-2}u(\ld^{-3}t,\ld^{-1}x)
 \label{scaling1}
\end{align}
 
 \noi
  is  a solution to   KdV  on $[0,\ld^3T]\times \T_\ld$,
where
  $\T_\ld = \R/(2\pi \ld \Z)$
denotes the dilated circle.
We note that, unlike the real line setting, 
the spatial domain is dilated under the scaling \eqref{scaling1}.

In the current modulated setting, 
we formally see that 
  $u$ is a solution to  the modulated KdV \eqref{kdv1} on $[0, T] \times \T$ with initial data $u_0$
  if and only if 
$u^\ld$ in \eqref{scaling1}
is a solution to the 
following scaled modulated KdV  on $[0, \ld^3 T] \times \T_\ld$:
\begin{align}
\dt u^{\ld} +    \dx^3 u^{\ld}\cdot \dt w^{\ld}  =\partial_x(u^{\ld} )^2
\label{kdv3}
\end{align}

\noi
with the scaled initial data
\begin{align}
u_0^\ld(x) = \ld^{-2} u_0(\ld^{-1}x), 
\label{scaling1a}
\end{align}

\noi
where the scaled modulation $w^\ld$ is given by
\begin{align}
w^{\ld}(t)= \ld^3w(\ld^{-3}t)
\label{scaling2} 
\end{align}

\noi
such that $\dt w^\ld(t) = \dt w(\ld^{-3} t)$.
Namely, the scaling \eqref{scaling2} on the modulation function
is $\dot W^{1, \infty}$-invariant.
See Lemma \ref{LEM:tri3}
for  equivalence of the unscaled and scaled modulated KdV equations.
In the following, we study the scaled modulated KdV \eqref{kdv3} on 
the dilated circle $\T_\ld$ with $\ld \ge 1$.

Before proceeding further, we 
introduce some notations.
Set
 $\Z_{\ld} = \Z/\ld$
and 
$\Z^*_{\ld} = \Z^*/\ld$. 
We define 
 the Fourier transform of a function $f$ on $\T_\ld$
by 
$$
\ft  f(n)=\int_{\T_\ld}f(x)e^{-i n x} \frac{d x}{2\pi}
$$

\noi
for $n\in\Z_{\ld}$.
Then, we have 
\begin{align}
\begin{split}
f(x)
& =
\frac{1}{\ld} \sum_{ n \in\Z_{\ld}}\ft  f(n)e^{i  nx}, \\
\int_{\T_\ld}|f(x)|^2 \frac{d x}{2\pi}
& =\frac{1}{\ld}\sum_{n\in\Z_{\ld}}|\ft  f(n)|^2,\\
\ft{fg}(n)
& =\frac{1}{\ld}
 \sum_{ \substack{n_1, n_2 \in \Z_\ld\\n = n_1+n_2}}
 \ft  f(n_1)\ft  g(n_2).
\end{split}
\label{FT3}
\end{align}

\noi
Given $s \in \R$, 
we define the non-homogeneous and homogeneous Sobolev spaces $H^{s}(\T_\ld)$ 
and $\dot H^{s}(\T_\ld)$ 
via the (semi-)norms: 
\begin{align}
\begin{split}
\|  f  \|^2_{H^{s}(\T_\ld)}
& =\frac{1}{\ld}
\sum_{n \in \Z_{\ld}}
\jb{n}^{2s} | \ft f(n) |^2, \\
\|  f  \|^2_{\dot H^{s}(\T_\ld)}
& =\frac{1}{\ld}
\sum_{n \in \Z_{\ld}}
|n|^{2s} | \ft f(n) |^2.
\end{split}
\label{FT4}
\end{align}

\noi
When $\ld = 1$, 
the $H^s(\T)$-norm and $\dot H^s(\T)$-norm
are equivalent in the current mean-zero setting.
When $\ld \gg1$, 
we still have 
$\|f \|_{H^s(\T_\ld)}\sim_\ld \|f\|_{\dot H^s(\T_\ld)}$
for a mean-zero function $f$ on $\T_\ld$
but the implicit constant diverges as $\ld \to \infty$.
In the following, 
we will  work with the homogeneous Sobolev spaces, 
since they enjoy a better scaling property; 
 noting from \eqref{scaling1a} that 
$\F_{\T_\ld}(u_0^\ld)(n) = \ld^{-1} \F_\T(u_0)(\ld n)$, $n \in \Z_\ld$, 
we have
\begin{align}
\|u_0^\ld\|_{\dot H^s(\T_\ld)}
=  \ld^{-\frac 32 - s} \|u_0\|_{\dot H^s(\T)}
\label{scaling3}
\end{align}

\noi
for any $s \in \R$, 
while, for the non-homogeneous Sobolev space,  we have 
\begin{align*}
\|u_0^\ld\|_{ H^s(\T_\ld)}
\le  \ld^{-\frac 32 - s} \|u_0\|_{ H^s(\T)}
\end{align*}

\noi
for $s \le 0$; see also Remark \ref{REM:scaling1}.
In particular, when $s > -\frac 32$ (namely, in the scaling subcritical regime),
we can make
the $\dot H^s(\T_\ld)$-norm of the scaled initial data $u_0^\ld$
small by choosing $\ld \gg_{\|u_0\|_{H^s(\T)}} 1$.

Now, we are ready to study the scaled modulated KdV \eqref{kdv3} on $\T_\ld$
with the scaled initial data $u_0^\ld$ in \eqref{scaling1a}.
The Duhamel formulation 
for the modulated interaction representation 
\begin{align}
\uu^\ld(t) = U^{w^\ld}(t)^{-1} u^\ld(t)
\label{ME0}
\end{align}

\noi
is given by 
\begin{align}
\uu^{\ld} (t)=u^\ld_0
+
 \int_0^t   U^{w^\ld}(t')^{-1}  \dx\big(  ( U^{w^\ld} (t') \uu^\ld(t'))^2   \big)dt' 
 \label{Dul3}
\end{align}

\noi
with the associated bilinear driver $X^\ld$ 
given by 
\begin{align}
X_{t,r}^{\ld} (f_1, f_2)
=\int_{r}^{t}
U^{w^\ld}(t')^{-1}
\dx \big( (U^{w^\ld}(t') f_1  )(
U^{w^\ld}(t') f_2)\big) dt'
\label{Xld1}
\end{align}

\noi
for functions $f_1$ and $  f_2$ on $\T_\ld$.
Then, local well-posedness of the scaled modulated KdV~\eqref{kdv3}
in $\dot H^s(\T_\ld)$
follows from Proposition \ref{PROP:main}
and the following lemma.

\begin{lemma}
\label{LEM:tri2}
Let $\ld \ge 1$.
Given $\rho \ge\frac 12$, $\frac 12 < \g < 1$, and $T> 0$, 
let  $w$ be $(\rho,\g)$-irregular on $[0, T]$ in the sense of Definition~\ref{DEF:ir}, 
and let $X^\ld$ be as in \eqref{Xld1}.
Suppose that $\rho \ge \frac 12$ and $s \in \R$
satisfy 
\eqref{reg1}.
Then, 
the driver $X^\ld$  in \eqref{Xld1}
belongs to $ \dot \cX^{s, \g}_2([0, \ld^3T] \times \T_\ld)$ defined in~\eqref{X1a}, 
satisfying the bound\textup{:}
 \begin{align}
 \| X^{\ld}_{t,r}  \|_{\cL_2 ( \dot H^s(\T_\ld)) }
 \les
 \ld^{ \frac32 +s -3\g} 
  \|\Phi^{w}\|_{\W_{T}^{\rho, \g}}
 |t-r|^{\g}
\label{D1}
\end{align}

\noi
for any $0\leq r<t \leq \ld^3T $.
%

\end{lemma}

\begin{proof}

Proceeding as in 
\eqref{K4} with \eqref{FT3} and \eqref{FT4}, we have 
\begin{align}
\begin{aligned}
& \|X^{\ld}_{t,r}(f_1,f_2) \|^2_{\dot H^s(\T_\ld)}
=  \ld^{-3} 
\sum_{n\in\Z_\ld^*} |n|^{2s+ 2}
\bigg |  \sum_{ \substack{n_1, n_2 \in \Z_\ld^*\\n = n_1+n_2}}
\Phi^{w^\ld}_{t,r}(\Xi_\KDV(\bar n) ) 
\ft f_1(n_1)  \ft f_2 (n_2)  \bigg|^2 \\
&\quad \leq \ld^{-1}
\bigg(\sup_{n_1 \in \Z_\ld^*}
\sum_{\substack{n, n_2 \in \Z_\ld^*\\n_1 = n - n_2}} 
|n|^{2s+ 2}
\frac{|\Phi^{w^\ld}_{t,r} (\Xi_\KDV(\bar n)) |^2}{|n_1|^{2s}|n_2|^{2s}} \bigg)
\|f_1\|_{\dot H^s(\T_\ld)}^2
\|f_2\|_{\dot H^s(\T_\ld)}^2, 
\end{aligned}
\label{D2}
\end{align}

\noi
where $\Phi^{w^\ld}_{t,r}$ and
$\Xi_\KDV(\bar n)$ 
are 
 as in \eqref{rho2} and \eqref{K3}, 
 respectively.
From \eqref{rho2}, 
\eqref{scaling2}, 
\eqref{rho1}, and~\eqref{K3a}, 
we have
\begin{align}
\begin{aligned}
|\Phi^{w^\ld}_{t,r} (\Xi_\KDV(\bar n)) |
& = \bigg| 
\int_{r}^{t}  e^{i  \Xi_\KDV (\bar n) w^{\ld}(t') } d t'
\bigg|
= \ld^3   \bigg|  \int_{\ld^{-3}r}^{\ld^{-3}t} 
e^{i  \ld^3  \Xi_\KDV(\bar n)  w(t') } dt'  \bigg|  \\
& \les 
\ld^{3- 3\rho -3\g} 
 \|\Phi^{w}  \|_{\W_{T}^{\rho, \g}}
 |t-r|^{\g} 
 |n n_1 n_2|^{-\rho}
\end{aligned}
\label{D3}
\end{align}

\noi
 for any $n, n_1, n_2\in \Z_\ld^*$,  satisfying $n = n_1+ n_2$.

Given $n, n_1, n_2 \in \Z_\ld^*$, 
let  $m  = \ld n$ and $m_j = \ld n_j$, $j = 1, 2$.
 Then, by 
proceeding as in the proof of Proposition \ref{PROP:kdv1}\,(i)
with 
\eqref{D2}
and  \eqref{D3}, 
 we have 
\begin{align*}
\|X_{t,r}^\ld \|_{\cL_2 (\dot H^s(\T_\ld))}
&\les
\ld^{\frac 52-3\rho -3\g}   \| \Phi^w \|_{\W_T^{\rho,\g}} 
  |t-r|^{\g} \\
  & \quad
  \times 
\bigg(\sup_{n_1 \in \Z_\ld^*}
\sum_{\substack{n, n_2 \in \Z_\ld^*\\n_1 = n - n_2}} 
|n|^{2- 4\rho}
 \frac{|n|^{2s+ 2\rho}}{|n_1|^{2s+ 2\rho} |n_2|^{2s + 2\rho}}\bigg)^\frac 12 \\
& \les 
 \ld^{\frac{3}{2}+s-3\g} 
 \| \Phi^w \|_{\W_T^{\rho,\g}} 
  |t-r|^{\g}\\
& \quad \times \bigg( \sup_{m_1 \in \Z^*}
\sum_{\substack{m, m_2 \in \Z^*\\m_1 = m - m_2}} 
|m|^{2- 4\rho}
 \frac{|m|^{2s+ 2\rho}}{|m_1|^{2s+ 2\rho} |m_2|^{2s+ 2\rho}}\bigg)^\frac{1}{2}\\
 &\les
 \ld^{\frac{3}{2}+s-3\g} 
 \| \Phi^w \|_{\W_T^{\rho,\g}} 
  |t-r|^{\g} 
\end{align*}

\noi
for any $0 \le r < t \le \ld^3 T$, 
provided that
\eqref{reg1} holds.
Hence,
from  Lemma \ref{LEM:OBS1}\,(i)
 and Remark~\ref{REM:OBS2}\,(ii), 
we conclude that 
 $X^\ld \in  \dot \cX^{s, \g}_2([0, \ld^3 T]\times \T_\ld)$ under~\eqref{reg1}.
\end{proof}

\begin{remark}\label{REM:scaling1}\rm

By repeating the proof for the non-homogeneous Sobolev space $H^s(\T_\ld)$
with \begin{align*}
\ld^{-1} \jb{m} \les\jb{n} \les \jb{m}, 
\end{align*}

\noi
we obtain
 \begin{align*}
 \| X^{\ld}_{t,r}  \|_{\cL_2 (  H^s(\T_\ld)) }
 \les
 \ld^{ \frac32 -s -3\g} 
  \|\Phi^{w}\|_{\W_{T}^{\rho, \g}}
 |t-r|^{\g}, 
\end{align*}

\noi
for $\ld \ge 1$, 
which gives a worse power of $\ld$ as compared to \eqref{D1},
yielding a worse regularity threshold
for global well-posedness on the  unscaled modulated KdV \eqref{kdv1}.

\end{remark}

The next lemma shows equivalence of the 
unscaled and scaled problems.
See Definition~\ref{DEF:con1}
for the definition of the class 
$\dot \cD_{w}^s(I \times \M)$.

\begin{lemma}
\label{LEM:tri3}

Let $\ld \ge 1$.
Given $\rho \ge\frac 12$,  $\frac 12 < \g < 1$, and $T> 0$, 
let  $w$ be $(\rho,\g)$-irregular on $[0, T]$ in the sense of Definition~\ref{DEF:ir}.
Fix $s \in \R$ satisfying \eqref{reg1}
and $u_0 \in H^s(\T)$.
Then,  $u\in \dot \cD_{w}^s( [0, \tau] \times \T)$ is a solution to the modulated KdV  \eqref{kdv1} 
on $[0, \tau ]\times \T$ 
with  $u|_{t = 0} = u_0$ 
for some $0 < \tau \le T$
if and only if the scaled function $u^\ld$ defined in \eqref{scaling1}
belongs to  $\dot \cD_{w^\ld}^s( [0, \ld^3\tau] \times \T_\ld)$
and is a   solution to the scaled modulated KdV~\eqref{kdv3}
on $[0, \ld^3\tau ]\times \T_\ld$ 
with the scaled initial data 
$u^\ld|_{t = 0} = u_0^{\ld}$ defined in~\eqref{scaling1a}.

\end{lemma}

\begin{proof}

We only show the forward direction.
Suppose that 
 $u\in \dot \cD_{w}^s( [0, \tau] \times \T) $ 
 is a solution to  the modulated KdV \eqref{kdv1}
 on $[0, \tau]\times \T$ 
with  $u|_{t = 0} = u_0$. 
From Definition \ref{DEF:con1}, 
the modulated interaction representation
$\uu(t)=  \uw(t)^{-1} u(t)$
belongs to $\CC^{\g}( [0,\tau]; \dot H^{s}(\T))$.
Then, 
a direct computation with 
 \eqref{scaling1} shows that 
 \[
 \uu^{\ld}(t)=
 U^{w^{\ld}} (t)^{-1}  u^{\ld} (t) \in \CC^\g([0,\ld^3\tau]; \dot H^{s}(\T_\ld)).
 \] 
 
\noi
Namely, 
$u^\ld \in \dot \cD_{w^\ld}^s( [0, \ld^3\tau] \times \T_\ld)$.

Now, we need to check that $u^\ld$ is indeed a solution 
to the scaled modulated KdV \eqref{kdv3}
on the time interval $[0, \ld^3\tau]$
in the sense that its modulated interaction representation 
$\uu^\ld$ satisfies~\eqref{Dul3}, 
where we interpret  the integral as a nonlinear Young integral.

Let  $f_1$ and $f_2$ be functions on $\T_\ld$.
Then, from 
\eqref{Xld1}
with 
\eqref{FT3}, 
\eqref{rho2}, 
\eqref{scaling2}, 
a change of variables, 
\eqref{K3}, 
and 
\eqref{K2}, we have
\begin{align}
\begin{aligned}
& \F_{\T_\ld}\big( X^{\ld}_{t,r} ( f _1,  f _2)\big)(n)
= \frac{in}{\ld} \sum_{ \sub {n_1,n_2 \in\Z^*_{\ld} \\ n=n_1+n_2}}
 \Phi^{w^\ld}_{t,r}(\Xi_\KDV (\bar n))
\ft  f _1(n_1) \ft  f _2(n_2) \\
&\quad = in \ld^{2} \sum_{\sub{m_1,m_2 \in \Z^*\\ \ld n = m_1+ m_2}}
 \Phi^{w}_{\ld^{-3}t, \ld^{-3}r}(\Xi_\KDV (\ld n, m_1, m_2))
\ft  f _1 (\ld^{-1} m_1)
\ft  f _2 (\ld^{-1} m_2)
 \\
&\quad = 
\ld^{-1}  \F_\T\big( X_{\ld^{-3}t, \ld^{-3}r} ( f _1^{\ld^{-1}},   f _2^{\ld^{-1}})\big) (\ld n)
\end{aligned}
\label{D6}
\end{align}

\noi
for any 
$n\in\Z_{\ld}^*$
and $0\leq r<t\leq \ld^3\tau$, 
where $m_j = \ld n_j$ and $f _j^{\ld^{-1}} (x) = \ld^2 f_j(\ld x)$ is a function on $\T$,  $j = 1, 2$.
Thus,  from \eqref{D6} with \eqref{scaling1}, we have 
\begin{align}
 X_{t,r}^{\ld} (\uu^{\ld}(r), \uu^{\ld}(r) )(x) 
= \ld^{-2} X_{\ld^{-3}t,\ld^{-3}r} ( \uu(\ld^{-3}r), \uu(\ld^{-3}r) )(\ld^{-1}x)
\label{D6a}
\end{align}

\noi
for $x \in \T_\ld$.
Since $\uu$ is a solution to \eqref{Dul3} 
on the time interval $[0, \tau]$
with $\ld = 1$, 
where we interpret
 the integral as a nonlinear Young integral, 
it follows from \eqref{Y2} in Proposition~\ref{PROP:young1} 
(adapted to the homogeneous Sobolev space) that 
\begin{align}
\| \uu(\ld^{-3}t) - \uu(\ld^{-3}r) 
- X_{\ld^{-3}t, \ld^{-3}r} ( \uu(\ld^{-3}r ) , \uu(\ld^{-3}r))  \|_{\dot H^s(\T)}
\les  |t-r|^{2\g}
\label{D6b}
\end{align}

\noi
for any $0\leq \ld^{-3}r<\ld^{-3}t\leq \tau$.
Then, from \eqref{D6b} with 
\eqref{scaling1}, \eqref{D6a}, and 
\eqref{scaling3}, we have
\begin{align*}
\|\uu^\ld(t) - \uu^\ld(r) - X^\ld_{t,r} (\uu^\ld(r), \uu^\ld (r))\|_{\dot H^s(\T_{\ld})}
\les \ld^{-\frac{3}{2}-s} |t-r|^{2\g}
\end{align*}

\noi
for any $0 \le r < t \le \ld^3 \tau$.
Let  $\Ta_{t, r} = X^\ld_{t,r} (\uu^\ld(r), \uu^\ld (r))$.
Noting that  $2\g > 1$, 
we then obtain 
\begin{align*}
\uu^\ld (t) - u_0^\ld  = \lim_{|\Pi([0,t])|\to 0} 
\sum^{n-1}_{j=0} \Theta_{t_j,t_{j+1}}
= \I^{X^\ld}(\uu^\ld)(t) 
\end{align*} 
 
 \noi
for any $0 < t \le \ld^3 \tau$,   where 
 the limit is in the sense of Lemma \ref{LEM:sew}\,(iii)
 and the second equality follows from the uniqueness
 of the nonlinear Young integral with the driver $X^\ld$.
This proves that 
 $\uu^\ld$ is a solution to \eqref{Dul3}
on $[0, \ld^3\tau]$.
\end{proof}

\subsection{Modulated $I$-KdV equation}
\label{SUBSEC:I2}

In this subsection, we go over  local well-posedness
of the following scaled modulated $I$-KdV
on the dilated circle $\T_\ld$:
\begin{align}
\dt I u^{\ld} +    \dx^3 I u^{\ld}
\cdot \dt  w^{\ld}
  =\partial_x I\big((u^{\ld} )^2\big). 
\label{kdv4}
\end{align}

\noi
By  writing
$\partial_x I\big((u^{\ld} )^2\big) = \partial_x I\big((I^{-1}Iu^{\ld} )(I^{-1}Iu^{\ld} )\big)$, 
we define the associated bilinear driver $Y^\ld$ by setting
\begin{align}
Y^\ld (f_1, f_2) = I X^{\ld} (I^{-1} f _1,  I^{-1} f _2)
\label{Xld2}
\end{align}

\noi
for functions $f_1$ and $f_2$ on $\T_\ld$, 
where $X^\ld$ is as in \eqref{Xld1}.

\begin{lemma}
\label{LEM:tri4}
Let $\ld \ge 1$.
Given $\rho \ge\frac 12$,  $\frac 12 < \g < 1$,  and $T> 0$, 
let  $w$ be $(\rho,\g)$-irregular on $[0, T]$ in the sense of Definition~\ref{DEF:ir}, 
and let $Y^\ld$ be as in \eqref{Xld2}.
Then, for any $s < 0$
satisfying~\eqref{reg1}, 
the driver $Y^\ld$  in \eqref{Xld2}
belongs to $ \cX^{0, \g}_2([0, \ld^3T] \times \T_\ld)$ defined in~\eqref{X1}, 
satisfying the bound\textup{:}
\begin{align}
 \| Y^{\ld}_{t,r}  \|_{\cL_2 ( L^2(\T_\ld)) }
 \les
 \ld^{ \frac32  -3\g} 
  \|\Phi^{w}\|_{\W_{T}^{\rho, \g}}
 |t-r|^{\g}
\label{IF1}
\end{align}

\noi
for any $0\leq r<t \leq \ld^3T$.

\end{lemma}

\begin{proof}

In view of Lemma  \ref{LEM:OBS1}\,(i), 
it suffices to prove \eqref{IF1}.
Moreover, by the interpolation lemma (see \cite[Lemma 12.1]{CKSTT04}),
 it suffices to prove \eqref{IF1} for $N = 1$. 

Proceeding as in \eqref{D2} with \eqref{Xld2}, we have 
 \begin{align}
\begin{aligned}
\|  & Y^{\ld}_{t,r}( f _1, f _2)  \|^2_{L^2(\T_\ld)}\\
&=
\ld^{-3} \sum_{ n\in\Z_\ld^*}
|n|^2 \bigg|  \sum_{ \substack{n_1, n_2 \in \Z_\ld^*\\n = n_1+n_2}}
 M (\bar n)
\Phi^{w^\ld}_{t,r}(\Xi_\KDV(\bar n) ) 
\ft f _1(n_1)
\ft  f _2 (n_2)
\bigg|^2 \\
& \leq \ld^{-1}
\bigg(\sup_{n_1 \in \Z_\ld^*}
\sum_{\substack{n, n_2 \in \Z_\ld^*\\n_1 = n - n_2}} 
|n|^2 | M (\bar n)|^2
|\Phi^{w^\ld}_{t,r} (\Xi_\KDV(\bar n)) |^2 \bigg)
\|f_1\|_{L^2(\T_\ld)}^2
\|f_2\|_{L^2(\T_\ld)}^2, 
\end{aligned}
\label{IF2}
\end{align}

\noi
where $ M (\bar n) =  M(n,n_1, n_2)$ is defined by 
\begin{align}
M (\bar n) =  M(n,n_1, n_2)=
\frac{m_{s, 1}(n) }
{m_{s, 1}(n_1) m_{s, 1}(n_2)}.
\label{IF3}
\end{align}

\noi
Then, 
by comparing \eqref{IF2} with \eqref{D2}, 
the bound 
\eqref{IF1} follows from a slight modification of  the proof of Lemma~\ref{LEM:tri2}
with 
\begin{align*}
 M (\bar n) \les \ld^{-s} \frac{|n|^{s}}{|n_1|^{s}|n_2|^{s}}
\end{align*}

\noi
for any $n, n_1, n_2 \in \Z_\ld^*$
with $n = n_1 + n_2$, 
which follows from  \eqref{IF3},  \eqref{Iop1}, and case by case analysis
(recall that $s < 0$).
\end{proof}

By noting that we have 
$ \ld^{ \frac32  -3\g} \le 1$
for any $\ld \ge 1$
in \eqref{IF1}, 
it follows from Lemma \ref{LEM:tri4} and 
Lemma \ref{LEM:OBS1}\,(i)
(more precisely, see \eqref{X1} and \eqref{X2})
that 
\begin{align}
\|Y^\ld\|_{\cX^{0, \g}_2([0, \ld^3T] \times \T_\ld)}
\les  
  \|\Phi^{w}\|_{\W_{T}^{\rho, \g}}, 
\label{I15x}
\end{align}

\noi
where the implicit constant is independent of $\ld \ge 1$.
Then, 
as a corollary to Lemma~\ref{LEM:tri4} (with~\eqref{I15x})
and Proposition \ref{PROP:main}, 
we obtain the following local well-posedness
result for the scaled modulated $I$-KdV \eqref{kdv4}.

\begin{corollary}
\label{COR:I2}
Let $\ld \ge 1$.
Given $\rho \ge\frac 12$,  $\frac 12 < \g < 1$, and $T> 0$, 
let  $w$ be $(\rho,\g)$-irregular
on $[0, T]$ 
 in the sense of Definition~\ref{DEF:ir}.
Fix  $s < 0$ satisfying \eqref{reg1}.
Then, given any 
$v_0 \in L^2(\T_\ld)$ and  $t_0 \in [0, \ld^3 T]$, 
there exist $C_0 > 0$, $\ta > 0$, both independent of $v_0$ and $t_0$, 
and a unique solution $Iu^\ld$
to 
the scaled modulated $I$-KdV equation~\eqref{kdv4} 
with $I u^\ld (t_0)= v_0$
on the time interval $[t_0, t_0 + \tau]\cap [0, \ld^3 T]$,  
belonging to the class\textup{:}
\begin{align*}
 \cD_{w^\ld}^0\big( ([t_0, t_0 + \tau] \cap [0, \ld^3 T]) \times \T_\ld\big), 
\end{align*}

\noi
where the local existence time 
$\tau = \tau \big(\| v_0 \|_{L^2(\T_\ld)}\big)> 0$ satisfies
\begin{align}
\tau  \ge C_0 \|Y^\ld\|_{\cX^{0,\g}_2([0, \ld^3 T]\times \T_\ld)}^{-\ta} \big( 1 + \| v_0 \|_{L^2(\T_\ld)}\big)^{-\ta}.
\label{ME0a}
\end{align}

\noi
Moreover, given any $0 < \al < \g$
with $\al + \g > 1$, 
there exists $C_\al > 0$ such that 
\begin{align*}
 \| I\uu^\ld \|_{\CC^\al ([t_0, t_0 + \tau];L^2(\T_\ld))}\le C_\al 
 \|v_0\|_{L^2(\T_\ld)}, 
\end{align*}

\noi
where $\uu^\ld$ denotes the modulated interaction representation of $u^\ld$
defined in  \eqref{ME0}.

\end{corollary}

\subsection{Commutator estimate}
\label{SUBSEC:I3}

In view of the blowup alternative \eqref{BA1}
with \eqref{I1} and \eqref{scaling3}, 
our main task is to control the growth
of the (scaled) modified energy:
\begin{align}
\| I u^\ld(t)\|_{L^2(\T_\ld)}^2 = \| I \uu^\ld(t)\|_{L^2(\T_\ld)}^2, 
\label{ME1}
\end{align}

\noi
where $\uu^\ld$ denotes the modulated interaction representation of $u^\ld$
defined in  \eqref{ME0}.
Recall from \eqref{GK0} that 
\begin{align}
\big\langle 
I \uu^\ld(r) , X^\ld_{t,r} ( I \uu^\ld(r) , I \uu^\ld(r) )
\big\rangle_{L^2(\T_\ld)}= 0.
\label{I4}
\end{align}

\noi
Using \eqref{I4}, 
we then have 
\begin{align}
\begin{split}
& \| I \uu^\ld(t) \|^2_{L^2(\T_\ld)} -  \| I \uu^\ld(r) \|^2_{L^2(\T_\ld)}\\
& \quad =2 \big\langle 
I \uu^\ld(r)  ,   \com_{t, r}(\uu^\ld(r), \uu^\ld(r))
\big\rangle_{L^2(\T_\ld)} 
+ R^\ld_{t,r}
\end{split}
\label{I2}
\end{align}

\noi
for any $t > r \ge 0$, 
where the commutator  $\com_{t, r}$ and the remainder $R^\ld_{t, r}$ are given by 
\begin{align}
\begin{split}
\com_{t, r}(f_1, f_2)
& =  I X^{\ld}_{t,r}(f_1, f_2 ) - X^{\ld}_{t,r}(I f_1, If_2), 
\\
R^\ld_{t,r}
& =\| I  \uu^\ld(t)-I\uu^\ld(r) \|_{L^2(\T_\ld)}^2\\
& \quad 
+
2
\big\langle 
I \uu^\ld(r) , I \uu^\ld(t) - I \uu^\ld(r) -  I X^{\ld}_{t,r}( \uu^\ld(r) , \uu^\ld(r) )
\big\rangle_{L^2(\T_\ld)}.
\end{split}
\label{I3}
\end{align}

\noi
In the next proposition, 
we establish a key commutator estimate, 
controlling the first term on the right-hand side of \eqref{I2}.
In particular,  note that the bound \eqref{I5}
on the commutator is controlled by a negative power of $N$ (see \eqref{I5a}),
yielding almost conservation of the modified energy 
$\| I \uu^\ld(t)\|^2_{L^2(\T_\ld)}$, 
which is the reason that this approach is called the method
of almost conservation laws.

\begin{proposition}
\label{PROP:com}

Let $\ld \ge 1$.
Given $\rho  >  \frac 12$,  $\frac 12 < \g < 1$,  and $T> 0$, 
let  $w$ be $(\rho,\g)$-irregular
on $[0, T]$ 
 in the sense of Definition~\ref{DEF:ir}, 
 and let $X^\ld$ be as in \eqref{Xld1}.
Fix  $s < 0$ satisfying~\eqref{reg1}.
Then, we have 
\begin{align}
\begin{aligned}
 \|\com_{t, r}(f_1, f_2) \|_{L^2(\T_\ld)}&  = \|IX^{\ld}_{t,r}( f _1, f _2) - X^{\ld}_{t,r} (I f _1,I f _2) \|_{L^2(\T_\ld)}\\
& \les
   K(\ld, N )\| \Phi^{w} \|_{\W_{T}^{\rho, \g}}
      |t-r|^{\g}
\| I  f _1 \|_{L^2(\T_\ld)}   \| I  f _2 \|_{L^2(\T_\ld)}  
\end{aligned}
\label{I5}
\end{align}

\noi
for any $0\leq r<t\leq \ld^3T$, $\ld \ge 1$, and $N \in \N$, 
where $ K(\ld, N )$ is given by 
\begin{align}
K(\ld, N )
= \begin{cases}
\ld^{3- 3\rho-3\g} N^{\frac 32 - 3\rho}, & \text{if $\frac 12 < \rho < \frac 32$}, \\
\ld^{- \frac 32 - 3\g} N^{-3} \sqrt{\log \jb{N} + \log \jb{\ld}}, & \text{if $\rho =  \frac 32$}, \\
\ld^{ \frac 32 - 2\rho - 3\g} N^{-2\rho}, & \text{if $ \rho >  \frac 32$}.
\end{cases}
\label{I5a}
\end{align}

\end{proposition}

\begin{proof}

From \eqref{Xld1} with \eqref{Ix2} and \eqref{FT3}, we have 
\begin{align}
\begin{aligned}
\|  & I X^{\ld}_{t,r}( f _1, f _2) - X^{\ld}_{t,r} (I  f_1, I  f_2) \|^2_{L^2(\T_\ld)}\\
&=
\ld^{-3} \sum_{ n\in\Z_{\ld}^*}
|n|^2 \bigg| \sum_{\substack{n_1, n_2 \in\Z_{\ld}^*\\n = n_1+n_2}}
M_N(\bar n)
\Phi^{w^\ld}_{t,r}(\Xi_\KDV(\bar n))
\ft  f _1(n_1) \ft  f _2 (n_2)
\bigg|^2 ,
\end{aligned}
\label{I6}
\end{align}

\noi
where $\Phi^{w^\ld}_{t, r}$ and 
$\Xi_\KDV(\bar n)$
are as in \eqref{rho2} 
 and \eqref{K3}, respectively, 
 and 
\begin{align}
M_N(\bar n)=  M_N(n,n_1, n_2)=  m_{s, N} (n) - m_{s, N} (n_1) m_{s, N}(n_2) .
\label{I6a}
\end{align}

\noi
By symmetry, assume $|n_1|\ge |n_2|$
in the following.
Then, given $n \in \Z_\ld^*$, 
we split  
\[\Si_n = \big\{ (n_1, n_2)\in (\Z_\ld^*)^2:
n = n_1+n_2,\, |n_1| \ge |n_2|\, \big\}
\]

\noi
 into four regions: 
\begin{align*}
\Si_n 
= \bigcup_{j=1,2,3,4} D_j(n), 
\end{align*}

\noi
where $D_j(n)$,  $j=1,\dots,4$,  is given by 
\begin{align}
\begin{aligned}
D_1(n)&=\Big\{ (n_1, n_2)\in \Si_n :\, |n_1| <  \frac N2, |n_2| <  \frac N2 \Big\}, \\
D_2(n)&=\Big\{(n_1, n_2)\in \Si_n:\, |n_1| \geq \frac N2, |n_2| <  \frac N2, |n|<  \frac N4 \Big\}, \\
D_3(n)&=\Big\{(n_1, n_2)\in \Si_n:\, |n_1| \geq \frac N2, |n_2| <  \frac N2, |n|\geq \frac N4 \Big\}, \\
D_4(n)&=\Big\{(n_1, n_2)\in \Si_n:\, |n_1| \geq \frac N2, |n_2| \geq \frac N2 \Big\}. 
\label{I8}
\end{aligned}
\end{align}

\noi
On $D_1(n)$ and $D_2(n)$, we have $M_N(\bar n) = 0$, 
since
$m_{s, N}(n)=m_{s, N}(n_1)=m_{s, N}(n_2)=1$. 

Proceeding as in \eqref{K4} and \eqref{KR0x}
with \eqref{FT4}, 
we have 
\begin{align*}
\text{LHS of \eqref{I6}}
&\les \sum_{j = 3}^4 
A_{\ld, N}^{(j)} \| I  f_1 \|_{L^2(\T_\ld)}^2 \| I  f_2 \|_{L^2(\T_\ld)}^2
\end{align*}

\noi
where $A_{\ld, N}^{(j)}$, $j =  3, 4$, is given by 
\begin{align}
A_{\ld, N}^{(3)}
& =
 \ld^{-1}\sup_{\substack{n\in \Z_\ld^*\\|n|\geq \frac N4}} |n|^2\sum_{D_3(n)}
 \frac{|\Phi^{w^\ld}_{t,r}(\Xi_\KDV(\bar n)) |^2  |M_N(\bar n)|^2 }{ |m_{s, N}(n_1)|^2|m_{s, N}(n_2)|^2},
\label{I11}\\
A_{\ld, N}^{(4)}
& =
\ld^{-1} 
\sup_{n_1 \in \Z_\ld^*} 
\sum_{n\in \Z_\ld^*} 
\ind_{D_4(n)}\cdot |n|^2 
\frac{ |\Phi^{w^\ld}_{t,r}(\Xi_\KDV(\bar n))|^2 |M_N(\bar n)|^2}{|m_{s, N}(n_1)|^2 |m_{s, N}(n_2)|^2}.
\label{I11a}
\end{align}

%
%

\medskip

\noi
$\bullet$ \textbf{Case 1:} on  $D_3$.\\
\indent
In this case, we have $m_{s, N}(n_2) = 1$
and $|n| \sim |n_1| \ges N \ges |n_2|$.
Then,  from \eqref{Iop1} and the mean value theorem, we have
\begin{align}
|M_N(\bar n)|
= |m_{s, N}(n)-m_{s, N}(n_1)|\les \frac{|n|^{s-1}|n_2|}{N^s}.
\label{I11c}
\end{align}

\noi
Then, from 
\eqref{I11}, \eqref{D3},   \eqref{I11c}, 
and a Riemann sum approximation,  we have 
\begin{align}
\begin{split}
A_{\ld, N}^{(3)}
&  \les 
\ld^{5-6\rho - 6\g}
  \|\Phi^{w} \|^2_{\W^{\rho,\g}_T} |t-r|^{2\g}
\sup_{|n|\geq \frac N4}|n|^{2s-2\rho}
\sum_{D_3(n)} \frac{1}{|n_1|^{2s+2\rho}|n_2|^{2\rho-2}}\\
&  \les 
\ld^{6-6\rho - 6\g}N^{-4\rho}
  \|\Phi^{w} \|^2_{\W^{\rho,\g}_T} |t-r|^{2\g}
  \sum_{\substack{m_2 \in \Z^*\\ |m_2|\les \ld N}} 
 \Big|\frac {m_2}\ld \Big|^{2-2\rho}\frac 1\ld\\
&  \les 
\ld^{6-6\rho - 6\g}N^{-4\rho}
  \|\Phi^{w} \|^2_{\W^{\rho,\g}_T} |t-r|^{2\g}
\int_{\ld^{-1}}^{2N} x^{2- 2\rho} dx\\
 &\les 
K(\ld, N)^2
 \| \Phi^{w} \|^2_{\W_{T}^{\rho, \g}}  
|t-r|^{2\g}
\end{split}
\label{I11b}
 \end{align}

 \noi
for any $0\leq r<t\leq \ld^3T$, 
where $K(\ld, N)$ is as in \eqref{I5a}. 
Here, the last step follows
from separately estimating 
the cases when $\rho < \frac 32$, 
$\rho = \frac 32$, 
and $\rho > \frac 32$.
This proves
 \eqref{I5} in this case.

\medskip

\noi
$\bullet$ \textbf{Case 2:} on  $D_4$.\\
\indent
From \eqref{I6a}, we have
\begin{align}
|M_N(\bar n)|\le |m_{s, N} (n)| + | m_{s, N} (n_1) m_{s, N}(n_2)|.
\label{I13x}
\end{align}

\noi
Then, from \eqref{I11a} with \eqref{D3}  and \eqref{I13x}, we have 
\begin{align}
\begin{split}
A_{\ld, N}^{(4)}
& \les 
\ld^{5-6\rho -6\g}N^{4s }
\| \Phi^{w} \|^2_{\W_{T}^{\rho, \g}} 
|t-r|^{2\g} \\
& \hphantom{XX}
\times
\sup_{n_1 \in \Z_\ld^*} 
\sum_{n\in \Z_\ld^*} 
\ind_{D_4(n)}\cdot
|n|^{2-2\rho}
\frac{  |m_{s, N}(n)|^2}{|n_1|^{2s + 2\rho}|n_2|^{2s + 2\rho}}\\
& \quad + 
\ld^{5-6\rho -6\g}
\| \Phi^{w} \|^2_{\W_{T}^{\rho, \g}} 
|t-r|^{2\g} \\
& \hphantom{XX}
\times
\sup_{n_1 \in \Z_\ld^*} 
\sum_{n\in \Z_\ld^*} 
\ind_{D_4(n)}\cdot|n|^{2-2\rho}
|n_1|^{- 2\rho}|n_2|^{- 2\rho}\\
 & =: \1 + \II.
\end{split}
\label{I13a}
\end{align}

Let us first estimate $\1$.

\noi
$\circ$ {\bf Subcase 2.a:}  $|n| \ll N$.\\
\indent
In this case, we have $m_{s, N}(n) = 1$.
Then, proceeding as in \eqref{I11b} with $|n_1|\sim |n_2|\ges N$, we have 
\begin{align}
\begin{split}
\1
&\les 
\ld^{5-6\rho -6\g}N^{-4\rho }
\| \Phi^{w} \|^2_{\W_{T}^{\rho, \g}} 
|t-r|^{2\g} 
\sum_{\substack{n \in \Z_\ld^*\\|n|\ll N}} |n|^{2-2\rho} \\
 &\les 
K(\ld, N)^2
 \| \Phi^{w} \|^2_{\W_{T}^{\rho, \g}}  
|t-r|^{2\g}
\end{split}
\label{I14}
\end{align}

 \noi
for any $0\leq r<t\leq \ld^3T$, 
thus yielding \eqref{I5}.

\medskip

\noi
$\circ$ {\bf Subcase 2.b:}  $|n| \ges N$.\\
\indent
In this case, we have   $|m_N(n)| \les  N^{-s} |n|^{s}$.
Without loss of generality, assume that $|n|\les |n_1|$.
We note that, for $\rho \ge \frac 32$, we have 
\begin{align}
\ld^{3- 3\rho-3\g} N^{\frac 32 - 3\rho}
\le
\ld^{ \frac 32 - 2\rho - 3\g} N^{-2\rho}.
\label{I5x}
\end{align}

\noi
Then, from \eqref{I13a}, a Riemann sum approximation, 
Lemma \ref{LEM:SUM}, 
and 
\eqref{I5x}, we have
\begin{align}
\1
&\les 
\ld^{5-6\rho -6\g}N^{2s }
\| \Phi^{w} \|^2_{\W_{T}^{\rho, \g}} 
|t-r|^{2\g}
\sup_{\substack{n_1 \in \Z_\ld^*\\|n_1|\ges N} }
\sum_{\substack{n \in \Z_\ld^*\\|n|, |n- n_1|\ges N}} 
\frac{1}{\jb{n}^{4\rho - 2}\jb{n - n_1}^{2s+2\rho}}
\notag \\
&\les 
\ld^{6-6\rho -6\g}N^{2s }
\| \Phi^{w} \|^2_{\W_{T}^{\rho, \g}} 
|t-r|^{2\g}  \notag \\
& \quad 
\times 
\sup_{\substack{n_1 \in \Z_\ld^*\\|n_1|\ges N} }
 \sum_{\substack{m \in \Z^*\\ |\frac m\ld |, |\frac m\ld - n_1|\ges  N}} 
\frac{1}{\jb{\frac{m}{\ld}}^{4\rho - 2}
\jb{ \frac m\ld - n_1}^{2s + 2\rho}}
 \frac 1\ld
\label{I16x}\\
&\les 
\ld^{6-6\rho -6\g}N^{2s }
\| \Phi^{w} \|^2_{\W_{T}^{\rho, \g}} 
|t-r|^{2\g}
\sup_{\substack{n_1 \in \Z_\ld^*\\|n_1|\ges N} }
\intt_{|\xi|, |\xi - n_1|\ges N} 
\frac{d\xi}{\jb{\xi}^{4\rho - 2}
\jb{ \xi - n_1}^{2s + 2\rho}} \notag \\
 &\les 
\ld^{6-6\rho -6\g}N^{3-6\rho }
 \| \Phi^{w} \|^2_{\W_{T}^{\rho, \g}}  
|t-r|^{2\g} \notag \\
 &\les 
K(\ld, N)^2
 \| \Phi^{w} \|^2_{\W_{T}^{\rho, \g}}  
|t-r|^{2\g}, \notag 
\end{align}

\noi
provided that 
  \eqref{K5b} holds.
First, suppose that $s>-\rho$.
In applying Lemma \ref{LEM:SUM} in the fourth step of \eqref{I16x}, 
we used 
\begin{align}
\begin{aligned}
& \intt_{|\xi|, |\xi - n_1|\ges N} 
\frac{d\xi}{\jb{\xi}^{4\rho - 2}
\jb{ \xi - n_1}^{2s + 2\rho}}\\
& \hphantom{XXXX}
\les N^{2- 2s - 6\rho + \al +\be}
\intt_{|\xi|, |\xi - n_1|\ges N} 
\frac{d\xi}{\jb{\xi}^{\al}
\jb{ \xi - n_1}^{\be}}
\\
& \hphantom{XXXX}
 \les N^{3- 2s - 6\rho}, 
\end{aligned}
\label{I16}
\end{align}

\noi
uniformly in $n_1\in \Z_\ld^*$, 
by taking some $0 \le \al, \be < 1$
such that $\al \le 4\rho - 2$
and $ \be \le 2s + 2\rho$
satisfying $\al + \be > 1$.
This choice is possible since $\rho>\frac{1}{2}$ and \eqref{K5b} give
$4\rho-2>0$, $2s+2\rho>0$, and
$(4\rho-2)+(2s+2\rho)>1$.

When $s=-\rho$, the condition \eqref{reg1}
implies $\rho>\frac{3}{4}$.
In this case, we simply use 
\begin{align}
\intt_{|\xi|,|\xi-n_1|\ges N}
\frac{d\xi}{\jb{\xi}^{4\rho-2}\jb{\xi - n_1}^0}
\le \intt_{|\xi|\ges N}
\frac{d\xi}{\jb{\xi}^{4\rho-2}}
\les N^{3-4\rho}
=N^{3-2s-6\rho}, 
\label{I16x1}
\end{align}

\noi
uniformly in $n_1\in \Z_\ld^*$.

\medskip

Next, we  estimate $\II$.

\noi
$\circ$ {\bf Subcase 2.c:}  $|n| \ll N$.\\
\indent
By substituting  $|n_1|, |n_2|\ges N$ in \eqref{I13a}, 
we see that \eqref{I14} also holds for $\II$ in this case.

\medskip

\noi
$\circ$ {\bf Subcase 2.d:}  $|n| \ges N$.\\
\indent
Without loss of generality, assume that $N \les |n| \les |n_1|$.
From \eqref{I13a}, 
 a Riemann sum approximation, 
Lemma \ref{LEM:SUM}, 
and 
\eqref{I5x}, 
we have
\begin{align*}
\II
&\les 
\ld^{5-6\rho -6\g}
\| \Phi^{w} \|^2_{\W_{T}^{\rho, \g}} 
|t-r|^{2\g} 
\sup_{\substack{n_1 \in \Z_\ld^*\\|n_1|\ges N} }
\sum_{\substack{n \in \Z_\ld^*\\|n|, |n- n_1|\ges N}} 
\frac{1}{\jb{n}^{4\rho - 2}\jb{n - n_1}^{2\rho}}\\
 &\les 
\ld^{6-6\rho -6\g}N^{3-6\rho }
 \| \Phi^{w} \|^2_{\W_{T}^{\rho, \g}}  
|t-r|^{2\g}\\
 &\les 
K(\ld, N)^2
 \| \Phi^{w} \|^2_{\W_{T}^{\rho, \g}}  
|t-r|^{2\g}, 
\end{align*}

\noi
provided that $\rho > \frac 12$, 
where, at the second inequality, we proceeded as in \eqref{I16}.
\end{proof}

\subsection{Global well-posedness on the circle}
\label{SUBSEC:I4}

We are now ready to put  together all the ingredients from the previous subsections
and prove
global well-posedness  
of the periodic modulated KdV~\eqref{kdv1} in negative Sobolev spaces
as  claimed in Theorem~\ref{THM:1}\,(ii).

Fix $u_0 \in H^s(\T)$ for some $s < 0$ (to be determined later).
Let 
 $u$ be a solution to  the modulated KdV~\eqref{kdv1} on $\T$ 
with $u|_{t = 0} = u_0$
whose existence for short times is guaranteed by Theorem \ref{THM:1}\,(i).
In view of the blowup alternative \eqref{BA1}, 
our goal is to show that 
\begin{align}
\sup_{0 \le t \le T} \|u(t)\|_{H^s(\T)}^2 \le C(T) < \infty
\label{ME2}
\end{align}

\noi
for  each finite $T > 0$.
In view of 
 \eqref{I1} (which also holds true on $\T_\ld$), \eqref{scaling1}, \eqref{scaling3}, and \eqref{ME1}, 
the bound \eqref{ME2} follows once we prove the following bound on the modified energy:
\begin{align*}
\sup_{0 \le t \le \ld^3 T} \|I\uu^\ld(t)\|_{L^2(\T_\ld)}^2 \le C(T, \ld ) < \infty
\end{align*}

\noi
for some $\ld \ge 1$.

From \eqref{I1} 
and \eqref{scaling3}, we have 
\begin{equation}
 \| I  u_0^{\ld} \|_{L^2(\T_\ld) }
 \les \ld^{-s-\frac 32} N^{-s} \|  u_0  \|_{H^{s}(\T)}.
 \label{IT1}
\end{equation}

\noi
Then, given small $\eps_0 > 0$, 
we can choose 
\begin{equation}
\ld \sim N^{-\frac{2s}{3+2s}}, 
 \label{IT2}
\end{equation}

\noi
such that 
\begin{align}
\| Iu_0^{\ld} \|^2 _{L^2(\T_\ld)}\le \eps_0, 
\label{IT2a}
\end{align}

\noi
provided that $s > -\frac 32$.
Here, 
 the implicit constant in \eqref{IT2} depends only on ${\| u_0 \|_{H^{s}(\T) }} $.

We now study the increment of the modified energy.
Let us first make preliminary analysis
on the remainder term
 $R^\ld_{t,r}$ defined in \eqref{I3}.
By rewriting \eqref{I2}, we have
\begin{align}
\begin{split}
R^\ld_{t,r} &=  \| I \uu^\ld(t) \|^2_{L^2(\T_\ld)} -  \| I \uu^\ld(r) \|^2_{L^2(\T_\ld)}\\
& \quad 
- 2 \big\langle 
I \uu^\ld(r)  ,   \com_{t, r}(\uu^\ld(r), \uu^\ld(r))
\big\rangle_{L^2(\T_\ld)} , 
\end{split}
\label{I2a}
\end{align}

\noi
where
$\com_{t, r}$ is  as in 
\eqref{I3}.
Then, 
from \eqref{dl1} and \eqref{I2a}
with 
\[  \com_{t_1, t_3}(\uu^\ld(t_3), \uu^\ld(t_3))
-    \com_{t_2, t_3}(\uu^\ld(t_3), \uu^\ld(t_3))
=   \com_{t_1, t_2}(\uu^\ld(t_3), \uu^\ld(t_3)), \]

%
\noi
which follows from \eqref{I3}, 
we have 
\begin{align}
\begin{aligned}
(\updl R^\ld)_{t_1, t_2, t_3} 
& = R^\ld_{t_1,t_3}-  R^\ld_{t_1, t_2}-   R^\ld_{t_2, t_3}\\
&=
2\big\langle 
I \uu^\ld(t_2) - I \uu^\ld(t_3)  ,   \com_{t_1, t_2}(\uu^\ld(t_2), \uu^\ld(t_2))
\big\rangle_{L^2(\T_\ld)} \\
& \quad + 2 \big\langle 
I \uu^\ld(t_3)  ,   \com_{t_1, t_2}(\uu^\ld(t_2)- \uu^\ld(t_3), \uu^\ld(t_2))
\big\rangle_{L^2(\T_\ld)}\\
& \quad + 2 \big\langle 
I \uu^\ld(t_3)  ,   \com_{t_1, t_2}(\uu^\ld(t_3), \uu^\ld(t_2)- \uu^\ld(t_3))
\big\rangle_{L^2(\T_\ld)}
\end{aligned}
\label{I2b}
\end{align}

\noi
for any  $t_1 >   t_2> t_3\ge 0$.

Fix an interval $J \subset [0, \ld^3T]$
and $0 < \al < \g$ such that $\al + \g > 1$.
Then, it follows from~\eqref{I2b}, Cauchy-Schwarz's inequality, 
and Proposition~\ref{PROP:com}
%
that  
\begin{align*}
\begin{aligned}
|(\updl R^\ld)_{t_1, t_2, t_3} |
&\les K(\ld, N)   \|  \Phi^{w}  \|_{\W_{T}^{\rho, \g}}   
|t_1-t_2|^{\g}
|t_2-t_3|^{\al}\|I \uu^\ld \|^3_{\CC^\al (J;L^2(\T_\ld))}
\end{aligned}
\end{align*}

\noi
for any 
$t_1 \ge   t_2\ge t_3\ge 0$
belonging to the interval $J$, 
where $K(\ld, N)$ is as in \eqref{I5a}.
Since $\al + \g > 1$, 
it follows from the sewing lemma 
(Lemma \ref{LEM:sew}) 
that 
\begin{align}
|  R^\ld_{t_1,t_2} | \les 
K(\ld, N)   \|  \Phi^{w}  \|_{\W_{T}^{\rho, \g}}   
|t_1-t_2|^{\al + \g}  \|I \uu^\ld \|^3_{\CC^\al (J;L^2(\T_\ld))}
\label{I2x}
\end{align}

\noi
for any 
$t_1 >   t_2\ge 0$
belonging to the interval $J$.

We now put everything together.
Let   $\tau \sim 1$ be  the local existence time
for the scaled modulated $I$-KdV \eqref{kdv4}
given by~\eqref{ME0a} in Corollary \ref{COR:I2}, 
corresponding to initial data $v_0$ with 
$\|v_0 \|_{L^2(\T_\ld)}^2
= 2\eps_0$.
Namely, if we have 
\begin{align}
 \|I\uu^\ld(t_0)\|_{L^2(\T_\ld)}^2 \le 2\eps_0
\label{I2c}
\end{align}

\noi
for some 
 $0 \le t_0 \le \ld^3T$, 
then Corollary \ref{COR:I2}
guarantees that the solution $Iu^\ld$ to the scaled modulated $I$-KdV~\eqref{kdv4} exists on 
the time interval $[t_0, t_0 + \tau]\cap [0, \ld^3 T]$, satisfying the bound:
\begin{align}
 \| I\uu^\ld \|_{\CC^\al ([t_0, t_0 + \tau];L^2(\T_\ld))}^2\le C_\al \eps_0.
\label{I2d}
\end{align}

\noi
for some $C_\al > 0$.
In particular, in view of~\eqref{IT2a} (which states that \eqref{I2c} is satisfied with $t_0 = 0$), 
we see that 
 the solution $Iu^\ld$  exists on $[0, \tau]$.

Suppose now that by choosing $N \gg 1$
(and thus $\ld = \ld(N) \gg1$
in view of \eqref{IT2}), 
we can make 
 $K(\ld, N)$  in~\eqref{I5a}
as small as we need;
see \eqref{IT0b} and \eqref{IT0d}.
Then,  
from 
\eqref{I2}, \eqref{IT2a}, 
Proposition~\ref{PROP:com}, 
 \eqref{I2x}, 
and \eqref{I2d} with $t_0 = 0$,  
 we have 
\begin{align}
\begin{split}
 \sup_{0\le t\le \tau} \| I \uu^\ld(t) \|^2_{L^2(\T_\ld)} 
&  \le \eps_0 + 
C_2 K(\ld, N)   \|  \Phi^{w}  \|_{\W_{T}^{\rho, \g}}   \\
& \le 2\eps_0.
\end{split}
\label{I2e}
\end{align}

\noi
In particular, the solution $Iu^\ld$ exists on $[\tau,  2\tau]$, 
satisfying 
\eqref{I2c} and \eqref{I2d}
with $t_0 = \tau$.
Then,  
from 
\eqref{I2}, \eqref{I2e} 
Proposition~\ref{PROP:com}, 
 \eqref{I2x}, and \eqref{I2d} with $t_0 = \tau$, we have 
\begin{align}
\begin{split}
 \sup_{\tau\le t\le2 \tau} \| I \uu^\ld(t) \|^2_{L^2(\T_\ld)} 
&  \le \eps_0 + 
2 C_2 K(\ld, N)   \|  \Phi^{w}  \|_{\W_{T}^{\rho, \g}}   \\
&  \le 2\eps_0, 
\end{split}
\label{I2f}
\end{align}

\noi
provided that 
we can make 
 $K(\ld, N)$  small by 
choosing $N \gg 1$
(and thus $\ld = \ld(N) \gg1$).
If~\eqref{I2f} holds, 
then 
the solution $Iu^\ld$ exists on $[2\tau, 3\tau]$, 
satisfying 
\eqref{I2c} and \eqref{I2d}
with $t_0 = 2\tau$, 
which allows us to iterate this process.

After $j$ steps, 
this iterative process shows that 
the solution $Iu^\ld$ exists on $[0, (j+1)\tau]$, 
satisfying~\eqref{I2c} and~\eqref{I2d}
 with $t_0 = j \tau$, and moreover, we have 
\begin{align}
 \sup_{0 \le t\le j \tau} \| I \uu^\ld(t) \|^2_{L^2(\T_\ld)} 
&  \le \eps_0 + 
j C_2 K(\ld, N)   \|  \Phi^{w}  \|_{\W_{T}^{\rho, \g}}.
\label{I2g}
\end{align}

\noi
Then, proceeding as in \eqref{I2f} with \eqref{I2g}, we obtain
\begin{align}
\begin{split}
 \sup_{0 \le t\le(j+1) \tau} \| I \uu^\ld(t) \|^2_{L^2(\T_\ld)} 
&  \le \eps_0 + 
(j+1) C_2 K(\ld, N)   \|  \Phi^{w}  \|_{\W_{T}^{\rho, \g}}   \\
&  \le 2\eps_0, 
\end{split}
\label{I2h}
\end{align}

\noi
provided that 
we can make 
 $K(\ld, N)$  small by 
choosing $N \gg 1$
(and thus $\ld = \ld(N) \gg1$
for $- \frac 32 < s < 0$).

Our goal is to construct the solution $Iu^\ld$ on the time interval $[0, \ld^3 T]$, 
for which we need to iterate this process $j \sim \ld^3 T$ times, 
since  $\tau \sim1$.
From \eqref{I2h}, 
we see that 
 it is possible 
to reach the target time $\ld^3 T$
(such that $\ld^3T$ is (at most) the doubling time
of the modified energy), 
provided that 
\begin{align}
\ld^3T\les  K(\ld, N)^{-1}
= \begin{cases}
\ld^{- 3+ 3\rho+3\g} N^{-\frac 32 + 3\rho}, & \text{if $\frac 12 < \rho < \frac 32$}, \\
\ld^{\frac 32 + 3\g} N^{3} (\log \jb{N} + \log \jb{\ld})^{-\frac 12}, & \text{if $\rho =  \frac 32$}, \\
\ld^{ -\frac 32 + 2\rho +3\g} N^{2\rho}, & \text{if $ \rho >  \frac 32$}.
\end{cases}
\label{IT0b}
\end{align}

\noi
In view of  \eqref{IT2}, 
we  can  guarantee \eqref{IT0b} by choosing $N \gg 1$, 
provided that 
\begin{align}
s > 
\begin{cases}
\max \big(\frac {3-6\rho }{6 - 4\g}, - \frac 32\big)& \text{if } \frac 12 < \rho < \frac 32, \\
- \frac 32, & \text{if }  \rho \ge \frac 32, 
\end{cases}
\label{IT0d}
\end{align}

\noi
where we used the fact that 
$\frac{-2\rho}{3 - 2\g}  < -\frac 32$ for 
 $\rho \ge \frac 32$ and
$ \frac 12< \g < 1$.
Hence, 
from~\eqref{reg1} and~\eqref{IT0d} with $ s > - \frac 32$, 
we conclude 
global well-posedness of the modulated KdV in $H^s(\T)$, 
provided that 
$\rho > \frac 12$ and $s\in \R$
satisfy \eqref{s1}.
This concludes the proof of Theorem~\ref{THM:1}\,(ii).

\subsection{Real line case}
\label{SUBSEC:I5}

We conclude this section by discussing the real line case.
In this case, we need to work with the non-homogeneous Sobolev spaces whose
scaling property is worse
(for example, compare \eqref{D3} and \eqref{RR4})
and thus we need to proceed with care.
In particular, we overcome the issue of the worse scaling properties
by changes of variables
(which are not available in the periodic case).
While we used the same notations from Subsections~\ref{SUBSEC:I1}-\ref{SUBSEC:I4}, 
 all the objects 
(equations, functions, and operators) are 
to be understood as those defined on the real line.


Let $s < 0$.
Then, we have 
\begin{align}
\|u_0^\ld\|_{ H^s(\R)}
\le  \ld^{-\frac 32 - s} \|u_0\|_{ H^s(\R)}, 
\label{RR1}
\end{align}

\noi
where $u_0^\ld$ is as in \eqref{scaling1a}.
As in the periodic case, 
our goal is to study the scaled modulated $I$-KdV \eqref{kdv4}
on the real line.
Given small $\eps_0 > 0$, it follows from 
\eqref{I1} (which also holds true on $\R$)
and  \eqref{RR1} that
\begin{equation*}
 \| I  u_0^{\ld} \|_{L^2(\R) }
 \les \ld^{-s-\frac 32} N^{-s} \|  u_0  \|_{H^{s}(\R)}
 \le \eps_0 , 
\end{equation*}

\noi
where the second step 
follows from choosing 
$\ld \sim N^{-\frac{2s}{3+2s}}$ as in \eqref{IT2}, 
provided that $s > -\frac 32$.
This allows us to reduce the study to the small data regime
as in the periodic case.

Given $\ld \ge 1$, let $X^\ld$ be the driver defined in \eqref{Xld1}, 
where we view all the operators and functions 
interpreted as those defined on $\R$.
Then, we have the following 
analogue of Lemma~\ref{LEM:tri2}
on the regularity of the driver $X^\ld$ 
whose proof follows
from a slight modification of the proof of 
Proposition \ref{PROP:KR1}\,(i).

\begin{lemma}
\label{LEM:Rtri2}
Let $\ld \ge 1$.
Given $\rho > \frac 12$, $\frac 12 < \g < 1$, and $T> 0$, 
let  $w$ be $(\rho,\g)$-irregular on $[0, T]$ in the sense of Definition~\ref{DEF:ir}, 
and let $X^\ld$ be as in \eqref{Xld1}.
Suppose that $\rho > \frac 12$ and $s <  0$
satisfy 
\eqref{regR1}.
Then, 
the driver $X^\ld$ 
belongs to $  \cX^{s, \g}_2([0, \ld^3T] \times \R)$ defined in~\eqref{X1}, 
satisfying the bound\textup{:}
 \begin{align*}
 \| X^{\ld}_{t,r}  \|_{\cL_2 ( H^s(\R)) }
 \les
 \ld^{ \frac32 -3\g} 
  \|\Phi^{w}\|_{\W_{T}^{\rho, \g}}
 |t-r|^{\g}
\end{align*}

\noi
for any $0\leq r<t \leq \ld^3T $.

\end{lemma}

Once we prove Lemma \ref{LEM:Rtri2}, 
an analogue of Lemma \ref{LEM:tri3} 
on  equivalence of the 
unscaled and scaled problems in the current real line case follows
from the same proof.

\begin{proof}[Proof of Lemma \ref{LEM:Rtri2}]
Proceeding as in 
 \eqref{D3}, 
we have
\begin{align}
\begin{aligned}
|\Phi^{w^\ld}_{t,r} (\Xi_\KDV(\bar \xi)) |
& = \bigg| 
\int_{r}^{t}  e^{i  \Xi_\KDV (\bar \xi) w^{\ld}(t') } d t'
\bigg|
= \ld^3   \bigg|  \int_{\ld^{-3}r}^{\ld^{-3}t} 
e^{i  \ld^3  \Xi_\KDV(\bar \xi)  w(t') } dt'  \bigg|  \\
& \les 
\ld^{3 -3\g} 
 \|\Phi^{w}  \|_{\W_{T}^{\rho, \g}}
 |t-r|^{\g} 
\jb{\ld^3 \xi \xi_1 \xi_2}^{-\rho}
\end{aligned}
\label{RR4}
\end{align}

\noi
 for any $\xi, \xi_1, \xi_2\in \R$ satisfying $\xi = \xi_1+ \xi_2$.

Proceeding as in \eqref{KR0a}, 
\eqref{KR0b}, and \eqref{KR2}, we have
\begin{align}
\|X_{t , r }^\ld\|^2_{\cL _2 (H^s (\R))}
 &
\le
\sup_{|\xi|< 1} 
  \intt_{|\xi_1|< 2} 
  \jb{\xi}^{ 2s }\xi^2
   \frac{|\Phi^{w^\ld}_{t,r} (\Xi_\KDV(\bar \xi)) |^2}{ \jb{\xi_1}^{2s}  \jb{\xi - \xi_1}^{2s}}  d\xi_1
\notag \\
& \quad + 
\sup_{|\xi_1|\ge 2}
\intt_{|\xi|< 1} \jb{\xi}^{2s}
\xi^2
\frac{|\Phi^{w^\ld}_{t,r} (\Xi_\KDV(\bar \xi)) |^2}{\jb{\xi_1}^{2s}\jb{\xi - \xi_1}^{2s}} d\xi
\notag \\ 
& \quad +
\sup_{|\xi_1|\ge \frac 12}
\intt_{\substack{|\xi|\ge 1\\|\xi-\xi_1|\ge \frac 12 }} 
\jb{\xi}^{2s}\xi^2
\frac{|\Phi^{w^\ld}_{t,r} (\Xi_\KDV(\bar \xi)) |^2}{\jb{\xi_1}^{2s}\jb{\xi - \xi_1}^{2s}} d\xi
\label{RR5} \\
& \quad +
2\sup_{|\xi|\ge 1} 
  \intt_{|\xi - \xi_1|< \frac 12} 
  \jb{\xi}^{ 2s }\xi^2
   \frac{|\Phi^{w^\ld}_{t,r} (\Xi_\KDV(\bar \xi)) |^2}{ \jb{\xi_1}^{2s}  \jb{\xi - \xi_1}^{2s}}  d\xi_1
\notag \\
& =:
A^\ld_1+ A^\ld_2 + A^\ld_3 + A^\ld_4.
\notag
\end{align}

We first estimate $A^\ld_1$.
In this case, we have $|\xi|, |\xi_1|, |\xi - \xi_1|\les 1$.
Without loss of generality, 
assume that  $|\xi| \les |\xi - \xi_1|$.
Then, from  \eqref{RR4}
and a change of variable $\z_1 = \ld^3 \xi^2 \xi_1$, we have
\begin{align}
\begin{split}
A^\ld_1 
& \les \ld^{6 -6\g} 
 \|\Phi^{w}  \|_{\W_{T}^{\rho, \g}}^2
 |t-r|^{2\g} 
\sup_{|\xi|< 1} 
\xi^2  \intt_{|\xi_1|< 2} 
\frac{1}{\jb{\ld^3 \xi^2 \xi_1 }^{2\rho}}
 d\xi_1\\
& \le \ld^{3 -6\g} 
 \|\Phi^{w}  \|_{\W_{T}^{\rho, \g}}^2
 |t-r|^{2\g} 
\sup_{|\xi|< 1} 
 \intt_\R
\frac{1}{\jb{\z_1 }^{2\rho}}
 d\z_1\\
& \les \ld^{3 -6\g} 
 \|\Phi^{w}  \|_{\W_{T}^{\rho, \g}}^2
 |t-r|^{2\g} 
 \end{split}
\label{RR6}
\end{align}

\noi
for any $0 \le r < t \le \ld^3T$, 
provided that  $\rho > \frac 12$.

From  \eqref{RR4}, 
 a change of variable $\z = \ld \xi$, 
and relabelling $\z_1 = \ld \xi_1$, 
 we have 
\begin{align}
\begin{split}
A^\ld_2 
& \les 
 \ld^{6 -6\g} 
 \|\Phi^{w}  \|_{\W_{T}^{\rho, \g}}^2
 |t-r|^{2\g} \\
 & 
 \quad 
 \times 
\sup_{|\xi_1|\ge 2}
\intt_{|\xi|< 1}
\xi^2
\frac{1}
{|\xi_1(\xi - \xi_1)|^{2s}\jb{\ld^3 \xi \xi_1 (\xi-\xi_1)}^{2\rho}} d\xi\\
& = 
 \ld^{3+4s -6\g} 
 \|\Phi^{w}  \|_{\W_{T}^{\rho, \g}}^2
 |t-r|^{2\g} \\
  & 
 \quad 
 \times 
\sup_{|\z_1|\ge 2\ld}
\intt_{|\z|< \ld}
\z^2
\frac1 
{|\z_1(\z - \z_1)|^{2s}\jb{\z\z_1(\z-\z_1)}^{2\rho}} d\z
\end{split}
\label{RR7}
\end{align}

\noi
for any $0 \le r < t \le \ld^3T$.
Denote by $A_2^{\ld, (1)}$
the contribution to $A^\ld_2$ from $|\z|\les 1$.
Then, 
proceeding as in \eqref{KR3c}, 
we have
\begin{align}
\begin{split}
A^{\ld, (1)}_2 
& \les  \ld^{3 + 4s -6\g} 
 \|\Phi^{w}  \|_{\W_{T}^{\rho, \g}}^2
 |t-r|^{2\g}\\
& \les  \ld^{3  -6\g} 
 \|\Phi^{w}  \|_{\W_{T}^{\rho, \g}}^2
 |t-r|^{2\g}
 \end{split}
\label{RR8}
\end{align}

\noi
for any $0 \le r < t \le \ld^3T$, 
provided that 
\begin{align*}
\rho > 0, \quad   -\frac 32 < s < 0, 
\quad \text{and} \quad  
s\ge -\rho.
\end{align*}

\noi
Let $A_2^{\ld, (2)} = A_2^{\ld}- A_2^{\ld, (1)}$.
Using $|\z_1| \sim |\z - \z_1|\ges \ld$, we have 
\begin{align}
\begin{split}
A^{\ld, (2)}_2 
& \les 
 \ld^{3 + 4s -6\g} 
 \|\Phi^{w}  \|_{\W_{T}^{\rho, \g}}^2
 |t-r|^{2\g} 
\sup_{|\z_1|\ge 2\ld}
\intt_{1\ll |\z|< \ld}
|\z|^{2- 2\rho}
\frac1 
{|\z_1(\z - \z_1)|^{2s+ 2\rho}} d\z\\
& \les 
 \ld^{3 -4\rho  -6\g} 
 \|\Phi^{w}  \|_{\W_{T}^{\rho, \g}}^2
 |t-r|^{2\g} 
\sup_{|\z_1|\ge 2\ld}
\intt_{1\ll |\z|< \ld}
|\z|^{2- 2\rho}
d\z\\
& \les  \ld^{3  -6\g} 
 \|\Phi^{w}  \|_{\W_{T}^{\rho, \g}}^2
 |t-r|^{2\g}
\end{split}
\label{RR10}
\end{align}

\noi
for any $0 \le r < t \le \ld^3T$, 
provided that $s + \rho \ge 0$
and $\rho \ge \frac 12$.

As for $A^\ld_3$, 
there is no low frequency issue.
Thus, by 
proceeding as in \eqref{K5} and \eqref{K5a} with~\eqref{RR4}, we have
\begin{align}
\begin{aligned}
A_3^\ld
& \les
\ld^{6 - 6\rho - 6\g}
\|\Phi^w\|^2_{\W_T^{\rho,\g}}|t-r |^{2\g}
\sup_{|\xi_1|\ge \frac 12}
\intt_{\substack{|\xi|\ge 1\\|\xi-\xi_1|\ge \frac 12 }} 
\frac{|\xi|^{2s+2-2\rho}}{|\xi_1(\xi - \xi_1)|^{2s+2\rho}} d\xi\\
& \les
\ld^{3 - 6\g}
\|\Phi^w\|^2_{\W_T^{\rho,\g}}|t-r |^{2\g}
\end{aligned}
\label{RR11}
\end{align}

\noi
for any $0 \le r <  t \le \ld^3T$, 
provided that
\eqref{K5b} holds.

As for $A_4^\ld$, 
we first note that 
$ |\xi_1|\sim |\xi| \ge 1$
and $\jb{\xi - \xi_1}\sim 1$
under 
 $|\xi|\ge1$ and $|\xi - \xi_1|< \frac12 $.
 Then, from~\eqref{RR4}, \eqref{KR1b},  and a change of variables
 $\z_2  = \ld^3 \xi^2 (\xi - \xi_1)$, we have 
 \begin{align}
\begin{aligned}
A_4^\ld
&\les \ld^{6 - 6\g}
\|\Phi^w\|^2_{\W_T^{\rho,\g}}
|t-r |^{2\g }    \sup_{|\xi|\ge 1}  \xi^2
  \intt_{|\xi - \xi_1|< \frac 12} 
\frac{ 1}{  \jb{\ld^3 \xi^2(\xi - \xi_1)}^{2\rho}} d \xi_1\\
&= \ld^{3 - 6\g}
\|\Phi^w\|^2_{\W_T^{\rho,\g}}
|t-r |^{2\g }  
   \sup_{|\xi|\ge 1}  
  \intt_{|\z_2|< \frac {\ld^3 \xi^2}2} 
 \frac{1}{\jb{\z_2}^{2\rho} }d \z_2\\
&  \les \ld^{3 - 6\g}\|\Phi^w\|^2_{\W_T^{\rho,\g}}
|t-r |^{2\g } 
\end{aligned}
\label{RR12}
\end{align}

\noi
for any $0 \le r <  t \le \ld^3T$, provided that $\rho > \frac 12$
(and any $s \in \R$).

Therefore, putting 
\eqref{RR5}, \eqref{RR6}, 
\eqref{RR7}, \eqref{RR8}, 
\eqref{RR10}, 
\eqref{RR11}, 
and \eqref{RR12}
together, we conclude 
from  Lemma \ref{LEM:OBS1}\,(i)
that 
 $X^\ld\in  \cX^{s, \g}_2([0, \ld^3 T]\times \R)$ under
 $s< 0$ and~\eqref{regR1}.
\end{proof}

We now consider the 
 scaled modulated $I$-KdV \eqref{kdv4}
 on $\R$.
Let $Y^\ld$ be the associated
driver
defined in \eqref{Xld2}, 
where we view all the operators and functions 
interpreted as those defined on $\R$.
Then, 
as a corollary to 
Lemma \ref{LEM:Rtri2}, we obtain the following 
regularity property
of the driver $Y^\ld$
and local well-posedness of 
the  scaled modulated $I$-KdV \eqref{kdv4}
on the real line.
See
 Lemma \ref{LEM:tri4}
and Corollary \ref{COR:I2}
for the corresponding statements
in the periodic setting
studied earlier.

\begin{corollary}
\label{COR:RI1}
Let $\ld \ge 1$.
Given $\rho > \frac 12$, $\frac 12 < \g < 1$, and $T> 0$, 
let  $w$ be $(\rho,\g)$-irregular on $[0, T]$ in the sense of Definition~\ref{DEF:ir}, 
and let $X^\ld$ be as in \eqref{Xld1}.
Suppose that $\rho >  \frac 12$ and $s <  0$
satisfy 
\eqref{regR1}.
Then, 
the driver $Y^\ld$  in \eqref{Xld2}
belongs to $  \cX^{0, \g}_2([0, \ld^3T] \times \R)$ defined in~\eqref{X1}, 
satisfying the bound\textup{:}
\begin{align}
 \| Y^{\ld}_{t,r}  \|_{\cL_2 ( L^2(\R)) }
 \les
 \ld^{ \frac32  -3\g} 
  \|\Phi^{w}\|_{\W_{T}^{\rho, \g}}
 |t-r|^{\g}
\label{RW1}
\end{align}

\noi
for any $0\leq r<t \leq \ld^3T$.

Moreover, given any 
$v_0 \in L^2(\R)$ and  $t_0 \in [0, \ld^3 T]$, 
there exist $C_0 > 0$, $\ta > 0$, both independent of $v_0$ and $t_0$, 
and a unique solution $Iu$
to 
the scaled modulated $I$-KdV equation~\eqref{kdv4} 
on the real line
with $I u^\ld (t_0)= v_0$
on the time interval $[t_0, t_0 + \tau]\cap [0, \ld^3 T]$,  
belonging to the class\textup{:}
\begin{align*}
 \cD_{w^\ld}^0\big( ([t_0, t_0 + \tau] \cap [0, \ld^3 T]) \times \R\big), 
\end{align*}

\noi
where the local existence time 
$\tau = \tau \big(\| v_0 \|_{L^2(\R)}\big)> 0$ satisfies
\begin{align*}
\tau  \ge C_0 \|Y^\ld\|_{\cX^{0,\g}_2([0, \ld^3 T]\times \R)}^{-\ta} \big( 1 + \| v_0 \|_{L^2(\R)}\big)^{-\ta}.
\end{align*}

\noi
Moreover, given any $0 < \al < \g$
with $\al + \g > 1$, 
there exists $C_\al > 0$ such that 
\begin{align*}
 \| I\uu^\ld \|_{\CC^\al ([t_0, t_0 + \tau];L^2(\R))}\le C_\al \|v_0\|_{L^2(\R)}.
\end{align*}

\end{corollary}

\begin{proof}

The bound \eqref{RW1}
follows from Lemma \ref{LEM:Rtri2}, 
 the interpolation lemma (see \cite[Lemma~12.1]{CKSTT04}),
 and 
the bound
\begin{align*}
 M (\bar \xi) \les  \frac{\jb{\xi}^{s}}{\jb{\xi_1}^{s}\jb{\xi_2}^{s}}
\end{align*}

\noi
for any
$\xi, \xi_1, \xi_2 \in \R$
with $\xi = \xi_1 + \xi_2$, 
where $ M (\bar \xi) =  M (\xi, \xi_1, \xi_2)$
is as in \eqref{IF3}.
 Then, from 
 Lemma  \ref{LEM:OBS1}\,(i), 
 we conclude that 
$Y^\ld \in  \cX^{0, \g}_2([0, \ld^3T] \times \R)$.

By noting that 
$ \ld^{ \frac32  -3\g} \le 1$
for any $\ld \ge 1$, 
the local well-posedness claim 
follows
from the regularity of the driver $Y^\ld$ and Proposition \ref{PROP:main}.
We omit details.
\end{proof}

Next, we present a key 
commutator estimate.
While we proceed as in the proof of 
Proposition \ref{PROP:com}, 
we provide details
since it involves changes of variables, 
unique to the current real line setting.

\begin{lemma}
\label{LEM:com2}

Let $\ld \ge 1$.
Given $\rho  >  \frac 12$,  $\frac 12 < \g < 1$,  and $T> 0$, 
let  $w$ be $(\rho,\g)$-irregular
on $[0, T]$ 
 in the sense of Definition~\ref{DEF:ir}, 
 and let $X^\ld$ be as in \eqref{Xld1}.
Fix  $s < 0$ satisfying~\eqref{regR1}.
Then, we have 
\begin{align}
\begin{aligned}
 \|\com_{t, r}(f_1, f_2) \|_{L^2(\R)}
 &  = \|IX^{\ld}_{t,r}( f _1, f _2) - X^{\ld}_{t,r} (I f _1,I f _2) \|_{L^2(\R)}\\
& \les
\wt    K(\ld, N )\| \Phi^{w} \|_{\W_{T}^{\rho, \g}}
      |t-r|^{\g}
\| I  f _1 \|_{L^2(\R)}   \| I  f _2 \|_{L^2(\R)}  
\end{aligned}
\label{RI1}
\end{align}

\noi
for any $0\leq r<t\leq \ld^3T$, $\ld \ge 1$, and $N \in \N$, 
where $ \wt K(\ld, N )$ is 
given by 
\begin{align}
\wt K(\ld, N )
= \begin{cases}
\ld^{3- 3\rho-3\g} N^{\frac 32 - 3\rho}, & \text{if $\frac 12 < \rho < \frac 32$}, \\
\ld^{- \frac 32 - 3\g} N^{-3} \sqrt{\log \jb{N} + \log \jb{\ld}}, & \text{if $\rho =  \frac 32$}, \\
\ld^{ - \frac 32 - 3\g } N^{-3}, & \text{if $ \rho >  \frac 32$}.
\end{cases}
\label{RI1a}
\end{align}

\end{lemma}

\begin{proof}

Proceeding as in 
the proof of Proposition \ref{PROP:com}, 
we have 
\begin{align*}
\text{LHS of \eqref{RI1}}
&\leq \bigg(\sum_{j = 3}^4 
B_{\ld, N}^{(j)} \| I  f_1 \|_{L^2(\R)}^2 \| I  f_2 \|_{L^2(\R)}^2\bigg)^\frac 12 , 
\end{align*}

\noi
where $B_{\ld, N}^{(j)}$, $j =  3, 4$, is given by 
\begin{align}
B_{\ld, N}^{(3)}
& =
\sup_{|\xi|\geq \frac N4} |\xi|^2
\int_{D_3(\xi)}
 \frac{|\Phi^{w^\ld}_{t,r}(\Xi_\KDV(\bar \xi)) |^2  |M_N(\bar \xi)|^2 }{ |m_{s, N}(\xi_1)|^2|m_{s, N}(\xi_2)|^2}
 d\xi_1,
\label{RI11}\\
B_{\ld, N}^{(4)}
& =
\sup_{\xi_1 \in \R}
\int_\R
\ind_{D_4(\xi)}\cdot |\xi|^2 
\frac{ |\Phi^{w^\ld}_{t,r}(\Xi_\KDV(\bar \xi))|^2 |M_N(\bar \xi)|^2}{|m_{s, N}(\xi_1)|^2 |m_{s, N}(\xi_2)|^2}
d\xi.
\label{RI11a}
\end{align}

\noi
Here, 
$M_N(\bar \xi) = M_N(\xi, \xi_1, \xi_2)$
and $D_j(\xi)$, $j =  3, 4$, 
are 
as in \eqref{I6a} and \eqref{I8}, respectively, 
with obvious modifications
to the current real line setting.

\medskip

\noi
$\bullet$ \textbf{Case 1:} on  $D_3$.\\
\indent
In this case, we have $|\xi|\sim |\xi_1|\ges N \ges |\xi_2|$.
Then, from~\eqref{RI11},  \eqref{RR4}, 
 \eqref{I11c}, 
and
a change of variables $\z_2 = \ld^3 \xi^2 \xi_2$, we have 
\begin{align}
\begin{split}
B_{\ld, N}^{(3)}
&  \les 
\ld^{6- 6\g}
  \|\Phi^{w} \|^2_{\W^{\rho,\g}_T} |t-r|^{2\g}
\sup_{|\xi|\geq \frac N4}
\int_{D_3(\xi)} \frac{|\xi_2|^2}{
\jb{\ld^3 \xi^2 \xi_2}^{2\rho}}\
d\xi_2\\
&  \les 
\ld^{-3- 6\g}
  \|\Phi^{w} \|^2_{\W^{\rho,\g}_T} |t-r|^{2\g}
\sup_{|\xi|\geq \frac N4}\xi^{-6}
\int_{|\z_2|\les \ld^3 \xi^2 N} \frac{|\z_2|^2}{
\jb{\z_2}^{2\rho}}\
d\z_2\\
 &\les 
\wt K(\ld, N)^2
 \| \Phi^{w} \|^2_{\W_{T}^{\rho, \g}}  
|t-r|^{2\g}
\end{split}
\label{RI12}
 \end{align}

 \noi
for any $0\leq r<t\leq \ld^3T$, 
where $\wt K(\ld, N)$ is as in \eqref{RI1a}. 
This proves
 \eqref{RI1} in this case.

\medskip

\noi
$\bullet$ \textbf{Case 2:} on  $D_4$.\\
\indent
From \eqref{RI11a} with \eqref{RR4}  and \eqref{I13x}, we have 
\begin{align}
\begin{split}
B_{\ld, N}^{(4)}
& \les 
\ld^{6- 6\g}N^{4s }
\| \Phi^{w} \|^2_{\W_{T}^{\rho, \g}} 
|t-r|^{2\g} \\
& \hphantom{XX}
\times
\sup_{\xi_1 \in \R}
\int_\R
\ind_{D_4(\xi)}\cdot
\frac{|\xi|^2  |m_{s, N}(\xi)|^2}{|\xi_1|^{2s}|\xi_2|^{2s}\jb{\ld^3 \xi \xi_1 \xi_2}^{2\rho}}d\xi \\
& \quad + 
\ld^{6 -6\g}
\| \Phi^{w} \|^2_{\W_{T}^{\rho, \g}} 
|t-r|^{2\g}
\sup_{\xi_1 \in \R}
\int_\R
\ind_{D_4(\xi)}
\frac{|\xi|^2 }{\jb{\ld^3 \xi \xi_1 \xi_2}^{2\rho}}d\xi \\
 & =: \1 + \II.
\end{split}
\label{RI13a}
\end{align}

Let us first estimate $\1$.


\noi
$\circ$ {\bf Subcase 2.a:}  $|\xi| \ll N$.\\
\indent
In this case, we have $m_{s, N}(\xi) = 1$
and  $|\xi_1|\sim |\xi_2|\ges N$.
Then, by 
a change of variables $\z = \ld^3 \xi_1^2 \xi$, we have 
\begin{align*}
\1
& \les 
\ld^{-3- 6\g}N^{ 4s}
 \| \Phi^{w} \|^2_{\W_{T}^{\rho, \g}}  
|t-r|^{2\g}
\sup_{|\xi_1|\ges N}|\xi_1|^{-6-4s}
\intt_{ |\z|\les \ld^3 \xi_1^2N } \frac{\z^{2}}{\jb{\z}^{2\rho}}d \z\\
 &\les 
\wt K(\ld, N)^2
 \| \Phi^{w} \|^2_{\W_{T}^{\rho, \g}}  
|t-r|^{2\g}
\end{align*}

\noi
for any $0\leq r<t\leq \ld^3T$, 
thus yielding \eqref{RI1}.
Here, 
in estimating the integral, 
we used the fact that 
(i)~$s+ \rho \ge 0$
when $\frac 12< \rho < \frac 32$
and 
(ii)~$s > -\frac 32$ when $\rho \ge \frac 32$.

\medskip

\noi
$\circ$ {\bf Subcase 2.b:}  $|\xi| \ges N$.\\
\indent
In this case, we have   $|m_N(\xi)| \les  N^{-s} |\xi|^{s}$.
Without loss of generality, assume that $|\xi|\les |\xi_1|$.
We note that, for $\rho \ge \frac 32$, we have 
\begin{align}
\ld^{3- 3\rho-3\g} N^{\frac 32 - 3\rho}
\le
\ld^{ - \frac 32  - 3\g} N^{-3}.
\label{RI14a}
\end{align}

\noi
Then, 
proceeding as in \eqref{I16x}
with 
 \eqref{RI13a}, 
\eqref{I16}, 
\eqref{I16x1}, 
and 
\eqref{RI14a} with \eqref{RI1a}, we have
\begin{align}
\begin{split}
\1
&\les 
\ld^{6-  6\rho  - 6\g}N^{2s }
\| \Phi^{w} \|^2_{\W_{T}^{\rho, \g}} 
|t-r|^{2\g} 
\sup_{\xi_1 \in \R}
\intt_{|\xi|\ges N}
\ind_{D_4(\xi)}\cdot
\frac{|\xi|^{2+2s-2\rho} }{|\xi_1|^{2s+2\rho}|\xi_2|^{2s+2\rho}}d\xi \\
&\les 
\ld^{6-  6\rho  - 6\g}N^{2s }
\| \Phi^{w} \|^2_{\W_{T}^{\rho, \g}} 
|t-r|^{2\g} 
\sup_{\xi_1 \in \R}
\intt_{|\xi|\ges N}
\ind_{D_4(\xi)}\cdot
\frac{1 }{\jb{\xi}^{4\rho - 2}\jb{\xi_2}^{2s+2\rho}}d\xi \\
 &\les 
\ld^{6-6\rho -6\g}N^{3-6\rho }
 \| \Phi^{w} \|^2_{\W_{T}^{\rho, \g}}  
|t-r|^{2\g}\\
 &\les 
\wt K(\ld, N)^2
 \| \Phi^{w} \|^2_{\W_{T}^{\rho, \g}}  
|t-r|^{2\g}, 
\end{split}
\label{RI14x}
\end{align}

\noi
provided that 
  \eqref{K5b} holds.

\medskip

Next, we  estimate $\II$.

\medskip

\noi
$\circ$ {\bf Subcase 2.c:}  $|\xi| \ll N$.\\
\indent
In this case, we have 
  $|\xi_1|\sim |\xi_2|\ges N$.
Then, proceeding as in \eqref{RI12}
with 
a change of variables $\z = \ld^3 \xi_1^2 \xi$, we have 
\begin{align*}
\II
& \les 
\ld^{-3- 6\g}
 \| \Phi^{w} \|^2_{\W_{T}^{\rho, \g}}  
|t-r|^{2\g}
\sup_{|\xi_1|\ges N}|\xi_1|^{-6}
\intt_{ |\z|\les \ld^3 \xi_1^2N } \frac{\z^{2}}{\jb{\z}^{2\rho}}d \z\\
 &\les 
\wt K(\ld, N)^2
 \| \Phi^{w} \|^2_{\W_{T}^{\rho, \g}}  
|t-r|^{2\g}
\end{align*}

\noi
for any $0\leq r<t\leq \ld^3T$, 
thus yielding \eqref{RI1}.

\medskip

\noi
$\circ$ {\bf Subcase 2.d:}  $|\xi| \ges N$.\\
\indent
Without loss of generality, assume that $N \les |\xi| \les |\xi_1|$.
Then, by noting
that $\jb{\ld^3 \xi \xi_1 \xi_2} \sim 
\ld^3 |\xi \xi_1 \xi_2|
\ges \ld^3 \xi^2 |\xi_2|$, 
we can repeat the computation in \eqref{RI14x}
and obtain \eqref{RI1}.
\end{proof}

We conclude this section by briefly discussing a proof of Theorem 
\ref{THM:2}\,(ii).
In view of 
\eqref{RR1}, Corollary \ref{COR:RI1}, 
and Lemma \ref{LEM:com2}
(compare \eqref{RI1a} with \eqref{I5a}), 
there is no difference from the periodic case, when $\frac 12< \rho \le \frac 32$, 
yielding 
global well-posedness of the modulated KdV \eqref{kdv1} on the real line
for the same range of $s$.
When $\rho > \frac 32$, 
by substituting $\rho = \frac 32$ in the last case of~\eqref{IT0b}, 
we obtain the condition $s > - \frac 3{3- 2\g}$ (which is smaller than $-\frac 32$ for $\frac 12 < \g < 1$).
In view of the regularity restriction $s > - \frac 32$
from the local well-posedness (see~\eqref{regR1}), 
we conclude that global well-posedness
holds for  $s > - \frac  32$.
This concludes the proof of Theorem~\ref{THM:2}\,(ii).

\section{Invariance of the white noise under the modulated KdV}
\label{SEC:WN}

In this section, we briefly discuss 
invariance of the white noise under
the modulated KdV, BO, ILW, and scaled ILW equations on the circle
(Theorems \ref{THM:white1} and \ref{THM:white2}).
For simplicity, we restrict our attention to the case
of the modulated KdV \eqref{kdv1} on $\T$.
As mentioned in Subsection \ref{SUBSEC:1.5}, 
Theorem \ref{THM:white1} follows once we prove the following proposition.

\begin{proposition}\label{PROP:white2}
Let $N \in \N$. The white noise $\mu$ defined in \eqref{white1}
is invariant under
the truncated modulated KdV dynamics \eqref{kdv1x}.
\end{proposition}

Given $N \in \N$, 
set
\begin{align*}
u_N = \P_N u^N\qquad \text{and}\qquad
u_N^\perp = \P_N^\perp u^N, 
\end{align*}

\noi
where  $\P_N^\perp = \Id - \P_N$.
Then, the truncated modulated KdV \eqref{kdv1x}
decouples into the low frequency nonlinear but finite-dimensional dynamics, 
satisfied by $u_N$:
\begin{equation}
\begin{cases}
\dt u_N+  \dx^3 u_N \cdot \dt w =\dx \P_N(u_N^2) \\
u_N(0) = \P_N u^N(0)
\end{cases}
\label{kdv1x2}
\end{equation}

\noi
and the high frequency linear dynamics, 
satisfied by $u_N^\perp$:
\begin{equation}
\begin{cases}
\dt u_N^\perp+  \dx^3 u_N^\perp \cdot \dt w = 0 \\
u_N^\perp(0) = \P_N^\perp u^N(0).
\end{cases}
\label{kdv1x3}
\end{equation}

\noi
Let $\rho > \frac 23$.
Then, by a slight modification of the proof of Theorem \ref{THM:1}
applied to the low frequency part \eqref{kdv1x2}
along with the conservation of the (truncated) $L^2$-norm, 
we see that 
the truncated modulated KdV \eqref{kdv1x}
is globally well-posed in $H^s(\T)$ for any $s \in \R$,
since the high frequency dynamics is linear
(which is globally well-posed in $\P_N^\perp H^s(\T)$ for any $s \in \R$).

Let us write the white noise $\mu$ as 
\begin{align}
 \mu 
  =  \mu^{N, \low} \otimes \mu^{N, \high}, 
\label{white1b}
\end{align}

\noi
where $\mu^{N, \low} = (\P_N)_* \mu$
and  $\mu^{N, \high} = (\P_N^\perp)_* \mu$
are the pushforward image
measures under the maps $\P_N$ and $\P_N^\perp$, 
respectively.
From the rotational invariance of complex-valued Gaussian random variables, we see that 
 $\mu^{N, \high}$ is invariant
under the high frequency linear dynamics~\eqref{kdv1x3}
whose action on each Fourier coefficient is given by 
$\F(u_N^\perp)(0, n) = g_n \mapsto  \F(u_N^\perp)(t, n) = e^{i n^3 w(t)}g_n$.

By Plancherel's theorem with $\ft u(-n) = \cj{\ft u(n)}$, 
we can write  the truncated white noise $\mu^{N, \low} = (\P_N)_*\mu$ as
\begin{align*}
d\mu^{N, \low}
& = Z_N^{-1}
\exp\bigg(- \frac 12 \int_\T (\P_N u)^2 dx_\T\bigg) d\P_N u\\
& = Z_N^{-1}
\exp\bigg(- \sum_{n = 1}^N |\ft u(n)|^2 \bigg) \prod_{n = 1}^N d\ft u(n), 
\end{align*}

\noi
where $d\ft u(n)$, $n = 1, \dots, N$, is the Lebesgue measure on $\C\cong \R^2$.
Arguing as in 
Proposition~\ref{PROP:L2}, 
we see that the truncated $L^2$-norm
$\|\P_N u\|_{L^2}$ 
(= the $\C^N$-Euclidean distance on the Fourier side)
is conserved under the flow of \eqref{kdv1x2}.
Hence, invariance of the truncated white noise 
$\mu^{N, \low}$ follows once we prove invariance of the Lebesgue measure
$\prod_{n = 1}^N d\ft u(n)$.

Very formally, we can write the (time-dependent) Hamiltonian for \eqref{kdv1x2}
as 
\begin{align*}
H_N(u) = \frac 12 \int_\T (\dx \P_N u)^2\dt w (t) dx_\T
+ \frac 13 \int_\T (\P_N u)^3dx_\T.
\end{align*}

\noi
Unfortunately, the expression above does not make sense
since $\dt w (t)$ is not well defined for fixed $t$. 
Similarly, by viewing \eqref{kdv1x2} on the Fourier side, 
one may try to verify invariance of the Lebesgue measure
$\prod_{n = 1}^N d\ft u(n)$ by directly computing
the divergence of the vector field; see, for example, 
\cite[Lemma 4.1]{NORS}.
However, the divergence involves $\dt w(t)$
and hence this approach fails as well.
(If we pretend that $\dt w(t)$ makes sense for fixed $t$, 
the divergence is indeed $0$ since the terms
with $\dt w(t)$ cancel out.)

A correct approach is to study the equations satisfied
by the modulated interaction representation.
Fix  $N \in \N$.
By 
applying $\uw(t)^{-1}$ to the truncated modulated KdV \eqref{kdv1x}, 
we obtain the following 
 truncated version of \eqref{kdv1a}:
\begin{align}
\dt \uu^N = \P_N\uw(t)^{-1} \dx\big( (\P_N \uw(t) \uu^N)^2\big).
\label{kdv1b}
\end{align}

\noi
Since the action of $\uw(t)^{-1}$ on each Fourier coefficient
is the multiplication by $e^{-in^3 w(t)}$, 
we see that the white noise is invariant under $\uw(t)^{-1}$, 
namely
\begin{align*}
\mu = (\uw(t)^{-1})_* \mu.
\end{align*}

\noi
See also Section 4 in \cite{OTz}.
Thus, Proposition \ref{PROP:white2}
follows once we prove invariance of the white noise~$\mu$
under \eqref{kdv1b}.

As before, we decouple the truncated modulated KdV dynamics \eqref{kdv1b}
 into
 the finite-dimensional nonlinear dynamics for the low frequency part
  $\uu_N = \P_{N} \uu^N$:
\begin{align}
\begin{cases}
\dt \uu_N = \P_N\uw(t)^{-1} \dx\big( ( \uw(t) \uu_N)^2\big)\\
\uu_N|_{t = 0} = \P_{N} u_0^\o, 
\end{cases}
\label{kdv1d}
\end{align}

\noi
and the constant  dynamics for the high frequency part 
$\uu_N^\perp = \P_N^\perp \uu^N$:
\begin{align}
\begin{cases}
\dt \uu_N^\perp = 0\\
\uu_N^\perp |_{t = 0} = \P_{N}^\perp u_0^\o.
\end{cases}
\label{kdv1e}
\end{align}

\noi
Thanks to \eqref{white1b}, 
invariance of the 
white noise $\mu = \mu^{N, \low} \otimes \mu^{N, \high}$
under \eqref{kdv1b}
follows once we prove invariance
of $\mu^{N, \low}$
(and $\mu^{N, \high}$) under \eqref{kdv1d}
(and \eqref{kdv1e}, respectively).
Obviously, $\mu^{N, \high}$ is invariant under 
the constant dynamics \eqref{kdv1e}.
In view of the $L^2$-conservation for~\eqref{kdv1d}
(which follows from the $L^2$-conservation for \eqref{kdv1x2}
and the $L^2$-unitarity of $\uw(t)^{-1}$), 
invariance of the truncated white noise $\mu^{N, \low}$ 
(and hence Proposition \ref{PROP:white2} and Theorem \ref{THM:white1})
follows
once we prove invariance of the Lebesgue measure
$\prod_{n = 1}^N d\ft \uu_N(n)$ under~\eqref{kdv1d}.

\begin{lemma}\label{LEM:div}

The Lebesgue measure $\prod_{n = 1}^N d\ft \uu_N(n)$ on $\C^N$
is invariant under the truncated dynamics \eqref{kdv1d} viewed on the Fourier side.

\end{lemma}

\begin{proof}

Let $p_n(t) =\Re  \ft \uu_N(t, n)$
and $q_n(t) = \Im \ft \uu_N( t, n)$ for $1\le  |n| \le N$.
Since $\uu_N$ is real-valued,  we have 
\begin{align}
p_{-n} = p_n \qquad \text{and} \qquad q_{-n} = - q_n.
\label{div1}
\end{align}

\noi
Then, \eqref{kdv1d} can be rewritten as 
\begin{align}
\dt p_n = P_n\qquad \text{and}\qquad
\dt q_n = Q_n 
\label{div2}
\end{align}

\noi
for $n = 1, \dots N$, 
where $P_n$ and $Q_n$ are defined by 
\begin{align}
\begin{split}
P_n 
& = \frac{in}{2}
\bigg[ \sum_{\substack{1\le |n_1|, |n_2|\le N\\n = n_1 + n_2}}
e^{i \Xi_\KDV (\bar n)w(t)} (p_{n_1} + i q_{n_1})(p_{n_2} + i q_{n_2})\\
& \hphantom{XXX}
-  \sum_{\substack{1\le |n_1|, |n_2|\le N\\n = n_1 + n_2}}
e^{- i \Xi_\KDV (\bar n)w(t)} (p_{n_1} - i q_{n_1})(p_{n_2} - i q_{n_2})\bigg], \\
Q_n 
& = \frac{n}{2}
\bigg[ \sum_{\substack{1\le |n_1|, |n_2|\le N\\n = n_1 + n_2}}
e^{i \Xi_\KDV (\bar n)w(t)} (p_{n_1} + i q_{n_1})(p_{n_2} + i q_{n_2})\\
& \hphantom{XXX}
+  \sum_{\substack{1\le |n_1|, |n_2|\le N\\n = n_1 + n_2}}
e^{- i \Xi_\KDV (\bar n)w(t)} (p_{n_1} - i q_{n_1})(p_{n_2} - i q_{n_2})\bigg], \\
\end{split}
\label{div3}
\end{align}

\noi
where $\Xi_\KDV (\bar n)$ is as in \eqref{K3a}.
From \eqref{div3} with \eqref{div1} (with $p_0 = q_0 = 0$), we have 
\begin{align}
\begin{split}
\dd_{p_n} P_n 
& = in \Big(e^{i \Xi_\KDV (n, -n, 2n)w(t)} (p_{2n} + i q_{2n})\\
& \hphantom{XXl}
 - e^{- i \Xi_\KDV (n, -n, 2n)w(t)} (p_{2n} - i q_{2n})\Big), \\
\dd_{q_n} Q_n 
& = n \Big(e^{i \Xi_\KDV (n, -n, 2n)w(t)}( -i) (p_{2n} + i q_{2n})\\
& \hphantom{XXl}
+ e^{- i \Xi_\KDV (n, -n, 2n)w(t)} i (p_{2n} - i q_{2n})\Big).
\end{split}
\label{div4}
\end{align}

\noi
Let $A(\bar p , \bar q) = A(p_1, \dots, p_N, q_1, \dots, q_N)$
be the vector field of the dynamics \eqref{div2} defined by 
\begin{align*}
A(\bar p , \bar q) = (P_1, \dots, P_N, Q_1, \dots, Q_N).
\end{align*}

\noi
Then, from \eqref{div4}, we have
\begin{align*}
\text{div}_{\bar p, \bar q}  A(\bar p , \bar q) 
= \sum_{n = 1}^N
(\dd_{p_n} P_n + \dd_{q_n} Q_n)
= 0.
\end{align*}

\noi
Therefore
from Liouville's theorem, 
we conclude that 
the Lebesgue measure $\prod_{n = 1}^N d\ft \uu_N(n)
= \prod_{n = 1}^N d p_n dq_n$ is invariant under the truncated dynamics \eqref{kdv1d}.
\end{proof}

\section{Stochastic modulated  KdV}
\label{SEC:SK}

\subsection{Local well-posedness}
In this subsection, we prove
local well-posedness of the stochastic modulated KdV \eqref{skdv1}
(Theorem~\ref{THM:skdv1}\,(i)).

Let us first give a precise meaning of the stochastic term $\Psi$ in~\eqref{sconv1}.
 Let  $W$ denote a cylindrical Wiener process on $L^2_0(\T)$:
\begin{align*}
W(t)
 = \sum_{n \in \Z^*} B_n (t) e_n, 
\end{align*}

\noi
where 
$\{ B_n \}_{n \in \Z^*}$ 
is defined by 
$B_n(t) = \jb{\xi, \ind_{[0, t]} \cdot e_n}_{t, x}$.
Here, $\jb{\cdot, \cdot}_{t, x}$ denotes 
the duality pairing on $ \R_+\times \T$.
As a result, 
we see that $\{ B_n \}_{n \in \Z^*}$ is a family of mutually independent complex-valued
Brownian motions conditioned  that $B_{-n} = \cj{B_n}$, $n \in \Z^*$. 
Note that we have, for any $n \in \Z^*$,  
 \[\text{Var}(B_n(t)) = \E\big[
 \jb{\xi, \ind_{[0, t]} \cdot e_n}_{t, x}\cj{\jb{\xi, \ind_{[0, t]} \cdot e_n}_{ t, x}}
 \big] = \|\ind_{[0, t]} \cdot e_n\|_{L^2_{t, x}}^2 = t.\]

\noi

\noi
Then, we can write the stochastic term $\Psi$ in \eqref{sconv1}
as 
\begin{align}
\Psi(t) = \int_0^t \uw(t')^{-1}\phi dW(t')
= \sum_{m,n\in\Z^*} \phi_{mn} I_{t,0}(m, n)e_m, 
\label{sconv2}
\end{align}

\noi
where $I_{t_1,t_2}(m,n)$, $m,n\in\Z^*$,
is the Wiener integral given by
\begin{align}
I_{t_1,t_2}(m,n)
= \int_{t_2}^{t_1} e^{-im^3w(t')}dB_n(t').
\label{sconv3}
\end{align}

\noi
Note that
$I_{t_1,t_2}(-m,-n)=\cj{I_{t_1,t_2}(m,n)}$,
$m,n\in\Z^*$.

\begin{lemma}\label{LEM:sto1}
Given $\rho >0$,  $\frac 12 < \g < 1$, and $T> 0$, 
let  $w$ be $(\rho,\g)$-irregular on $[0, T]$ in the sense of Definition~\ref{DEF:ir}.
Given $s \in \R$, 
let $\phi \in \HS(L^2_0(\T); H^s_0(\T))$.
Then, 
given any $0<  \al < \frac 12$, 
$\Psi$ defined in \eqref{sconv2}
almost surely belongs to 
 $ \CC^\al([0, T]; H^s(\T))$.

\end{lemma}

\begin{proof}

By writing 
\begin{align}
\phi(e_n) = \sum_{m \in \Z^*} \phi_{mn} e_m, \quad n \in \Z^*,
\label{sconv3a} 
\end{align}

\noi
the Hilbert-Schmidt norm of $\phi$ is given by 
\begin{align}
\| \phi\|_{\HS(L^2_0; H^s_0)}
= \bigg(\sum_{n \in \Z^*}\|\phi(e_n)\|_{H^s_0}^2\bigg)^\frac 12 
= \bigg(\sum_{n, m  \in \Z^*}\jb{m}^{2s} |\phi_{mn}|^2 \bigg)^\frac 12 .
\label{sconv4}
\end{align}

\noi
\noi
Then, from the equivalence of the $L^p(\O)$-norms 
of a Gaussian random variable, \eqref{sconv2}, 
\eqref{sconv3a}, 
\eqref{sconv3}, 
\eqref{sconv4}, 
and the Wiener isometry, we have
\begin{align}
\begin{aligned}
\big\| \|\Psi(t_1) - \Psi(t_2)\|_{H^s}\big\|_{L^p(\O)}
& \sim_p  \big\| \|\Psi(t_1) - \Psi(t_2)\|_{H^s}\big\|_{L^2(\O)}\\
& = \bigg(\sum_{m \in \Z^*} \jb{m}^{2s}
\E\bigg[\Big|\sum_{n \in \Z^*} I_{t_1, t_2}(m, n)\phi_{mn}\Big|^2\bigg]\bigg)^\frac 12 \\
& =  \bigg(\sum_{n, m \in \Z^*} \jb{m}^{2s}|\phi_{mn}|^2
\E\big[|I_{t_1, t_2}(m, n)|^2\big]\bigg)^\frac 12 \\
& \sim \|\phi\|_{\HS(L^2_0; H^s_0)}|t_1 - t_2|^\frac 12 
\end{aligned}
\label{sconv5}
\end{align}

\noi
for any $t_1 \ge t_2 \ge 0$ and any finite $p \ge 1$.
Then, given $0 < \al < \frac 12$, 
by applying
the Kolmogorov continuity criterion (\cite[Theorem 2.1]{Baldi})
to \eqref{sconv5}
with  $p > \frac{2}{1-2\al}$, 
 we conclude that 
 $\Psi \in \CC^\al([0, T]; H^s(\T))$, 
 almost surely.
 \end{proof}

We now present 
a proof of Theorem \ref{THM:skdv1}\,(i).

\begin{proof}[Proof of Theorem \ref{THM:skdv1}\,(i)]
Fix $\al = \frac 12 - \eps$ for some small $\eps > 0$
such that $\al + \g > 1$.
Define the map 
$\G^\Psi  = \G^\Psi_{X, u_0}$ on $\cC^\al([0, \tau]; H^s(\T))$
by setting
\begin{align*}
\G^\Psi(\uu)(t) 
& = u_0 + \I^X(\uu)(t) + \Psi(t)\\
& = \G(\uu) + \Psi(t), 
\end{align*}

\noi
where $\G(\uu)$ is as in \eqref{YD3}.
Then, 
by setting 
an almost surely finite radius $R = 2 \big(\|u_0\|_{H^s}+ \|\Psi\|_{\CC^\al_T H^s_x}\big) $
and choosing 
an almost surely positive local existence time
$\tau = 
\tau\big(R, \|X\|_{\cX^{s, \g}_k(T)}\big) 
= 
\tau\big(\|u_0\|_{H^s}, 
\|\Psi\|_{\CC^\al_T H^s_x}, 
\|X\|_{\cX^{s, \g}_k(T)}\big)$
sufficiently small, 
it follows from~\eqref{YD5a} and \eqref{YD14} that 
\begin{align*}
 \|\G^\Psi(\uu)\|_{\CC^\al_\tau H^s_x}
&  \le \|u_0\|_{H^s} 
+ \|\Psi\|_{\CC^\al_T H^s_x} +  C \tau^{\g - \al}
\|X\|_{\cX^{s, \g}_k(T)}
(1+R) ^{k-1}
 \|\uu\|_{\CC^\al_{\tau } H^s_x} \\
& \le R
\end{align*}

\noi
and 
\begin{align*}
 \|\G^\Psi(\uu^1)- \G^\Psi(\uu^2)\|_{\cC^\al _{\tau} H^s_x}
&   \le
C
\tau^{\g - \al}\|X\|_{\cX^{s, \g}_k(T)} (1+ R)^{k-1}\|\uu^1- \uu^2 \|_{\CC^\al _\tau H^s_x}\\
& \le \frac 12 \|\uu^1- \uu^2 \|_{\CC^\al _\tau H^s_x}
\end{align*}

\noi
for any  $\uu, \uu^1, \uu^2
\in B_R\subset 
\cC^\al([0, \tau]; H^s(\M))$, 
from which
we conclude 
that $\G^\Psi$ is a contraction on $B_R$.
This proves 
 local well-posedness
of the stochastic modulated KdV \eqref{skdv1}
on $\T$.
\end{proof}

\subsection{Global dynamics}

In this subsection, we briefly 
go over a proof of Theorem~\ref{THM:skdv1}\,(ii).
Since the argument is identical to that in \cite{OQS}, 
we only
sketch the main ideas.

Given $N \in \N$, 
consider 
the following truncated stochastic modulated KdV:
\begin{equation}
\label{KdV4}
\begin{cases}
\dt u^N+  \dx^3 u^N \cdot \dt w =\dx \P_N\big((\P_Nu^N)^2\big) +  \xi\\
u^N|_{t = 0} = u_0^\o, 
\end{cases}
\end{equation}

\noi
where $u_0^\o$ is the white noise given in \eqref{series1}.
In fact,
as in Section~\ref{SEC:WN}, 
 it is more convenient to work at the level of the modulated interaction representation
 $\uu(t) = \uw(t)^{-1} u(t)$.
Then,~\eqref{KdV4} becomes
\begin{equation}
\label{KdV4a}
\begin{cases}
\dt \uu^N = \P_N\uw(t)^{-1} \dx\big( (\P_N \uw(t) \uu^N)^2\big) +  \uw(t)^{-1} \xi\\
\uu^N|_{t = 0} = u_0^\o.
\end{cases}
\end{equation}

\noi
As in Section~\ref{SEC:WN}, the truncated dynamics \eqref{KdV4a}
 decouples into 
 the finite-dimensional nonlinear dynamics for the low frequency part 
 $\uu_N = \P_{N} \uu^N$:
\begin{align}
\begin{cases}
\dt \uu_N = \P_N\uw(t)^{-1} \dx\big( ( \uw(t) \uu_N)^2\big)
+ \P_N  \uw(t)^{-1} \xi\\
\uu_N|_{t = 0} = \P_{N} u_0^\o, 
\end{cases}
\label{KdV5}
\end{align}

\noi
and the linear dynamics for the high frequency part 
$\uu_N^\perp = \P_N^\perp \uu^N$:
\begin{align}
\begin{cases}
\dt \uu_N^\perp =  \P_N^\perp  \uw(t)^{-1} \xi \\
\uu_N^\perp |_{t = 0} = \P_{N}^\perp u_0^\o.
\end{cases}
\label{KdV6}
\end{align}

\noi
By a minor modification of \cite[Lemma 2.1]{OQS}, 
we see  that both \eqref{KdV5} and \eqref{KdV6}
are globally well-posed, 
 which in particular implies that \eqref{KdV4a} is globally well-posed.
For $t_1 \ge t_2 \ge 0$,
we denote by 
 $\Phi^{N, \low}_{t_1, t_2}$
 and $\Phi^{N, \high}_{t_1, t_2}$  
 the solution maps for \eqref{KdV5} and \eqref{KdV6}
sending data $\varphi$  at time $t_2$ to the solutions 
$\Phi^{N, \low}_{t_1, t_2}\varphi$
and $\Phi^{N, \high}_{t_1, t_2}\varphi$ at time $t_1$.
We let $P^{N,\low}_{t_1, t_2}$
and $P^{N,\high}_{t_1, t_2}$ denote the transition semigroups
for \eqref{KdV5} and \eqref{KdV6}, respectively, 
defined as in \eqref{trans1}, where the expectation 
is taken over the noise restricted to the time interval $[t_2, t_1]$.
We also use 
$\Phi^{N}_{t_1, t_2}$
and $P^{N}_{t_1, t_2}$
to denote the solution map and the transition semigroup for the truncated
stochastic modulated KdV \eqref{KdV4a}.

Given $\s \ge 0$, let $\mu_\s$ be the white noise measure (with variance $\s$)
as in \eqref{mu1}.
Then, as in~\eqref{white1b}, we can write $\mu_\s$ as
\begin{align*}
 \mu_\s 
 & =  \mu_\s^{N, \low} \otimes \mu_\s^{N, \high}, 
\end{align*}

\noi
where $\mu_\s^{N, \low} = (\P_N)_* \mu_\s$ and 
$\mu_\s^{N, \high} = (\P_N^\perp)_* \mu_\s$
are 
the pushforward image measures of $\mu_\s$ under
$\P_N$ and $\P_N^\perp$, respectively.
The high frequency dynamics \eqref{KdV6}
is linear and it is easy to verify 
that 
\begin{align}
(P^{N, \high}_{t_1, t_2})^* \mu^{N, \high}_{1+t_2}
=   \mu^{N, \high}_{1+t_1 }.
\label{mu3}
\end{align}

\noi
See \eqref{evo2} (and \eqref{evo1})
for the definition of the adjoint $(P^{N, \high}_{t_1, t_2})^*$.
By writing it on the Fourier side, 
we see that the low frequency dynamics \eqref{KdV5}
is nothing but a finite-dimensional system of SDEs,
which can be viewed as the superposition of 
the finite-dimensional modulated KdV dynamics:
\begin{align}
\dt \uu_N = \P_N\uw(t)^{-1} \dx\big( ( \uw(t) \uu_N)^2\big)
\label{KdV7}
\end{align}

\noi
and the linear stochastic dynamics:
\begin{align}
\dt \uu_N = 
 \P_N  \uw(t)^{-1} \xi.
\label{KdV8}
\end{align}

\noi
As we saw in Section \ref{SEC:WN}, 
the truncated white noise $\mu_\s^{N, \low}$ (with any variance $\s$)
is invariant under  the former dynamics \eqref{KdV7}, 
while 
the latter dynamics~\eqref{KdV8}
increases the variance of the white noise initial data by the length of the time interval under consideration.
Then, in view of the Lie-Trotter product formula \eqref{tro}, 
we see that 
\begin{align}
(P^{N, \low}_{t_1, t_2})^* \mu^{N, \low}_{1+t_2}
=   \mu^{N, \low}_{1+t_1 }.
\label{mu4}
\end{align}


\noi
Putting \eqref{mu3} and \eqref{mu4} together
with the invariance of the white noise $\mu_\s$ (with any variance $\s$)
under $\uw(t)$, 
we then obtain the following proposition.

\begin{proposition}\label{PROP:finite}
Let $N \in \NB$. 
Then, for any $t_1 \ge t_2 \ge 0$, we have 
\begin{align*}
(P^{N}_{t_1, t_2})^* \mu_{1+t_2}
=   \mu_{1+t_1 }, 
\end{align*}

\noi
where 
  $P^{N}_{t_1, t_2}$
is  the transition semigroup for the truncated dynamics \eqref{KdV4a}.
Namely, $\{\mu_{1+t}\}_{t \in \R_+}$
is an evolution system of measures
for the truncated dynamics \eqref{KdV4a}
and also for the truncated
stochastic modulated KdV equation \eqref{KdV4}.
\end{proposition}

In \cite{OQS}, 
the authors proved an analogous claim
for the truncated stochastic (unmodulated) KdV 
by studying 
 the  Kolmogorov forward equation 
(= the Fokker-Planck equation)
 for the distribution 
 of the Fourier coefficients
 of a solution $u_N(t)$ to the low frequency dynamics
 and showing that the truncated white noise
 is the unique  solution to the 
  Kolmogorov forward equation;
 see \cite[Section 2]{OQS}.
We point out that for such an argument, 
we need to make use of the computations in the proof of Lemma \ref{LEM:div}.
Namely, in the current modulated setting, 
we need to work at the level of the modulated interaction representation.

Once we have
Proposition \ref{PROP:finite} at hand, 
 a variant of 
 Bourgain's invariant measure argument~\cite{BO94}
 yields
 the following 
probabilistic  uniform (in $N$) growth bound 
on the solutions to the truncated dynamics \eqref{KdV4a}.

\begin{proposition}\label{PROP:main1}
Given $\rho>\frac23$ and $\frac12<\g<1$,  let $w$ be
$(\rho,\g)$-irregular on $\R_+$.
Let $s<-\frac12$ satisfy~\eqref{reg3a}.
Given $N \in \NB$, let $\uu^N$ be the global solution to the truncated dynamics~\eqref{KdV4a}
with the mean-zero white noise initial data $u_0^\o$ in \eqref{series1}.
Then, 
given any  $T\gg1$ and $0 < \eps\ll 1$, there exists a set  $\Omega_{T,\eps}(N)$ such that 
$\mathbb{P}(\Omega_{T,\eps}(N)^c)<\eps$ and
  \begin{align}
    \sup_{t\in [0,T]}\|\uu^N(t)\|_{H^s}&\le  C\big( \| \Phi^w \|_{\W_T^{\rho,\g}} \big)
 \sqrt {\log \frac 1 \eps}\sqrt{T \log T}
\label{growth1}
  \end{align}

\noi
on $\Omega_{T,\eps}(N)$, 
where $C:[0,\infty)\to(0,\infty)$ is a non-decreasing function, 
independent of $N \in \NB$, $T\gg1$, and $\eps \ll 1$.


\end{proposition}

Proposition~\ref{PROP:main1}
follows from a straightforward modification of the proof of 
 \cite[Proposition~4.1]{OQS} and thus we omit details.
The probabilistic uniform (in $N$) growth bound
\eqref{growth1}
together with the following lemma
on convergence of the Galerkin approximation 
yields
the almost sure global well-posedness
claimed in 
Theorem \ref{THM:skdv1}\,(ii).
 Finally, \eqref{th1}
 follows from Proposition \ref{PROP:finite}
 and 
the following lemma
(with $\phi = \Id$).

\begin{lemma}[Galerkin approximation]\label{LEM:Ga1}
Given $\rho > \frac 12$,  $\frac12< \g < 1$, and $T> 0$, 
let  $w$ be $(\rho,\g)$-irregular on $[0, T]$ in the sense of Definition~\ref{DEF:ir}.
Suppose that  $s \in \R$ satisfies~\eqref{reg3a}
and fix 
 $u_0 \in H^s(\T)$
and  $\phi \in \HS(L^2_0(\T); H^s_0(\T))$.
Then, as $N \to \infty$, 
the solution $\uu^N$
to the following truncated dynamics\textup{:}
\begin{equation}
\label{KdV9}
\begin{cases}
\dt \uu^N = \P_N\uw(t)^{-1} \dx\big( (\P_N \uw(t) \uu^N)^2\big) +  \uw(t)^{-1}\phi \xi\\
\uu^N|_{t = 0} = u_0 
\end{cases}
\end{equation}

\noi
converges to 
the modulated interaction $\uu(t) = \uw(t)^{-1}u(t)$
of the solution $u$
to 
 the stochastic modulated KdV equation~\eqref{skdv1}
 with $u|_{t = 0} = u_0$
in $\cC^\g([0, \tau]; H^s(\T))$, 
where $\tau > 0$ denotes the almost surely positive local existence time
from Theorem \ref{THM:skdv1}\,(i)
\textup{(}but 
possibly  made  smaller by a multiplicative constant\textup{)}.

\end{lemma}

It is worthwhile to note that while we have $\uu, \uu^N 
\in \cC^\al([0, \tau]; H^s(\T))$
only for   $  \al < \frac 12$, 
the difference 
 $\uu - \uu^N$ 
 belongs to a better space
 $\cC^\g([0, \tau]; H^s(\T))$ (with $\g > \frac 12$).

\begin{proof}[Proof of Lemma \ref{LEM:Ga1}]

Let $X$ and $X^N$ be as in 
\eqref{K1} and 
\eqref{XN},   respectively.
Then, 
from Proposition \ref{PROP:kdv1}\,(iv), 
we see that 
$X^N$ converges to  $X$ 
in  $ \cX^{s, \g}_2([0, T]\times \T)$ as $N \to \infty$.

Note that $\uu$ and $\uu^N$ satisfy the following equations:
\begin{align*}
\uu(t) & = u_0 + \I^{X}(\uu)(t) + \Psi, \\
\uu^N(t) & = u_0 + \I^{X^N}(\uu^N)(t) + \Psi, 
\end{align*}

\noi
where $\Psi$ is as in \eqref{sconv1} and \eqref{sconv2}.
In particular, $\uu - \uu^N$ satisfies \eqref{YDE3}.
Then, by repeating the proof of 
Proposition \ref{PROP:main}\,(iv), 
we obtain the desired convergence.
\end{proof}

\begin{ackno}\rm
The authors would like to thank Damiano Greco and Shao Liu
for their comments on Section 7.
K.C., G.L., and T.O.~were supported by the European Research Council (grant no.~864138 ``SingStochDispDyn").
M.G.~was partly supported by the UK Research and Innovation via the grant 
(grant no.~EP/Z534328/1 ``StochFields'')  and by the Simons Collaboration on Probabilistic Paths to QFT.
G.L. was also supported by the NSFC (grant no.~12501181)
and   by the Fundamental Research Funds for the Central Universities. 
T.O.~was also supported by the EPSRC 
Mathematical Sciences
Small Grant  (grant no.~EP/Y033507/1)
and acknowledges support from  
the NSFC (grant no.~W2531005).

\end{ackno}

\end{document}